\documentclass[a4paper]{article}
\pdfoutput=1
\usepackage{latexsym}
\usepackage{amsmath}% \eqref
\usepackage{amssymb}% \mathbb
\usepackage{amsthm}

\usepackage{placeins}% for \FloatBarrier

\usepackage[svgnames]{xcolor}% color names for links
\usepackage[pstarrows]{pict2e}% no slope restrictions

\newtheorem{theorem}{Fact}

\usepackage{microtype}

\usepackage{parskip}% zero Parindent and non-zero Parskip.
\newcommand{\jjterm}[1]{\emph{#1}\marginpar{\small\emph{#1}}}
\newcommand{\jjtermMP}[1]{\marginpar{\small\emph{#1}}}
\newcommand{\eqq}[1]{``{#1}''}%% ``quote''
\newcommand{\jjseqref}[1]{\href{http://oeis.org/#1}{#1}}
\newcommand{\Lmap}[2]{\texttt{#1} $\mapsto$ \texttt{#2}}
\newcommand{\CID}[1]{\texttt{#1}}
\newcommand{\Ref}[1]{\ref{#1} on page \pageref{#1}}
\newcommand{\Tile}[1]{\Theta_{#1}}
\newcommand{\NS}[1]{#1}
\newcommand{\abs}[0]{\mathrm{abs}}
\newcommand{\adeg}[0]{^{\circ}}

\usepackage{graphicx}
\graphicspath{ {./figures/} }

\usepackage{float}%% see \jjfile{float.pdf}
\floatstyle{boxed}%% caption below
\restylefloat{figure}%% let choice take effect

\makeatletter
\renewcommand{\fps@figure}{h!tbp}
\renewcommand{\fps@table}{h!tbp}
\makeatother

\usepackage{ifpdf}
\ifpdf
\usepackage[hyperindex,pdftex,breaklinks=true,backref=page,colorlinks=true,linkcolor=DarkBlue,citecolor=DarkBlue,urlcolor=DarkBlue]{hyperref}%
\hypersetup{%
pdfauthor={Joerg Arndt},%
pdftitle={Plane-filling curves on all uniform grids},%
pdfsubject={plane-filling curves},%
pdfkeywords={plane-filling curves, uniform tilings},%
bookmarksopen,% bookmark tree open
bookmarksopenlevel={0},% chapter
}
\else
\usepackage[hypertex,hyperindex,breaklinks=true,backref=page,colorlinks=true,linkcolor=DarkBlue,citecolor=DarkBlue,urlcolor=DarkBlue]{hyperref}%
\fi

\makeatletter
\@addtoreset{equation}{section}

\@addtoreset{figure}{subsection}

\makeatother
\begin{document}
\bibliographystyle{plain}
\title{Plane-filling curves on all uniform grids}
\author{J\"{o}rg Arndt,  \texttt{<arndt@jjj.de>}\\
 Technische Hochschule N\"{u}rnberg}
%
%\date{\today}
%
%\date{July 8, 2016}%% xxx date, for arxiv
%
\date{May 31, 2017}%% xxx date, for Computers in Scientific Discovery 8 [CSD8]
% August 23-26, 2017, Mons, Belgium
%
\pagestyle{plain}
\maketitle

%%%%%%%%%%%%%%%%%%%%%%%%%%%%%%%%%%%%%%%%%%%%%%%%%%%%%%%%%%%%
\begin{abstract}\noindent
%
%% see file abstract.txt
%
We describe a search for plane-filling curves
traversing all edges of a grid once.
The curves are given by Lindenmayer systems with only
one non-constant letter.
%  which we call \jjterm{simple} L-systems.
%
All such curves for small orders on three grids have been found.
For all uniform grids we show how curves traversing all points once
can be obtained from the curves found.
Curves traversing all edges once are described
for the four uniform grids where they exist.
%
%A non-commutative product of curves is defined
%using a composition of L-systems,
%and also a method for the division of a curve into parts.
%
\ifpdf
%\\ \quad% same layout through formats
\else
\\ \textbf{WARNING:} Only the pdf and postscript formats contain the images.
\fi
\end{abstract}

%%%%%%%%%%%%%%%%%%%%%%%
{% xxx layout
%\small
%\vspace*{-2em}
%\baselineskip3.8mm
%\baselineskip0.95em
\baselineskip4mm
\tableofcontents
}
%%%%%%%%%%%%%%%%%%%%%%%

\clearpage% xxx

%%%%%%%%%%%%%%%%%%%%%%%%%%%%%%%%%%%%%%%%%%%%%%%%%%%%%%%%%%%%
%%%%%%%%%%%%%%%%%%%%%%%%%%%%%%%%%%%%%%%%%%%%%%%%%%%%%%%%%%%%
\section{Introduction}%\label{sect:intro}

%%%%%%%%%%%%%%%%%%%%%%%%%%%%%%%%%%%%%%%%%%%%%%%%%%%%%%%%%%%%
\subsection{Self-avoiding edge-covering curves on a grid}
%\the\textwidth

%%%%%%%%%%%%%%%%%%%%%%%%%%
%
\begin{figure}[h!tbp]
\ifpdf
\begin{center}
{\includegraphics*[width=30mm, viewport={110 300 310 500}]{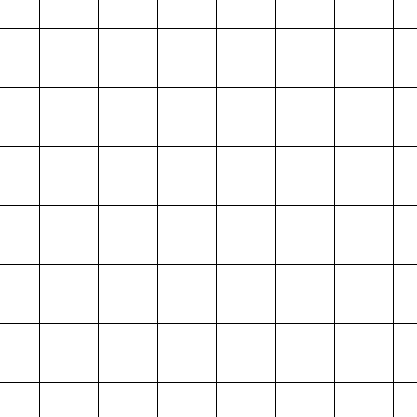}}%
\hspace*{2mm}% layout
{\includegraphics*[width=30mm, viewport={200 390 240 430}]{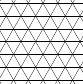}}%
\hspace*{2mm}% layout
{\includegraphics*[width=30mm, viewport={200 405 300 505}]{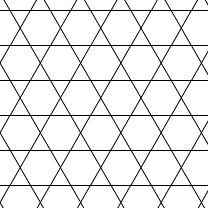}}%
\hspace*{2mm}% layout
{\includegraphics*[width=30mm, viewport={200 400 230 430}]{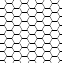}}%
\end{center}
\else
\verb+{see pdf for image}+
\fi
\caption{\label{fig:simple-grids}
From left to right: square grid, triangular grid, tri-hexagonal grid, and hexagonal (honeycomb) grid.}
\end{figure}
%
%%%%%%%%%%%%%%%%%%%%%%%%%%

%
We are interested in self-avoiding, plane-filling curves on certain grids.
A curve is \jjterm{self-avoiding} if it neither crosses itself
nor has an edge that is traversed twice.
It is \jjterm{edge-covering} if it traverses all edges of some grid.

%%%%%%%%%%%%%%%%%%%%%%%%%%
% stringsubst 4 F  F F+F+F-F-F + + - -  | tail -1 | ./bin 4 2 0 0 0.15 > tmp-pic.tex && make dotex # R5-dragon
\begin{figure}[h!tbp]
\ifpdf
\begin{center}
{\includegraphics*[width=80mm, viewport={50 430 500 750}]{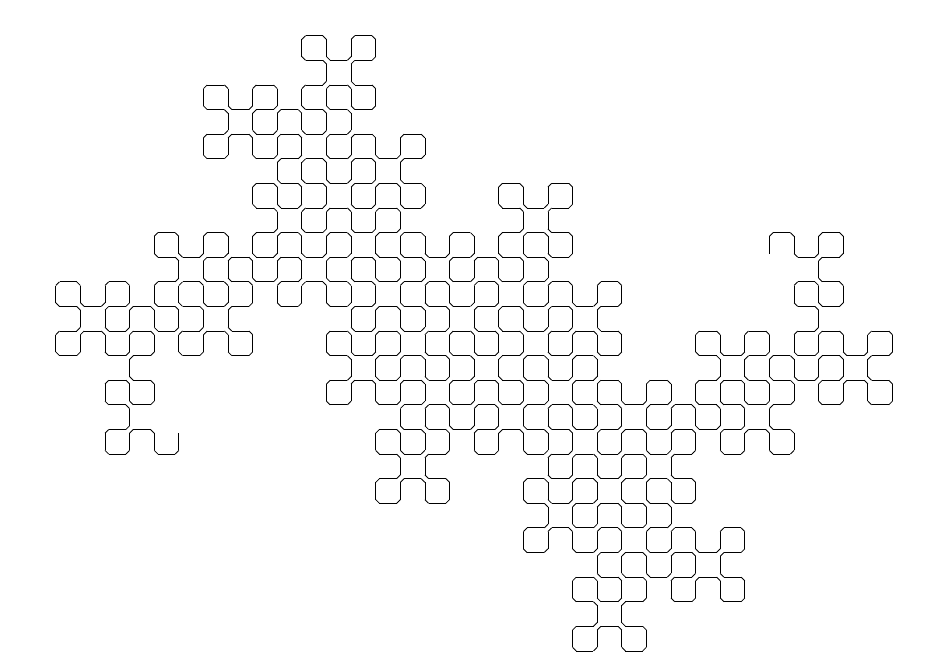}}
\end{center}
\else
\verb+{see pdf for image}+
\fi
\caption{\label{fig:r5-dragon}
The R5-dragon, a curve on the square grid (\CID{R5-1}).}
\end{figure}
%
%%%%%%%%%%%%%%%%%%%%%%%%%%

%%%%%%%%%%%%%%%%%%%%%%%%%%
% stringsubst 6 F  F F+F-F + + - -  | tail -1 | ./bin 3 2 0 0 0.15 > tmp-pic.tex && make dotex # terdragon
\begin{figure}[h!tbp]
\ifpdf
\begin{center}
{\includegraphics*[width=80mm, viewport={50 490 500 740}]{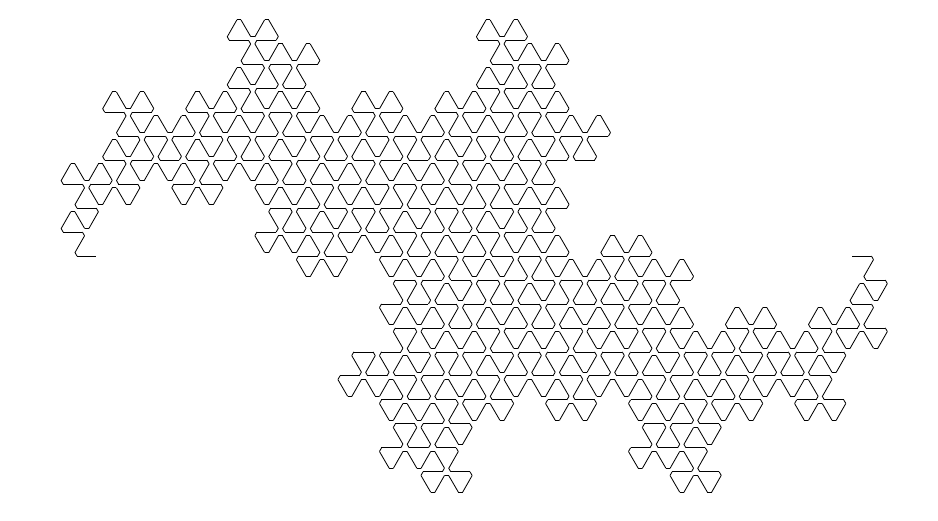}}
\end{center}
\else
\verb+{see pdf for image}+
\fi
\caption{\label{fig:terdragon}
The terdragon, a curve on the triangular grid (\CID{R3-1}).}
\end{figure}
%
%%%%%%%%%%%%%%%%%%%%%%%%%%

%%%%%%%%%%%%%%%%%%%%%%%%%%
% stringsubst 3 F  F F+F+F+F--F--F+F   + + - - | tail -1 | ./bin 6 2 0 0 0.15 > tmp-pic.tex && make dotex # R7-b-1
\begin{figure}[h!tbp]
\ifpdf
\begin{center}
{\includegraphics*[width=80mm, viewport={70 430 480 740}]{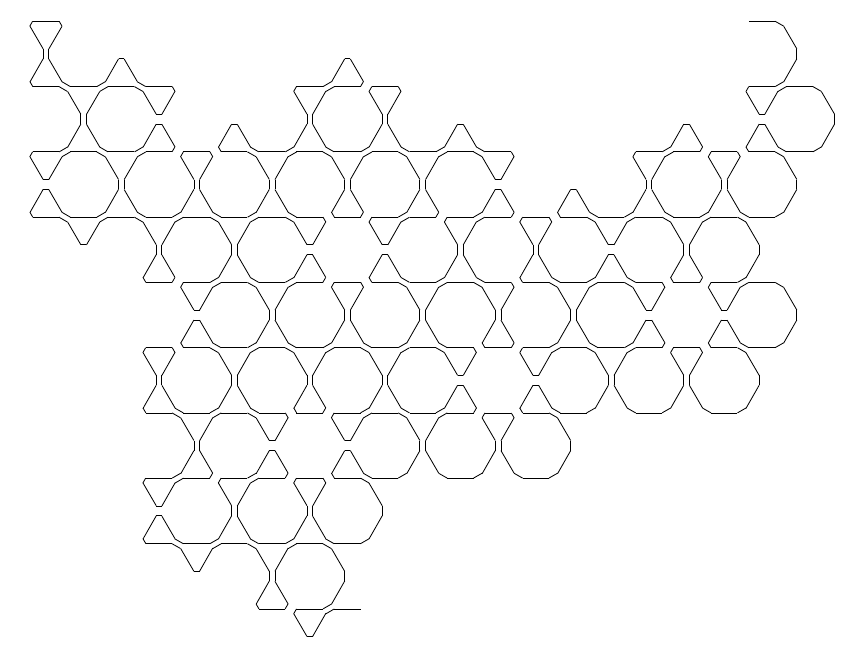}}
\end{center}
\else
\verb+{see pdf for image}+
\fi
\caption{\label{fig:r7-b-curve}
Ventrella's curve, a curve on the tri-hexagonal grid (\CID{R7-1}).}
\end{figure}
%
%%%%%%%%%%%%%%%%%%%%%%%%%%

The curves traverse three grids:
the square grid (Figure~\ref{fig:r5-dragon}),
the triangular grid (Figure~\ref{fig:terdragon}),
and the tri-hexagonal grid (Figure~\ref{fig:r7-b-curve}).
The ID of each curve such as \CID{R7-1} is
explained in section \ref{sect:lsys-and-id}.
The turns are slightly rounded in these figures to make
the curves visually apparent.

For edge-covering curves to exist,
the grid must have an even number of incident edges at each point.
Otherwise a dead end is produced after the point is traversed sufficiently often.
This rules out the hexagonal (honeycomb) grid as every point has incidence 3.

%%%%%%%%%%%%%%%%%%%%%%%%%%%%%%%%%%%%%%%%%%%%%%%%%%%%%%%%%%%%
\subsection{Description via simple Lindenmayer systems}%\label{sect:}
%

%% definition
A \jjterm{Lindenmayer system}
(L-system, see \cite[Section~1.2, pp.~3ff]{algorithmic-beauty-of-plants})
 is a triple $(\Omega, A, P)$
where $\Omega$ is an alphabet,
$A$ a word over $\Omega$ (called the \jjterm{axiom}),
and $P$ a set of maps from letters $\in{}\Omega$
to words over $\Omega$ that contains one map for each letter.

The word that a letter is mapped to is called the
\jjterm{production} of the letter.
If the map for a letter is the identity,
we call the letter a \jjterm{constant} of the L-system.

We specify curves by L-systems interpreted as a sequence of (unit-length) edges and turns.
The curves can be rendered via turtle-graphics,
see \cite[Section~1.3, pp.~6ff]{algorithmic-beauty-of-plants}.
%
%We start at an arbitrary point and keep track of
The initial position and direction are arbitrary.
Letters are interpreted as \eqq{draw a unit stroke in the current direction},
\texttt{+} and \texttt{-} as turns by $\pm$ a fixed angle $\phi$
(set to either $60\adeg$, $90\adeg$, or $120\adeg$).
We will also use the constant letter \texttt{0}
for turns by $0\adeg$ (non-turns).

As an example we take the L-system with
alphabet $\Omega=\{$\texttt{L, R, +, -}$\}$,
axiom $A=$~\texttt{L},
and the maps $P=\{$\Lmap{L}{L+R}, \Lmap{R}{L-R}, \Lmap{+}{+}, \Lmap{-}{-}$\}$.
%
%Here \texttt{L+R} is the production of \texttt{L}
There are two constants, \texttt{+} and \texttt{-}, in this L-system.

Starting with the axiom and repeatedly applying the maps $\in{}P$,
we obtain the words shown in Figure~\ref{fig:lsys-iterates}.
We call the word obtained after the maps were applied $n$ times ($n\geq{}0$)
the $n$th \jjterm{iterate} of the L-system.

%%%%%%%%%%%%%%%%%%%%%%%%%%
% for t in $(seq 0 6) ; do; stringsubst $t L L L+R R L-R + + - - | tail -1; done
\begin{figure}[h!tbp]
\begin{center}
\begin{verbatim}
0: L
1: L+R
2: L+R+L-R
3: L+R+L-R+L+R-L-R
4: L+R+L-R+L+R-L-R+L+R+L-R-L+R-L-R
5: L+R+L-R+L+R-L-R+L+R+L-R-L+R-L-R+L+R+L-R+L+R-L-R-L+R+L-R-L+R-L-R
\end{verbatim}
\end{center}
%6: L+R+L-R+L+R-L-R+L+R+L-R-L+R-L-R+L+R+L-R+L+R-L-R-L+R+L-R-L+R-L-R+L+R+L-R+L+R-L-R+L+R+L-R-L+R-L-R-L+R+L-R+L+R-L-R-L+R+L-R-L+R-L-R
\caption{\label{fig:lsys-iterates}
The $n$th iterates ($n\in{}\{0,1,2,3,4,5\}$) of a Lindenmayer system.}
\end{figure}
%
%%%%%%%%%%%%%%%%%%%%%%%%%%

%%%%%%%%%%%%%%%%%%%%%%%%%%
% stringsubst 10 -L  L L+R    R L-R + + - -  | tail -1 | ./bin 4 2 0 0 0.20 > tmp-pic.tex && make dotex # Heighway-Harter dragon
\begin{figure}[h!tbp]
\ifpdf
\begin{center}
{\includegraphics*[width=90mm, viewport={60 450 500 750}]{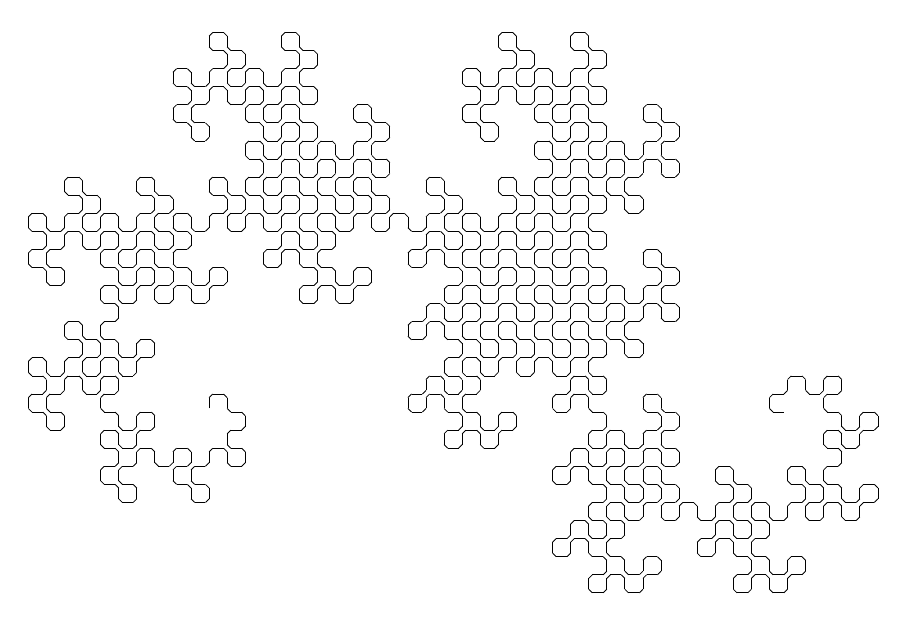}}
\end{center}
\else
\verb+{see pdf for image}+
\fi
\caption{\label{fig:hh-dragon}
The Heighway-Harter dragon.}
\end{figure}
%
%%%%%%%%%%%%%%%%%%%%%%%%%%

%
The curve shown in Figure~\ref{fig:hh-dragon}
found in the 1960s by John Heighway
(sometimes called \jjterm{Heighway-Harter dragon})
corresponds to the L-system just described:
we interpret \texttt{+} and \texttt{-} as turns by $\pm{}90\adeg$ and
both \texttt{L} and \texttt{R} as edges of unit length.
The curve shown corresponds to the $10$th iterate of the L-system.
A description in terms of the paper-folding sequence
is given in \cite[Chapter~5, pp.~152ff]{AutoSeq}.

The so-called \jjterm{terdragon} shown in Figure~\ref{fig:terdragon} corresponds to the L-system
with axiom \texttt{F} and maps
%\texttt{F} $\mapsto$ \texttt{F+F-F}, \texttt{+} $\mapsto$ \texttt{+}, and \texttt{-} $\mapsto$ \texttt{-}
\Lmap{F}{F+F-F}, \Lmap{+}{+}, and \Lmap{-}{-}
where the turns are by $120\adeg$ and \texttt{F} is an edge.
Both curves were described by Chandler Davis and Donald Knuth \cite{davis-knuth-old} in 1970
(an extended version of the paper is reprinted in \cite{davis-knuth-new}).
%

%%%%%%%%%%%%%%%%%%%%%%%%%%
%
%% Render with thick lines: lnth *= 3.0; :
% stringsubst 0 -F_0F_+F_0F_-F_-F_+F_  _ _  F F0F+F0F-F-F+F   0 0 + + - - | tail -1 | ./bin 3 3 1 > tmp-pic.tex && make dotex # R7-1
% stringsubst 1 -F_0F_+F_0F_-F_-F_+F_  _ _  F F0F+F0F-F-F+F   0 0 + + - - | tail -1 | ./bin 3 3 1 > tmp-pic.tex && make dotex # R7-1
%
\begin{figure}[h!tbp]
\ifpdf
\begin{center}
{\includegraphics*[width=50mm, viewport={100 360 450 700}]{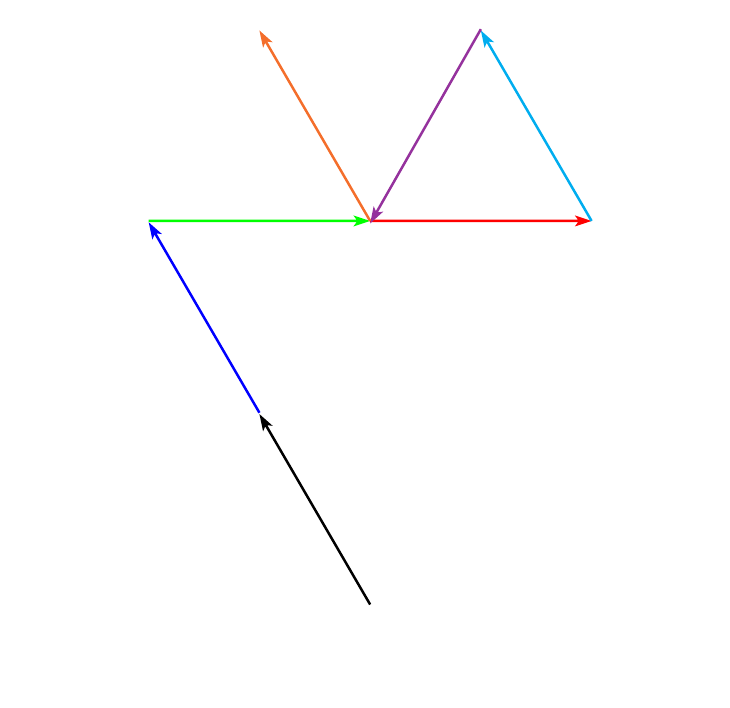}}%
\qquad% layout
{\includegraphics*[width=50mm, viewport={110 350 450 720}]{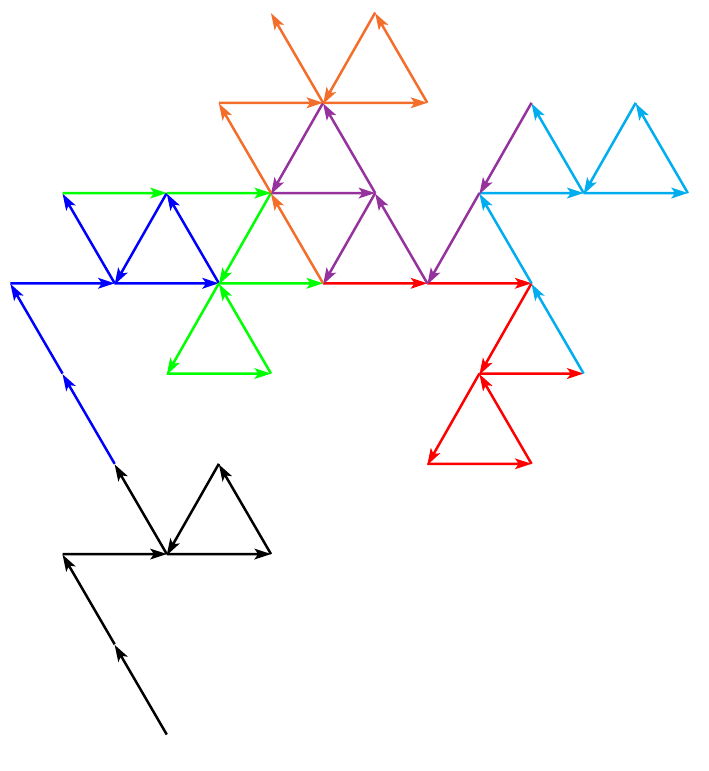}}
\end{center}
\else
\verb+{see pdf for image}+
\fi
\caption{\label{fig:iterate-1-2-decomp}
First iterate (motif) and second iterate of a curve of order 7 (\CID{R7-1}).}
\end{figure}
%
%%%%%%%%%%%%%%%%%%%%%%%%%%

%%%%%%%%%%%%%%%%%%%%%%%%%%
%% with lnth *= 2.0;  // thicker lines
% stringsubst 2 -F_0F_+F_0F_-F_-F_+F_  _ _  F F0F+F0F-F-F+F   0 0 + + - - | tail -1 | ./bin 3 3 0 > tmp-pic.tex && make dotex # R7-1
% stringsubst 3 -F_0F_+F_0F_-F_-F_+F_  _ _  F F0F+F0F-F-F+F   0 0 + + - - | tail -1 | ./bin 3 3 0 > tmp-pic.tex && make dotex # R7-1
% stringsubst 4 -F_0F_+F_0F_-F_-F_+F_  _ _  F F0F+F0F-F-F+F   0 0 + + - - | tail -1 | ./bin 3 3 0 > tmp-pic.tex && make dotex # R7-1
%
\begin{figure}[h!tbp]
\ifpdf
\begin{center}
\includegraphics*[width=130mm, page=1, viewport={0 410 600 740}]{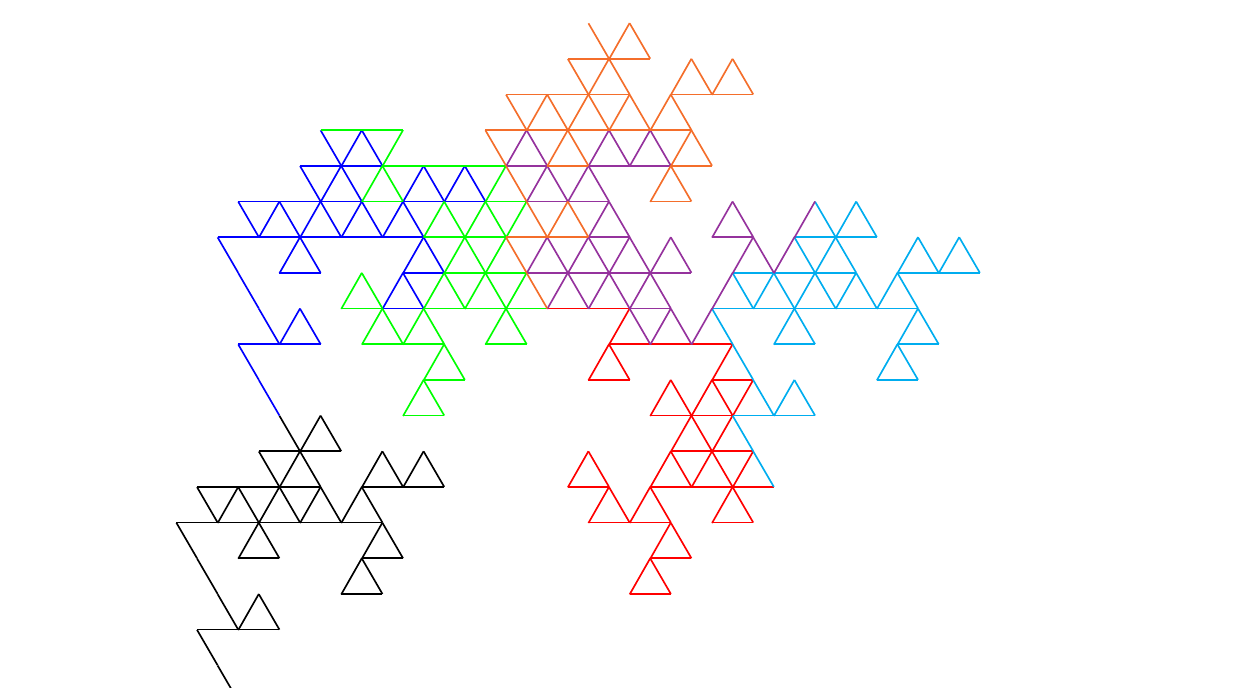}
\includegraphics*[width=130mm, page=2, viewport={0 450 600 740}]{iterate-decomp-3-5.pdf}
\includegraphics*[width=130mm, page=3, viewport={0 450 600 740}]{iterate-decomp-3-5.pdf}
\end{center}
\else
\verb+{see pdf for image}+
\fi
\caption{\label{fig:iterate-3-4-5-decomp}
Third, fourth, and fifth iterate of the curve \CID{R7-1}.}
\end{figure}
%
%%%%%%%%%%%%%%%%%%%%%%%%%%

We call an L-system \jjterm{simple} if it has just one non-constant letter.
The terdragon corresponds to a simple L-system but the Heighway-Harter dragon does not.
Only curves with simple L-systems are considered for the search to keep the search space manageable.

For simple L-systems we always use \texttt{F} for the only non-constant letter.
The axiom (\texttt{F}) and the maps for the constant letters
(\Lmap{+}{+}, \Lmap{-}{-}) will be omitted.
The \jjterm{order} $R$ of a curve is the number of \texttt{F}s
in the production of \texttt{F} by the (simple) L-system.
For example, the terdragon (Figure~\ref{fig:terdragon}) has order $R=3$
(\Lmap{F}{F+F-F})
and the curve shown in Figure~\ref{fig:r5-dragon} has order $R=5$
(\Lmap{F}{F+F+F-F-F}).

\begin{figure}[h!tbp]
\ifpdf
\begin{center}
{\includegraphics*[width=35mm, viewport={160 360 420 750}]{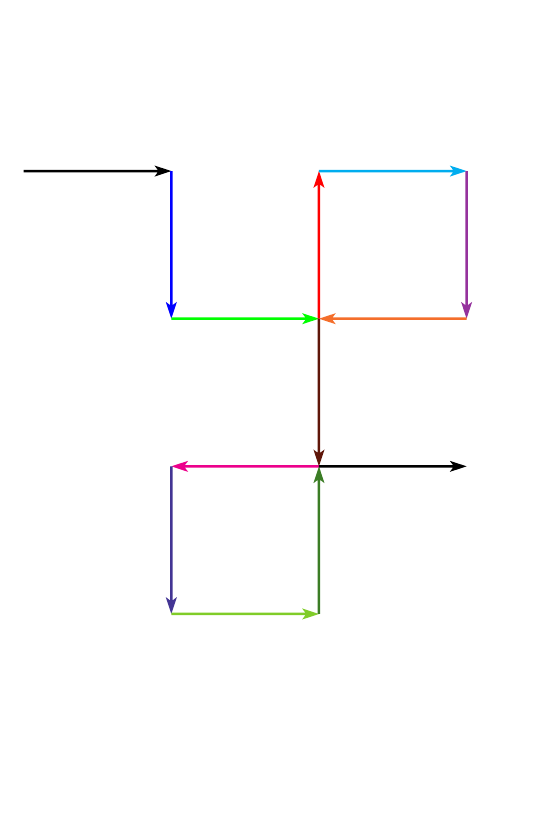}}%
{\includegraphics*[width=90mm, viewport={40 380 490 750}]{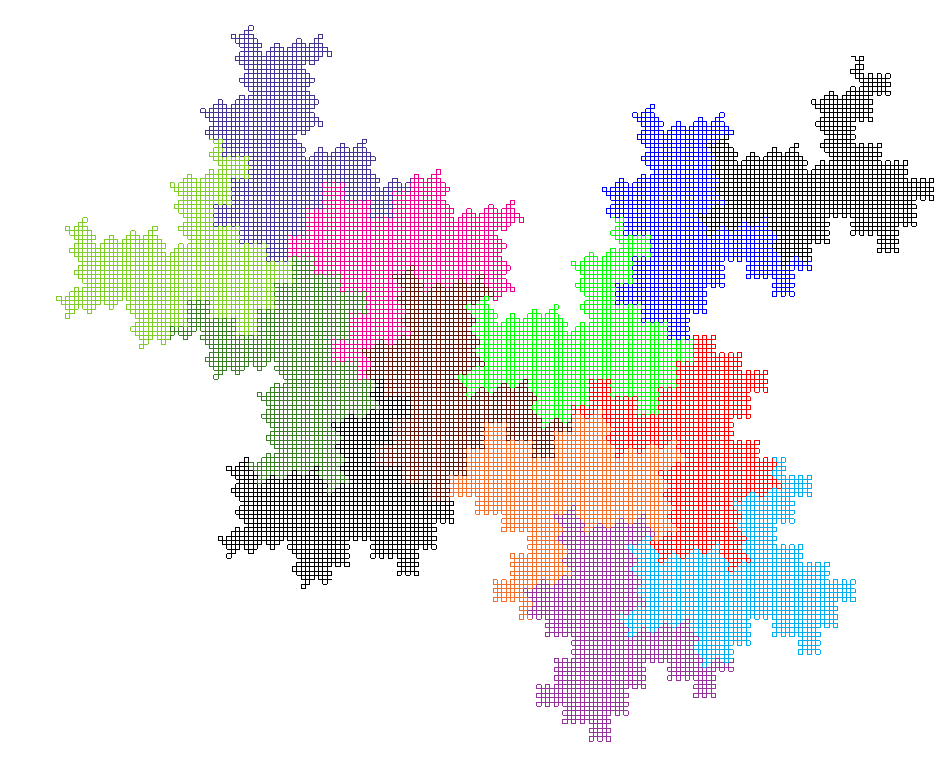}}
\end{center}
\else
\verb+{see pdf for image}+
\fi
\caption{\label{fig:r13-q-3-asymmm}
First iterate (motif) and fourth iterate of a curve of order 13 (\CID{R13-3}).}
\end{figure}

%
% \xxx{term "approximation" (Helmberg)}
We call the curve corresponding to the $n$th iterate of an L-system
the $n$th \jjterm{iterate of the curve}.
We call the first iterate the \jjterm{motif} of the curve.
Iterate 0 corresponds to a single edge and
iterate $n$ is obtained from iterate $n-1$ by replacing
every edge by the motive.
Figures~\ref{fig:iterate-1-2-decomp} and \ref{fig:iterate-3-4-5-decomp}
show the iterates 1 through 5 of a curve of order 7 on the triangular grid with
L-system \Lmap{F}{F0F+F0F-F-F+F}.
Here the constant letter \texttt{0} stands for \eqq{no turn}.
%and order 7.
%
The motif and fourth iterate of the curve of order 13
on the square grid with
L-system \Lmap{F}{F+F-F-F+F+F+F-F+F-F-F-F+F}
are shown in Figure~\ref{fig:r13-q-3-asymmm}.

Note that Davis and Knuth use \eqq{order} for what we call \eqq{iterate}
\cite[Figure~21, p.~600]{davis-knuth-new}.

It is our goal to find \emph{all} self-avoiding edge-covering
curves with simple L-systems of small orders on the three grids mentioned.

%%%%%%%%%%%%%%%%%%%%%%%%%%
%
%% Render with thick lines: lnth *= 3.0; :
% stringsubst 0 F_+F_-F_-F_+F_+F_+F_-F_+F_-F_-F_-F_+F_ _ _ F F+F-F-F+F+F+F-F+F-F-F-F+F + + - - | tail -1 | ./bin 4 3 1 > tmp-pic.tex && make dotex # R13-3
%
% Render with normal lines:
% stringsubst 3 F_+F_-F_-F_+F_+F_+F_-F_+F_-F_-F_-F_+F_ _ _ F F+F-F-F+F+F+F-F+F-F-F-F+F + + - - | tail -1 | ./bin 4 3 0 > tmp-pic.tex && make dotex # R13-3
%
%
%%%%%%%%%%%%%%%%%%%%%%%%%%

%%%%%%%%%%%%%%%%%%%%%%%%%%%%%%%%%%%%%%%%%%%%%%%%%%%%%%%%%%%%
\subsection{Some known plane-filling curves}%\label{sect:}

%%%%%%%%%%%%%%%%%%%%%%%%%%
%% with   lnth *= 2.0;  // thicker lines
% stringsubst 3 L  L LtRtL-t-RtLtR+t+LtRtL  R RtLtR+t+LtRtL-t-RtLtR  t t  + +  - - | tail -1 | tr -d RL |./bin 4 3 0 > tmp-pic.tex && make dotex # Peano curve
%
% stringsubst 4 L  L +Rt-LtL-tR+  R -Lt+RtR+tL-  + + - - t t | tr -d LR | tail -1 | ./bin 4 3 0 > tmp-pic.tex && make dotex # Hilbert curve
%
\begin{figure}[h!tbp]
\ifpdf
\begin{center}
{\includegraphics*[width=55mm, viewport={50 330 500 740}]{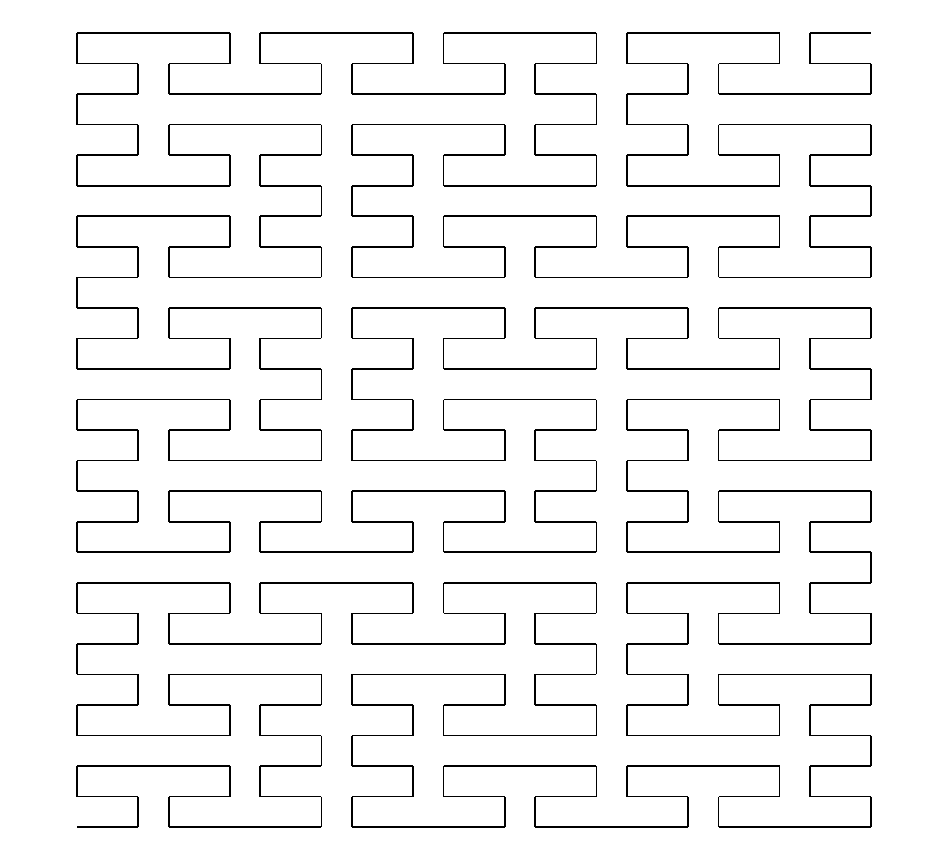}}%
{\includegraphics*[width=55mm, viewport={50 340 500 740}]{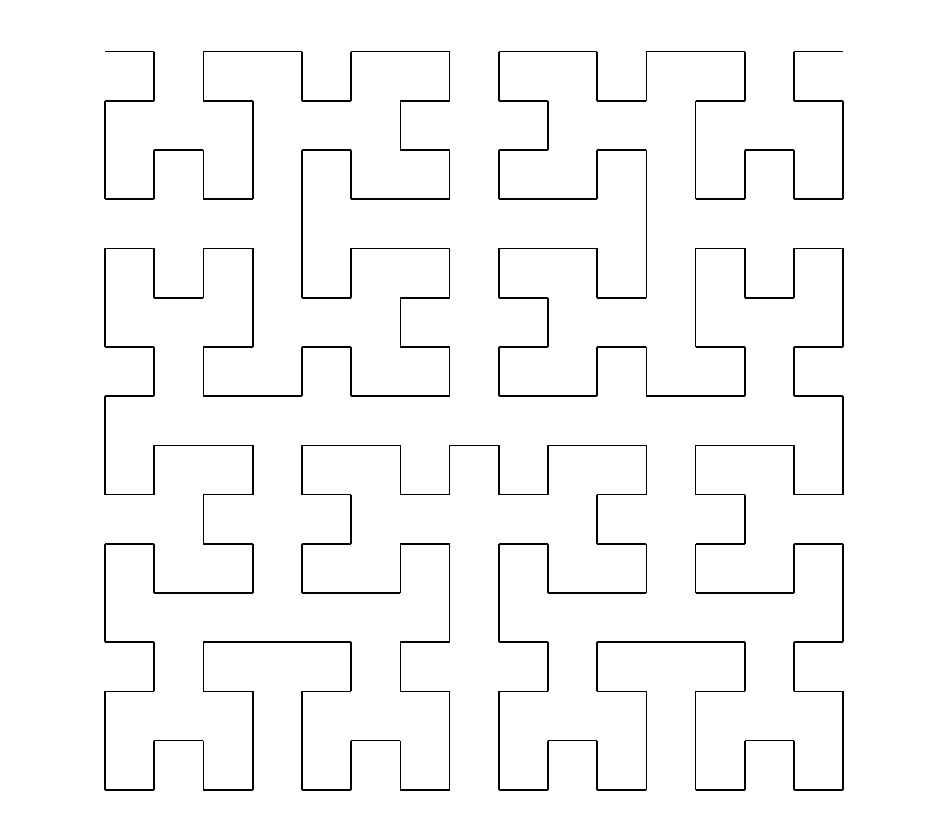}}
\end{center}
\else
\verb+{see pdf for image}+
\fi
\caption{\label{fig:peano-hilbert}
The Peano curve (left) and the Hilbert curve (right).}
\end{figure}
%
%%%%%%%%%%%%%%%%%%%%%%%%%%

%%%%%%%%%%%%%%%%%%%%%%%%%%
%% with   lnth *= 2.0;  // thicker lines
% stringsubst 3 L   L L+R++R-L--LL-R+   R -L+RR++R+L--L-R   + +   - - | tail -1 | ./bin 6 3 0 > tmp-pic.tex && make dotex # Gosper's "flowsnake" curve
%
% stringsubst 3 L   L L+R++R-L--LL-R+   R -L+RR++R+L--L-R   + +   - - | tail -1 | sed 's/L/t-t+/g; s/R/-t+t/g;' | ./bin 6 3 0 > tmp-pic.tex && make dotex
%
\begin{figure}[h!tbp]
\ifpdf
\begin{center}
{\includegraphics*[width=55mm, viewport={40 300 500 740}]{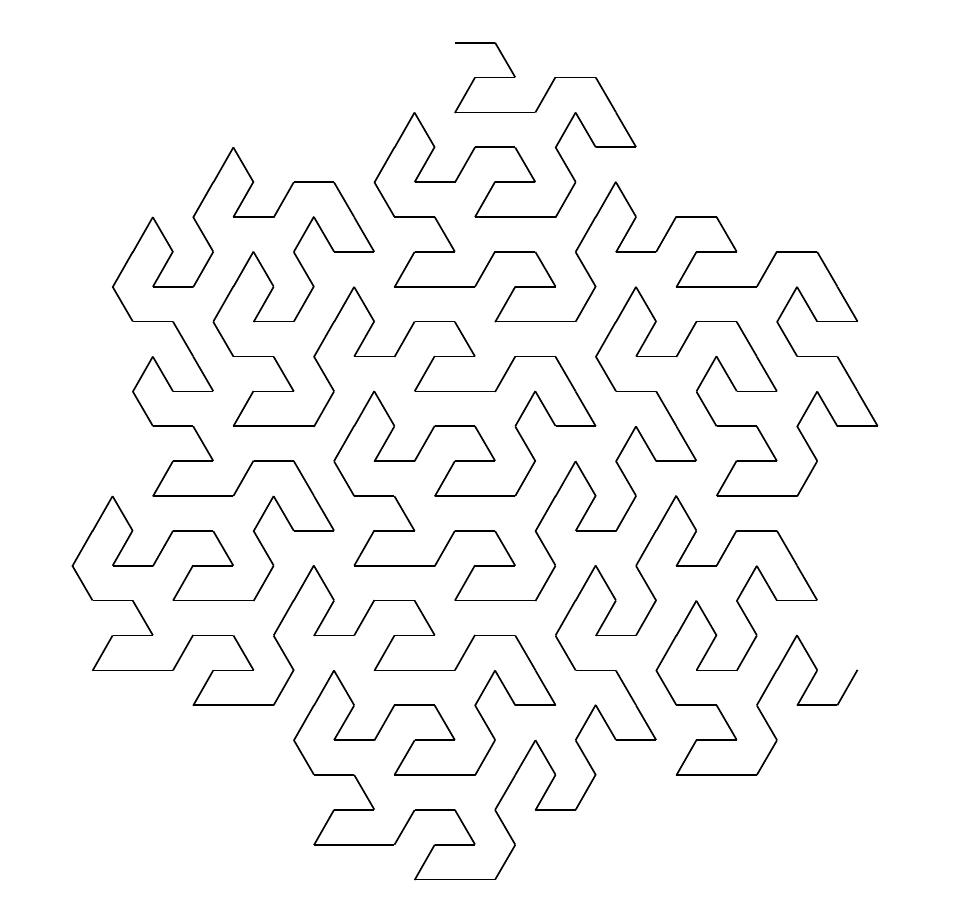}}%
{\includegraphics*[width=55mm, viewport={40 320 500 740}]{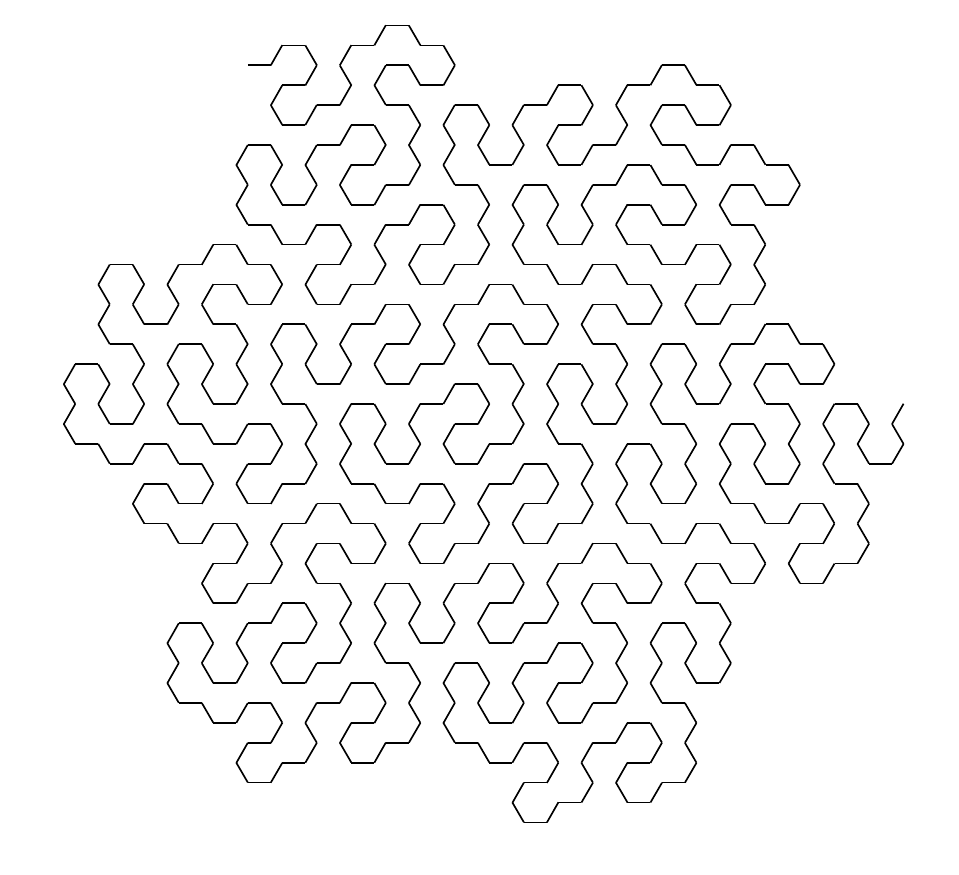}}
\end{center}
\else
\verb+{see pdf for image}+
\fi
\caption{\label{fig:gosper-flowsnake}
Two renderings of Gosper's flowsnake.}
\end{figure}
%
%%%%%%%%%%%%%%%%%%%%%%%%%%

%The two earliest examples of plane-filling curves were described
%by Giuseppe Peano in 1890 \cite{peano-curve-1890}
%(also see \cite[Chapter~3, pp.~31ff]{sagan}
%and \cite{cole-peano-gray3} for the construction of $n$-dimensional
%Peano curves via the Gray code for base 3)
%and
%David Hilbert in 1891 \cite{hilbert-hilbert}, see Figure~\ref{fig:peano-hilbert}.
%

The earliest example of a plane-filling curve was described\jjtermMP{Peano curve}
by Giuseppe Peano in 1890 \cite{peano-curve-1890}
(also see \cite[Chapter~3, pp.~31ff]{sagan}
and \cite{cole-peano-gray3} for the construction of $n$-dimensional
Peano curves via the Gray code for base 3).
In 1891 David Hilbert\jjtermMP{Hilbert curve}
described the curve now named after him \cite{hilbert-hilbert}.
Both curves are shown in Figure~\ref{fig:peano-hilbert}.

In the 1970s Bill Gosper discovered his \jjterm{flowsnake} curve
shown in Figure~\ref{fig:gosper-flowsnake} \cite{gardner-monster-curves}.
%% Bill's email:
%% Subject: Re: Question about your "flowsnake" curve
%% Date: Sun, 1 Nov 2015 09:47:51 -0800
%% Message-ID: <CAA-4O0Evk2j7KV6fzwxWyZyyTajjEzrWaEujWjRG6sxVNoYvzw@mail.gmail.com>
%

Gosper's curve \ref{fig:gosper-flowsnake}
can be given by the axiom \texttt{L}
and non-constant maps
\Lmap{L}{L+R++R-L--LL-R+} and
\Lmap{R}{-L+RR++R+L--L-R}
where the turns are by $\pm{}60$ degrees.
Both \texttt{L} and \texttt{R} correspond to edges.
%
% sed 's/L/t-t+/g; s/R/-t+t/g;'
The rendering on the right is obtained via
post-processing replacing
all \texttt{L} by \texttt{L-L+}
and
all \texttt{R} by \texttt{-R+R}.

More curves resembling Gosper's flowsnake were given by
\NS{Jin} \NS{Akiyama},
\NS{Hiroshi} \NS{Fukuda},
\NS{Hiro} \NS{Ito}, and
\NS{Gisaku} \NS{Nakamura}
in 2007 \cite{akiyama-fukuda}.

%%%%%%%%%%%%%%%%%%%%%%%%%%%
%% stringsubst 4 L  L -RL+L+R-L  R R+L-R-RL+  + +  - -  | tail -1 | ./bin 4 3 0 > tmp-pic.tex && make dotex # Mandelbrot, The Fractal Geometry of Nature, Plate 49, p.49
%%
%\begin{figure}[h!tbp]
%%
%\ifpdf
%\begin{center}
%\includegraphics*[width=85mm, viewport={0 330 590 740}]{mandelbrot-5.pdf}
%\end{center}
%\else
%\verb+{see pdf for image}+
%\fi
%\caption{\label{fig:mandelbrot-5}
%A curve on the square grid given by Mandelbrot.}
%\end{figure}
%%
%%%%%%%%%%%%%%%%%%%%%%%%%%%
%
%
%%%%%%%%%%%%%%%%%%%%%%%%%%%
%% stringsubst 2 L L -RR+L+L-R-RL-R+LL+R+L-RLL+R+LR-R-L+L+R-R-LL R RR+L+L-R-R+L+LR-L-RRL+R-L-RR-L+RL+L+R-R-LL+ + + - - | tail -1 | ./bin 4 2 0 0 0.15 > tmp-pic.tex && make dotex # E-curve
%%
%\begin{figure}[h!tbp]
%%
%\ifpdf
%\begin{center}
%\includegraphics*[width=85mm, viewport={0 330 590 740}]{e-curve.pdf}
%\end{center}
%\else
%\verb+{see pdf for image}+
%\fi
%\caption{\label{fig:e-curve}
%The E-Curve on the square grid given by \NS{McKenna}.}
%\end{figure}
%%
%%%%%%%%%%%%%%%%%%%%%%%%%%%

%%%%%%%%%%%%%%%%%%%%%%%%%%
%% with   lnth *= 2.0;  // thicker lines
% stringsubst 4 L  L -RL+L+R-L  R R+L-R-RL+  + +  - -  | tail -1 | ./bin 4 3 0 > tmp-pic.tex && make dotex # Mandelbrot, The Fractal Geometry of Nature, Plate 49, p.49
%
% stringsubst 2 L L -RR+L+L-R-RL-R+LL+R+L-RLL+R+LR-R-L+L+R-R-LL R RR+L+L-R-R+L+LR-L-RRL+R-L-RR-L+RL+L+R-R-LL+ + + - - | tail -1 | ./bin 4 3 0 > tmp-pic.tex && make dotex # E-curve
%
\begin{figure}[h!tbp]
\ifpdf
\begin{center}
{\includegraphics*[width=55mm, viewport={60 310 485 730}]{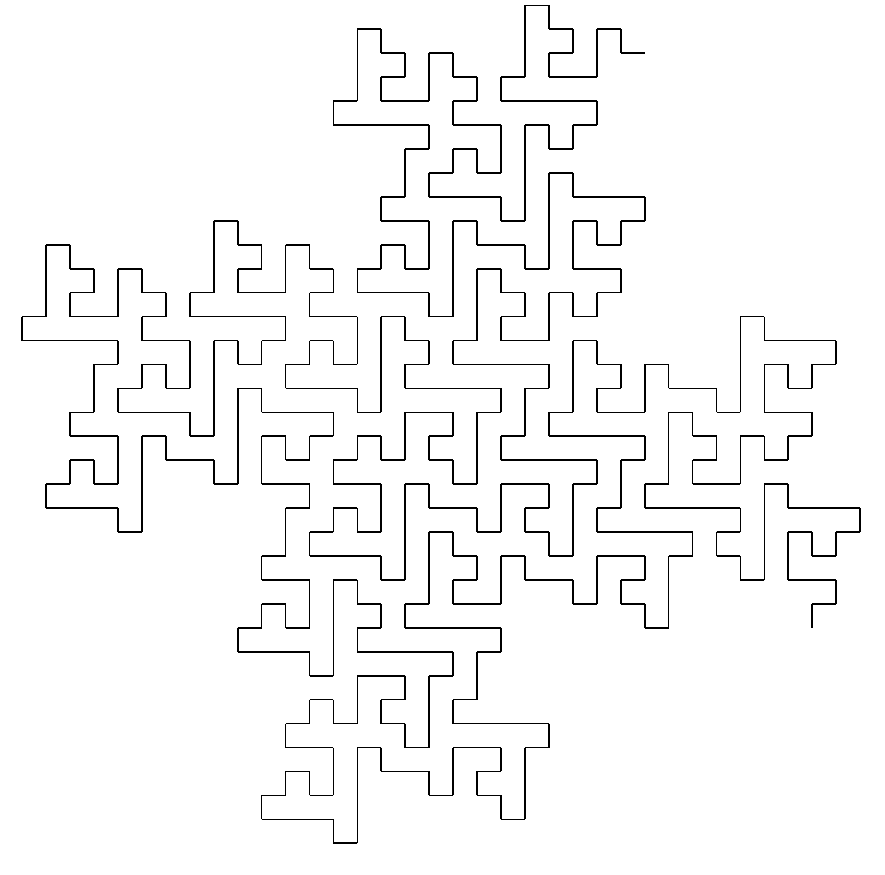}}%
{\includegraphics*[width=55mm, viewport={40 310 510 730}]{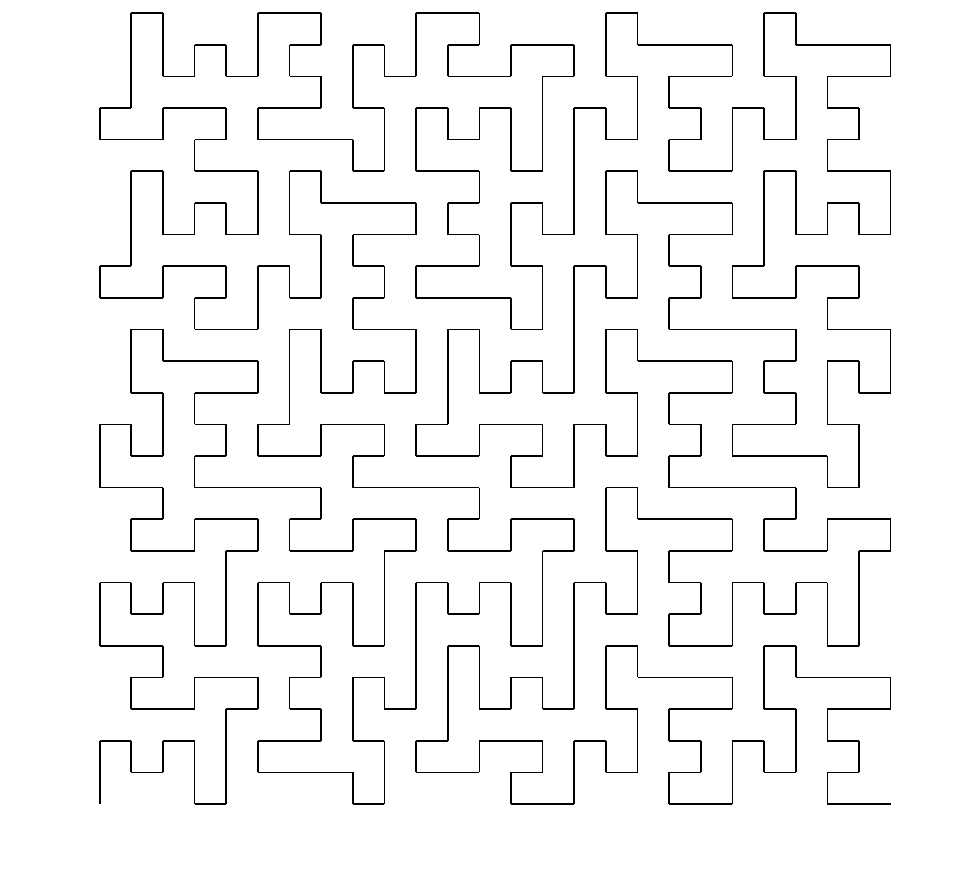}}
\end{center}
\else
\verb+{see pdf for image}+
\fi
\caption{\label{fig:square-curves}
A curve on the square grid given by Mandelbrot (left) and
the E-Curve given by \NS{McKenna} (right).}
\end{figure}
%
%%%%%%%%%%%%%%%%%%%%%%%%%%

Benoit Mandelbrot \cite[Plate~49, p.~49]{mandelbrot-fractal-geom-nature} shows a
curve with the L-system with axiom
\texttt{L} and (non-constant) maps
\Lmap{L}{-RL+L+R-L},
\Lmap{R}{R+L-R-RL+},
%see Figure~\ref{fig:mandelbrot-5}.
see Figure~\ref{fig:square-curves} (left).

Douglas M.\ \NS{McKenna} gave what he calls the \jjterm{E-Curve} in
1994 \cite[p.~60]{mckenna-e-curve},
%see Figure~\ref{fig:e-curve}.
shown in Figure~\ref{fig:square-curves} (right).
It can be generated by the L-system with maps
\begin{center}
\Lmap{L}{-RR+L+L-R-RL-R+LL+R+L-RLL+R+LR-R-L+L+R-R-LL}
\hphantom{.}and\\
\Lmap{R}{RR+L+L-R-R+L+LR-L-RRL+R-L-RR-L+RL+L+R-R-LL+}.
\hphantom{and}
\end{center}

None of these curves corresponds to a simple L-system
or is edge-covering.
They do traverse every point of the grid once.
We will call such curves \jjterm{point-covering}.

The already mentioned curve shown in Figure~\ref{fig:r5-dragon}
first appeared in \cite{dekking-1981} and \cite{dekking-folds-3}
 and later in
\cite[Section~1.31.5, p.~95ff]{fxtbook} (therein called \jjterm{R5-dragon})
and \cite[p.~84]{ventrella}.

%%%%%%%%%%%%%%%%%%%%%%%%%%%
%% stringsubst 5 -F  F F+F0F-F + + - - 0 0  | tail -1 | ./bin 3 2 0 0 0.15 > tmp-pic.tex && make dotex # R4-crab-dragon
%%
%\begin{figure}[h!tbp]
%%
%\ifpdf
%\begin{center}
%\fbox{\includegraphics*[width=130mm, viewport={0 380 600 750}]{r4-crab-dragon.pdf}}
%\end{center}
%\else
%\verb+{see pdf for image}+
%\fi
%\caption{\label{fig:r4-crab-dragon}
%A curve of order 4 on the triangular grid (\CID{R4-1}).}
%\end{figure}
%%
%%%%%%%%%%%%%%%%%%%%%%%%%%%
%
%
%%%%%%%%%%%%%%%%%%%%%%%%%%
%
%% stringsubst 5 -F  F F+F0F-F + + - - 0 0  | tail -1 | ./bin 3 2 0 0 0.15 > tmp-pic.tex && make dotex # R4-crab-dragon
%% stringsubst 1 1 1 1221 2 2332 3 3113 | tail -1 | ./bin 3 3 1 1 > tmp-pic.tex && make dotex #
%
%% Render with thick lines: lnth *= 3.0; and with \Huge letter-size:
% stringsubst 0 _1_2_2_1 _ _ 1 1221 2 2332 3 3113 | tail -1 | ./bin 3 3 1 1 > tmp-pic.tex && make dotex #
% stringsubst 1 _1_2_2_1 _ _ 1 1221 2 2332 3 3113 | tail -1 | ./bin 3 3 1 1 > tmp-pic.tex && make dotex #
%
%% with lnth *= 2.0;  // thicker lines
% stringsubst 6 _2_3_3_2 _ _ 1 1221 2 2332 3 3113 | tail -1 | ./bin 3 3 0 > tmp-pic.tex && make dotex #
%
\begin{figure}[h!tbp]
\ifpdf
\begin{center}
{\includegraphics*[width=20mm, viewport={220 380 360 690}]{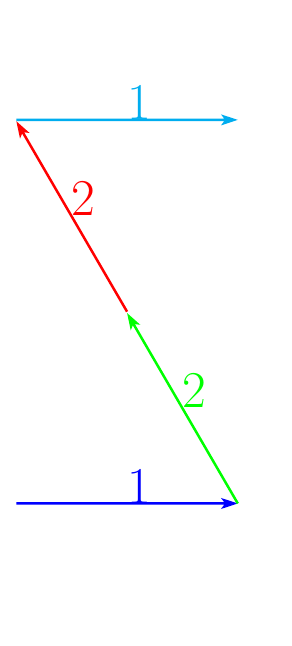}}%
{\includegraphics*[width=40mm, viewport={90 360 420 720}]{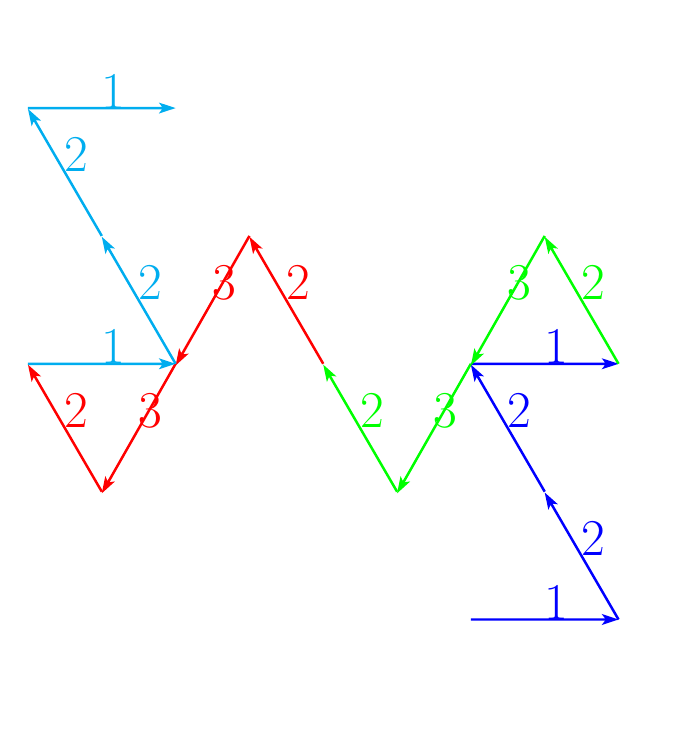}}%
%\fbox{\includegraphics*[width=60mm, viewport={100 150 450 750}]{r4-crab-dragon-it7.pdf}}
{\includegraphics*[width=65mm, viewport={60 370 500 740}]{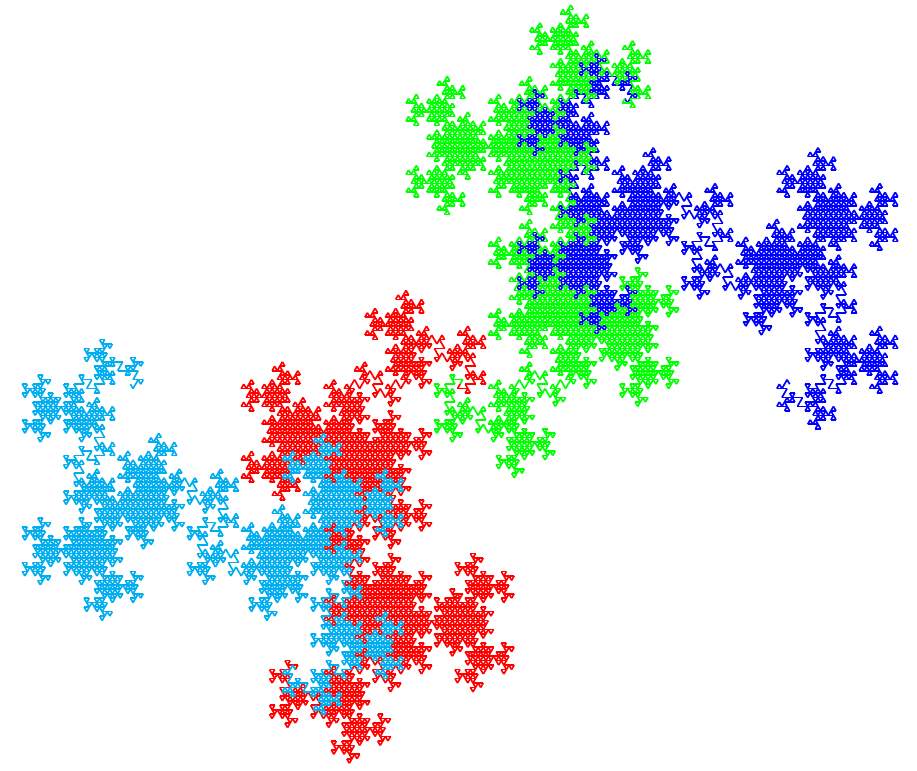}}
\end{center}
\else
\verb+{see pdf for image}+
\fi
\caption{\label{fig:r4-crab-dragon}
First, second, and seventh iterate of a curve of order 4 on the triangular grid (\CID{R4-1}).
The numbers (left and middle) correspond to the directions of the edges.}
\end{figure}
%
%%%%%%%%%%%%%%%%%%%%%%%%%%
%
A lesser known curve on the triangular grid with L-system \Lmap{F}{F+F0F-F},
shown in Figure~\ref{fig:r4-crab-dragon},
is called the \jjterm{crab} by Gilbert Helmberg \cite[Section~1.1.5.2, p.~24ff]{helmberg}.
Note that Helmberg calls \eqq{generator} what we call \eqq{motif}.

A great number of plane-filling curves are given in
Jeffrey Ventrella's 2012 book \eqq{Brain Filling Curves} \cite{ventrella}.
Some of the curves correspond to simple L-systems.
He especially gives the first curve on the tri-hexagonal grid
(page 105 in \cite{ventrella},
see Figure~\ref{fig:r7-b-curve})
and an asymmetric curve on the triangular grid
(top of page 107 in \cite{ventrella},
see Figures~\ref{fig:iterate-1-2-decomp} and \ref{fig:iterate-3-4-5-decomp}).
The latter curve has already been found by Bill Gosper in the 1970s,
but there seems to be no publication prior to Ventrella describing it.
%% Photo of the plot: http://gosper.org/IMG_0245.JPG
%% also see http://www.tweedledum.com/rwg/tril7.htm
%% attempt to recreate image (close, but not exactly):
%% stringsubst 4 F-F-F F F0F+F0F-F-F+F 0 0 + + - - | tail -1 | tr +- -+ | ./bin 3 1 0 > tmp-pic.tex && make dotex

We will only find edge-covering curves in our search.
Methods to obtain point-covering curves from
the ones found are given in section \Ref{sect:PC-curves}.
% like the curves of Hilbert, Peano, Gosper, and Mandelbrot

%%%%%%%%%%%%%%%%%%%%%%%%%%%%%%%%%%%%%%%%%%%%%%%%%%%%%%%%%%%%
\subsection{Other descriptions using L-systems}%\label{sect:}

%%%%%%%%%%%%%%%%%%%%%%%%%%
%% with const char *ssz = "\\Huge";
%% with lnth *= 10.0;  // Very thick lines
% echo [1][2][3] | ./bin 3 3 1 1 > tmp-pic.tex && make dotex
% echo [1][2][3][4] | ./bin 4 3 1 1 > tmp-pic.tex && make dotex
% echo [1][2][3][4][5][6] | ./bin 6 3 1 1 > tmp-pic.tex && make dotex
%
\begin{figure}[h!tbp]
\ifpdf
\begin{center}
%{\includegraphics*[width=40mm, viewport={100 450 490 750}]{directions-3.pdf}}%
%{\includegraphics*[width=40mm, viewport={100 450 490 750}]{directions-4.pdf}}%
%{\includegraphics*[width=40mm, viewport={100 450 490 750}]{directions-6.pdf}}
%% better image(s) by Julia Handl (667261 bytes for this alone!):
%{\includegraphics*[width=120mm, viewport={0 00 400 90}]{curve-search-pfeile.pdf}}% artefact!
%{\includegraphics*[width=120mm, viewport={0 00 400 90}]{curve-search-pfeile-big.pdf}}
%% smaller (115kB):
{\includegraphics*[width=120mm, viewport={0 00 400 90}]{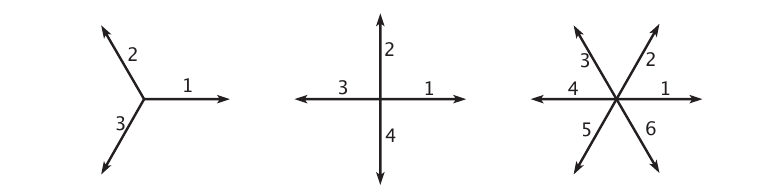}}
\end{center}
\else
\verb+{see pdf for image}+
\fi
\caption{\label{fig:directions}
Numbering for the directions in the
triangular grid (left),
square grid (middle),
and the hexagonal and tri-hexagonal grid (right).
}
\end{figure}
%
%%%%%%%%%%%%%%%%%%%%%%%%%%

%%%%%%%%%%%%%%%%%%%%%%%%
\paragraph{Description by directed edges.}
The curves with simple L-systems
can also be described by L-systems specifying the directions
of the edges in the drawing.
There are 3 possible directions on the triangular grid.
We assign the letters \texttt{1} (right), \texttt{2} (upward left), and \texttt{3} (downward left)
to them, as shown in Figure~\ref{fig:directions}.
The word for the motif (left image in Figure~\ref{fig:r4-crab-dragon}),
here \texttt{1221}, is the production of its leftmost letter,
giving the map \Lmap{1}{1221}.
The other maps corresponding to the remaining rotations
are \Lmap{2}{2332} and \Lmap{3}{3113}.
The map for $k\geq{2}$ is obtained by incrementing all numbers in the
map for $k-1$, taking the result modulo 3.
The image in the middle of Figure~\ref{fig:r4-crab-dragon}
shows the second iterate with a coloration that
allows to identify the corresponding parts in the motif.

%% stringsubst 4 1 1 12321 2 23432 3 34143 4 41214 | tail -1 | ./bin 4 3 1 > tmp-pic.tex && make dotex #
The L-system for the R5-dragon shown in
Figure~\ref{fig:r5-dragon} is \Lmap{F}{F+F+F-F-F}
and the maps for the edge directions are
\Lmap{1}{12321},
\Lmap{2}{23432},
\Lmap{3}{34143}, and
\Lmap{4}{41214}.
The interpretations as directions are
\texttt{1} (right),
\texttt{2} (up),
\texttt{3} (left), and
\texttt{4} (down).
%

%% for R13-15 see https://oeis.org/A265671
%
% \Lmap{F}{F+F0F0F-F-F+F0F+F+F-F0F-F}
% \Lmap{+}{+00--+0++-0-+}
% \Lmap{0}{+00--+0++-0-0}
% \Lmap{-}{+00--+0++-0--}
%
% \Lmap{1}{1222131123221}
% \Lmap{2}{2333212231332}
% \Lmap{3}{3111323312113}
%
% stringsubst 0 1_2_2_2_1_3_1_1_2_3_2_2_1 _ _ 1 1222131123221 2 2333212231332 3 3111323312113 | tail -1 | ./bin 3 3 1 1 > tmp-pic.tex && make dotex # R13-15
% r13-15-edge-subst-iterate1.pdf
%
% stringsubst 1 1_2_2_2_1_3_1_1_2_3_2_2_1 _ _ 1 1222131123221 2 2333212231332 3 3111323312113 | tail -1 | ./bin 3 3 1 1 > tmp-pic.tex && make dotex # R13-15
% r13-15-edge-subst-iterate2.pdf
%
% stringsubst 2 1_2_2_2_1_3_1_1_2_3_2_2_1 _ _ 1 1222131123221 2 2333212231332 3 3111323312113 | tail -1 | ./bin 3 3 0 0 > tmp-pic.tex && make dotex # R13-15
% r13-15-edge-subst-iterate3.pdf
%
% stringsubst 3 1_2_2_2_1_3_1_1_2_3_2_2_1 _ _ 1 1222131123221 2 2333212231332 3 3111323312113 | tail -1 | ./bin 3 3 0 0 > tmp-pic.tex && make dotex # R13-15
% r13-15-edge-subst-iterate4.pdf
% sequence data:
% stringsubst 3 1 1 1222131123221 2 2333212231332 3 3111323312113 | tail -1
%
%% for Gosper flowsnake, see https://oeis.org/A229214 (Arie Bos):
% stringsubst 2 1 1 1243116 2 1224532 3 3465332 4 3446154 5 5621554 6 5662316 | tail -1 | ./bin 6 3 1 1 > tmp-pic.tex && make dotex #
% stringsubst 1 _1_2_4_3_1_1_6 _ _ 1 1243116 2 1224532 3 3465332 4 3446154 5 5621554 6 5662316 | tail -1 | ./bin 6 3 1 1 > tmp-pic.tex && make dotex #

The maps of these L-systems by definition commute with the map $\sigma$ that
replaces $k$ by $k+1\bmod{n}$ where $n$ is the number of directions.
Dekking's \jjterm{folding morphisms} \cite[Section~2]{dekking-tiles-final}
commute with $\sigma\,\tau$ where $\tau$ is the map
that reverses the productions.
These curves can also be described by L-systems
with two non-constant maps for, say, \texttt{L} and \texttt{R} as follows.
The production of \texttt{L} starts with \texttt{L} and the letters \texttt{L} and \texttt{R}
appear alternatingly, with turns (\texttt{+} and \texttt{-}) between them.
The production of \texttt{R} is obtained from that of \texttt{L}
by reversing and swapping \texttt{L} with \texttt{R} and \texttt{+} with \texttt{-}.
For example, the maps \Lmap{L}{L+R-L-R} and \Lmap{R}{L+R+L-R}
correspond to the curve shown in figure \Ref{fig:alt-paperfold-carousel}.

The intersection of the sets of Dekking's curves and the curves with simple
L-systems consists of the curves with two-fold symmetry.

%%%%%%%%%%%%%%%%%%%%%%%%
\paragraph{Description by turns.}
Yet another way to describe curves with simple L-systems is
by an L-system for the succession of turns.
The maps for the turns $T$
(where $T$ $\in\{$ \texttt{+}, \texttt{-}, \texttt{0} $\}$
and \texttt{0} denotes a turn by $0\adeg$, a non-turn)
are the turns in the production of \texttt{F} (in the simple L-system)
followed by $T$.
Edges are drawn at the start and after each turn.
For example, the curve \CID{R4-1} in Figure~\ref{fig:r4-crab-dragon}
has the L-system \Lmap{F}{F+F0F-F},
the sequence of turns is \texttt{+0-},
and we get the maps
\Lmap{+}{+0-+},  \Lmap{-}{+0--}, and \Lmap{0}{+0-0}.
%
% stringsubst 1 +  + +0-+  0 +0-0  - +0-- | sed 's/\(.\)/F\1/g;' | tail -1 | ./bin 3 3 1 1 > tmp-pic.tex && make dotex # crab
%
For the R5-dragon with map \Lmap{F}{F+F+F-F-F},
we get the maps \Lmap{+}{++--+} and \Lmap{-}{++---}.
%
% terdragon: \Lmap{+}{+-+} and \Lmap{-}{+--};
% R5-dragon: \Lmap{+}{++--+} and \Lmap{-}{++---};
% stringsubst 6 +  + +-+  - +-- | sed 's/\(.\)/F\1/g;' | tail -1 | ./bin 3 2 0 0 0.15 > tmp-pic.tex && make dotex # terdragon
% stringsubst 4 + + ++--+ - ++--- | sed 's/\(.\)/F\1/g;' | tail -1 | ./bin 4 2 0 0 0.15 > tmp-pic.tex && make dotex # 5-curve
%
The choice of the axiom does not matter as the maps are identical except for the last turn.
This description is equivalent to the method of looking up the $n$th turn
($n\geq{}1$) from the lowest non-zero digit in the base-$R$ expansion of $n$
indicated in \cite[Section~1.31.5, pp.~95ff]{fxtbook}.
% in fxt-dir:  ee $(fh bit-dragon)

%%%%%%%%%%%%%%%%%%%%%%%%
\paragraph{L-systems for the Hilbert curve.}
The following descriptions of the Peano and Hilbert curves by L-systems are slightly more
complicated than for other curves shown.
The Peano curve can be described by the L-system with axiom
\texttt{L} and (non-constant) maps
\Lmap{L}{LtRtL-t-RtLtR+t+LtRtL},
\Lmap{R}{RtLtR+t+LtRtL-t-RtLtR}.
Here only the letter \texttt{t} corresponds to an edge,
both \texttt{L} and \texttt{R} are ignored in the drawing process.
For the Hilbert curve a possible L-system with axiom \texttt{L}
and (non-constant) maps
\Lmap{L}{+Rt-LtL-tR+},
\Lmap{R}{-Lt+RtR+tL-}
can be used (only \texttt{t} corresponds to an edge).

Again, the descriptions by strokes and turns is by no means the only way to specify these curves.
For the Hilbert curve the L-system with axiom \texttt{A}
and non-constant maps
\Lmap{A}{DrAuAlC},
\Lmap{B}{ClBdBrD},
\Lmap{C}{BdClCuA},
\Lmap{D}{AuDrDdB}
can be used.
The drawing is obtained by ignoring the capital letters and using
\texttt{r}, \texttt{d}, \texttt{l}, and \texttt{u}
as moves respectively
\eqq{right}, \eqq{down}, \eqq{left}, and \eqq{up},
see \cite[Figure~1.31-C, p.~85]{fxtbook}.
Bader \cite[Section~3.3, pp.~37ff]{bader-sfc}
gives a similar description for the Peano curve.

%%% Emacs:
%%% Local Variables:
%%% mode: latex
%%% MyRelDir: "."
%%% TeX-master: "arndt-curve-search.tex"
%%% dvi-file: "arndt-curve-search"
%%% makefile-dir: "./"
%%% frame-title-format: "CURVE-SEARCH (intro)"
%%% End:

%%%%%%%%%%%%%%%%%%%%%%%%%%%%%%%%%%%%%%%%%%%%%%%%%%%%%%%%%%%%
%%%%%%%%%%%%%%%%%%%%%%%%%%%%%%%%%%%%%%%%%%%%%%%%%%%%%%%%%%%%
\section{The search}%\label{sect:search}

%\xxx{Intro text}

%%%%%%%%%%%%%%%%%%%%%%%%%%%%%%%%%%%%%%%%%%%%%%%%%%%%%%%%%%%%
\subsection{Conditions for a curve to be self-avoiding and edge-covering}\label{sect:conditions}

Michel Dekking's 2012 paper \eqq{\NS{Paperfolding} morphisms, \NS{planefilling} curves, and fractal tiles}
\cite{dekking-tiles-final} has been crucial for this project.

Let $X$ be the production of \texttt{F} in a simple L-system of order $R$.
Let $C_n$ be the $n$th iterate of the curve corresponding to the L-system
and $C$ be the set of all iterates $C_n$ of the curve $n\geq{}0$.
We say $C$ is self-avoiding or edge-covering
if every curve in $C$ has the respective property.
%

%%%%%%%%%%%%%%%%%%%%%%%%%%%%%%
\subsubsection{Sufficient conditions}%\label{sect:}
%

%%%%%%%%%%%%%%%%%%%%%%%%%%
%
%% Render with thick lines: lnth *= 3.0; :
% stringsubst 1 _F+_F+_F+_F  _ _ F F+F-F-F+F+F+F-F+F-F-F-F+F   + + - - | tail -1 | ./bin 4 3 1 > tmp-pic.tex && make dotex # R13-3
% stringsubst 1 _F-_F-_F-_F  _ _ F F+F-F-F+F+F+F-F+F-F-F-F+F   + + - - | tail -1 | ./bin 4 3 1 > tmp-pic.tex && make dotex # R13-3
%
\begin{figure}[h!tbp]
\ifpdf
\begin{center}
\includegraphics*[width=48mm, viewport={40 390 470 720}]{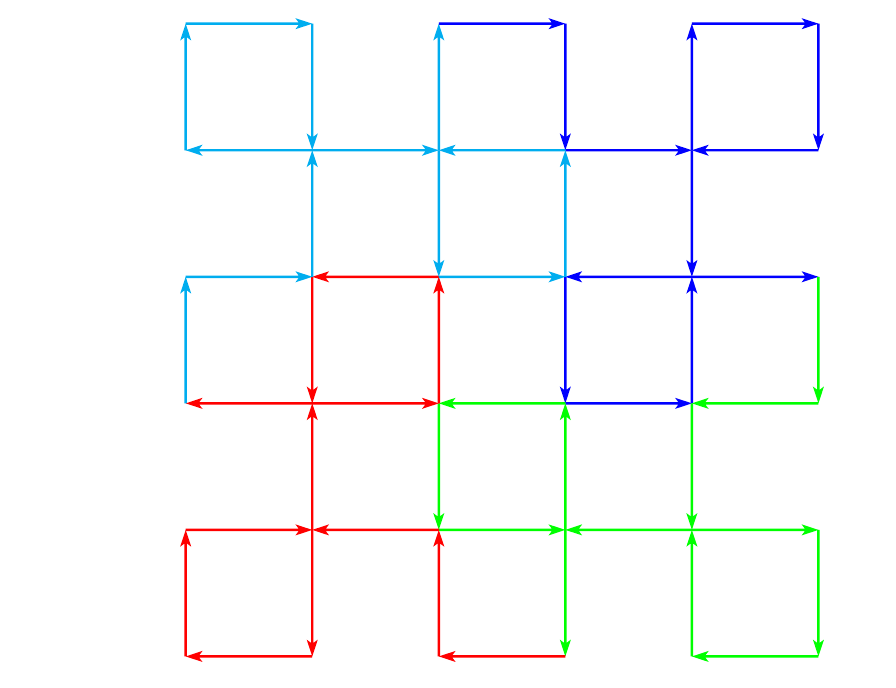}%
\includegraphics*[width=48mm, viewport={40 380 450 720}]{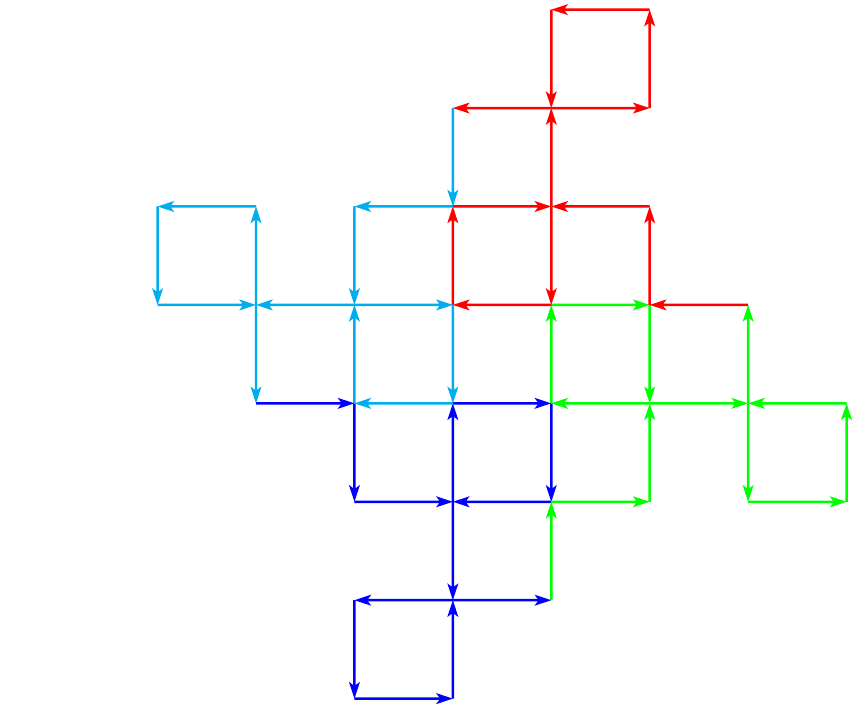}
\end{center}
\else
\verb+{see pdf for image}+
\fi
\caption{\label{fig:r13-q-3-tiles-it1}
%Tiles $\Tile{+1}$ (top) and $\Tile{-1}$ (bottom) for the curve of order 13 on the square grid
Tiles $\Tile{+1}$ (left) and $\Tile{-1}$ (right) for the curve of order 13 on the square grid
with L-system \Lmap{F}{F+F-F-F+F+F+F-F+F-F-F-F+F} (\CID{R13-3}).}
\end{figure}
%
%%%%%%%%%%%%%%%%%%%%%%%%%%

%%%%%%%%%%%%%%%%%%%%%%%%%%
%
%
%% Render with thick lines: lnth *= 3.0; :
% stringsubst 1 _F+_F+_F  _ _  F F+F0F0F-F-F+F0F+F+F-F0F-F   0 0 + + - - | tail -1 | ./bin 3 3 1 > tmp-pic.tex && make dotex # R13-15
% stringsubst 1 _F-_F-_F  _ _  F F+F0F0F-F-F+F0F+F+F-F0F-F   0 0 + + - - | tail -1 | ./bin 3 3 1 > tmp-pic.tex && make dotex # R13-15
%
\begin{figure}[h!tbp]
\ifpdf
\begin{center}
\includegraphics*[width=46mm, viewport={120 440 450 710}]{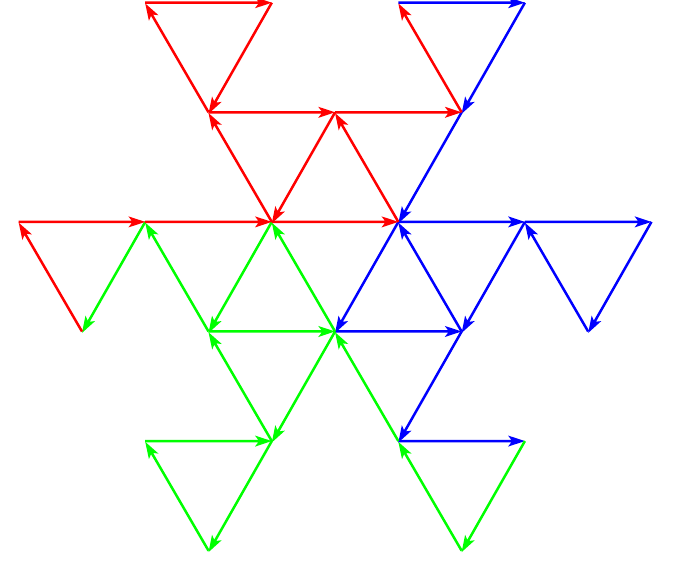}%
\includegraphics*[width=46mm, viewport={120 420 450 710}]{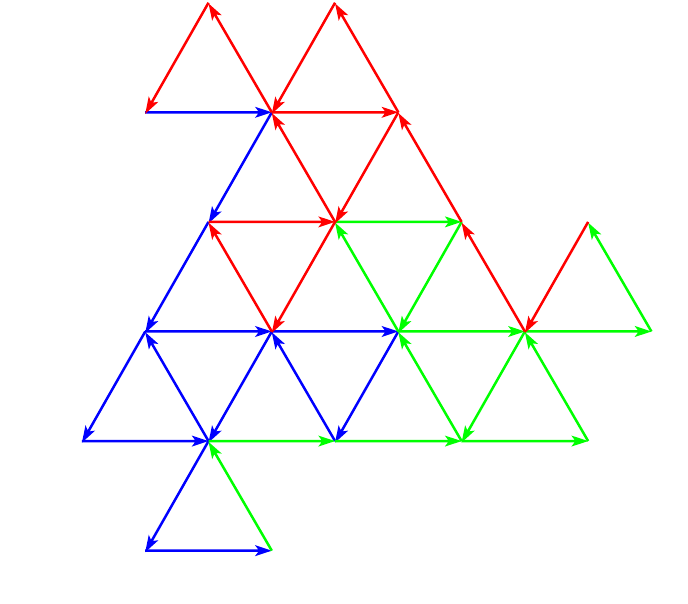}
\end{center}
\else
\verb+{see pdf for image}+
\fi
\caption{\label{fig:r13-t-15-tiles-it1}
%Tiles $\Tile{+1}$ (top) and $\Tile{-1}$ (bottom) for the curve of order 13 on the triangular grid
Tiles $\Tile{+1}$ (left) and $\Tile{-1}$ (right) for the curve of order 13 on the triangular grid
with L-system \Lmap{F}{F+F0F0F-F-F+F0F+F+F-F0F-F} (\CID{R13-15}).}
\end{figure}
%
%%%%%%%%%%%%%%%%%%%%%%%%%%

We need the concept of a \jjterm{tile} for the following facts.
Here we define them only for the square and the triangular grid,
the tiles for the tri-hexagonal grid are described in section \Ref{sect:tri-hex}.

For the square grid, let $\Tile{+1}$ be the (closed) curve corresponding to
the first iterate of the map of the L-system with axiom \texttt{F+F+F+F}
and $\Tile{-1}$ for the axiom \texttt{F-F-F-F}.
Here the turns are by $\phi=90\adeg$, see Figure~\ref{fig:r13-q-3-tiles-it1} for an example.
For the triangular grid, the respective axioms are \texttt{F+F+F} and \texttt{F-F-F},
and the turns are by $\phi=120\adeg$, see Figure~\ref{fig:r13-t-15-tiles-it1}.

The tiles of edge-covering curves do indeed tile the grid: infinitely many disjoint translations
of them do cover all edges of the grid, see section \ref{sect:curve-self-sim}.

\begin{theorem}[Tiles-SA]\label{thm:tiles-sa}
$C$ is self-avoiding if and only if both tiles $\Tile{+1}$ and $\Tile{-1}$ are self-avoiding.
\end{theorem}
Dekking's proof for the curves described by folding morphisms
 \cite[Theorem~1, p.~24]{dekking-tiles-final}
can be adapted for the curves with simple L-systems
[Dekking, personal communication, March 2016].
Note that for Dekking's curves on the square grid,
it suffices that one of the two tiles is self-avoiding

We call a tile edge-covering if all edges
in its interior are traversed once.

\begin{theorem}[Tiles-Fill]\label{thm:tiles-fill}
$C$ is edge-covering if and only if
both tiles $\Tile{+1}$ and $\Tile{-1}$ are edge-covering.
\end{theorem}
The proof for Theorem~2 in \cite[p.~27]{dekking-tiles-final}
can be modified to show this
[Dekking, personal communication, March 2016].

What is called tile (for a self-avoiding, edge-covering curve) here
is called a (maximally simple) $\Theta$-loop by
Dekking \cite[Definition~3, p.~24]{dekking-tiles-final}.

%%%%%%%%%%%%%%%%%%%%%%%%%%%%%%
\subsubsection{Necessary conditions}%\label{sect:}
%

%%%%%%%%%%
The motif must obviously be self-avoiding:
\begin{theorem}[Obv]\label{thm:obv}
For $C$ to be self-avoiding, $C_1$ must be self-avoiding.
\end{theorem}

%%%%%%%%%%
\begin{theorem}[Turn]\label{thm:turn}
For $C$ to be self-avoiding and edge-covering,
the net rotation of the curve $C_1$ must be zero.
\end{theorem}
That is, the number of \texttt{+} and \texttt{-} in $X$ must be equal.
If this condition does not hold,
the tiles cannot be self-avoiding.

%%%%%%%%%%
\begin{theorem}[Dist]\label{thm:dist}
For $C$ to be self-avoiding and edge-covering,
the squared distance between the start and the endpoint
of the curve must be equal to $R$.
%the number of \texttt{F} in $X$.
\end{theorem}
If this condition does not hold,
the tiles cannot be both edge-covering.

Also, if the dimension of the limiting
curve is $\log{R}/\log{d} = 2\, \log{R}/\log{d^2}$
where $d$ is the distance (see \cite[Section~1.1.4, pp.~28ff]{peitgen-science-fract-img}),
we need $R=d^2$ for dimension $2$.

For the square grid, the possible orders are the
numbers of the form $x^2+y^2$, sequence \jjseqref{A001481} in \cite{OEIS},
otherwise numbers of the form $x^2+x\,y+y^2$
(equivalently, numbers of the form $3\,x^2+y^2$), sequence \jjseqref{A003136} in \cite{OEIS}.

The curves on the square grid must turn at every point,
otherwise two dead ends would be created.
As we do not allow turns as first or last letters in $X$,
we will find only curves of odd orders, sequence \jjseqref{A057653} in \cite{OEIS}.
%

%For even orders, pairs of maps of the form
%\Lmap{L}{f(L)} and \Lmap{R}{g(R)}
%where the word \texttt{g(R)} is obtained from \texttt{f(L)}
%by reversing the word and swapping
%all \texttt{+} and \texttt{-}
%and all \texttt{L} and \texttt{R}
%correspond to curves.
%%
%For example, the maps
%\Lmap{L}{L+R} and \Lmap{R}{L-R}
%(that is, \texttt{f(L)}=\texttt{L+R} and \texttt{g(R)}=\texttt{L-R})
%for the Heighway-Harter dragon shown in Figure~\ref{fig:hh-dragon}
%satisfy the condition.
%%
%Allowing such rules for odd orders
%would have given additional curves.
% whose tiles are less symmetric.

%%%%%%%%%%%%
\subsection{The shape of a curve}

%%%%%%%%%%%%%%%%%%%%%%%%%%
%
% Render with thick lines: lnth *= 3.0; :
% stringsubst 1 F  F F0F+F+F-F-F0F + + - - 0 0  | tail -1 | ./bin 3 2 1 0 0.10 > tmp-pic.tex && make dotex # R7-2
% stringsubst 1 F  F F+F-F-F+F+F-F + + - - 0 0  | tail -1 | ./bin 3 2 1 0 0.10 > tmp-pic.tex && make dotex # R7-5
%
\begin{figure}[h!tbp]
\ifpdf
\begin{center}
\includegraphics*[height=30mm, viewport={170 560 450 710}]{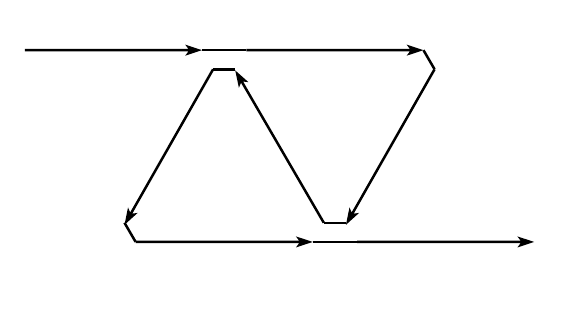}%
\includegraphics*[height=30mm, viewport={170 560 450 710}]{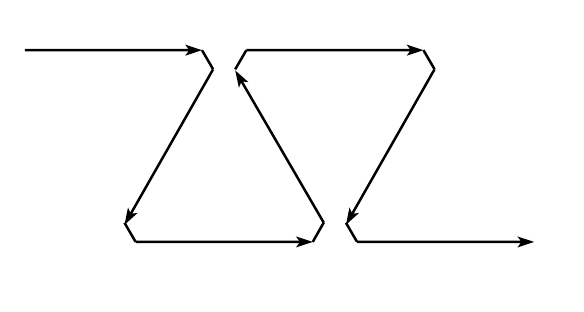}
\end{center}
\else
\verb+{see pdf for image}+
\fi
\caption{\label{fig:shape-7}
Two different curves of order 7 with same shape (\CID{R7-2} and \CID{R7-5}).}
\end{figure}
%
%%%%%%%%%%%%%%%%%%%%%%%%%%

%
The \jjterm{shape} of a curve is the set of grid edges traversed in the first iterate.
%% corrected 2018-July-03
%The \jjterm{shape} of a curve is the set of grid points traversed in the first iterate.
%
Different curves can have the same shape, as shown in Figure~\ref{fig:shape-7}.
The L-systems are respectively
\Lmap{F}{F0F+F+F-F-F0F} and
\Lmap{F}{F+F-F-F+F+F-F}.
Here the letter \texttt{0} (interpreted as \eqq{no turn}) is used
to keep the maps for one order of fixed length.
% only in the triangular grid
These two curves were already given
in \cite[Figures~1.31-P and 1.31-Q, pp.~97-98]{fxtbook}.
%and also in \cite[p.~91ff]{ventrella}.

We consider two curves to be of the same shape
whenever any transformation of the symmetry
of the underlying grid (rotations and flips) maps one shape into the other.
If two curves have the same shape, we call them \jjterm{similar}.

Rendering only one curve for each shape makes
the results (visually) much more manageable,
as for higher orders more and more curves
tend to have the same shape.
%

%%%%%%%%%%%%%%%%%%%%%%%%%%%%%%
\subsection{Structure of the program for searching}%\label{sect:}

The program consists of the following parts:
generation of the L-systems,
testing of the corresponding curves,
detection of similarity to shapes seen so far,
and detection of symmetries.

The generation of the L-systems is equivalent to counting in base 2 or 3:
For the square grid the allowed turns are by $+90\adeg$ and $-90\adeg$.
For the triangular grid turns by $+120\adeg$, $-120\adeg$, and $0\adeg$ are possible.
For the tri-hexagonal, two values again suffice, see section \ref{sect:tri-hex}.

Testing for self-avoidance and the edge-covering property uses the conditions of
section \ref{sect:conditions} in some order where the criteria cheaply testable
are checked before the more costly ones.
Condition \ref{thm:turn} [Turn] is tested before the curve
is even computed (as a sequence of points),
then \ref{thm:dist} [Dist], followed by
\ref{thm:obv} [Obv],
\ref{thm:tiles-sa} [Tiles-SA], and
\ref{thm:tiles-fill} [Tiles-Fill].

Every curve found is assigned a positive ID.
The triple $($grid, order, ID$)$ identifies a curve.

The shape of the curve is computed as a sorted sequence of points without duplicates.
%also the shapes of certain transformations of the grid,
%such as reflections by the $x$ and $y$ axis.
%
For each curve found, the shape is compared to the shapes found so far.  If a
match is detected, the new curve is marked a duplicate shape.
% otherwise a new shape has been found.
%
If a curve is a duplicate, the lowest ID of a curve
of the same shape and the transformations to match the shapes
are documented.

The symmetries (if any) of each curve are computed and documented.

The generation of all L-systems is very simple but certainly not very elegant.
Still, after careful optimization, the tests and the generation of the L-systems
each use about half of the CPU time, and approximately one million curves
per second are checked.
% // job started 2013-December-24:
% ## from htop (2013-December-28 (18:26)):
% ## PID USER     PRI  NI  VIRT   RES   SHR S CPU% MEM%   TIME+  Command
% ## 27772 jj      20   0 1247M 1236M   956 R 108.  7.7 101h41:34 ./abin 31 3
% 487383361197 curves were tested.  cct = 116861  sct = 1848
% ./abin 31 3  489524.19s user 15.76s system 99% cpu 136:01:16.85 total
% 487383361197/489524.19 == 995,626 per second
% 487383361197/116861 == 4,170,624 candidates per curve
% 487383361197/1848 == 263,735,585 candidates per shape
%
% start:  2015-January-12 (13:18)  tmp-r53-q.txt
% end:  53M -rw-r--r--  1 arndt users  53M Jan 19 00:05 tmp-r53-q.txt
%  ==> about 6.5 days = 156 hours
% 793747044784 curves were tested.  cct = 200869  sct = 2467
%  793747044784/(156*3600.) == 1,413,367 per second
%
Due to the combinatorial explosion, searches for large orders
would take very long.
The largest orders searched and their approximate duration were
$R=53$ for the square grid (7 days),
$R=31$ for the triangular grid (6 days), and
$R=61$ for the tri-hexagonal grid (76 days).
%$53$ for the square grid (taking 156 hours),
%$31$ for the triangular grid (taking 136 hours), and
%$61$ for the tri-hexagonal grid (taking 1815 hours).
% ? 156/24.
% 6.50000000000000
% ? 136/24.
% 5.66666666666667
% ? 1815/24.
% 75.6250000000000
%
%
Parallelization of the existing program is not planned,
instead a completely different method for searching will be used.

%%%%%%%%%%%%%%%%%%%%%%%%%%%%%%
\subsection{Format of the files specifying the curves}%\label{sect:}

%%%%%%%%%%%%%%%%%%%%%%%%%%
\begin{figure}[h!tbp]
%%  ./scbin 17 4 0 c4/search-r17-q-curves.txt | ./shorten-output.pl
{\small
\begin{verbatim}
F F+F+F-F+F-F-F-F+F-F+F+F+F-F+F-F-F  R17-1  # # symm-dr
F F+F+F-F-F+F+F+F-F+F-F-F-F+F+F-F-F  R17-2  # # symm-dr
F F+F+F-F-F+F+F+F-F-F+F+F-F-F-F+F-F  R17-3  #
F F+F+F-F-F-F+F+F+F-F+F+F-F-F-F+F-F  R17-4  # # symm-r ## same = 1 P R
F F+F-F+F+F+F-F-F+F+F-F-F-F+F+F-F-F  R17-5  # ## same = 3 R X
F F+F-F+F+F+F-F-F+F-F+F+F-F-F-F+F-F  R17-6  # # symm-dr
F F+F-F+F+F+F-F-F+F-F-F-F+F+F+F-F-F  R17-7  # # symm-r ## same = 1 P R
F F+F-F+F+F+F-F-F-F+F+F+F-F-F-F+F-F  R17-8  # # symm-dr ## same = 1 P R
F F+F-F+F+F-F+F+F+F-F-F-F+F-F-F+F-F  R17-9  # # symm-dr ## same = 1 P R
F F+F-F+F+F-F+F+F-F-F-F+F+F-F-F-F+F  R17-10  #
F F+F-F+F+F-F+F-F+F+F-F-F-F+F-F-F+F  R17-11  #
F F+F-F-F+F-F-F-F+F+F-F+F-F+F+F-F+F  R17-12  # ## same = 11 Z T
F F+F-F-F-F+F+F-F-F-F+F+F-F+F+F-F+F  R17-13  # ## same = 10 Z T
\end{verbatim}
}
%% cf.   ./scbin 9 4 0 c4/search-r09-q-curves.txt | ./shorten-output.pl
%% F F+F-F-F-F+F+F+F-F  R9-1  # # dragon # symm-dmrqz
%
\caption{\label{fig:output-17-q}
Descriptions of the curves of order 17 on the square grid.}
\end{figure}
%
%%%%%%%%%%%%%%%%%%%%%%%%%%

Figure~\ref{fig:output-17-q} shows the output format for
the curves of order 17 on the square grid.

The descriptions of the curves found
are available for download
at \url{http://jjj.de/3frac/}
in the file \texttt{short-lsys.tar.xz} (of size 4.3 MB)
which unpacks into a directory named
\texttt{short-lsys/} of size 136 MB.
Therein the subdirectories
%\texttt{c3/}, \texttt{c4/}, and \texttt{c6/}
%% after correction 2018-July-03:
\texttt{x3/}, \texttt{x4/}, and \texttt{x6/}
contain respectively the descriptions of
the curves in the
triangular, square, and tri-hexagonal grids.
The file for an order, say $R=17$,
is named \texttt{lsys-r17.txt}.
%% see short-lsys.tar.xz  (4.3M, unpacked 136M)
%% Created with command
%%   zsh shorten-all.zsh

%%%%%%%%%%%%
\subsubsection{L-system and ID}\label{sect:lsys-and-id}
Each line starts with the non-constant map for the L-system:
\texttt{F X} where \texttt{X} is the production of \texttt{F}.
No listed L-system has a map with initial turn \texttt{-}:
these can be omitted because such a map is either
the reversal of some map with last turn \texttt{-}
or can be obtained from a map with \texttt{+}
both as initial and final turn by swapping \texttt{-} and \texttt{+}.

The following string (like \CID{R17-6}) is of the form \texttt{R<order>-<ID>}.
The order of the curves is constant within a file and \texttt{ID} identifies a
curve among all curves of the same order on the same grid.

%%%%%%%%%%%%
\subsubsection{Information about symmetries of the curve}

%%%%%%%%%%%%%%%%%%%%%%%%%%
% stringsubst 1 F  F F0F+F-F-F0F+F+F-F   0 0 + + - - | tail -1 | ./bin 3 2 0 0 0.15 > tmp-pic.tex && make dotex # R9-2 # symm-m
%
%% with lnth *= 2.0;  // thicker lines
%  stringsubst 1 F  F F0F+F-F-F0F+F+F-F   0 0 + + - - | tail -1 | ./bin 3 2 0 0 0.15 > tmp-pic.tex && make dotex # R9-2 # symm-m
%  stringsubst 1 F  F F0F0F+F0F-F-F+F-F0F0F+F+F-F+F-F   0 0 + + - - | tail -1 | ./bin 3 2 0 0 0.15 > tmp-pic.tex && make dotex # R16-1 # symm-m
%  stringsubst 1 F  F F0F0F+F+F-F-F0F-F0F-F+F+F0F+F-F   0 0 + + - - | tail -1 | ./bin 3 2 0 0 0.15 > tmp-pic.tex && make dotex # R16-4 # symm-m
%
\begin{figure}[h!tbp]
\ifpdf
\begin{center}
%\fbox{\includegraphics*[width=100mm, viewport={0 390 600 740}]{r9-t-2-manta.pdf}}
%
{\includegraphics*[width=30mm, viewport={160 500 450 700}]{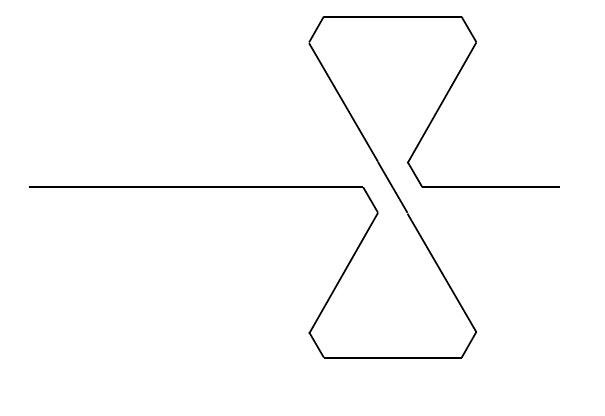}}%
\qquad% layout
\qquad% layout
{\includegraphics*[width=30mm, viewport={60 410 450 710}]{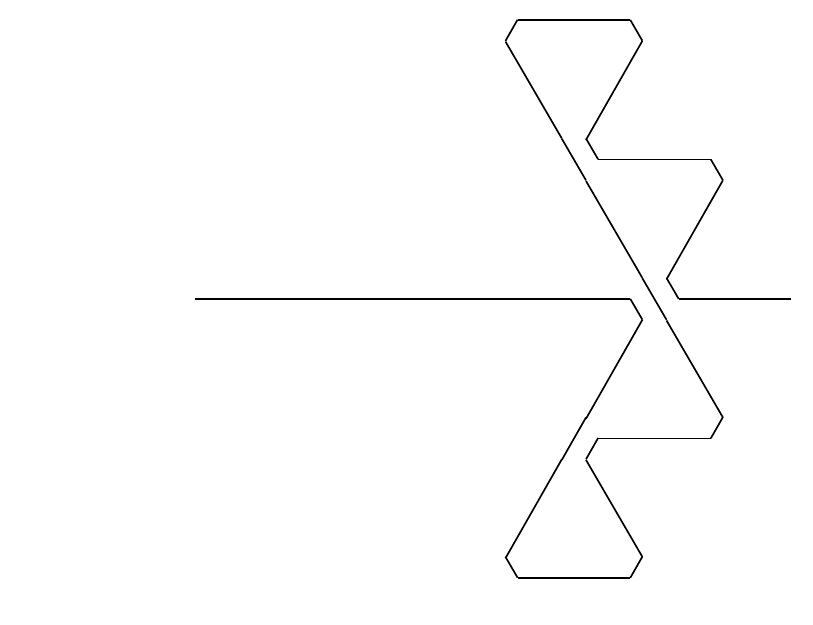}}%
\qquad% layout
{\includegraphics*[width=30mm, viewport={160 520 450 705}]{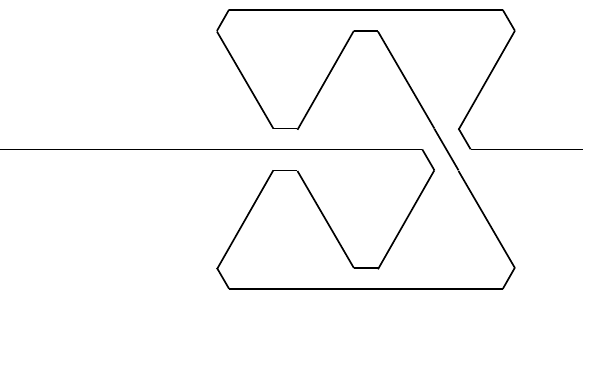}}
\end{center}
\else
\verb+{see pdf for image}+
\fi
\caption{\label{fig:r9-t-2-manta}
Motifs of curves on the triangular grid with shape invariant to reflection of the $y$-axis
(\CID{R9-2}, \CID{R16-1}, and \CID{R16-4}).}
\end{figure}
%
%%%%%%%%%%%%%%%%%%%%%%%%%%
%

%%%%%%%%%%%%%%%%%%%%%%%%%%
%
%% with lnth *= 2.0;  // thicker lines
% stringsubst 1 F  F F+F-F+F+F-F-F+F-F-F-F+F+F-F+F-F-F+F-F-F-F+F+F+F-F+F+F-F-F+F+F-F+F+F+F-F-F+F-F+F-F+F-F-F-F+F+F+F-F   + + - - | tail -1 | ./bin 4 2 0 0 0.15 > tmp-pic.tex && make dotex # R49-14164 # symm-m
%
% stringsubst 1 F  F F+F-F+F+F-F-F-F+F+F-F-F-F+F+F-F-F+F-F-F-F+F+F+F-F+F+F-F-F+F+F+F-F-F+F+F+F-F-F+F-F+F-F-F-F+F+F+F-F   + + - - | tail -1 | ./bin 4 2 0 0 0.15 > tmp-pic.tex && make dotex # R49-15059 # symm-m
%
\begin{figure}[h!tbp]
\ifpdf
\begin{center}
{\includegraphics*[width=32mm, viewport={90 280 500 730}]{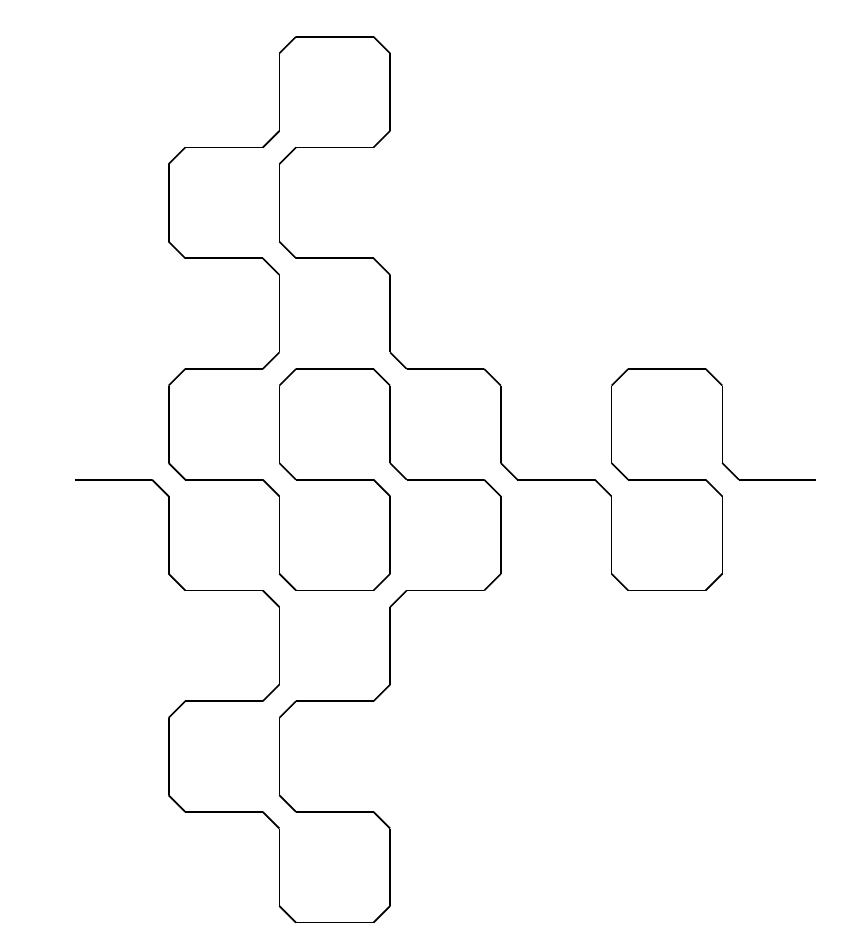}}%
\qquad% layout
\qquad% layout
{\includegraphics*[width=40mm, viewport={90 370 500 730}]{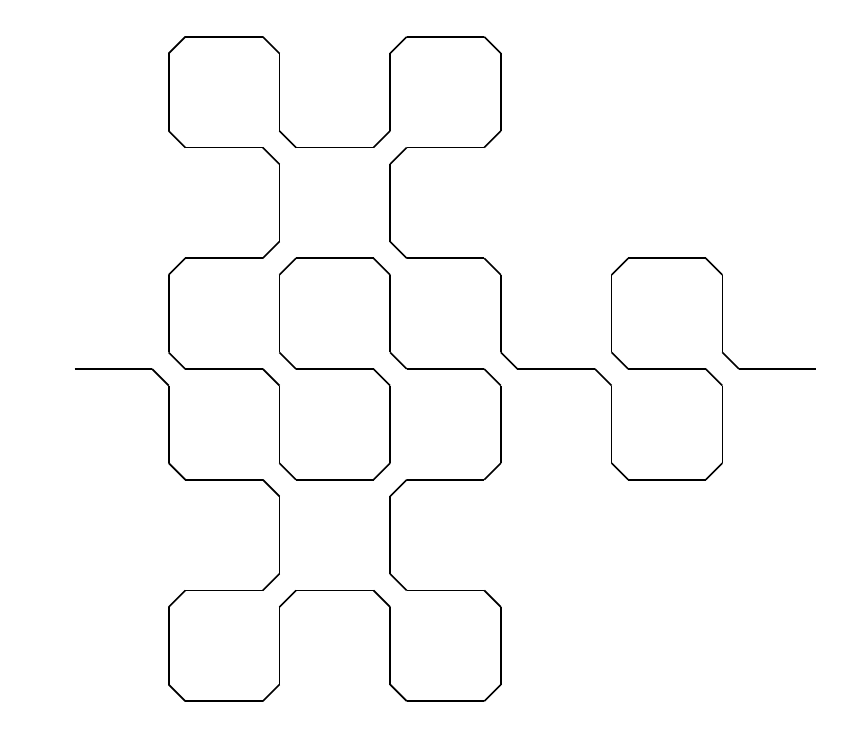}}
\end{center}
\else
\verb+{see pdf for image}+
\fi
\caption{\label{fig:r49-q-manta}
Motifs of two curves of order 49 on the square grid
invariant to reflection of the $y$-axis (\CID{R49-14164} and \CID{R49-15059}).}
%% and not of square shape at the same time
\end{figure}
%
%%%%%%%%%%%%%%%%%%%%%%%%%%
%

%
Information about symmetry and similarity of the shape is given
after the character~\texttt{\#}.
%
%The word \texttt{dragon} marks curves where \texttt{X}
%is invariant under the map reversing the word and swapping \texttt{-} and \texttt{+}.
%These curves do have a 2-fold rotational symmetry about the midpoint.
%%
%Curves that have this point symmetry but not the mentioned property
%for \texttt{X} are marked \texttt{dragon-shape}.

If a curve has any non-trivial symmetries then these are
specified as letters following \texttt{symm-}.
%
%\xxx{Define rev/swap map?}
A \texttt{d} indicates that the word \texttt{X}
is invariant under the map reversing it and swapping \texttt{-} and \texttt{+}.
These curves do have a 2-fold (rotational) symmetry.
An \texttt{r} indicates 2-fold symmetry,
so whenever there is a \texttt{d} there is an \texttt{r}.
%% grep symm-d c3/*txt | grep -v symm-dr

%
In the following it is assumed that the first edge
goes in positive $x$-direction.
The letter \texttt{m} stands for invariance under reflection of the $y$-axis.
Curves with this symmetry are shown
in Figures~\ref{fig:r9-t-2-manta} (orders 9 and 16) and \ref{fig:r49-q-manta} (order 49).
Curves with just this symmetry on the triangular grid exist for orders $R=n^2$
where $n\pmod{3}\not\equiv{}2$,
and on the square grid for $R=n^2$ with $n$ odd.
%\xxx{More cond. on square grid?}
% triangular grid: ...because we must have (n^2-n) % 6 == 0
% ? vector(10,n,[n^2,(n^2-n)%6])
% [[1, 0], [4, 2], [9, 0], [16, 0], [25, 2], [36, 0], [49, 0], [64, 2], [81, 0], [100, 0]]

The specifications \texttt{symm-mrqz} (or \texttt{symm-dmrqz})
are given for curves whose shapes are perfect squares,
these exists only on the square grid and for orders that are squares.

%{\scriptsize
%\begin{verbatim}
% file src/shape.h:
%
%        A_ = new iVec[np+nh];  // original Shape
%        for (ulong j=0; j<np; ++j)  A_[j] = P[j];
%
%        // x_mirror() etc. in iVec.h:  q_mirror() := "flip q-axis"
%
%//        Amx_ = new iVec[np+nh];  // x-mirrored reversed Shape
%//        for (ulong j=0; j<np; ++j)  Amx_[j] = (P[np-1] - P[j]).x_mirror();
%
%        Amy_ = new iVec[np+nh];  // y-mirrored Shape
%        for (ulong j=0; j<np; ++j)  Amy_[j] = P[j].y_mirror();
%
%        Ar_ = new iVec[np+nh];  // reversed Shape (curve with reversed directions)
%        for (ulong j=0; j<np; ++j)  Ar_[j] = P[np-1] - P[j];
%
%        Arm_ = new iVec[np+nh];  // y-mirrored reversed Shape
%        for (ulong j=0; j<np; ++j)  Arm_[j] = Ar_[j].y_mirror();
%
%\end{verbatim}
%}
%
%{\scriptsize
%\begin{verbatim}
% file src/shape.h:
%
%    cout << "symm-";  // following comments copied from above
%    // L-system:
%    if ( tstb(symm_, 0) )  cout << "d";  // d (anti-symmetric L-sys, proper Dragon)
%//    if ( tstb(symm_, 1) )  cout << "?";  // ? (symmetric L-sys, never?)
%    // overall shape:
%//    if ( tstb(symm_, 12) )  cout << "x";  // x (self == reversed x-mirrored)
%    if ( tstb(symm_, 13) )  cout << "m";  // m (self == y-mirrored, Manta)
%    if ( tstb(symm_, 14) )  cout << "r";  // r (self == Reversed, dragon-shape == point-symmetric shape)
%    if ( tstb(symm_, 15) )  cout << "q";  // q (y-mirrored == reversed, e.g. for sQuares)
%    if ( tstb(symm_, 16) )  cout << "z";  // z (self == mirrored reversed)
%
%\end{verbatim}
%}

%%%%%%%%%%%%
\subsubsection{Information about similarities to shapes of other curves}

If a curve has a shape similar (via any transformation of the grid)
to any curve found earlier then this is documented as
%\texttt{\#\#}
\texttt{same = ID} where \texttt{ID} refers to the earlier curve.
Among multiple possible choices of similar curves the one
with the lowest ID is taken.

The letters following specify which transformations
map the shape of the new curve into the shape of the earlier curve.
Here
\texttt{P} is put if the two shapes are identical (without any transformation),
\texttt{M} for reflection of the $y$-axis,
\texttt{R} for reversing directions of all edges,
\texttt{Z} for reversing directions and reflection of the $y$-axis,
\texttt{T} if the turns of the curves after reversing one L-system are identical,
and
\texttt{X} if the turns of the curves after reversing one L-system and swapping signs are identical.
% conjugates under Dekking's operation \pi( \tau( . ))
%
For shapes with non-trivial symmetries more than one letter may be given.

%{\small
%\begin{verbatim}
%  file src/shape.h:
%
%    if ( tstb(qm, 10) )  cout << " P";  // (shape proper)
%//    if ( tstb(qm, 11) )  cout << " ?";  // (shape x-mirrored)
%    if ( tstb(qm, 12) )  cout << " M";  // (shape y-mirrored)
%    if ( tstb(qm, 15) )  cout << " R";  // (shape reversed)
%    if ( tstb(qm, 17) )  cout << " Z";  // (shape mirrored reversed)
%
%    if ( tstb(qm, 0) )  cout << " T";  // (turns reversed)
%    if ( tstb(qm, 1) )  cout << " X";  // (turns reversed negated)
%
%\end{verbatim}
%}

%%%%%%%%%%%%%%%%%%%%%%%%%%%%%%
\subsection{Numbers of shapes found}\label{sect:stats}

%\xxx{Layout}
We now give the numbers $n$ of shapes found for the orders $R$
(up to the search limits) where any curve exists as lists of entries $R$\,:\,$n$.
\\
Triangular grid, sequence \jjseqref{A234434} in \cite{OEIS}:\\
%{3:1, 4:1, 7:3, 9:5, 12:10, 13:15, 16:17, 19:71, 21:212, 25:184, 27:543, 28:842, 31:1848}.
%% corrected 2018-July-03
{3:1, 4:1, 7:3, 9:5, 12:10, 13:15, 16:17, 19:71, 21:213, 25:184, 27:549, 28:845, 31:1850}.
\\
Square grid, sequence \jjseqref{A265685} in \cite{OEIS}:\\
{5:1, 9:1, 13:4, 17:6, 25:33, 29:39, 37:164, 41:335, 49:603, 53:2467}.
\\
Tri-hexagonal grid, sequence \jjseqref{A265686} in \cite{OEIS}:\\
{7:1, 13:3, 19:7, 25:10, 31:63, 37:157, 43:456, 49:1830, 61:8538}.

The number of curves is much greater than the number of shapes.
For example, for order $R=53$ on the square grid
there are 2467 shapes and 401738 curves,
so about 162 curves share one shape on average.
The corresponding file lists the 200869 curves where the map for \texttt{F} does not start
with \text{F-F}.
These are counted by adding the number of entries ending in
\text{F-F} (131111) to the number (69758) of entries beginning and ending in
\text{F+F} (which after swapping signs give yet unseen L-systems).
%% wc -l c4/search-r53-q-curves.txt # 200869
%% grep -F 'F-F  ' c4/search-r53-q-curves.txt | wc -l # 131111
%% grep -F 'F+F  ' c4/search-r53-q-curves.txt | grep -F ' F+F' | wc -l # 69758

%For the triangular and the square grid the search\xxx{Could drop for layout.}
%for curves whose L-system is invariant under the map
%reversing the production of \texttt{F} and swapping
%signs (these curves have a 2-fold rotational symmetry)
%was done for a few more orders:
%%
%\\
%Triangular grid:
%{36:226, 37:200, 39:559, 43:547, 48:1271, 49:1903}.
%\\
%Square grid:
%{61:183, 65:523, 73:534, 81:634}.

% grep  'same' c?/search-r*drag*.txt | wc -l
% 156887
% grep -v 'same' c?/search-r*drag*.txt | wc -l
% 6580

%%% Emacs:
%%% Local Variables:
%%% mode: latex
%%% MyRelDir: "."
%%% TeX-master: "arndt-curve-search.tex"
%%% dvi-file: "arndt-curve-search"
%%% makefile-dir: "./"
%%% frame-title-format: "CURVE-SEARCH (search)"
%%% End:

\FloatBarrier

%%%%%%%%%%%%%%%%%%%%%%%%%%%%%%%%%%%%%%%%%%%%%%%%%%%%%%%%%%%%
%%%%%%%%%%%%%%%%%%%%%%%%%%%%%%%%%%%%%%%%%%%%%%%%%%%%%%%%%%%%
\section{Properties of curves and tiles}%\label{sect:prop}

%%%%%%%%%%%%%%%%%%%%%%%%%%%%%%%
\subsection{Self-similarity, symmetries, and tiling property}\label{sect:curve-self-sim}

%%%%%%%%%%%%%%%%%%%%%%%%%%
% stringsubst 3 F_+F_0F_0F_-F_-F_+F_0F_+F_+F_-F_0F_-F_  _ _  F F+F0F0F-F-F+F0F+F+F-F0F-F   0 0 + + - - | tail -1 | ./bin 3 3 0 > tmp-pic.tex && make dotex # R13-15
% stringsubst 3 RF_+F_0F_0F_-F_-F_+F_0F_+F_+F_-F_0F_-F_ R R _ _ F F+F0F0F-F-F+F0F+F+F-F0F-F 0 0 + + - - | tail -1 | sed 's/+/++/g; s/-/--/g; s/R/+/;' | ./bin 6 3 0 > tmp-pic.tex && make dotex # R13-15
\begin{figure}[h!tbp]
\ifpdf
\begin{center}
{\includegraphics*[width=110mm, viewport={60 270 490 750}]{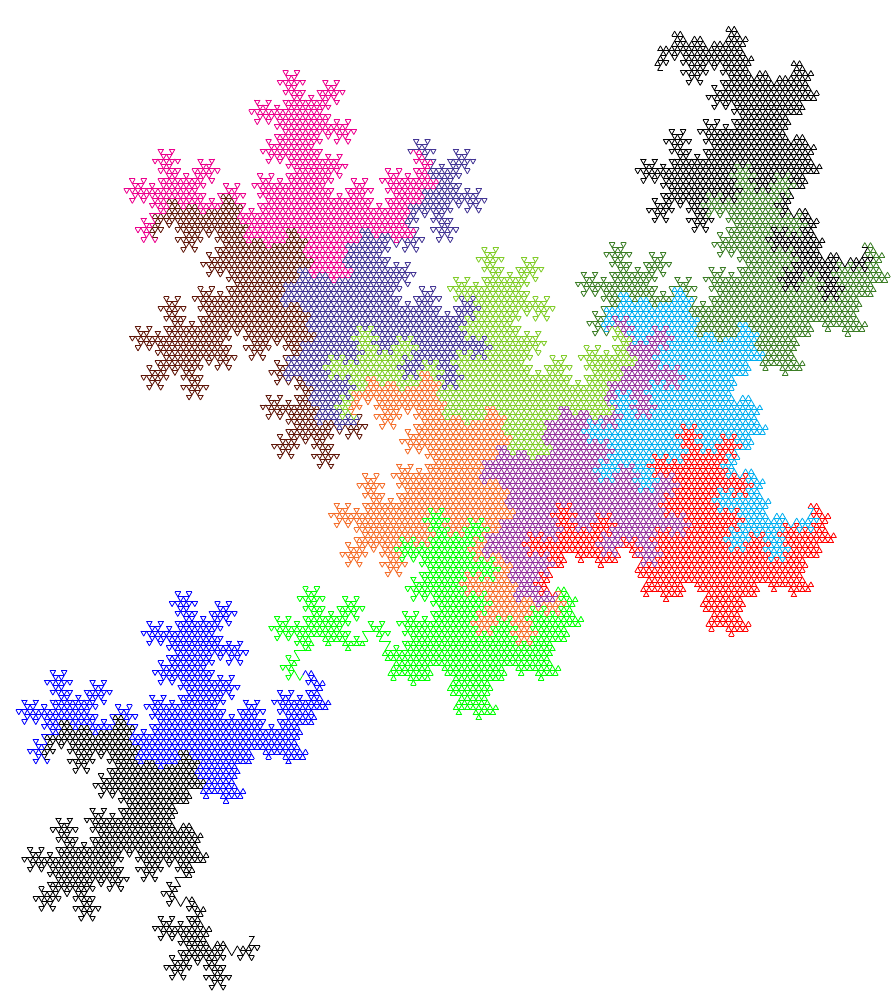}}
%{\includegraphics*[width=130mm, viewport={60 430 490 740}]{r13-t-15-decomp-rot.pdf}}
\end{center}
\else
\verb+{see pdf for image}+
\fi
\caption{\label{fig:r13-t-15-decomp}
Self-similarity of the order-13 curve \CID{R13-15} on the triangular grid.}
\end{figure}
%
%%%%%%%%%%%%%%%%%%%%%%%%%%

%%%%%%%%%%%%%%%%%%%%%%%%%%
%% with lnth *= 2.0;  // thicker lines
% stringsubst 3 _F+_F+_F  _ _  F F+F0F0F-F-F+F0F+F+F-F0F-F   0 0 + + - - | tail -1 | ./bin 3 2 0 0 0.33 > tmp-pic.tex && make dotex # R13-15
% stringsubst 3 _F-_F-_F  _ _  F F+F0F0F-F-F+F0F+F+F-F0F-F   0 0 + + - - | tail -1 | ./bin 3 2 0 0 0.33 > tmp-pic.tex && make dotex # R13-15
%
\begin{figure}[h!tbp]
\ifpdf
\begin{center}
\includegraphics*[width=65mm, viewport={70 290 500 750}]{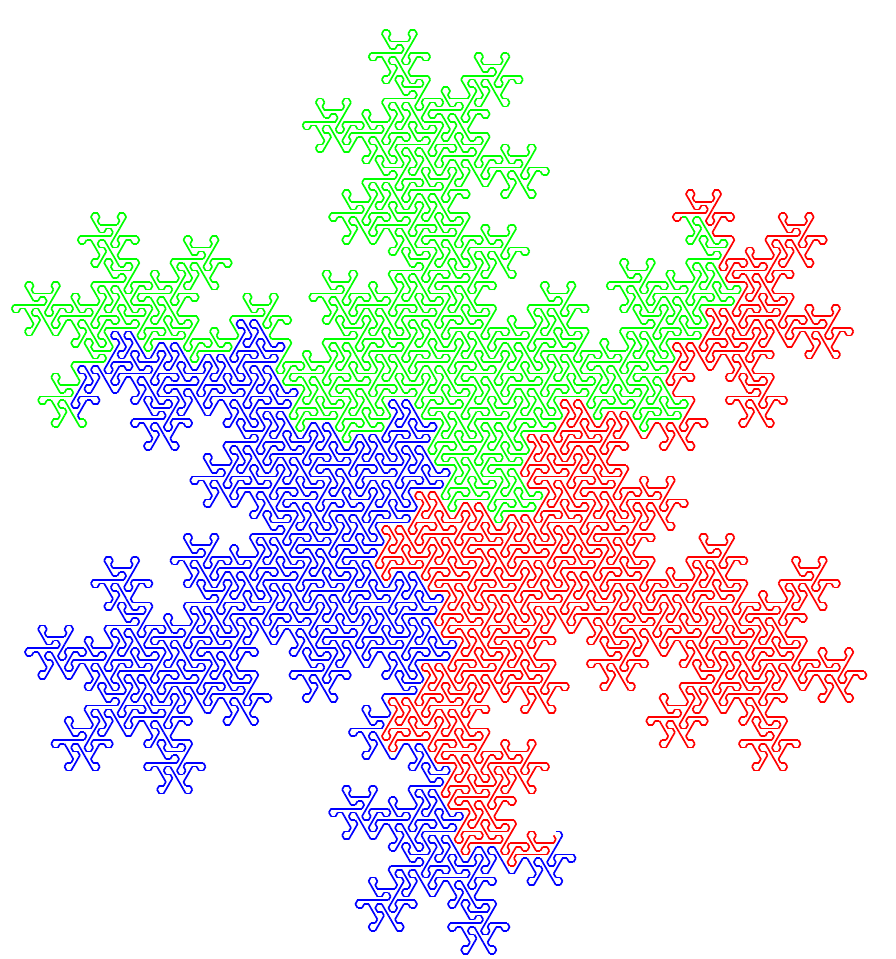}%
\includegraphics*[width=65mm, viewport={70 320 500 750}]{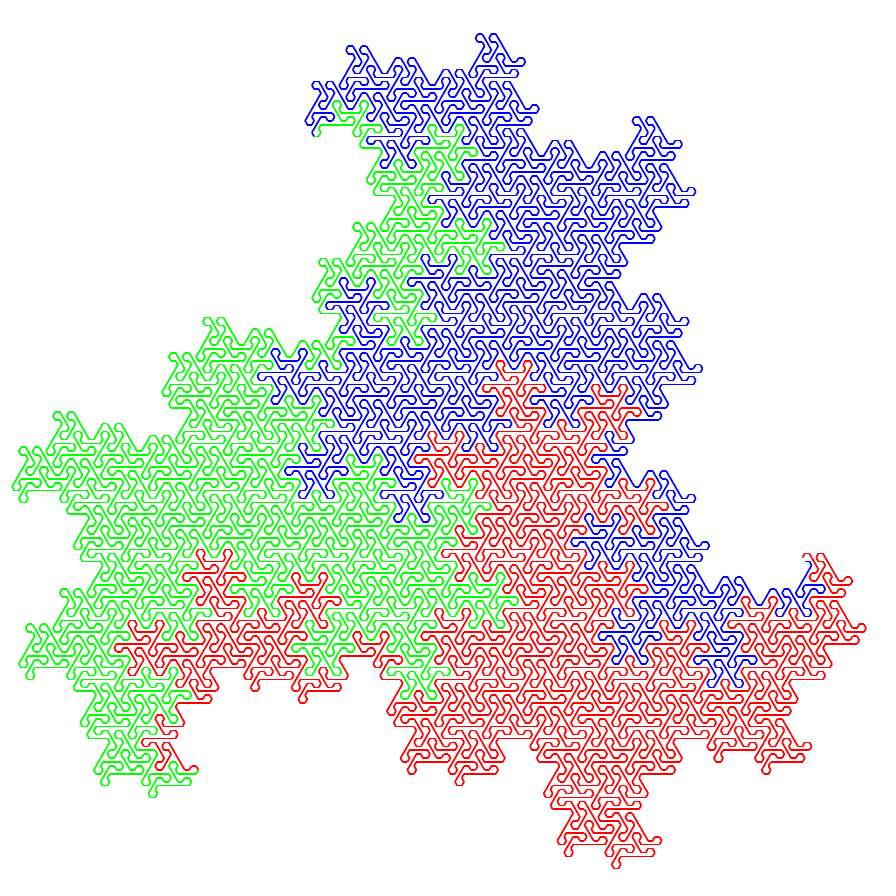}
\end{center}
\else
\verb+{see pdf for image}+
\fi
\caption{\label{fig:r13-t-15-tiles}
The two tiles $\Tile{+3}$ and $\Tile{-3}$ of the curve in the previous figure (\CID{R13-15}).}
\end{figure}
%
%%%%%%%%%%%%%%%%%%%%%%%%%%

%%%%%%%%%%%%%%%%%%%%%%%%%%
% stringsubst 4 [A-[A-[A-[A-[A-[A]B]B]B]B]B]B  A F  B ++F+F[+F+F]-F+F+F+F+F  F F+F+F+F-F-F-F  [ [ ] ] + + - - | tail -1 | ./bin 6 3 0 > tmp-pic.tex && make dotex # R13-t-15 tile-plus, tiling by border
\begin{figure}[h!tbp]
\ifpdf
\begin{center}
\includegraphics*[width=90mm, viewport={60 330 500 750}]{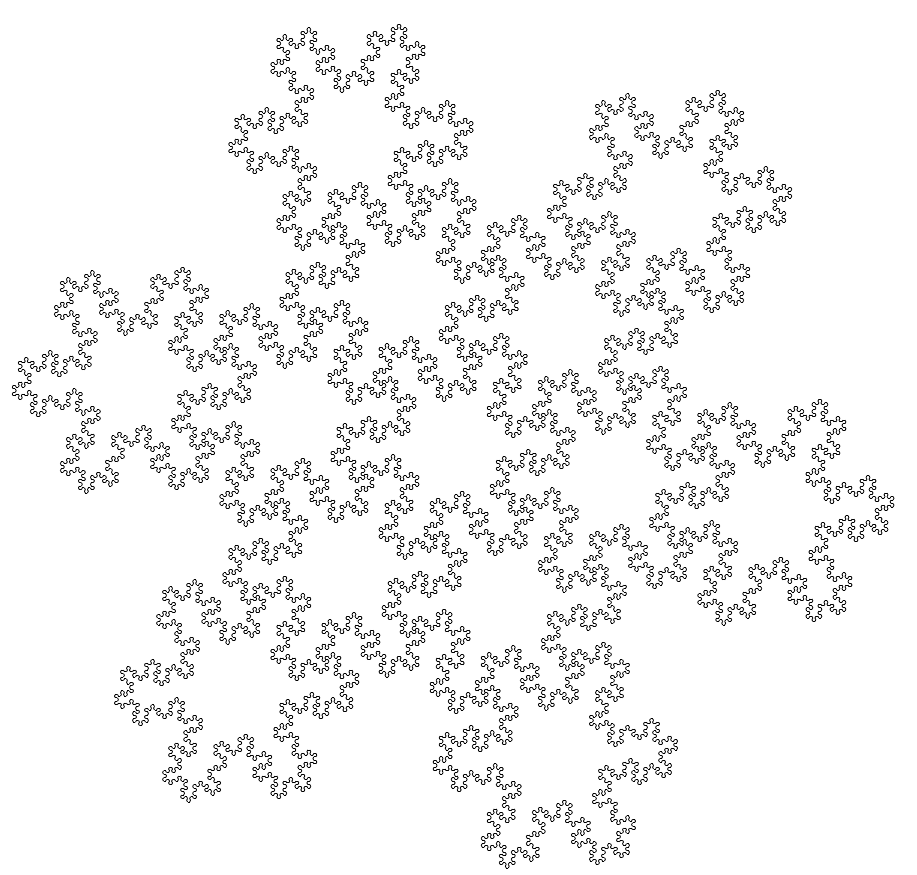}
%\fbox{\includegraphics*[width=108mm, viewport={60 330 500 750}, angle=12]{r13-t-15-plus-tiling-by-border.pdf}}
\end{center}
\else
\verb+{see pdf for image}+
\fi
\caption{\label{fig:r13-t-15-plus-tiling}
The shape $\Tile{+\infty}$ of the tile at the left in the previous figure,
decomposed into 13 smaller rotated copies of itself.
}
\end{figure}
%
%%%%%%%%%%%%%%%%%%%%%%%%%%

All curves are self-similar by construction:
every curve of order $R$ can be decomposed into $R$
disjoint rotated copies of itself.
Figure~\ref{fig:r13-t-15-decomp} shows this for a curve
of order 13 with L-system \Lmap{F}{F+F0F0F-F-F+F0F+F+F-F0F-F}
on the triangular grid.
%
% cf. Figures~\ref{fig:iterate-1-2-decomp} and \ref{fig:iterate-3-4-5-decomp}

Let $\Tile{+k}$ and $\Tile{-k}$ be the tiles for the $k$th iterate of a curve,
and $\Tile{+\infty}$ and $\Tile{-\infty}$ the limiting shape of the tiles.
The tiles $\Tile{+3}$ and $\Tile{-3}$ of this curve are
shown in Figure~\ref{fig:r13-t-15-tiles}
(for $\Tile{+1}$ and $\Tile{-1}$ see Figure~\Ref{fig:r13-t-15-tiles-it1}).
The tile can be decomposed into 13 small copies of itself as
shown in Figure~\ref{fig:r13-t-15-plus-tiling}.

%%%%%%%%%%%%%%%%%%%%%%%%%%
%% with lnth *= 2.0;  // thicker lines
% stringsubst 7 _[A-[A-[A-[A]B]B]B]B  A F  B _++F+F+F  F F+F-F  _ _ [ [ ] ] + + - - | tail -1 | ./bin 4 3 0 > tmp-pic.tex && make dotex # 5-tiling by border
%
% stringsubst 6 X+X+X+X  X F+F-F+F-F-F+F-F-F+F-F-F-F+F+F+F-F+F+F-F+F+F-F+F-F  F F+F-F + + - - | tail -1 | ./bin 4 3 0 > tmp-pic.tex && make dotex # 5-tiling using tile of R25-76
%
\begin{figure}[h!tbp]
\ifpdf
\begin{center}
{\includegraphics*[width=63mm, viewport={60 310 500 740}]{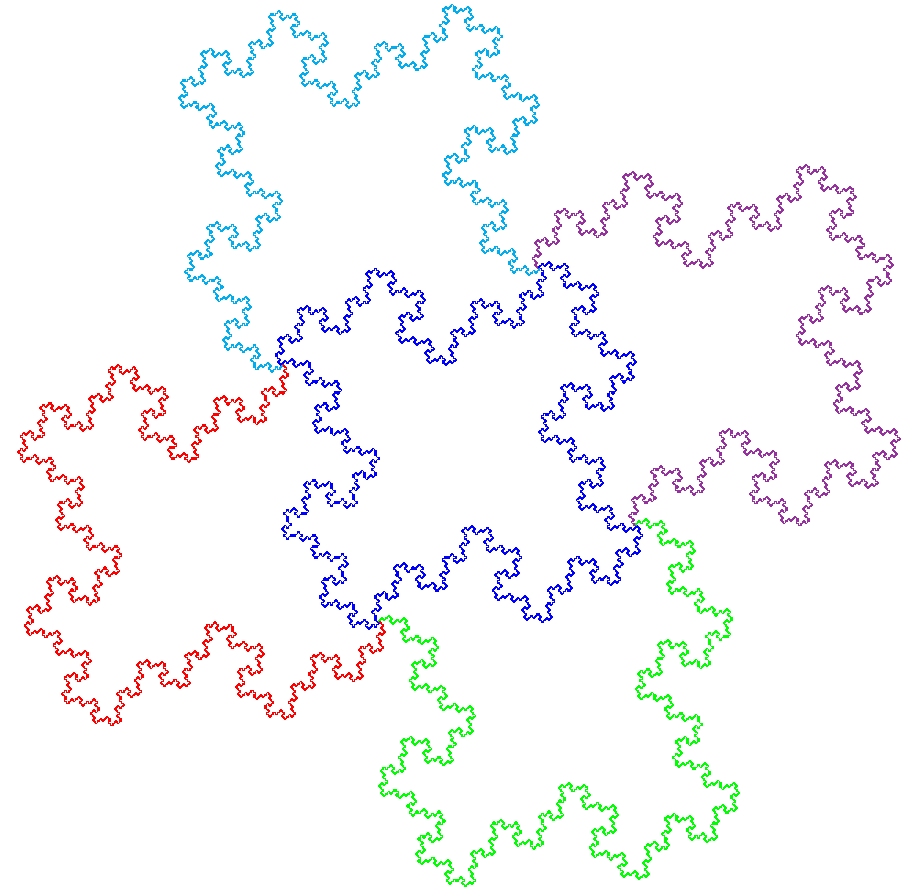}}%
{\includegraphics*[width=63mm, viewport={60 310 500 740}]{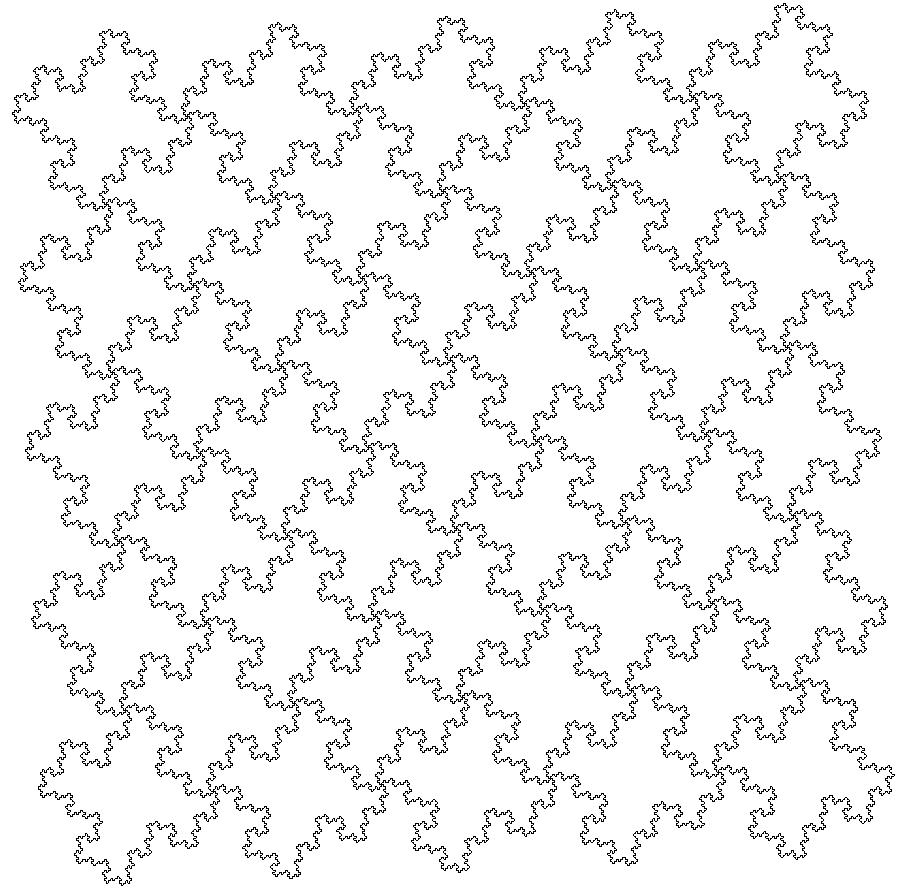}}
\end{center}
\else
\verb+{see pdf for image}+
\fi
\caption{\label{fig:r5-q-tile-tiling}
Decomposition of the tile $\Tile{+\infty}$ of the R5-dragon (\CID{R5-1}) into rotated copies of itself (left)
and the tiling of the plane by such tiles (right).}
\end{figure}
%
%%%%%%%%%%%%%%%%%%%%%%%%%%
%% cf. 5-tiling-by-border*

Such a statement is true for all tiles:
the tiles of a curve of order $R$ can always
be decomposed into $R$ scaled (by $1/\sqrt{R}$)
rotated copies of itself.
Figure~\ref{fig:r5-q-tile-tiling} (left) shows this for the R5-dragon, the
tiling of the plane by translations of the tile is shown on the right.
The curves tile the plane by translations and rotations.

Our tiles are always connected,
and the tilings are lattice tilings (the positions of the tiles are on a lattice).
%  often disk-like, and
% rotational symmetry

%%%%%%%%%%%%%%%%%%%%%%%%%%
% stringsubst 2 F_+F_+F_-F_-F_-F_+F_+F_-F_-F_-F_+F_+F_-F_-F_-F_+F_+F_+F_-F_+F_+F_-F_+F_-F_  _ _  F F+F+F-F-F-F+F+F-F-F-F+F+F-F-F-F+F+F+F-F+F+F-F+F-F   + + - - | tail -1 | ./bin 4 2 0 0 0.3 > tmp-pic.tex && make dotex # R25-46
%
\begin{figure}[h!tbp]
\ifpdf
\begin{center}
{\includegraphics*[width=90mm, viewport={60 350 500 740}]{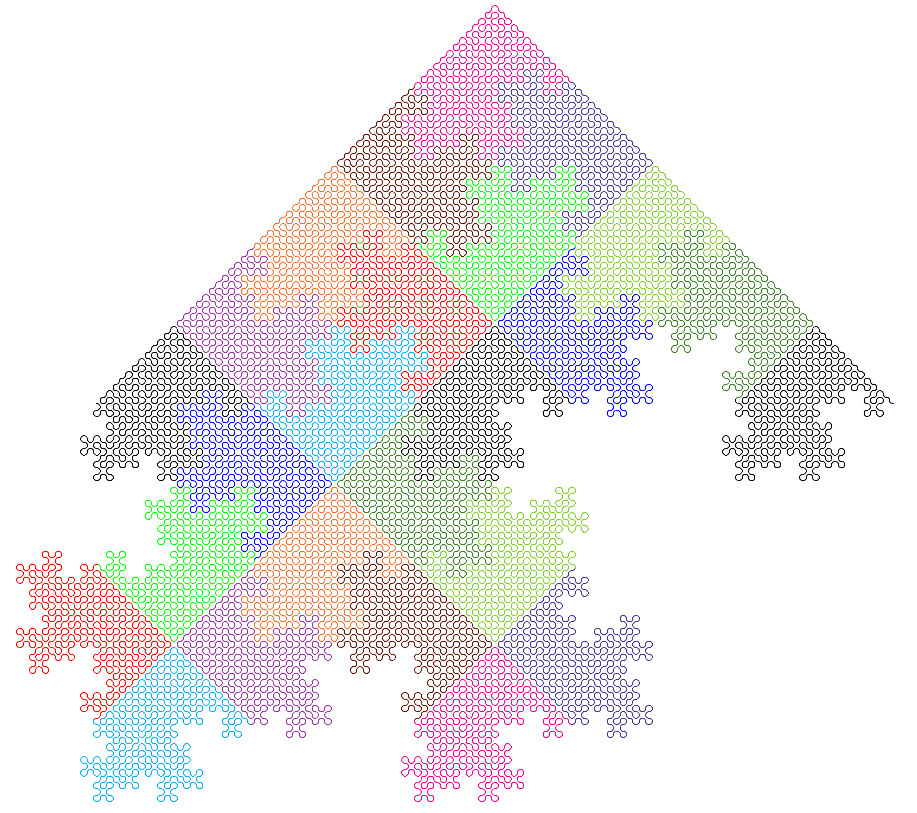}}
\end{center}
\else
\verb+{see pdf for image}+
\fi
\caption{\label{fig:r25-q-46-decomp}
Self-similarity of a curve of order 25 on the square grid (\CID{R25-46}).}
%The map of the L-system is {\small\Lmap{F}{F+F+F-F-F-F+F+F-F-F-F+F+F-F-F-F+F+F+F-F+F+F-F+F-F}}.}
\end{figure}
%
%%%%%%%%%%%%%%%%%%%%%%%%%%

The tiles for the curves on the square grid
always have 4-fold rotational symmetry.
Figure~\ref{fig:r25-q-46-decomp} shows the self-similarity
of the curve \CID{R25-46} on the square grid.
%with map \Lmap{F}{F+F+F-F-F-F+F+F-F-F-F+F+F-F-F-F+F+F+F-F+F+F-F+F-F}.
The tiles appear as groups of four smaller curves.

%%%%%%%%%%%%%%%%%%%%%%%%%%
%% with lnth *= 5.0;  // very thick lines
% stringsubst 1 L  L L+R+L-R+L-R-L-R+L-R+L+R-L  R R+L-R-L+R-L+R+L+R-L+R-L-R   + + - - | tail -1 | ./bin 4 3 1 > tmp-pic.tex && make dotex # R13-rl
%
%% with lnth *= 2.0;  // thicker lines
% stringsubst 3 -L_+R_+L_+R _ _ L L+R+L-R+L-R-L-R+L-R+L+R-L R R+L-R-L+R-L+R+L+R-L+R-L-R + + - - | tail -1 | ./bin 4 2 0 > tmp-pic.tex && make dotex # R13-rl
%
\begin{figure}[h!tbp]
\ifpdf
\begin{center}
{\includegraphics*[width=40mm, viewport={110 420 500 740}]{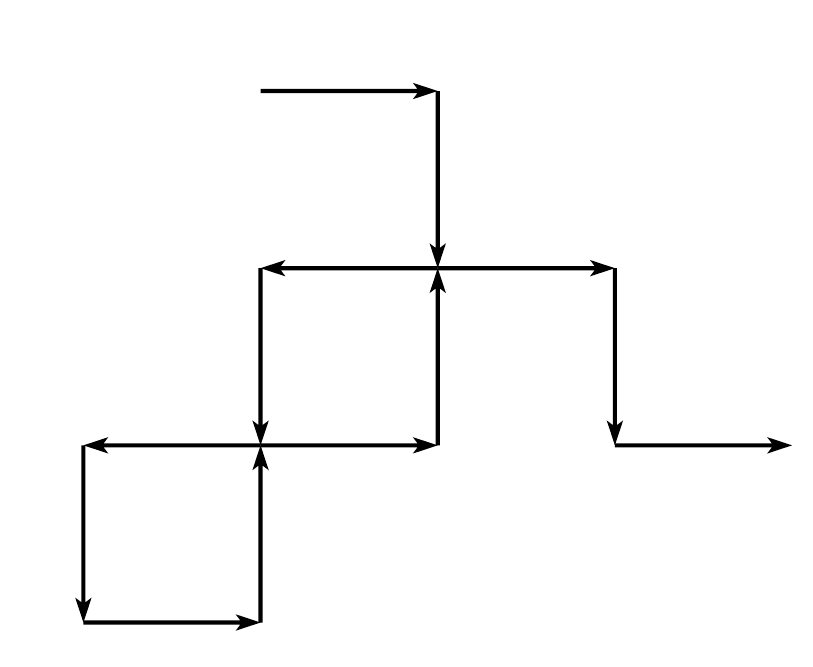}}%
{\includegraphics*[width=85mm, viewport={60 530 500 740}]{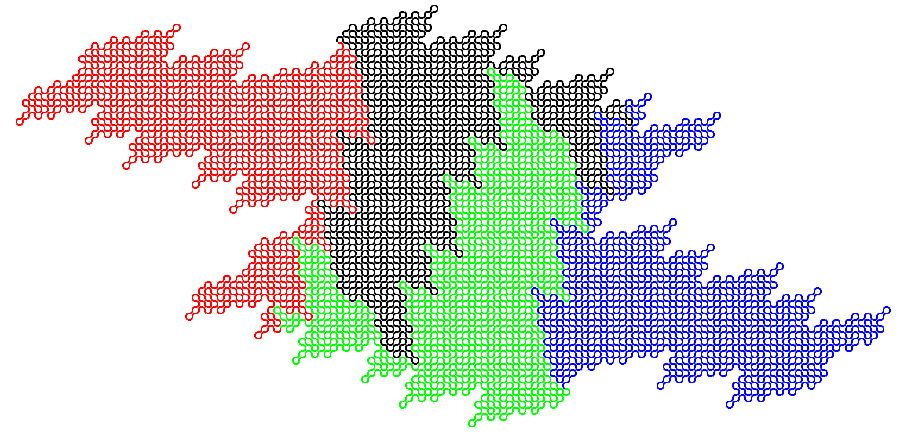}}
\end{center}
\else
\verb+{see pdf for image}+
\fi
\caption{\label{fig:r13-rl-tile}
Motif of a curve of order 13 with an L-system with two non-constant letters
\texttt{L} and \texttt{R} (left)
and third iterate of the tile with axiom \texttt{L+R+L+R} (right).
}
\end{figure}
%
%%%%%%%%%%%%%%%%%%%%%%%%%%

While not found in our search, curves on the square grid with tiles
having only 2-fold rotational symmetry exist.
Figure~\ref{fig:r13-rl-tile} shows the motif of
such a curve of order 13 generated by an L-system with two non-constant letters
\texttt{L} and \texttt{R} (left)
and the tile from the axiom \texttt{L+R+L+R} (right).
The maps for the curve are
\Lmap{L}{L+R+L-R+L-R-L-R+L-R+L+R-L}
and
\Lmap{R}{R+L-R-L+R-L+R+L+R-L+R-L-R},
the production for \texttt{R}
can be obtained from that for \texttt{L}
by reversing and swapping all \texttt{+} and \texttt{-}
and all \texttt{L} and \texttt{R}.

On both grids the two tiles $\Tile{+j}$ and $\Tile{-j}$
for $j\geq{}1$ (and $j=\infty$)
have different shapes in general
(the only exception are tiles that are perfect squares).
If the curve has 2-fold rotational symmetry,
then both tiles are identical.
%(indicated be the letter \texttt{r} in the file)
%% Check:
%% stringsubst 5 F_+F_+F_+F _ _ F F+F+F-F-F + + - - | tail -1 | ./bin 4 3 0 > tmp-pic.tex && make dotex # R5-1
%% stringsubst 5 F_-F_-F_-F _ _ F F+F+F-F-F + + - - | tail -1 | ./bin 4 3 0 > tmp-pic.tex && make dotex # R5-1
%%
%% stringsubst 7 F_+F_+F_  _ _ F F+F-F + + - - | tail -1 | ./bin 3 3 0 > tmp-pic.tex && make dotex # R5-1
%% stringsubst 7 F_-F_-F_  _ _ F F+F-F + + - - | tail -1 | ./bin 3 3 0 > tmp-pic.tex && make dotex # R5-1
%

%%%%%%%%%%%%%%%%%%%%%%%%%%
%% with lnth *= 2.0;  // thicker lines
% stringsubst 2 _F-_F-_F _ _  F F0F+F0F+F-F0F+F-F-F0F-F+F0F0F+F+F-F0F-F+F0F-F-F+F   0 0 + + - - | tail -1 | ./bin 3 2 0 > tmp-pic.tex && make dotex # R25-247
% stringsubst 2 _F-_F-_F _ _  F F0F+F0F+F-F0F+F-F-F0F-F+F0F-F+F+F0F0F-F+F0F-F-F+F   0 0 + + - - | tail -1 | ./bin 3 2 0 > tmp-pic.tex && make dotex # R25-248
%
\begin{figure}[h!tbp]
\ifpdf
\begin{center}
{\includegraphics*[width=63mm, viewport={60 330 500 750}]{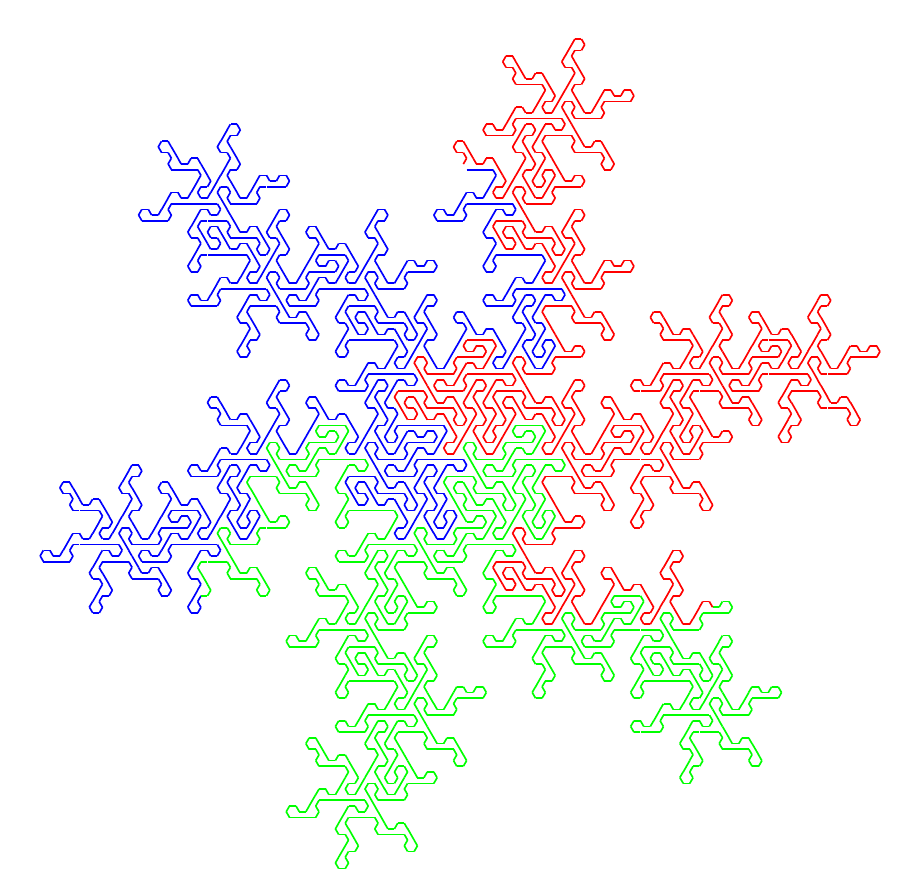}}%
{\includegraphics*[width=63mm, viewport={60 330 500 750}]{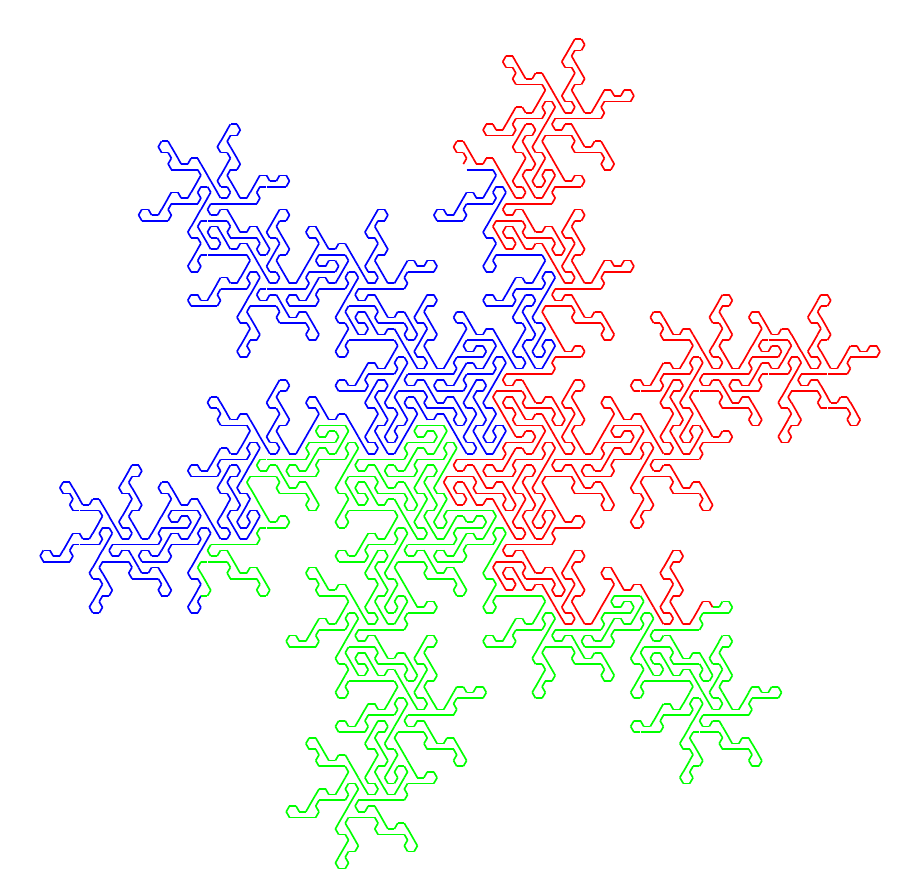}}
\end{center}
\else
\verb+{see pdf for image}+
\fi
\caption{\label{fig:r25-t-tiles-same-shape}
Two curves with different shapes but same shape of the tiles $\Tile{-k}$ for all $k$,
shown are the tiles $\Tile{-2}$ (\CID{R25-247} and \CID{R25-248}).}
\end{figure}
%
%%%%%%%%%%%%%%%%%%%%%%%%%%

Note that curves of different shapes can lead to the same shape of tiles,
see Figure~\ref{fig:r25-t-tiles-same-shape} for an example
with curves of order 25 on the triangular grid.
%

%%%%%%%%%%%%%%%%%%%%%%%%%%
%% with lnth *= 2.0;  // thicker lines
%
%% the tile
% stringsubst 2 _F+_F+_F _ _ F F+F0F0F-F-F+F0F+F+F-F0F-F 0 0 + + - - | tail -1 | sed 's/+/++++/g; s/-/----/g;' | ./bin 12 3 0 > tmp-pic.tex && make dotex # R13-15
%
%% elongated lozenges:
%% stringsubst 2 _F+_F+_F _ _ F F+F0F0F-F-F+F0F+F+F-F0F-F 0 0 + + - - | tail -1 | sed 's/+/++++/g; s/-/----/g; s/F/[F-F+]-F+F/g;' | ./bin 12 3 0 > tmp-pic.tex && make dotex # R13-15
% better (rotation OK):
% stringsubst 2 R_F+_F+_F R R _ _ F F+F0F0F-F-F+F0F+F+F-F0F-F 0 0 + + - - | tail -1 | sed 's/+/++++++++/g; s/-/--------/g; s/F/[F--F++]--F++F/g; s/R/+/;' | ./bin 24 3 0 > tmp-pic.tex && make dotex # R13-15
%
%% full lozenges (cubes):
% stringsubst 2 _F+_F+_F _ _ F F+F0F0F-F-F+F0F+F+F-F0F-F 0 0 + + - - | tail -1 | sed 's/+/++++/g; s/-/----/g; s/F/[+F--F+]-F++F-/g;' | ./bin 12 3 0 > tmp-pic.tex && make dotex # R13-15
%
%% full lozenges, no center line (hexes):
%% stringsubst 2 _F+_F+_F _ _ F F+F0F0F-F-F+F0F+F+F-F0F-F 0 0 + + - - | tail -1 | sed 's/+/++++/g; s/-/----/g; s/F/[++F--F]--F++F/g;' | ./bin 12 3 0 > tmp-pic.tex && make dotex # R13-15
% better (no overlap):
%% stringsubst 2 _F+_F+_F _ _ F F+F0F0F-F-F+F0F+F+F-F0F-F 0 0 + + - - | tail -1 | sed 's/+/++++/g; s/-/----/g; s/F/--F++F/g;' | ./bin 12 3 0 > tmp-pic.tex && make dotex # R13-15
%% still better (rotation OK):
% stringsubst 2 R_F+_F+_F R R _ _ F F+F0F0F-F-F+F0F+F+F-F0F-F 0 0 + + - - | tail -1 | sed 's/+/++++/g; s/-/----/g; s/F/--F++F/g; s/R/+/;' | ./bin 12 3 0 > tmp-pic.tex && make dotex # R13-15
%
\begin{figure}[h!tbp]
\ifpdf
\begin{center}
{\includegraphics*[width=50mm, viewport={70 330 480 750}]{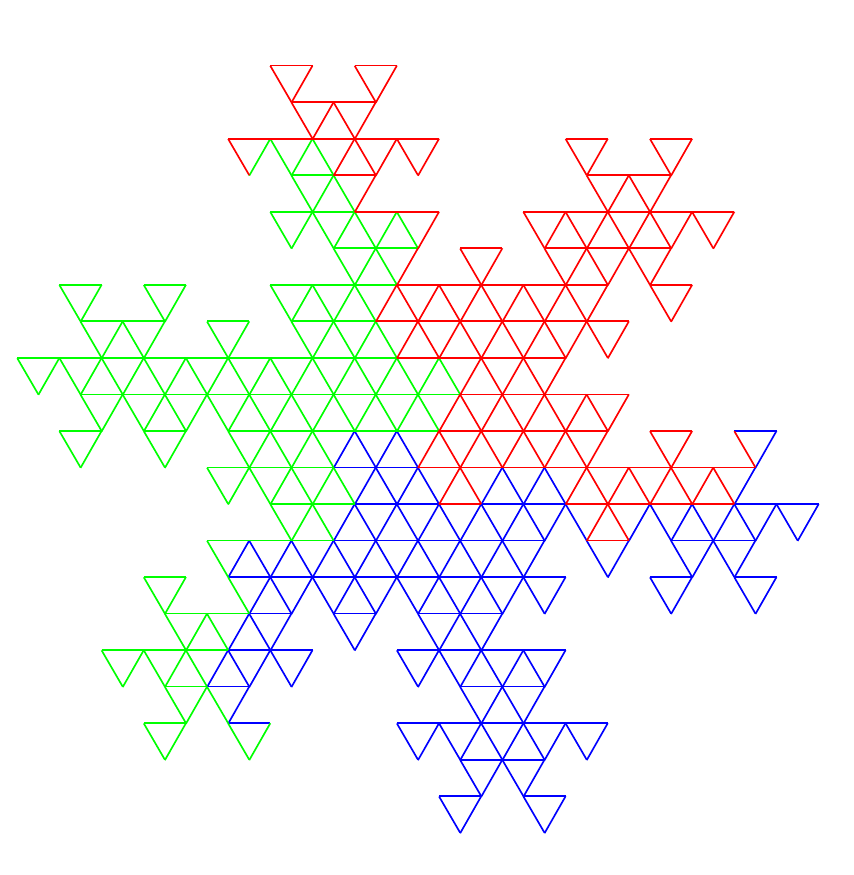}}%
{\includegraphics*[width=50mm, viewport={70 330 480 750}]{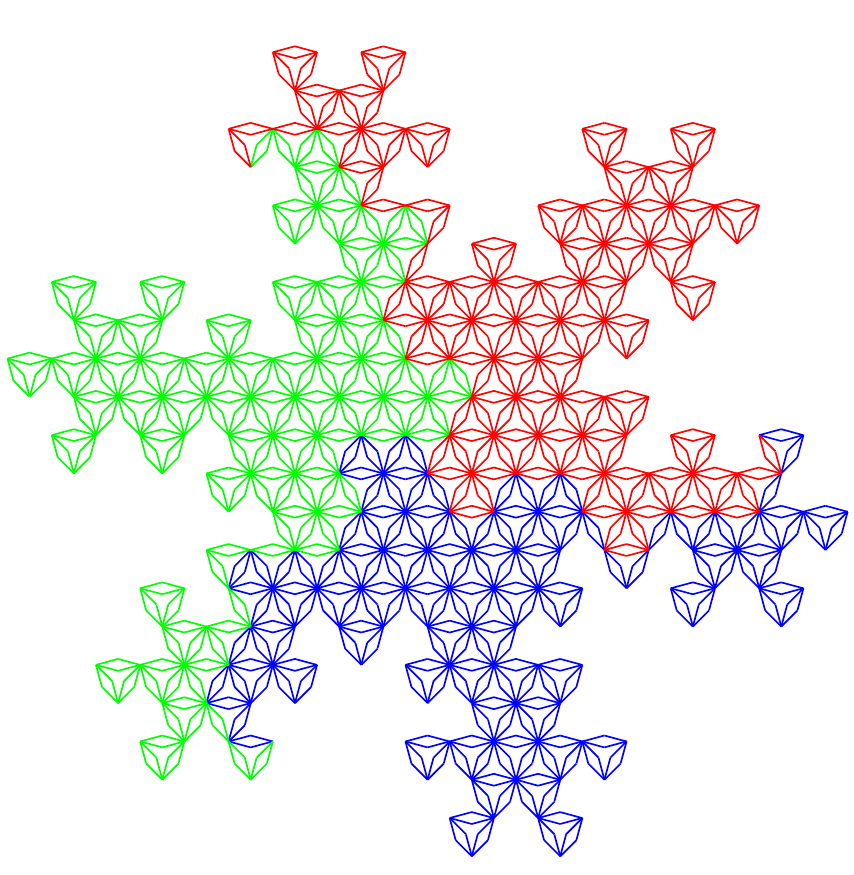}}
{\includegraphics*[width=50mm, viewport={70 330 480 750}]{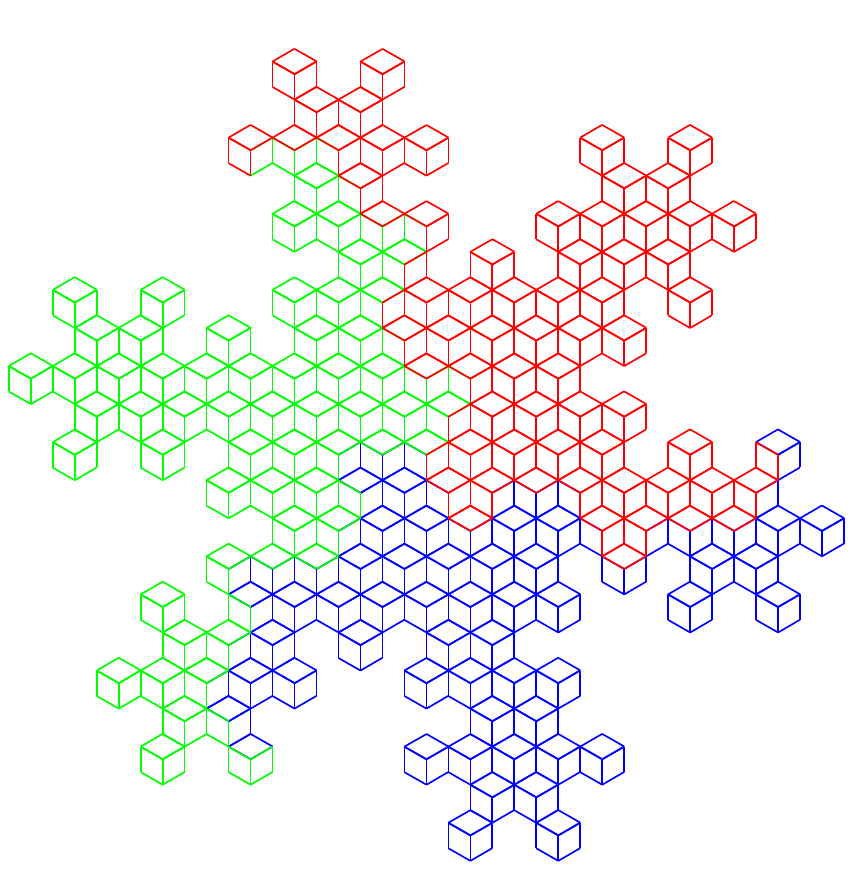}}%
{\includegraphics*[width=50mm, viewport={70 320 480 740}]{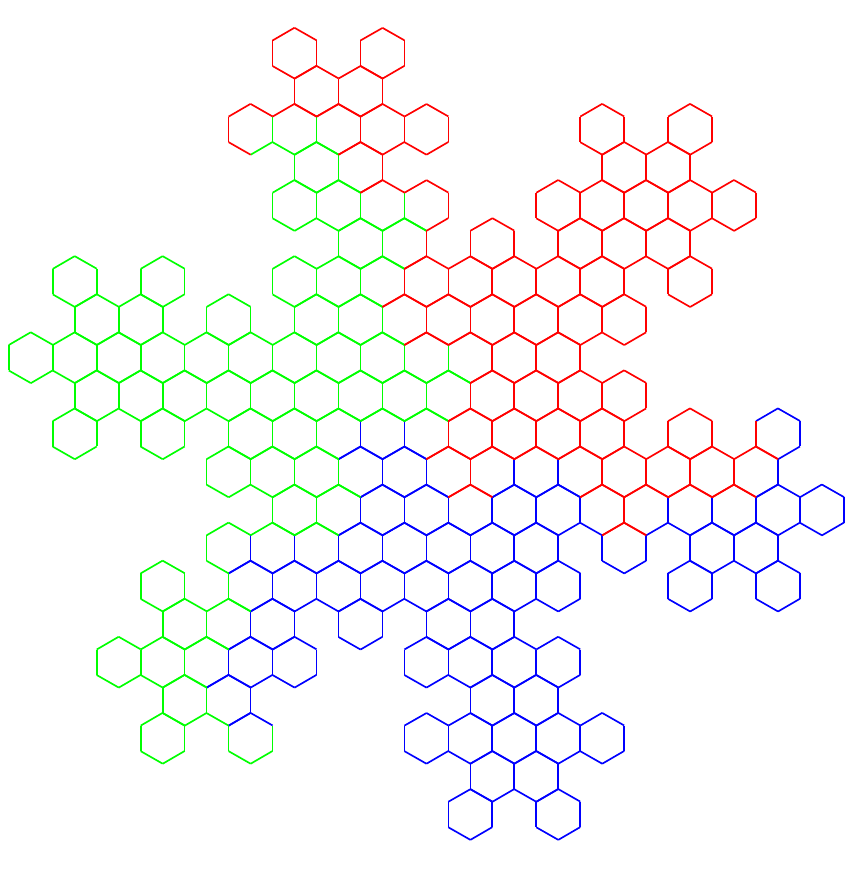}}
\end{center}
\else
\verb+{see pdf for image}+
\fi
\caption{\label{fig:r13-t-15-tile-6fold-symm}
Rendering a tile with 6-fold symmetry, see text.}
\end{figure}
%
%%%%%%%%%%%%%%%%%%%%%%%%%%

The tiles for the curves on the triangular grid always have 3-fold rotational symmetry.
Some tiles appear to have 6-fold rotational symmetry
%(see Figure~\ref{fig:r25-t-tiles-same-shape})
but this is only true in the limit.
With some artistic license such tiles
can be rendered with 6-fold symmetry,
we use $\Tile{+2}$ of \CID{R13-15}
(for $\Tile{+3}$ see Figure~\ref{fig:r13-t-15-tiles}).
The upper left of Figure~\ref{fig:r13-t-15-tile-6fold-symm}
shows that the tile does not actually have a 6-fold symmetry.
In the rendering to the right
each line has been replaced with a thin lozenge.
Choosing $60\adeg$ for the acute angle of the lozenges
leads to the rendering in the lower left, an apparent
arrangement of cubes.
If we keep just the outline of the cubes
(lower right), we obtain the arrangement of
hexagons corresponding to a 6-symmetric tile
of the hexagonal grid.

%%%%%%%%%%%%%%%%%%%%%%%%%%
%% with lnth *= 2.0;  // thicker lines
%
%% for 5-tile:
% stringsubst 3 F+_F+_F+_F _ _ F F+F+F-F-F 0 0 + + - - | tail -1 | ./bin 4 3 0 > tmp-pic.tex && make dotex
%
%% thin lozenges:
% stringsubst 3 RF+_F+_F+_F R R _ _ F F+F+F-F-F 0 0 + + - - | tail -1 | sed 's/+/++++/g; s/-/----/g; s/F/[F--F++]--F++F/g; s/R/+/;' | ./bin 16 3 0 > tmp-pic.tex && make dotex
%
%% full lozenges:
% stringsubst 3 RF+_F+_F+_F R R _ _ F F+F+F-F-F 0 0 + + - - | tail -1 | sed 's/+/++/g; s/-/--/g; s/F/[F--F++]--F++F/g; s/R/+/;' | ./bin 8 3 0 > tmp-pic.tex && make dotex
%
%% as tile shape:
% stringsubst 3 RF+_F+_F+_F R R _ _ F F+F+F-F-F 0 0 + + - - | tail -1 | sed 's/+/++/g; s/-/--/g; s/F/--F++F/g; s/R/+/;' | ./bin 8 3 0 > tmp-pic.tex && make dotex
%
\begin{figure}[h!tbp]
\ifpdf
\begin{center}
{\includegraphics*[width=50mm, viewport={60 315 490 745}]{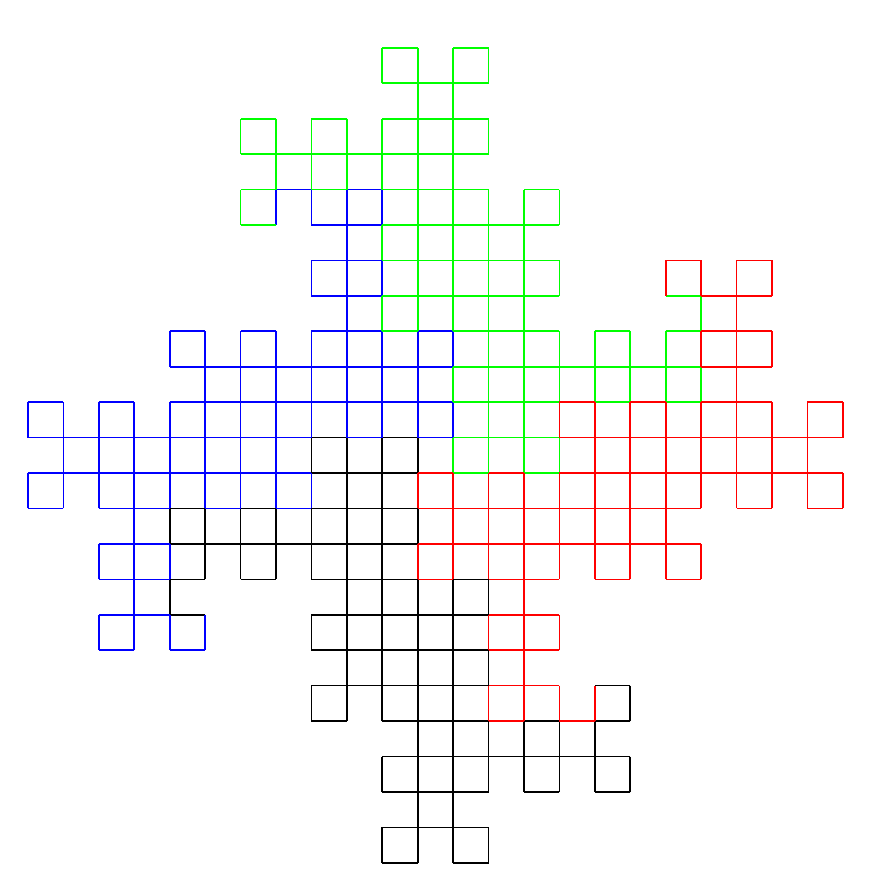}}%
{\includegraphics*[width=50mm, viewport={60 310 490 745}]{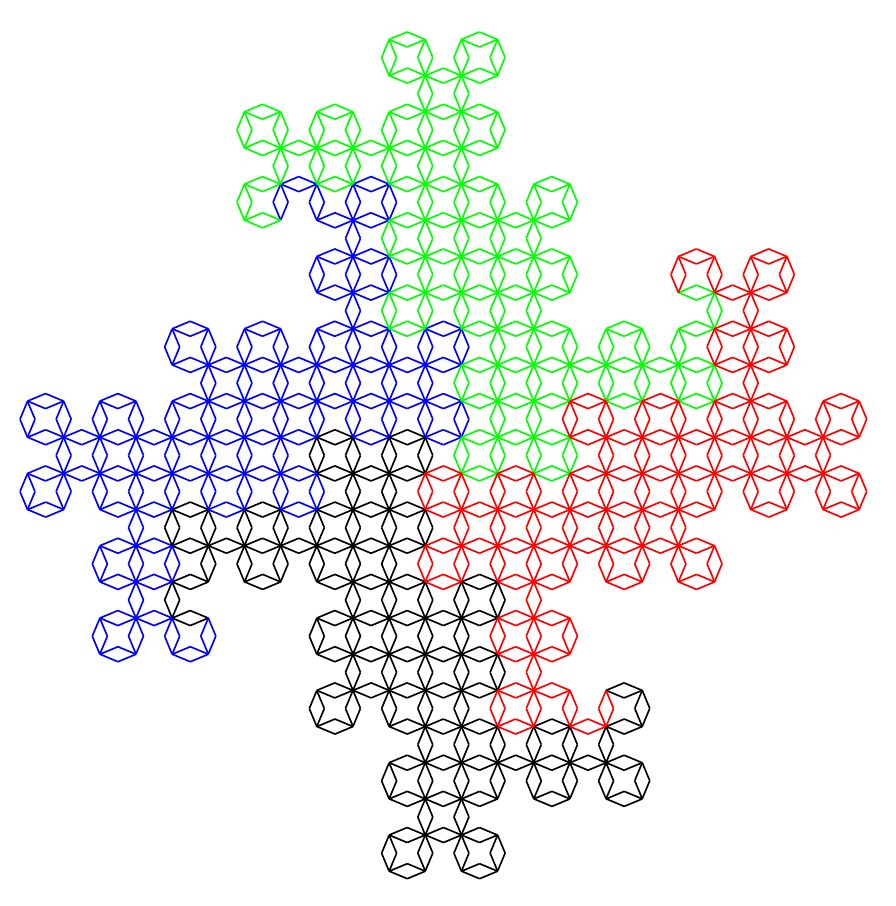}}
{\includegraphics*[width=50mm, viewport={60 320 490 740}]{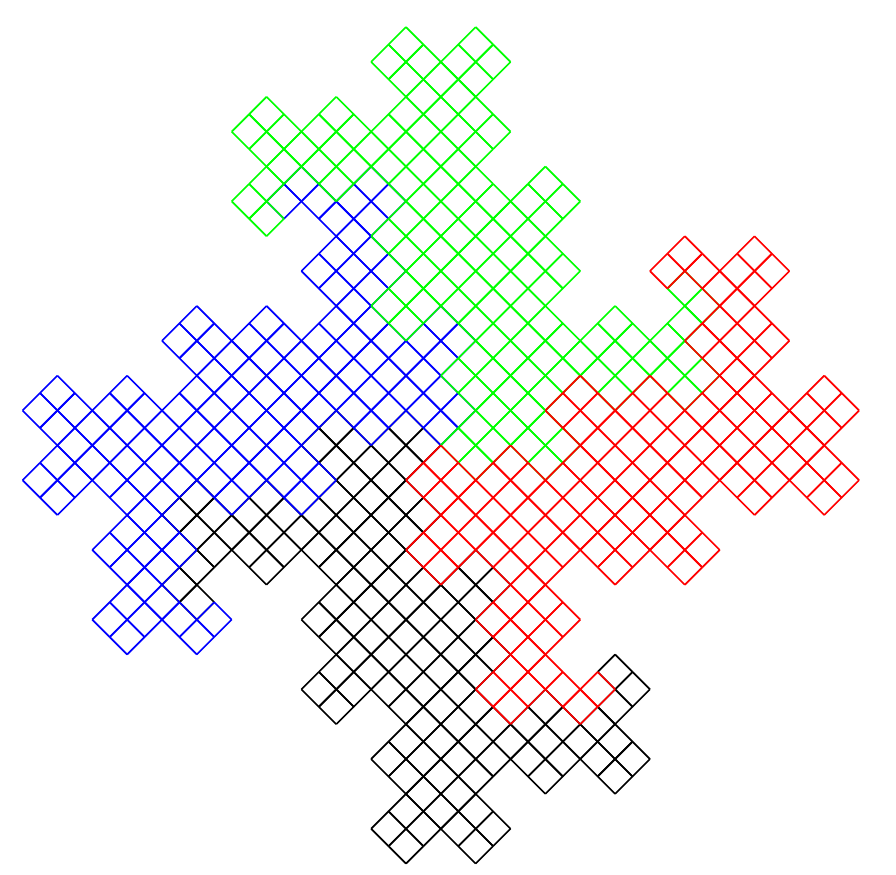}}%
{\includegraphics*[width=50mm, viewport={60 320 490 740}]{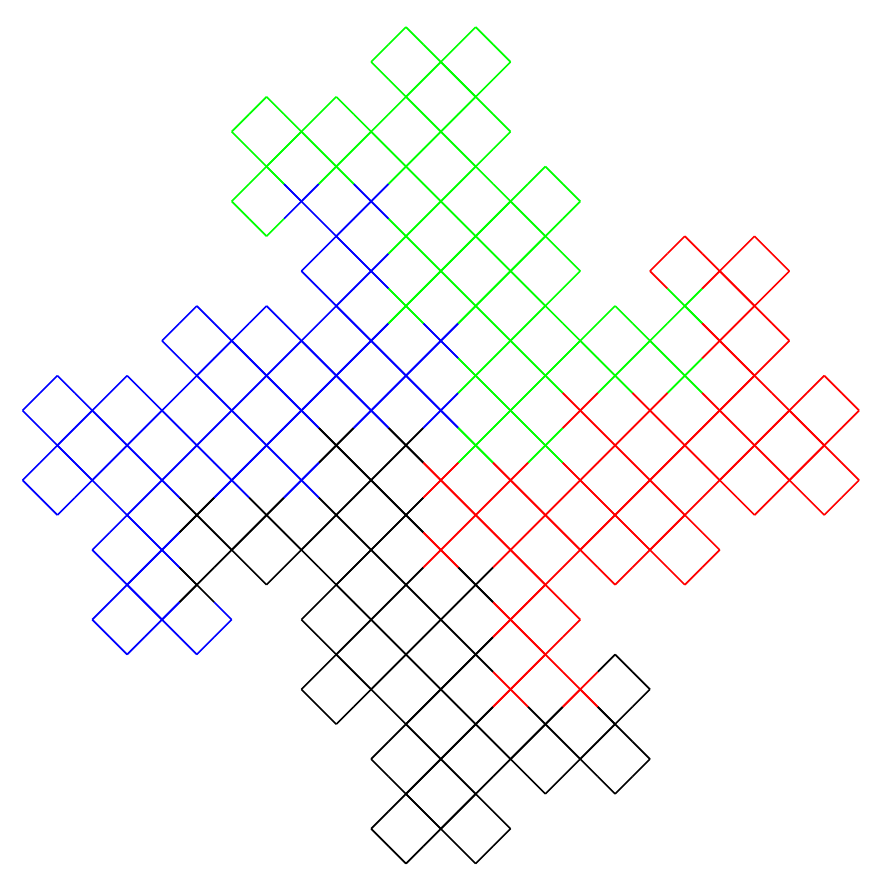}}
\end{center}
\else
\verb+{see pdf for image}+
\fi
\caption{\label{fig:r5-q-tile-rendering}
Rendering a tile on the square grid, see text.}
\end{figure}
%
%%%%%%%%%%%%%%%%%%%%%%%%%%

The analogous rendering for the square grid
is shown in Figure~\ref{fig:r5-q-tile-rendering}
for the tile $\Tile{+3}$ of the R5-dragon.
Here the full \eqq{lozenges} shown in the lower left
are squares that are rotated by $45\adeg$ against the original grid.

\FloatBarrier
%\clearpage% xxx
%%%%%%%%%%%%%%%%%%%%%%%%%%%%%%%
%\subsection{Carousels on the triangular grid}%\label{sect:}
\subsection{Certain arrangements of curves}%\label{sect:}

%%%%%%%%%%%%%%%%%%%%%%%%%%
% stringsubst 4 [_F+_F+_F]-[_F+_F+_F]-[_F+_F+_F] _ _ [ [ ] ] F F0F+F0F-F-F+F 0 0 + + - - | tail -1 | ./bin 3 3 0 > tmp-pic.tex && make dotex # R7-t-1 tiling
%
%% stringsubst 4 [F]+_[G]+_[F]+_[G]+_[F]+_[G] _ _ [ [ ] ] F F0F++F0F--F--F++F  G G--G++G++G0G--G0G 0 0 + + - - | tail -1 | ./bin 6 3 0 > tmp-pic.tex && make dotex # carousel
%% colors as above:
% stringsubst 4 _[F]+__[G]---_[F]+__[G]---_[F]+__[G] _ _ [ [ ] ] F F0F++F0F--F--F++F G G--G++G++G0G--G0G 0 0 + + - - | tail -1 | ./bin 6 3 0 > tmp-pic.tex && make dotex # carousel
%
\begin{figure}[h!tbp]
\ifpdf
\begin{center}
{\includegraphics*[width=63mm, viewport={60 330 520 740}]{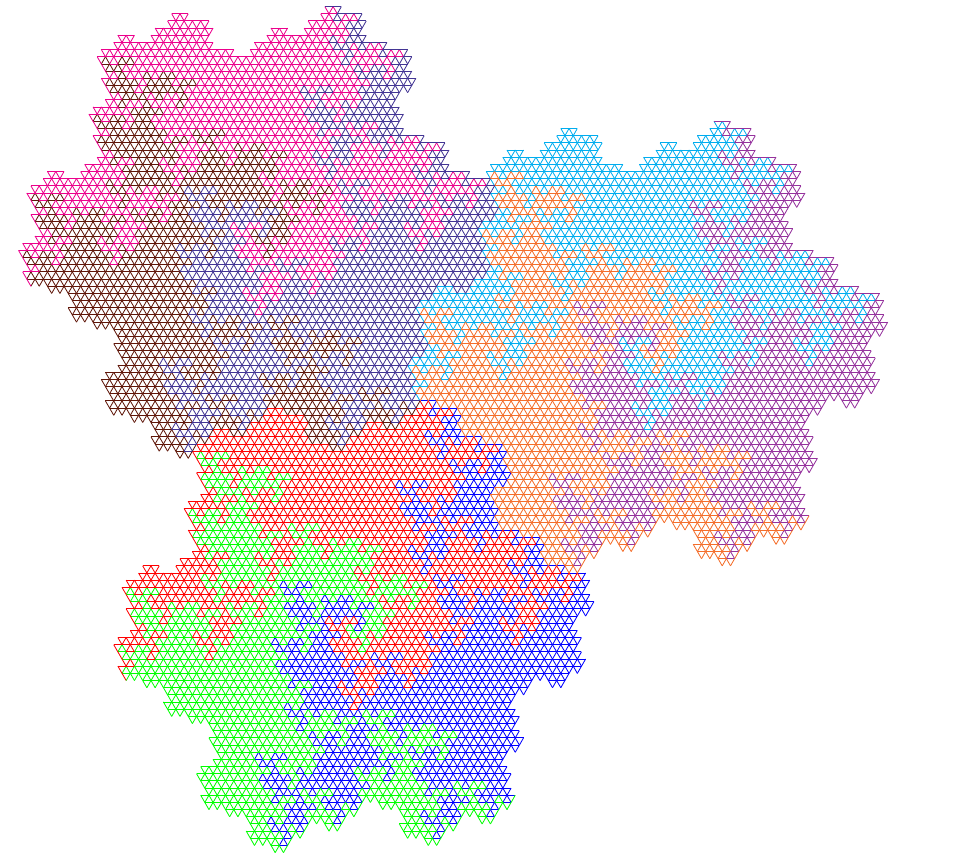}}
{\includegraphics*[width=63mm, viewport={60 320 520 740}]{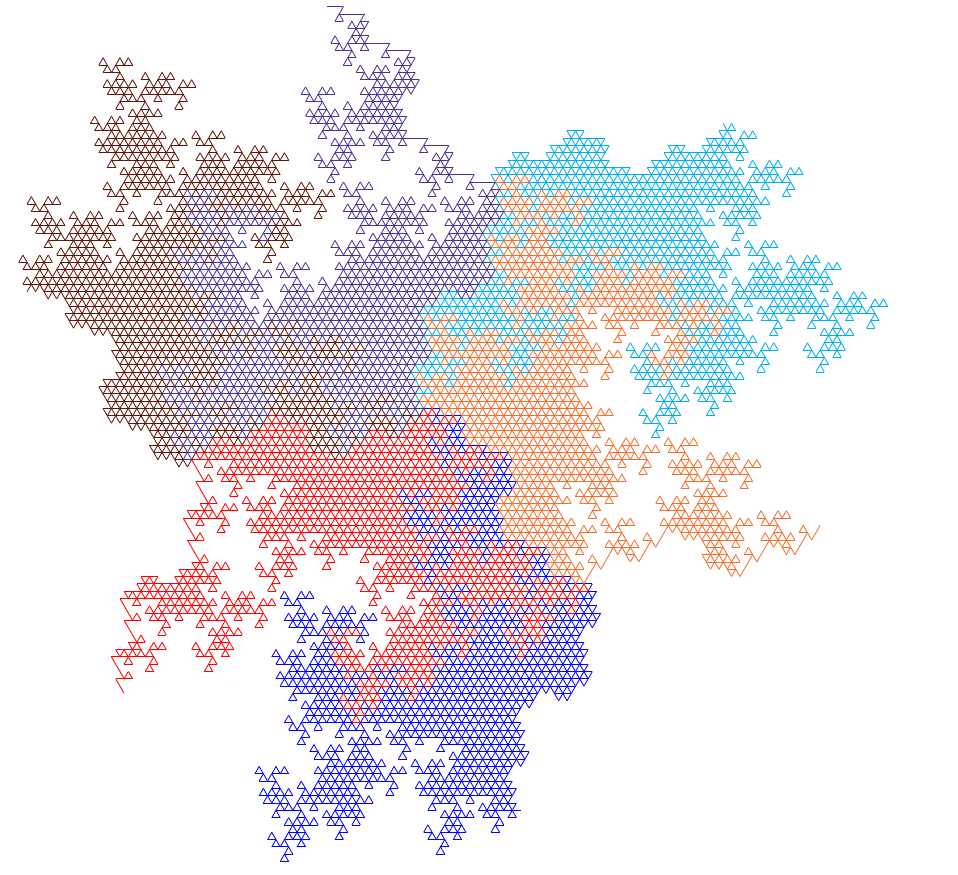}}
\end{center}
\else
\verb+{see pdf for image}+
\fi
\caption{\label{fig:r7-t-1-carousel}
Three tiles of the curve \CID{R7-1} in the tiling of the plane (left)
and the arrangement of the curves meeting at the center (right).}
%Analogue of a carousel for the unsymmetric curve \CID{R7-1}.}
\end{figure}
%
%%%%%%%%%%%%%%%%%%%%%%%%%%

In the tilings of the triangular grid there are points where six curves end.
These curves give an arrangement with 3-fold symmetry as shown for
the curve \CID{R7-1} in Figure~\ref{fig:r7-t-1-carousel}.
%For curves without 2-fold symmetry, the analogue of the carousel
%can be picked at the positions just described,
%here every second curve must be drawn reversed,
%see Figure~\ref{fig:r7-t-1-carousel} for an example using the curve \CID{R7-1}.

%%%%%%%%%%%%%%%%%%%%%%%%%%
% stringsubst 7 _F+t+F+t+F_-t-F+t+F+t+F_-t-F+t+F+t+F+t _ _ F F+t+F-t-F t t + + - - | tail -1 | ./bin 6 3 0 > tmp-pic.tex && make dotex
%% smaller version:
% stringsubst 6 _F+F+F_-F+F+F_-F+F+F _ _ F F+F-F + + - - | tail -1 | ./bin 3 3 0 > tmp-pic.tex && make dotex
%
% stringsubst 7 _F+t+_F+t+_F-t-_F+t+_F+t+_F-t-_F+t+_F+t+_F+t _ _ F F+t+F-t-F t t + + - - | tail -1 | ./bin 6 3 0 > tmp-pic.tex && make dotex
%% smaller version:
% stringsubst 6 _F+_F+_F-_F+_F+_F-_F+_F+_F _ _ F F+F-F + + - - | tail -1 | ./bin 3 3 0 > tmp-pic.tex && make dotex
%
\begin{figure}[h!tbp]
\ifpdf
\begin{center}
%
%{\includegraphics*[width=64mm, viewport={60 290 520 740}]{terdragon-tile-2.pdf}}%
%{\includegraphics*[width=64mm, viewport={60 290 520 740}]{terdragon-carousel.pdf}}
%% smaller version:
{\includegraphics*[width=64mm, viewport={60 330 520 740}]{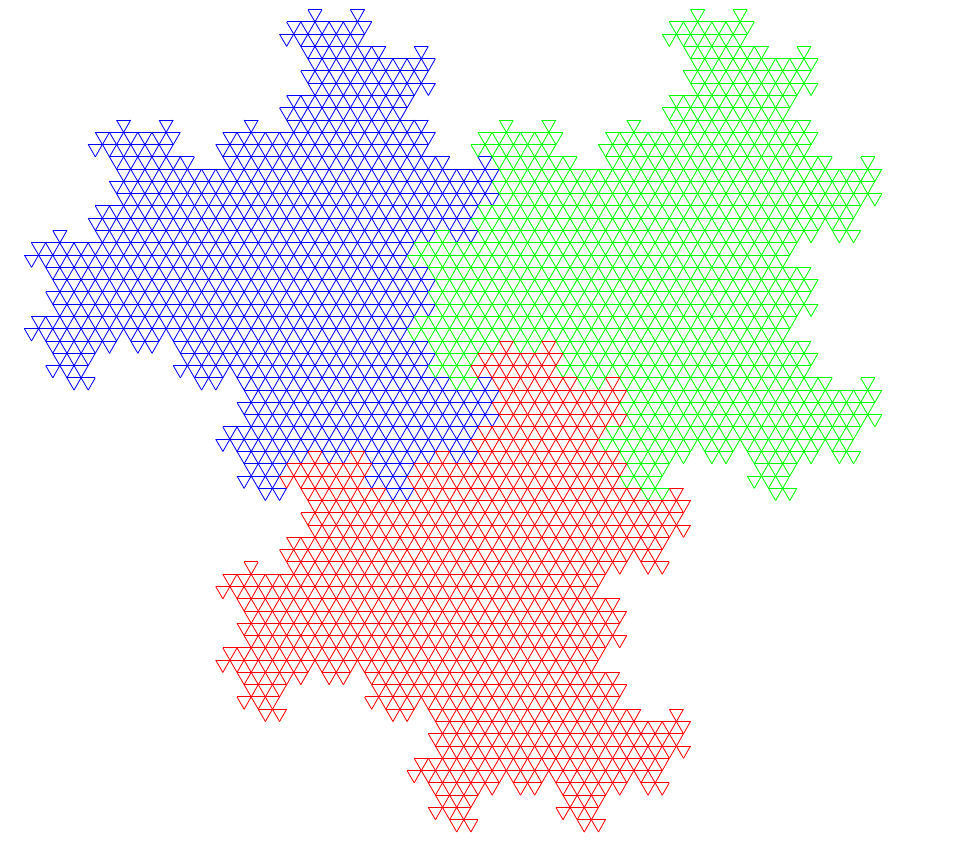}}%
{\includegraphics*[width=64mm, viewport={60 330 520 740}]{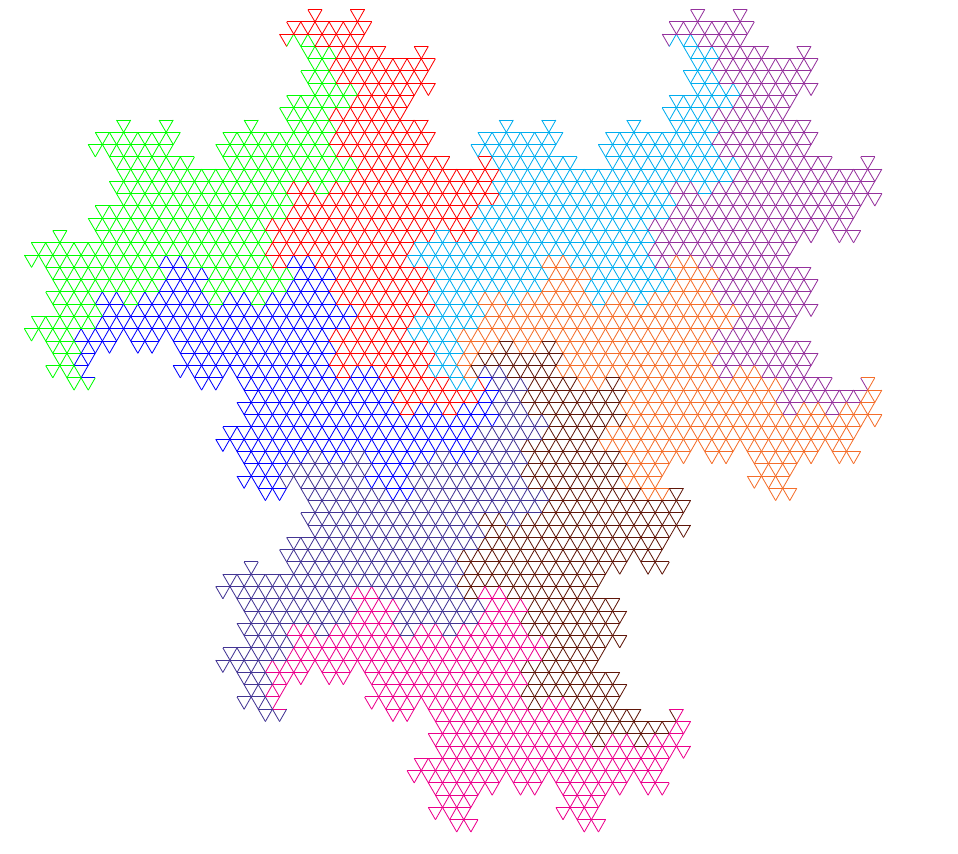}}
\end{center}
\else
\verb+{see pdf for image}+
\fi
\caption{\label{fig:terdragon-carousel}
%The tile $\Tile{+7}$ of the terdragon (top) and its shape $\Tile{+\infty}$ (bottom).}
Three tiles of the terdragon in the tiling of the plane (left).
The six curves meeting at the center give an arrangement with 6-fold symmetry (right).}
%The decomposition into terdragons shows a carousel (right, inner six terdragons).}
\end{figure}
%
%%%%%%%%%%%%%%%%%%%%%%%%%%

%%%%%%%%%%%%%%%%%%%%%%%%%%%
%% with  lnth *= 2.0;  // thicker lines
% stringsubst 4 [L]+[_L]+[_L]+[_L] _ _ [ [ ] ] L L+R-L-R R L+R+L-R + + - - | tail -1 | ./bin 4 3 0 > tmp-pic.tex && make dotex #
%
% stringsubst 4 [R]+[_R]+[_R]+[_R] _ _ [ [ ] ] L L+R-L-R R L+R+L-R + + - - | tail -1 | ./bin 4 3 0 0 0.15 > tmp-pic.tex && make dotex #
%
\begin{figure}[h!tbp]
\ifpdf
\begin{center}
{\includegraphics*[width=50mm, viewport={50 310 520 740}]{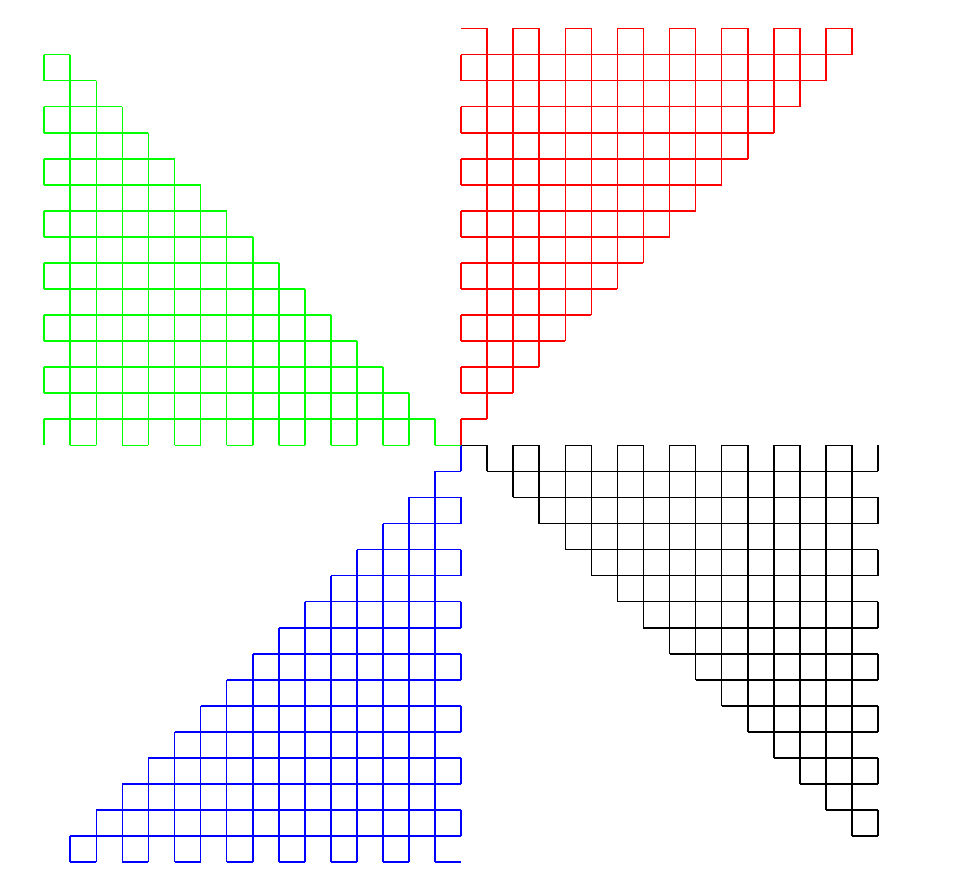}}%
{\includegraphics*[width=52mm, viewport={10 310 490 750}]{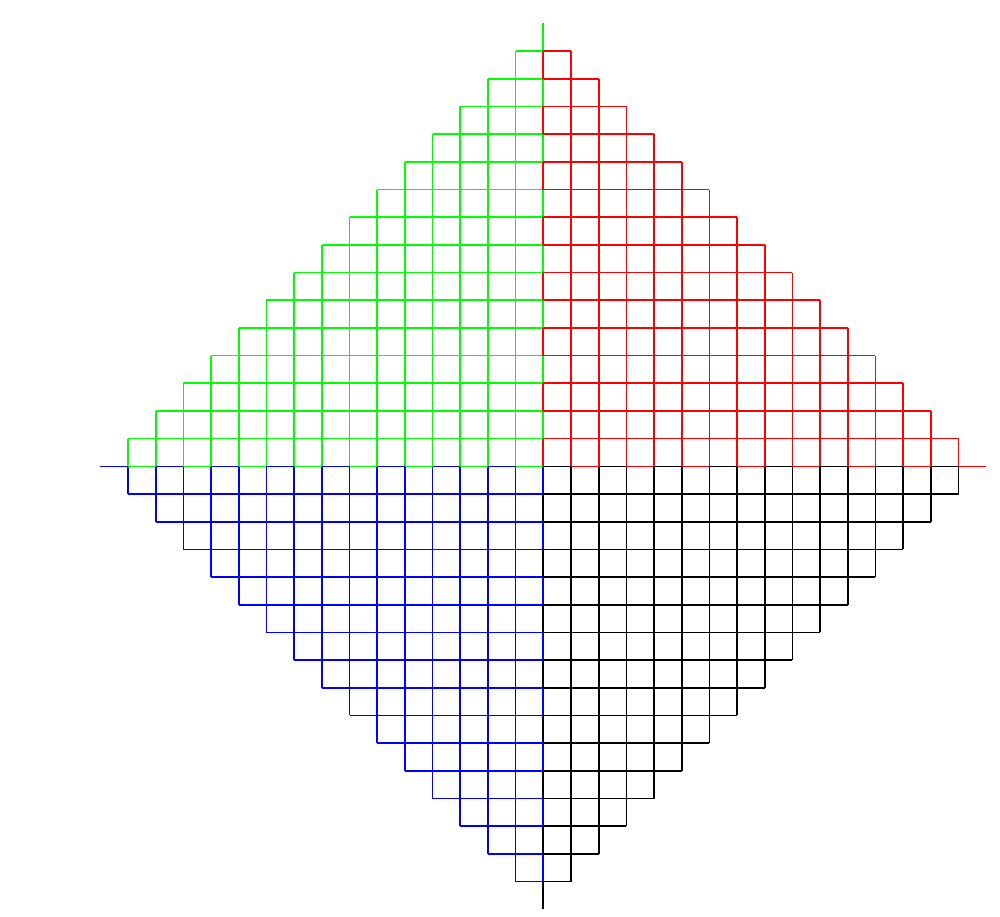}}
\end{center}
\else
\verb+{see pdf for image}+
\fi
\caption{\label{fig:alt-paperfold-carousel}
Only one of the arrangements of the curve
with maps \Lmap{L}{L+R-L-R} and \Lmap{R}{L+R+L-R}
contains a neighborhood of the origin.}
\end{figure}
%
%%%%%%%%%%%%%%%%%%%%%%%%%%%

%\ref{fig:terdragon-tile}
%
For curves with 2-fold symmetry, such as the terdragon,
the arrangement has 6-fold symmetry
as shown in Figure~\ref{fig:terdragon-carousel}.
This arrangement is also shown in \cite[Figure~19, p.~596]{davis-knuth-new},
also see \cite[Section~8, pp.~29-33]{dekking-tiles-final}
where these arrangements are called \eqq{carousels},
and curves for which the arrangements fill the plane
are called \eqq{perfect}.
An example of a curve that is not perfect
is shown in Figure~\ref{fig:alt-paperfold-carousel},
called \eqq{alternate dragon curve} in \cite[Figure~10, p.~590]{davis-knuth-new}.

\begin{figure}[h!tbp]
\ifpdf
\begin{center}
{\includegraphics*[width=52mm, viewport={70 160 490 740}]{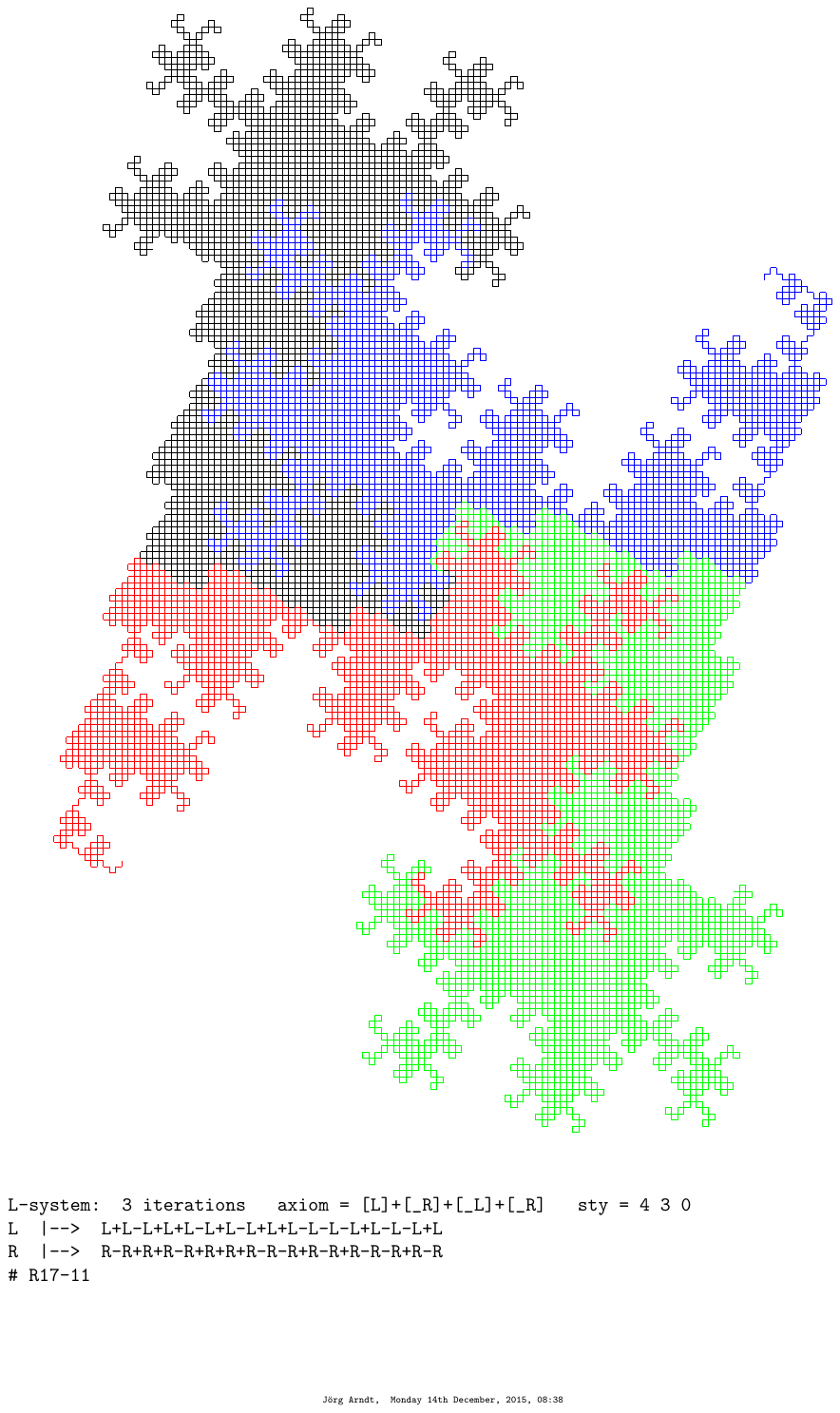}}%
{\includegraphics*[width=72mm, viewport={60 310 500 740}]{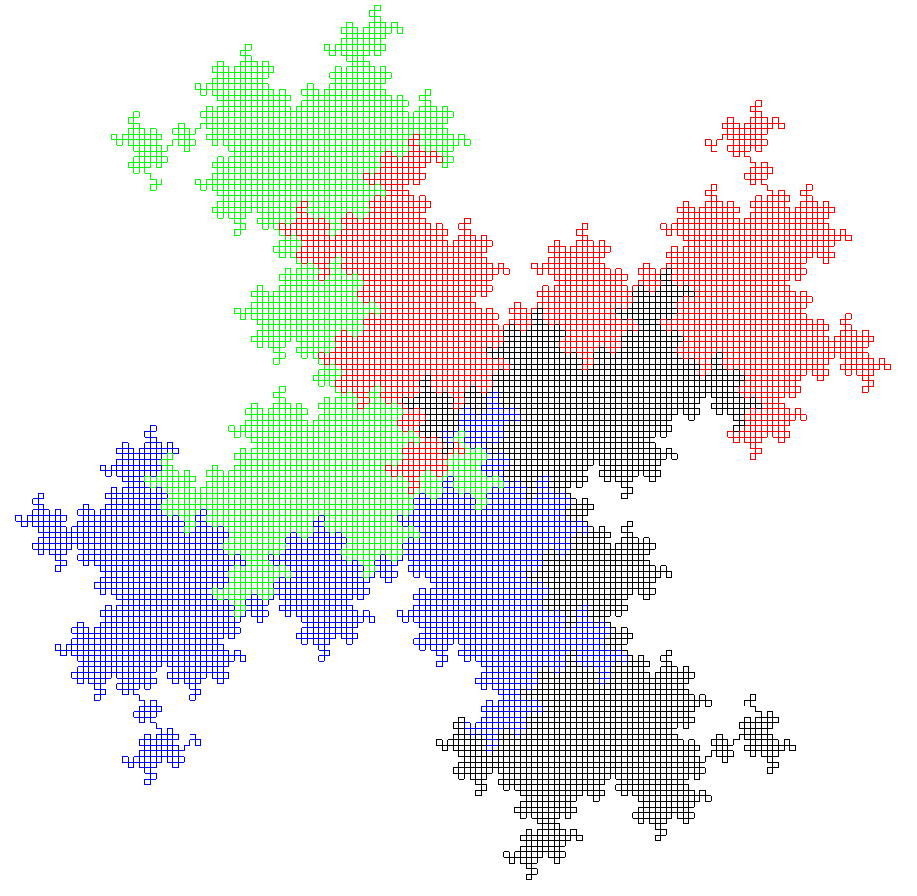}}
\end{center}
\else
\verb+{see pdf for image}+
\fi
\caption{\label{fig:r17-q-carousel}
The 2-symmetric arrangement of curves \CID{R17-11} on the square grid (left)
and the 4-symmetric arrangement of 2-symmetric curves \CID{R17-6} (right).}
\end{figure}
%
%%%%%%%%%%%%%%%%%%%%%%%%%%%

For curves on the square grid
these arrangements have a 2-fold symmetry
and a 4-fold symmetry for curves with 2-fold symmetry,
as shown in Figure~\ref{fig:r17-q-carousel}
for the curves \CID{R17-11} and \CID{R17-6}.
% Cf. Figures~\Ref{fig:R17-11} and \Ref{fig:R17-6}.

%
%
%%%%%%%%%%%%%%%%%%%%%%%%%%%
% stringsubst 3 [L]+[_L]+[_L]+[_L] _ _ [ [ ] ] L L+R+L-R+L-R-L-R+L-R+L+R-L R R+L-R-L+R-L+R+L+R-L+R-L-R + + - - | tail -1 |  ./bin 4 3 0 > tmp-pic.tex && make dotex # R13-rl
%
% stringsubst 3 [R]+[_R]+[_R]+[_R] _ _ [ [ ] ] L L+R+L-R+L-R-L-R+L-R+L+R-L R R+L-R-L+R-L+R+L+R-L+R-L-R + + - - | tail -1 |  ./bin 4 3 0 > tmp-pic.tex && make dotex # R13-rl
%
\begin{figure}[h!tbp]
\ifpdf
\begin{center}
{\includegraphics*[width=52mm, viewport={60 310 490 740}]{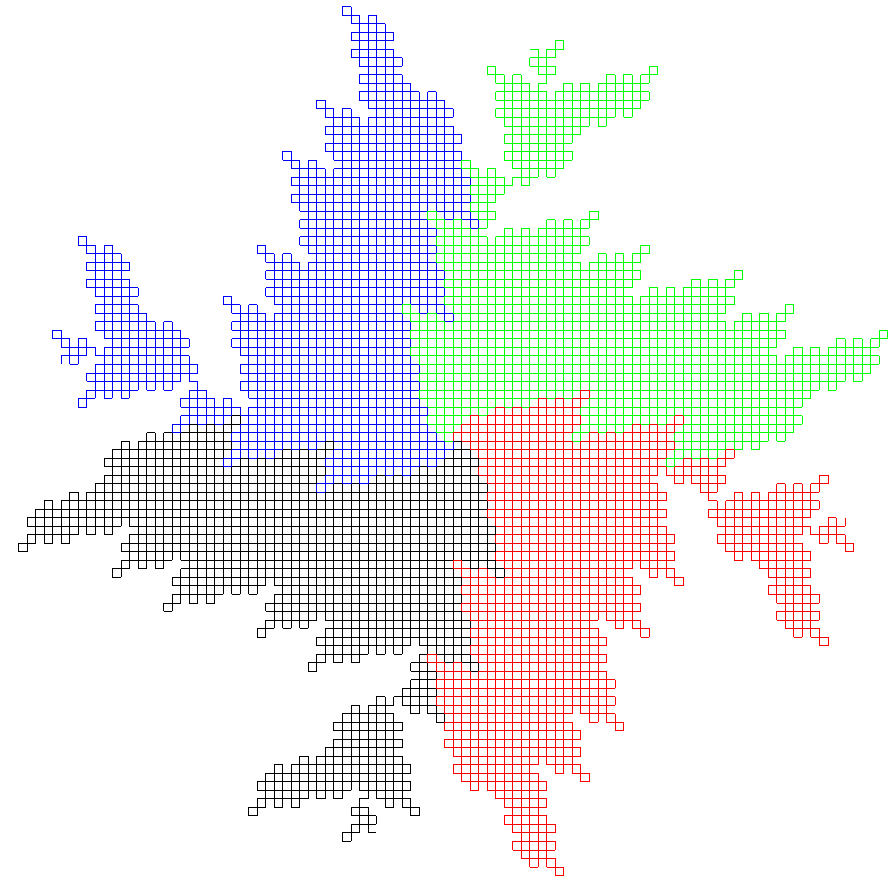}}%
\qquad% layout
{\includegraphics*[width=52mm, viewport={60 310 490 740}]{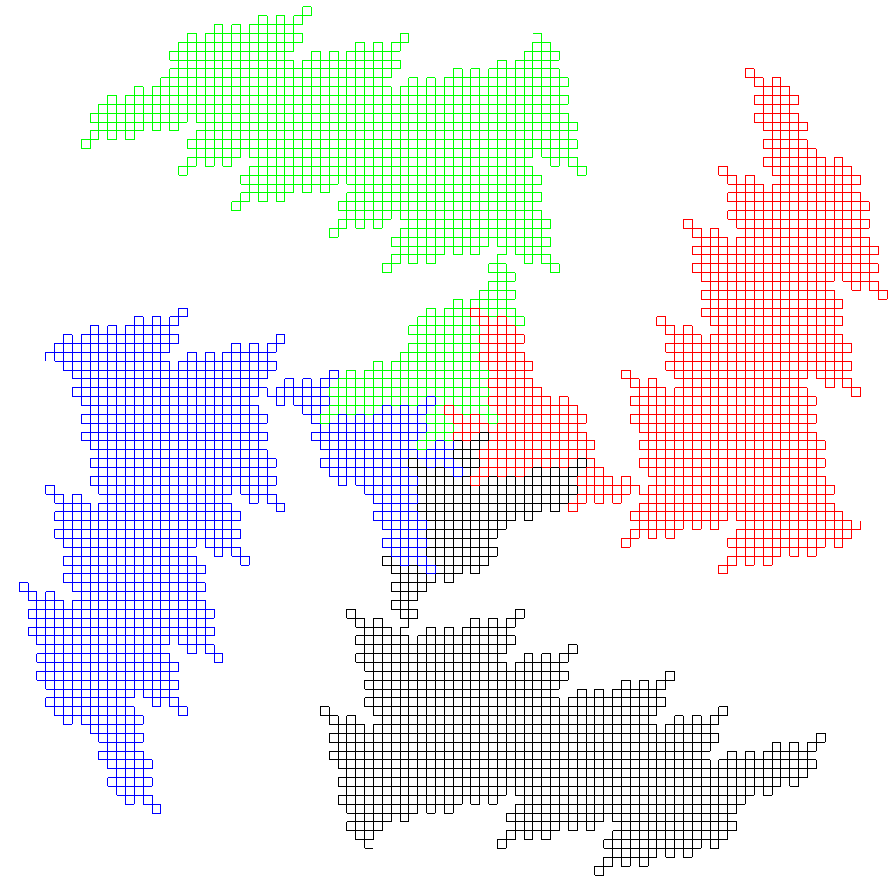}}
\end{center}
\else
\verb+{see pdf for image}+
\fi
\caption{\label{fig:r13-rl-carousel}
The two 4-symmetric arrangements of curves of order 13 on the square grid
generated by an L-system with two non-constant letters.}
\end{figure}
%
%%%%%%%%%%%%%%%%%%%%%%%%%%%

Unsymmetric curves with arrangements of 4-fold symmetry exist.
The arrangements for the curve already shown in Figure~\Ref{fig:r13-rl-tile}
are shown in Figure~\ref{fig:r13-rl-carousel}, the curves on the right
are flipped versions of those on the left.
This happens for curves with L-systems with two non-constant letters,
say \texttt{L} and \texttt{R},
and maps \Lmap{L}{f(L)}, \Lmap{L}{g(L)}
where the letters
 \texttt{L} and \texttt{R} alternate in \texttt{f(L)}
and \texttt{g(R)} is the reversed word \texttt{f(L)}
with letters \texttt{L} and \texttt{R} swapped.
If such a curve has 2-fold symmetry,
a curve of the same shape with simple L-system exists.
%% By identifying both \texttt{L} and \texttt{R} with \texttt{F}.
%% May need to re-arrange the curve into a dragon!

% stringsubst 7 [F]+_[F]+_[F]+_[F]+_[F]+_[F] _ _ [ [ ] ] F F++F--F + + - - | tail -1 | ./bin 6 3 0 > tmp-pic.tex && make dotex # carousel
% stringsubst 6 [F]+_[F]+_[F]+_[F]+_[F]+_[F] _ _ [ [ ] ] F F++F0F--F + + - - 0 0 | tail -1 | ./bin 6 3 0 > tmp-pic.tex && make dotex #
%% alt carousels:
% stringsubst 5 F[++F]+_F[++F]+_F[++F]+_F[++F]+_F[++F]+_F[++F] _ _ [ [ ] ] F F++F--F 0 0 + + - - | tail -1 | ./bin 6 3 0 0 0.1 > tmp-pic.tex && make dotex #
% stringsubst 5 F[++F]+_F[++F]+_F[++F]+_F[++F]+_F[++F]+_F[++F] _ _ [ [ ] ] F F++F0F--F 0 0 + + - - | tail -1 | ./bin 6 3 0 0 0.1 > tmp-pic.tex && make dotex #

% For \ref{fig:R12-4}:
% stringsubst 3 [_F]+[_F]+[_F]+[_F]+[_F]+[_F] _ _ [ [ ] ] F F0F++F++F--F--F0F++F++F--F--F0F 0 0 + + - - | tail -1 | ./bin 6 3 0 > tmp-pic.tex && make dotex # R12-4 # symm-dr
% stringsubst 3 [_F]+[_F]+[_F]+[_F]+[_F]+[_F] _ _ [ [ ] ] F F++F0F++F--F++F0F0F--F++F--F0F--F 0 0 + + - - | tail -1 | ./bin 6 3 0 > tmp-pic.tex && make dotex #  R13-17 (thin-3-curve) # symm-dr

%% NON-carousel, Dekking p.32:
%% tile:
% stringsubst 4 L+_R+_L+_R _ _ L L+R-L-R R L+R+L-R + + - - | tail -1 | ./bin 4 3 0 0 0.15 > tmp-pic.tex && make dotex #

%\clearpage% xxx
%%%%%%%%%%%%%%%%%%%%%%%%%%%%%%%
%%%%%%%%%%%%%%%%%%%%%%%%%%%%%%%
\subsection{Tiles and complex numeration systems}%\label{sect:complnumsys}

%%%%%%%%%%%%%%%%%%%%%%%%%%
%
\begin{figure}[h!tbp]
\ifpdf
\begin{center}
{\includegraphics*[width=70mm, viewport={0 0 968 950}]{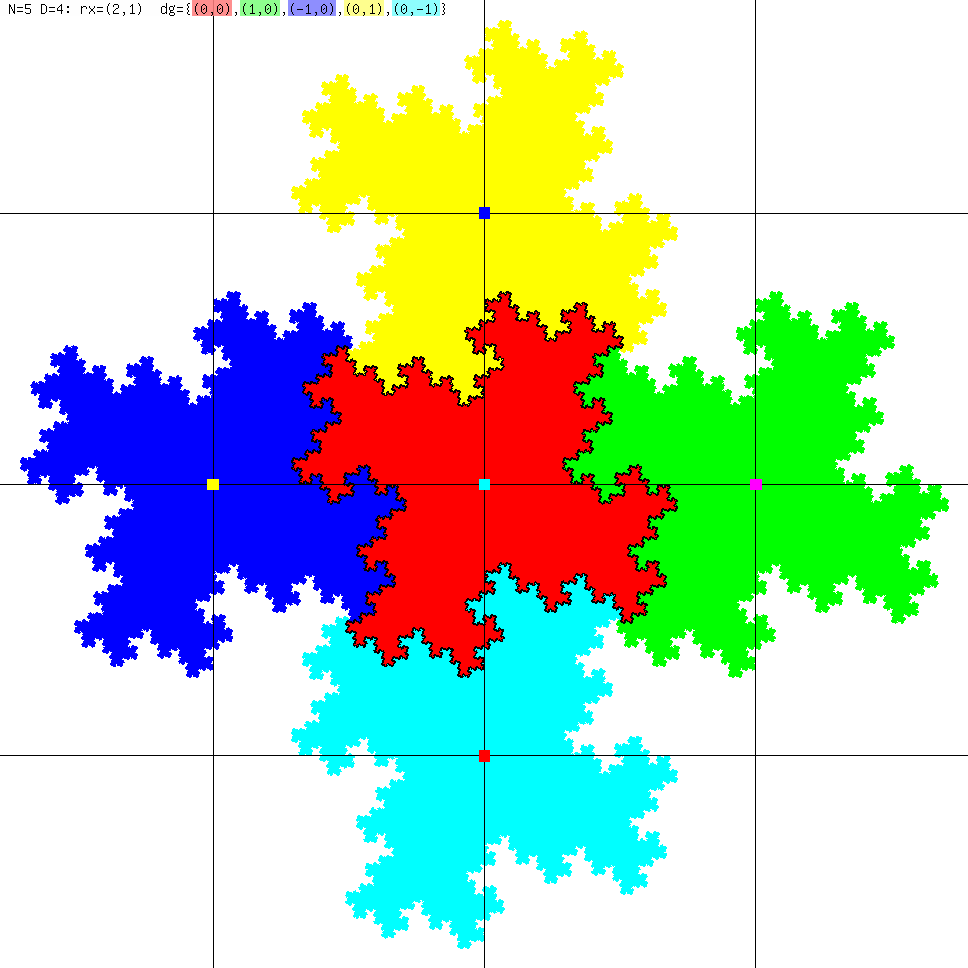}}
\end{center}
\else
\verb+{see pdf for image}+
\fi
\caption{\label{fig:r5-numsys}
The fundamental region for the complex numeration system
with base $2+i$ and digit set $\{0, +1, -1, +i, -i\}$ where $i = \sqrt{-1}$.
The set is scaled up by the factor $\sqrt{5}$ to let the digits (small squares) lie in the five subsets.
}
\end{figure}
%
%%%%%%%%%%%%%%%%%%%%%%%%%%

%%%%%%%%%%%%%%%%%%%%%%%%%%
%
% Render with thick lines: lnth *= 3.0; :
% stringsubst 1 F+_F+_F+_F _ _ F F+F+F-F-F + + - - | tail -1 | sed 's/F/FFF/g; s/+/+t+/g; s/-/-t-/g; s/^/t+/g;' | ./bin 8 3 0 > tmp-pic.tex && make dotex # R5-1
%
% stringsubst 5 F+_F+_F+_F _ _ F F+F+F-F-F + + - - | tail -1 | ./bin 4 3 0 > tmp-pic.tex && make dotex # R5-1
%
%
\begin{figure}[h!tbp]
\ifpdf
\begin{center}
{\includegraphics*[width=50mm, viewport={60 280 520 740}]{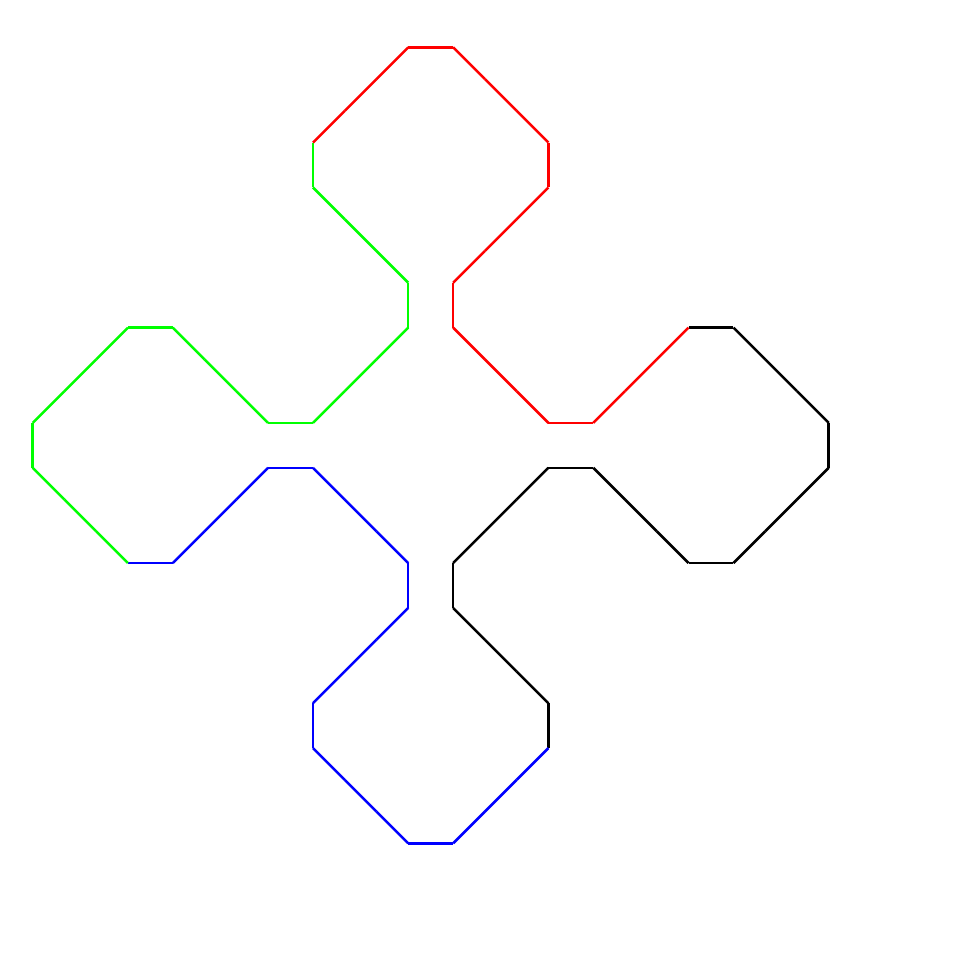}}%
{\includegraphics*[width=60mm, viewport={60 310 490 745}]{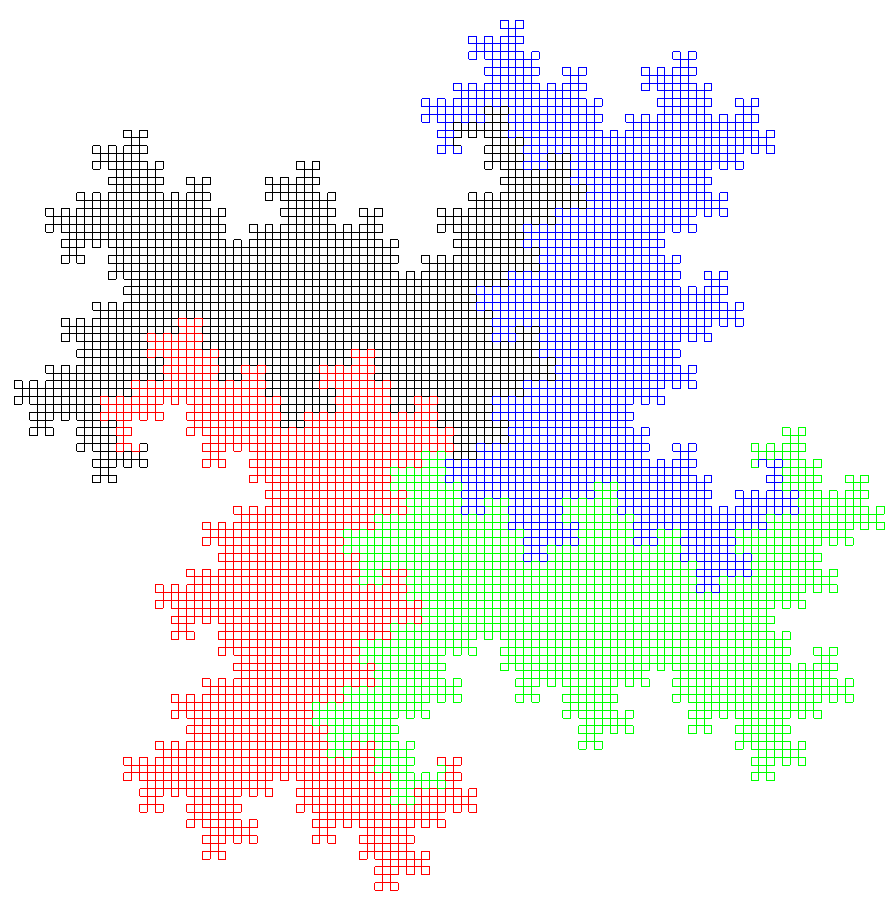}}
\end{center}
\else
\verb+{see pdf for image}+
\fi
\caption{\label{fig:r5-q-1-tile-plus}
The tiles $\Tile{+1}$ (left)  $\Tile{+5}$ (right) of the curve \CID{R5-1} on the square grid.}
\end{figure}
%
%%%%%%%%%%%%%%%%%%%%%%%%%%

A \jjterm{numeration system}
is a pair $(B,\,D)$ where
$B$ is the base and $D$ the set of digits.
The expansion of a number in a numeration system
has the form $\sum_{k=-\infty}^{m}{d_k\,B^k}$ where $d_k\in{}D, m\in\mathbb{N}$.
We call
the part left of the radix point ($\sum_{k=0}^{m}{d_k\,B^k}$)
the \jjterm{integer part}
and the part right of the radix point ($\sum_{k=-\infty}^{-1}{d_k\,B^k}$)
the \jjterm{fractional part}.

The well-known positional numeration systems
are defined by an integer base $B\geq{}2$
and a digit set $D$, usually $D=\{0, 1, \ldots, B-1\}$.
We call a numeration system \jjterm{integral}
if each integer has a unique representation
as $\sum_{k=0}^{m}{d_k\,B^k}$
(a representation without fractional part).
The positional numeration systems
with digits $D=\{0, 1, \ldots, B-1\}$ are integral.
All numbers can be represented in
an integral numeration system \cite[Section~7.4, pp.~257ff]{alg-comb-words}.

We consider complex numeration systems
where both the base $B$ and the digits in $D$
are in general complex.
Let $\omega_k = \exp(2\pi{}i/k)$ where $i=\sqrt{-1}$
be a primitive complex $k$th root of unity,
we will only use the values
$\omega_3 = (-1+i\,\sqrt{3})/2$, $\omega_4 = i$, and $\omega_6 = (1+i\,\sqrt{3})/2$.
The integers of interest are
the Gaussian integers
(numbers $x+\omega_4\,y=x+i\,y$ where $x,y\in{}\mathbb{Z}$,
points on the square grid)
and the Eisenstein integers
(numbers $x+\omega_3\,y$ where $x,y\in{}\mathbb{Z}$,
points on the triangular grid).
We require that both the base and all digits are such integers.

Following Knuth \cite[Section~4.1]{DEK2}
we call the set of numbers of the
form  $\sum_{k=-\infty}^{-1}{d_k\,B^k}$
(that is numbers of the form $0.*$ where $*$ is any sequence of digits)
the \jjterm{fundamental region} of the numeration system.

An example for the square grid is the numeration system with
base $B=2+i$ and digit set D=$\{0, +1, -1, +i, -i\}$
whose fundamental region is shown in Figure~\ref{fig:r5-numsys}.
The set consists of five subsets that are
similar to the whole set:
the subsets shown in red, green, blue, yellow, and cyan
are respectively the subsets of numbers of the form
$0.(0)*$, $0.(+1)*$, $0.(-1)*$, $0.(+i)*$ $0.(-i)*$.
The positions of the values of the digits
are indicated by the small squares inside the subsets
(for this purpose the set is shown scaled up by $\sqrt{\abs(B)}=\sqrt{5}$).

The tiles of our curves correspond to complex numeration systems.
The shape of the tile $\Tile{+5}$ of the curve \CID{R5-1} on the square grid
(with L-system \Lmap{F}{F+F+F-F-F}) shown at the right of
Figure~\ref{fig:r5-q-1-tile-plus} matches (up to rotation)
the form of the fundamental region shown in \ref{fig:r5-numsys}.
The values of the digits correspond to the positions of
the squares surrounded by the tile $\Tile{+1}$
in the left of Figure~\ref{fig:r5-q-1-tile-plus}.

For every tile we can find a complex numeration system
and identify the tile with the fundamental region as follows.
Let $R$ be the order of the underlying curve.
There exists a complex base $B$ such that $\abs(B)^2=R$
and a set of digits $\{d_0, d_1, \ldots, d_{R-1}\}$
corresponding to the squares surrounded by $\Tile{+1}$
such that the shape of the tile ($\Tile{+\infty}$)
matches the fundamental region of the numeration system.
We take the digit $0$ for the center of the tile where possible:
on the square grid the orders are of the form $R=4k+1$ and there always is a central digit,
on the triangular grid there is a central digit for orders of the form $R=3k+1$,
otherwise (orders $R=3k$) we choose one of the three digits around the center.
The subdivision into sets for numbers of the forms $0.d_j*$ for $0\leq{}j<R$
reproduces the tiling of the fundamental region
(and of $\Tile{+\infty}$) into $R$ smaller rotated copies of itself.
The same holds for $\Tile{-1}$ and $\Tile{-\infty}$,
giving a different numeration system in general.

The digits form a \jjterm{complete residue system}
($D$ contains exactly one representative of each coset modulo $B$).
Note that in general there is more than one choice for the base.

%%%%%%%%%%%%%%%%%%%%%%%%%%
%
% Render with thick lines: lnth *= 3.0; :
%
% stringsubst 1 F+_F+_F+_F _ _ F F+F+F-F+F+F-F+F-F-F+F-F-F + + - - | tail -1 | sed 's/F/FFF/g; s/+/+t+/g; s/-/-t-/g; s/^/t+/g;' | ./bin 8 3 0 > tmp-pic.tex && make dotex # R13-1
%
% stringsubst 1 F+_F+_F+_F _ _ F F+F-F-F-F+F+F-F-F+F+F+F-F + + - - | tail -1 | sed 's/F/FFF/g; s/+/+t+/g; s/-/-t-/g; s/^/t+/g;' | ./bin 8 3 0 > tmp-pic.tex && make dotex # R13-4
%
\begin{figure}[h!tbp]
\ifpdf
\begin{center}
{\includegraphics*[width=35mm, viewport={70 330 520 740}]{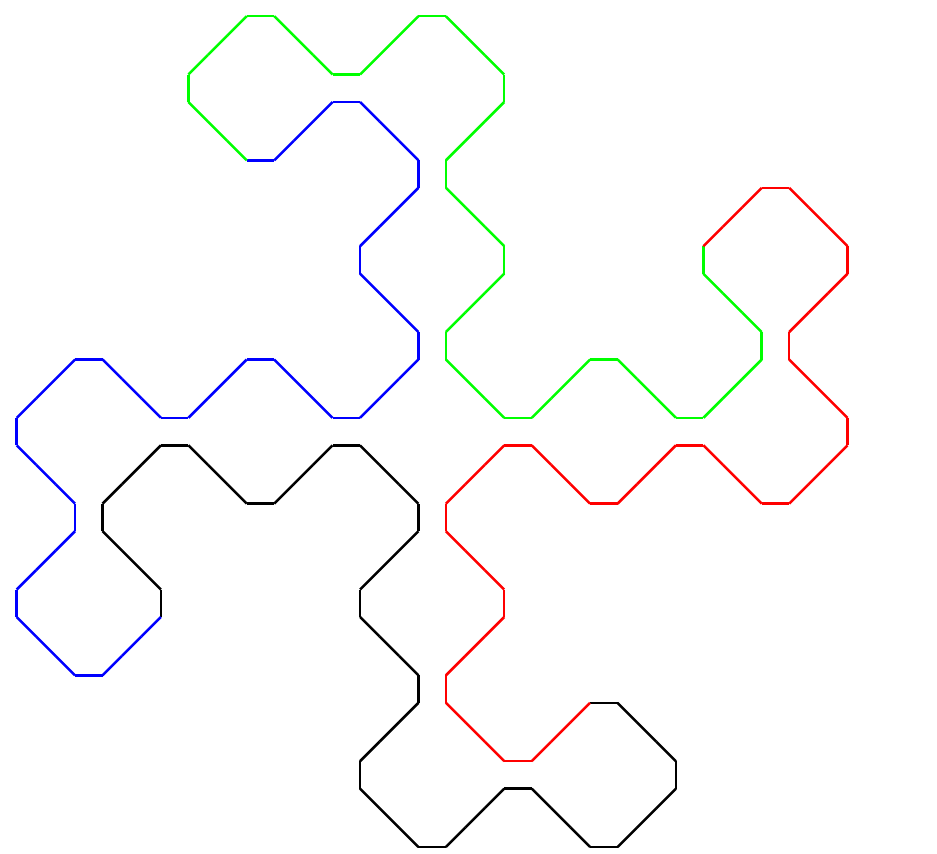}}%
{\includegraphics*[width=35mm, viewport={10 310 480 740}]{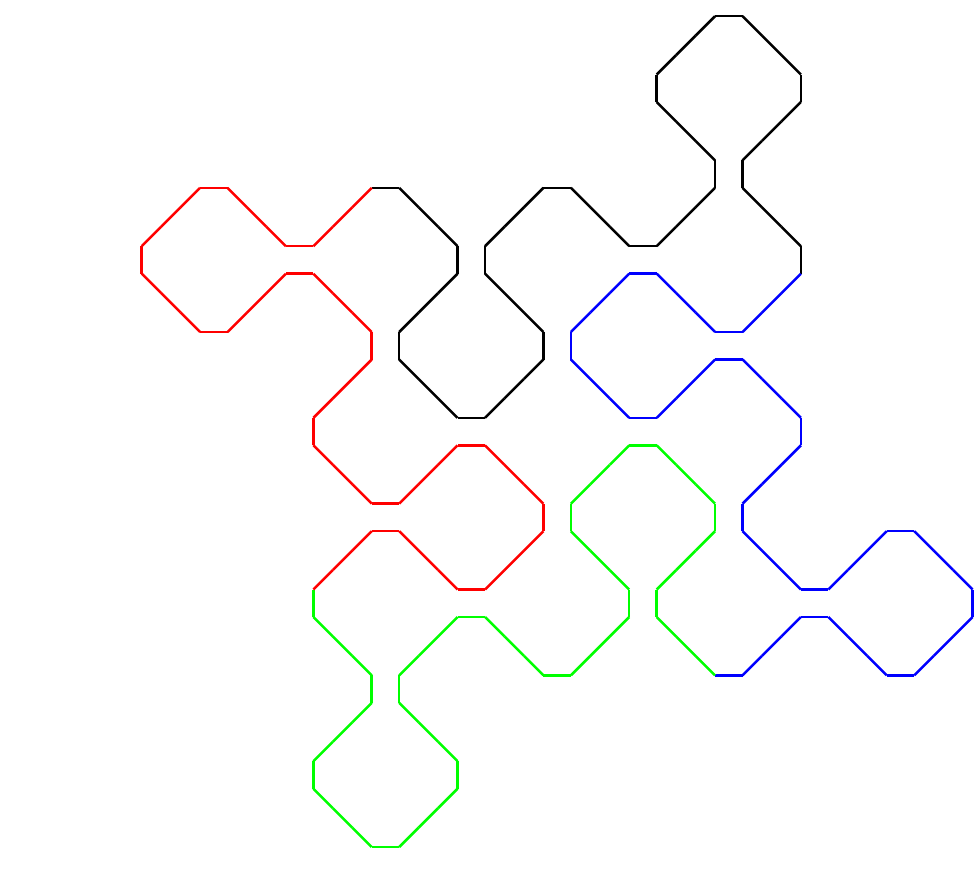}}
\end{center}
\else
\verb+{see pdf for image}+
\fi
\caption{\label{fig:r13-q-4-tile-plus}
The tiles $\Tile{+1}$ of the curves \CID{R13-1} (left) and \CID{R13-4} (right) on the square grid.}
\end{figure}
%
%%%%%%%%%%%%%%%%%%%%%%%%%%

%%%%%%%%%%%%%%%%%%%%%%%%%%
% stringsubst 3 F+_F+_F+_F _ _ F F+F+F-F+F+F-F+F-F-F+F-F-F + + - - | tail -1 | sed 's/^/+/g;' | ./bin 4 3 0 > tmp-pic.tex && make dotex # R13-1
%
\begin{figure}[h!tbp]
\ifpdf
\begin{center}
{\includegraphics*[width=62mm, viewport={10 15 960 950}]{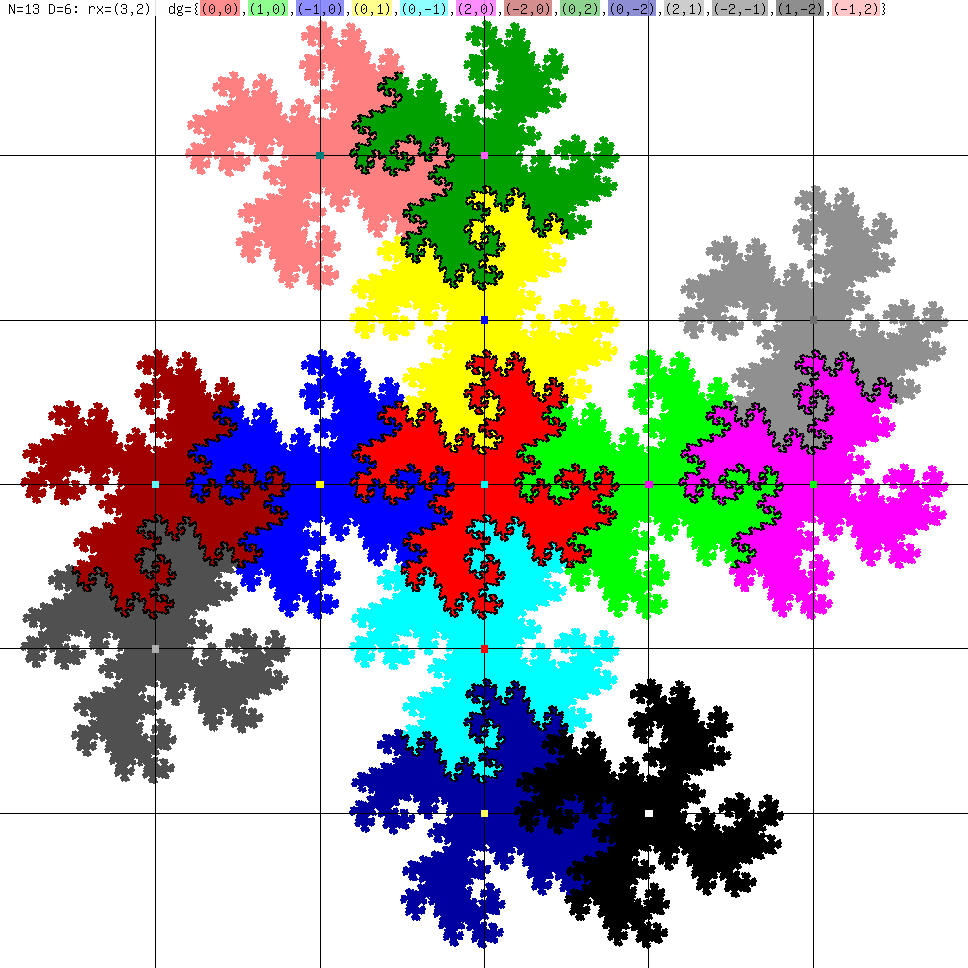}}%
{\includegraphics*[width=62mm, viewport={65 315 490 740}]{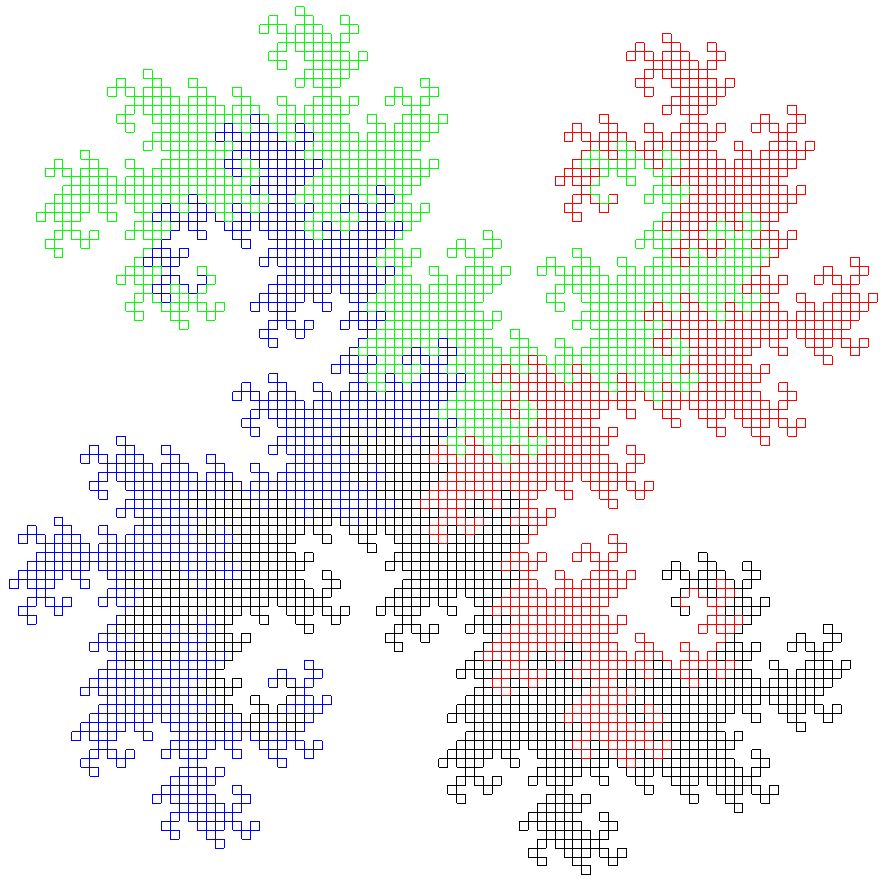}}
\end{center}
\else
\verb+{see pdf for image}+
\fi
\caption{\label{fig:r13-q-1-numsys}
Fundamental region of a numeration system (left, scaled up by $\sqrt{13}$) and
the tile $\Tile{+3}$ of the curve \CID{R13-1} on the square grid (right).}
\end{figure}
%
%%%%%%%%%%%%%%%%%%%%%%%%%%

%%%%%%%%%%%%%%%%%%%%%%%%%%
% stringsubst 3 F+_F+_F+_F _ _ F F+F-F-F-F+F+F-F-F+F+F+F-F + + - - | tail -1 | sed 's/^/+/g;' | ./bin 4 3 0 > tmp-pic.tex && make dotex # R13-4
%
\begin{figure}[h!tbp]
\ifpdf
\begin{center}
{\includegraphics*[width=62mm, viewport={10 15 960 950}]{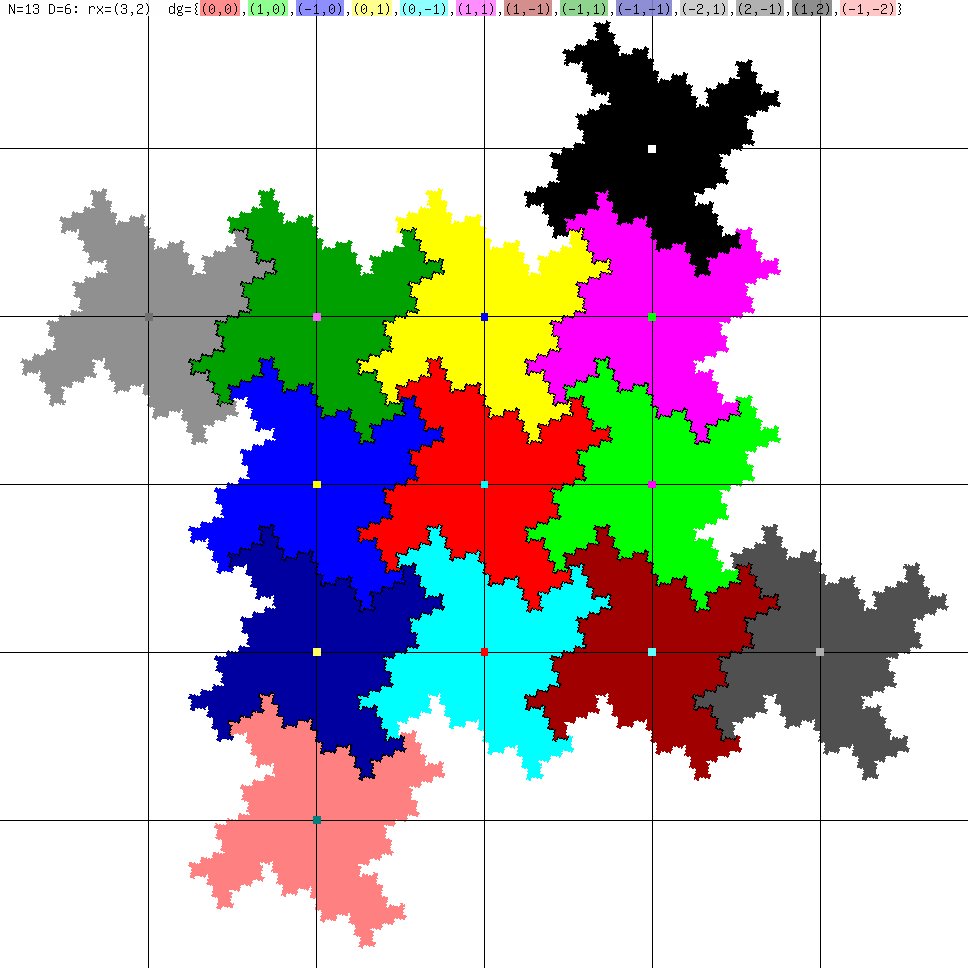}}%
{\includegraphics*[width=62mm, viewport={65 315 490 740}]{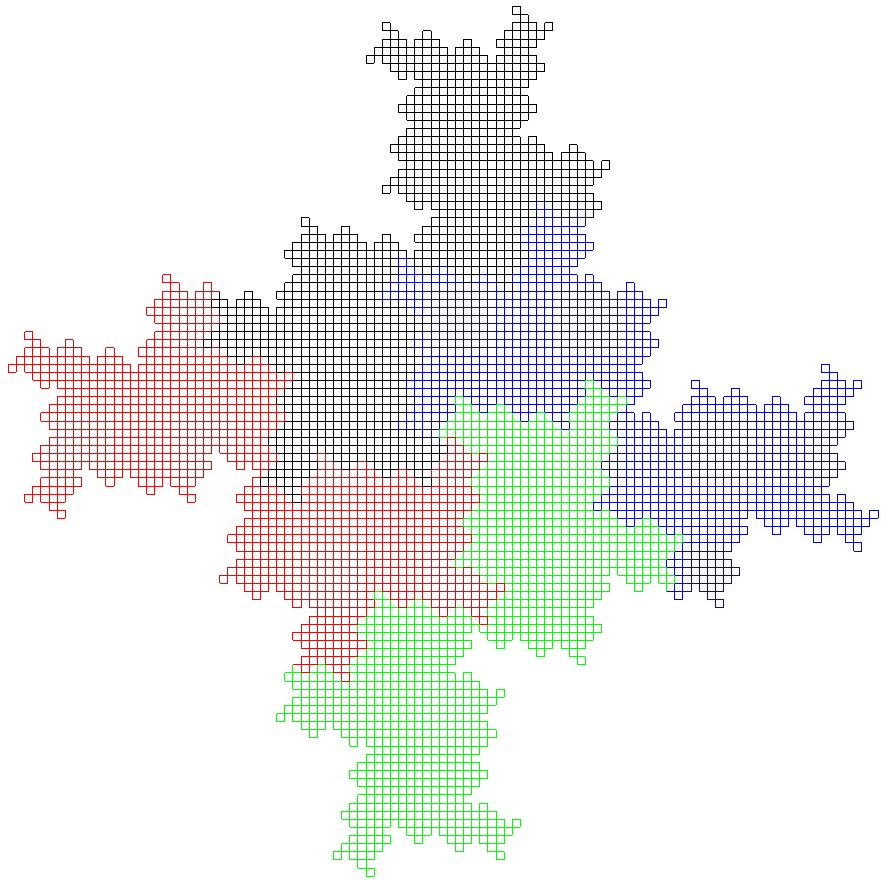}}
\end{center}
\else
\verb+{see pdf for image}+
\fi
\caption{\label{fig:r13-q-4-numsys}
Fundamental region of a numeration system (left, scaled up by $\sqrt{13}$) and
the tile $\Tile{+3}$ of the curve \CID{R13-4} on the square grid (right).}
\end{figure}
%
%%%%%%%%%%%%%%%%%%%%%%%%%%

We give two more examples for the square grid.

The tiles $\Tile{+1}$ for the order-13 curves
with maps
 \Lmap{F}{F+F+F-F+F+F-F+F-F-F+F-F-F}
and
 \Lmap{F}{F+F-F-F-F+F+F-F-F+F+F+F-F}
are shown in Figure~\ref{fig:r13-q-4-tile-plus}.
The 13 squares inside the tiles give the digit sets,
$D=\{0,\, \pm{}1,\, \pm{}i,\, \pm{}2,\, \pm{}2i,\,
  +2+i,\, -2-i,\, +1-2i,\, -1+2i\}$
for the curve \CID{R13-1}
and
$D=\{0,\, \pm{}1,\, \pm{}i,\, \pm{}1+\pm{}i,\,
  -2+i,\, +2-i,\, +1+2i,\, -1-2i\}$
for the curve \CID{R13-4}.
The base $B=3+2i$ can be chosen for both numeration systems.
The fundamental regions of the numeration systems
are shown next to the tiles $\Tile{+3}$
in Figures~\ref{fig:r13-q-1-numsys} and \ref{fig:r13-q-4-numsys}.

%%%%%%%%%%%%%%%%%%%%%%%%%%
%% with  lnth *= 2.0;  // thicker lines
% stringsubst 7 +_F++_F++_F  F F++F--F  _ _  + +  - - | tail -1 | ./bin 6 2 0 > tmp-pic.tex && make dotex # terdragon tile
%
\begin{figure}[h!tbp]
\ifpdf
\begin{center}
{\includegraphics*[width=62mm, viewport={0 0 968 950}]{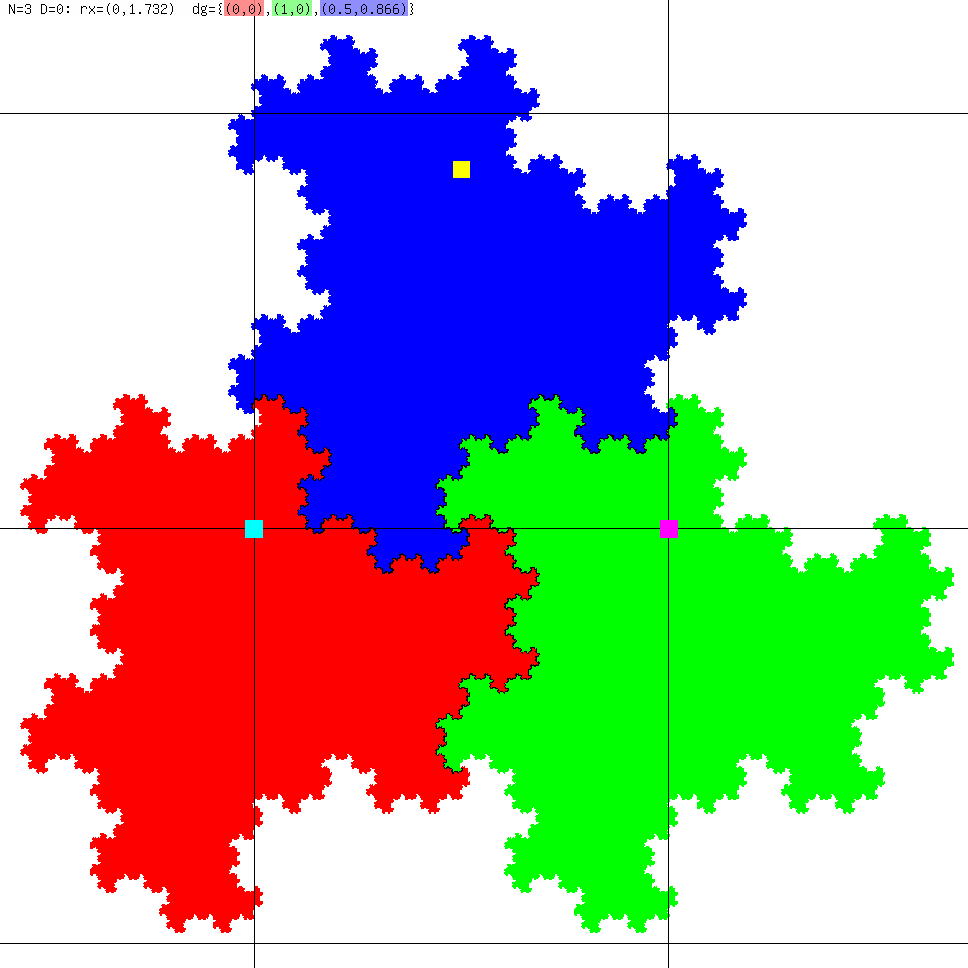}}%
{\includegraphics*[width=65mm, viewport={60 320 500 740}]{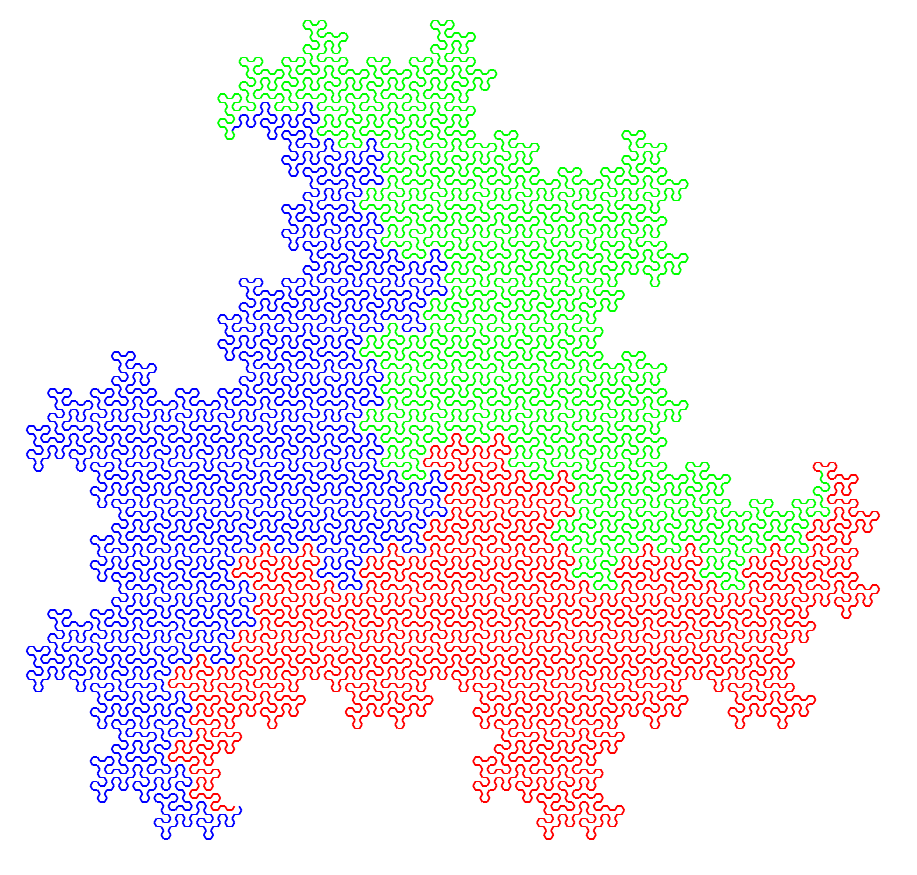}}
\end{center}
\else
\verb+{see pdf for image}+
\fi
\caption{\label{fig:terdragon-numsys-and-tile}
The fundamental region  (left) for the complex numeration system
with base $\sqrt{3}\,i$ and digit set $\{0, +1, \omega_6\}$ where $\omega_6 = \exp(+2\pi{}i/6) = (1+i\,\sqrt{3})/2$
(the set is scaled up by the factor $\sqrt{3}$ to let the digits (small squares) lie in the three subsets)
and the tile $\Tile{+7}$ for the terdragon (right).
}
\end{figure}
%
%%%%%%%%%%%%%%%%%%%%%%%%%%

The fundamental region for the numeration system with
base $B=\sqrt{3}\,i$ and digit set D=$\{0, +1, \omega_6\}$
(the number $\omega_6$ is an Eisenstein integer)
is shown in Figure~\ref{fig:terdragon-numsys-and-tile}.
The set consists of the three subsets
shown in red, green, and blue,
respectively for numbers of the form
$0.(0)*$, $0.(1)*$, and $0.(\omega_6)*$.

Helmberg \cite[Section~4, pp.~374ff]{helmberg-crab-fractal}
shows that the tile for the curve in Figure~\ref{fig:r4-crab-dragon}
corresponds to the numeration systems with $B=-2$
and digits $D=\{0, 1, \omega_3, \omega_3^2\}$.
%% file eisenstein-tile.png

The numeration systems mentioned so far are all integral.
A sufficient condition for that is
that a non-empty neighborhood of the origin (the digit $0$)
is eventually contained in the tile.
%

%This always holds if all polygons of the grid\xxx{Danger!}
%sharing a corner with the polygon at the origin
%are contained in the tile.
%%
%The are four such squares on the square grid and
%six such triangles on the triangular grid.
%%
%By construction this is true for all
%curves on the square grid.
%%

Note that for integral numeration systems there
is neither a separation into real and imaginary parts
nor is a sign needed.

For the triangular grid, tiles exists which
do not contain a neighborhood of zero,
we now give one example.
%

%%%%%%%%%%%%%%%%%%%%%%%%%%
%
% Render with thick lines: lnth *= 3.0; :
%
% stringsubst 1 _F-_F-_F _ _ F F0F+F0F-F-F+F 0 0 + + - - | tail -1 | ./bin 3 3 1 > tmp-pic.tex && make dotex # R07-t-1 tile-minus
%
% stringsubst 1 R_F-_F-_F R R _ _ F F0F+F0F-F-F+F 0 0 + + - - | tail -1 | sed 's/+/++++/g; s/-/----/g; s/F/++F--F/g; s/R/-/;' | ./bin 12 3 0 > tmp-pic.tex && make dotex # R07-t-1 tile-minus
%
\begin{figure}[h!tbp]
\ifpdf
\begin{center}
{\includegraphics*[width=41mm, viewport={150 450 450 730}]{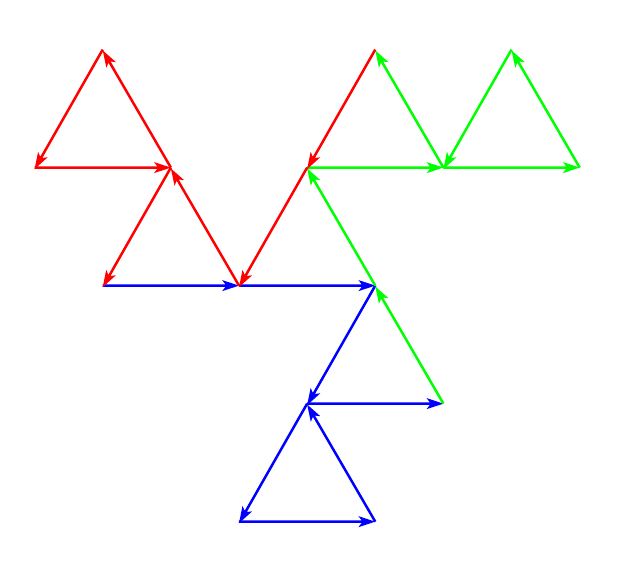}}%
{\includegraphics*[width=41mm, viewport={50 390 470 730}]{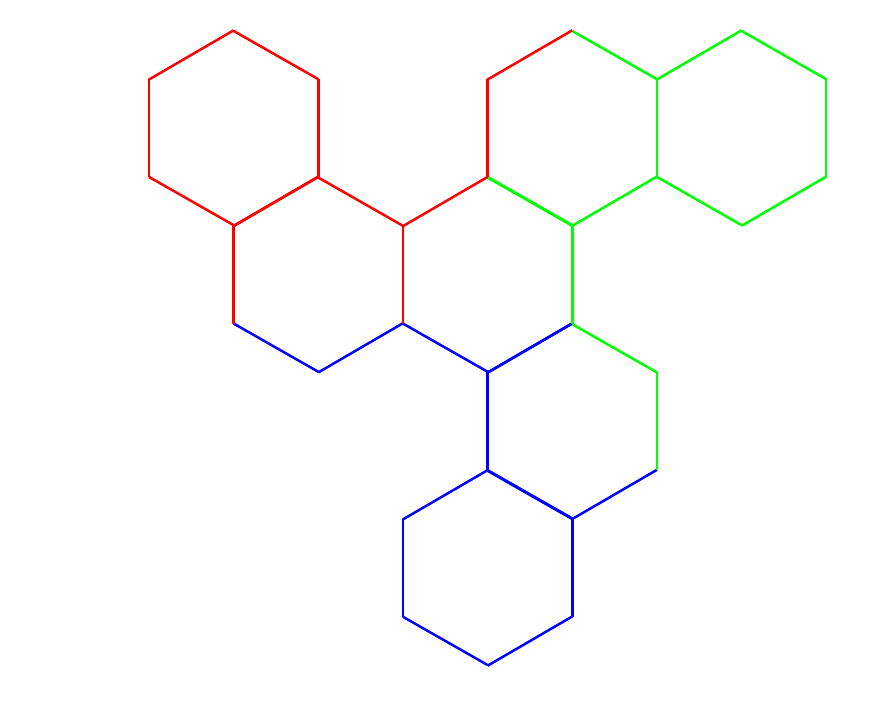}}%
\end{center}
\else
\verb+{see pdf for image}+
\fi
\caption{\label{fig:r7-t-1-tile-minus}
The tile $\Tile{-1}$ of the curve \CID{R7-1} (left),
the same rendered as in Figure~\ref{fig:r13-t-15-tile-6fold-symm} (right).}
\end{figure}
%
%%%%%%%%%%%%%%%%%%%%%%%%%%
%% alt hex rendering:
% stringsubst 4 _F-_F-_F _ _ F F+F-F-F0F+F0F 0 0 - - + + | tail -1 | sed 's/+/++/g; s/-/--/g; s/F/+F-F/g;' | tr +- -+ | ./bin 6 2 0 0 0.15 > tmp-pic.tex && make dotex

%%%%%%%%%%%%%%%%%%%%%%%%%%
% stringsubst 5 _F-_F-_F _ _ F F0F+F0F-F-F+F 0 0 + + - - | tail -1 | ./bin 3 3 0 > tmp-pic.tex && make dotex # R07-t-1 tile-minus
%
\begin{figure}[h!tbp]
\ifpdf
\begin{center}
{\includegraphics*[width=65mm, viewport={10 30 980 950}]{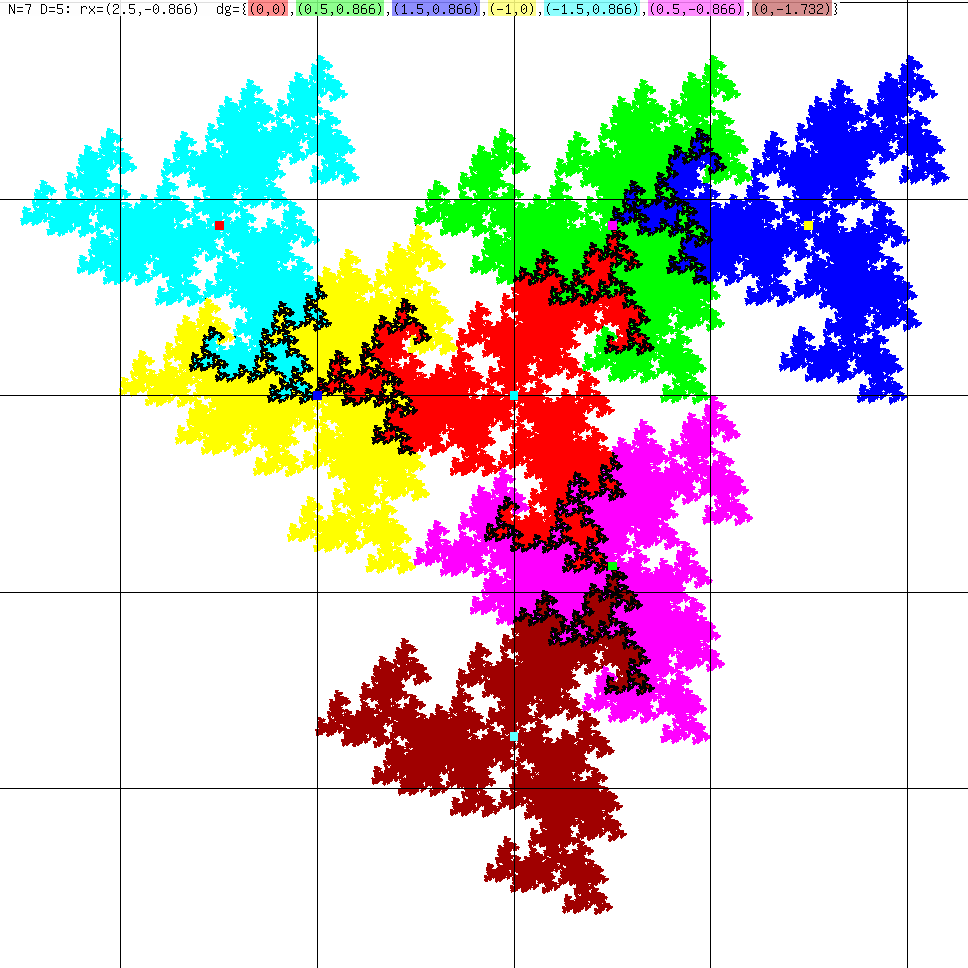}}%
{\includegraphics*[width=64mm, viewport={55 320 490 740}]{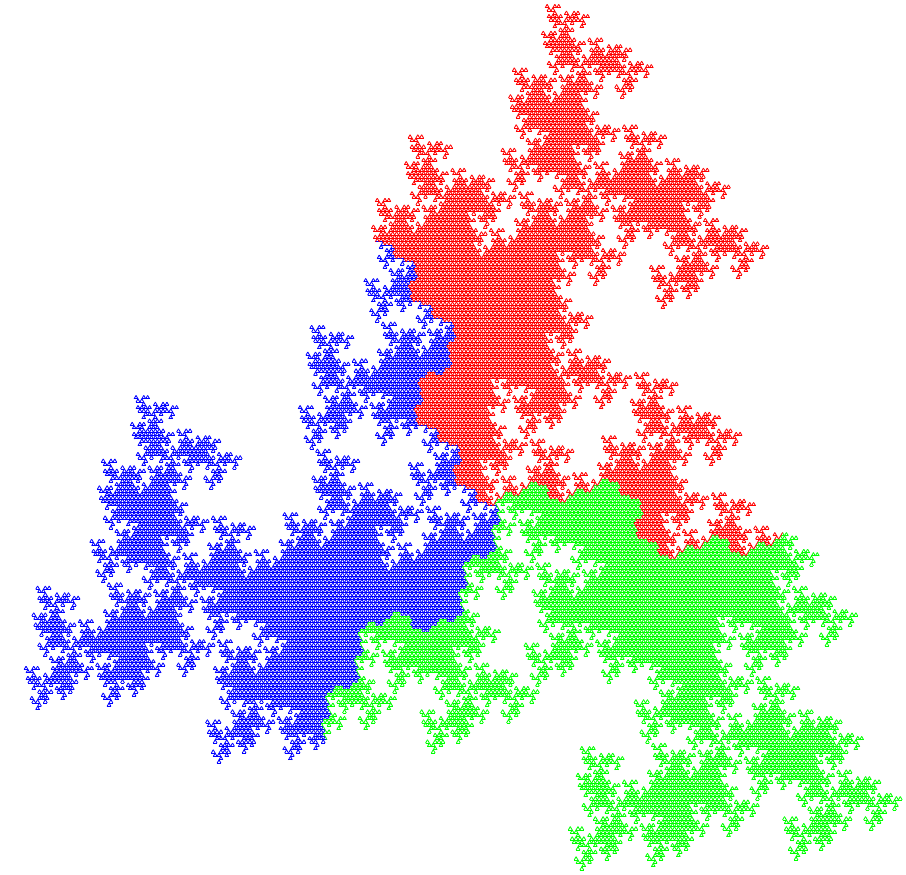}}
\end{center}
\else
\verb+{see pdf for image}+
\fi
\caption{\label{fig:r7-numsys-tile}
Fundamental region of a numeration system (left, scaled up by $\sqrt{7}$) and
fifth iterate $\Tile{-5}$ of the tile for the curve \CID{R7-1} (right).
}
\end{figure}
%
%%%%%%%%%%%%%%%%%%%%%%%%%%

A complex numeration system corresponding to the
tile $\Tile{-1}$ of the curve \CID{R7-1} on the triangular grid
can be obtained as follows.
The $R=7$ digits are points on the grid that are enclosed by the tile $\Tile{-1}$,
Figure~\ref{fig:r7-t-1-tile-minus} gives two renderings of the tile.
We choose the center as digit $0$,
the other digits are
 $\omega_6$ and $1+\omega_6$ (upper right arm),
 $-1$ and $-1+\omega_3$ (upper left arm),
and
$-\omega_3$ and $-\omega_3-\omega_6$ (lower arm).

The choice of base $B=$
% $2\,\omega_3+\omega_6=$ $\frac{1}{2}\,(-1+i\,3\,\sqrt{3})$
%% this one is directly suggested by the blue curve
 $2 - \omega_3=$ $\frac{1}{2}\,(5-i\,\sqrt{3})$
($\abs(B)^2=7$ and $B$ lies on the triangular lattice)
leads to the fundamental region
as shown left in Figure~\ref{fig:r7-numsys-tile},
after rotation it matches the tile $\Tile{-5}$ on the right.
% Check with Pari
% ? w(k)=exp(2*Pi*I/k)
% (k)->exp(2*Pi*I/k)
% ? 2*w(3)+w(6)
% -0.500000000000000 + 2.59807621135332*I
% ? 1./2*(-1+I*3*sqrt(3))
% -0.500000000000000 + 2.59807621135332*I
%
This tile (and fundamental region)
does not contain a neighborhood of zero.
%
% xxx identify numeration system with negated base

Tiles with 6-fold symmetry necessarily contain all
digits adjacent to $0$, so the corresponding
numeration systems are always integral.
%
%
% file r13-b-1-tile.png
% file radix-dragon.inc.cc:
% //Duda(13,+2);  Cx dg[] = { 0, +1, w6, w3, -1, -w6, -w3,  2* +1, 2* w6, 2* w3, 2* -1, 2* -w6, 2*  -w3 };
% // OK are: -7, -5, -2, +2, +5, +7
For example,
the (filled out) shape shown in Figure~\ref{fig:r13-t-15-plus-tiling},
after rotation such that 5 of the smaller tiles are on the real axis,
corresponds to the complex base $1+i\,\sqrt{12}$
and digit set $D=\{0, \pm{}1, \pm{}2, \pm{}\omega_6, \pm{}2\omega_6, \pm{}\omega_3, \pm{}2\omega_3\}$,
an integral numeration system.

%%%%%%%%%%%%%%%%%%%%%%%%%%
%
% Render with thick lines: lnth *= 3.0; :
% stringsubst 1 R_F-_F-_F R R _ '' F F+F0F+F-F+F0F0F+F-F+F-F+F0F0F0F-F+F-F+F-F0F0F-F+F-F0F-F 0 0 + + - - | tail -1 | sed 's/+/++++/g; s/-/----/g; s/F/++F--F/g; s/R/-/;' | ./bin 12 3 0 > tmp-pic.tex && make dotex # R28-10790 # symm-dr
% stringsubst 2 R_F-_F-_F R R _ '' F F+F0F+F-F+F0F0F+F-F+F-F+F0F0F0F-F+F-F+F-F0F0F-F+F-F0F-F 0 0 + + - - | tail -1 | sed 's/+/++++/g; s/-/----/g; s/F/++F--F/g; s/R/-/;' | ./bin 12 3 0 > tmp-pic.tex && make dotex # R28-10790 # symm-dr
%
\begin{figure}[h!tbp]
\ifpdf
\begin{center}
{\includegraphics*[width=60mm, viewport={70 330 480 740}]{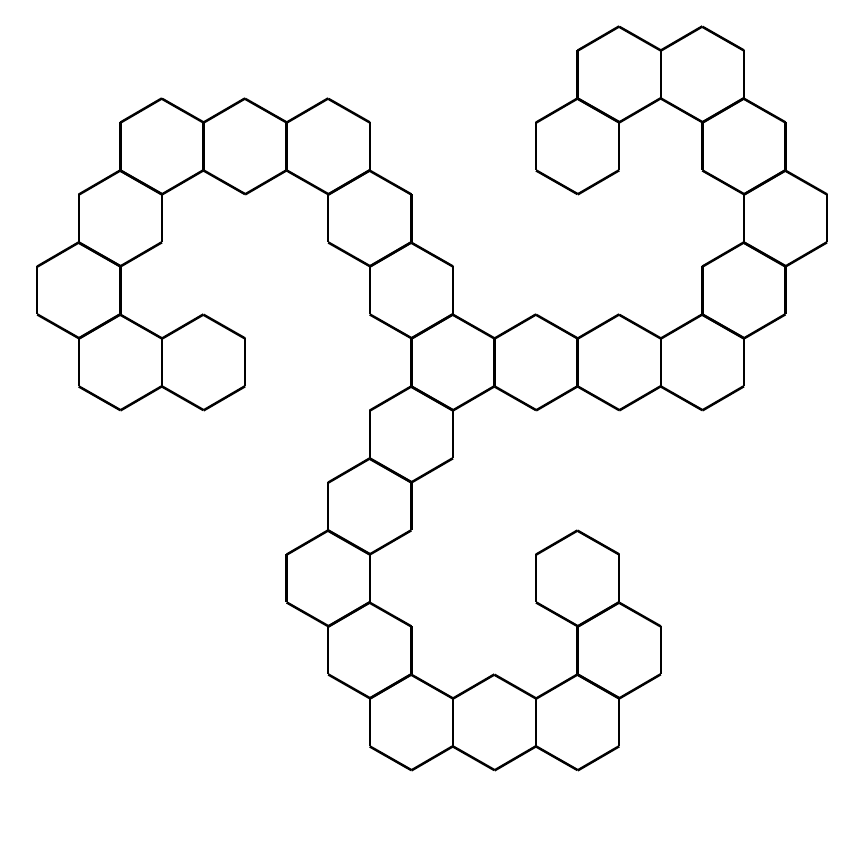}}%
{\includegraphics*[width=67mm, viewport={60 285 500 740}]{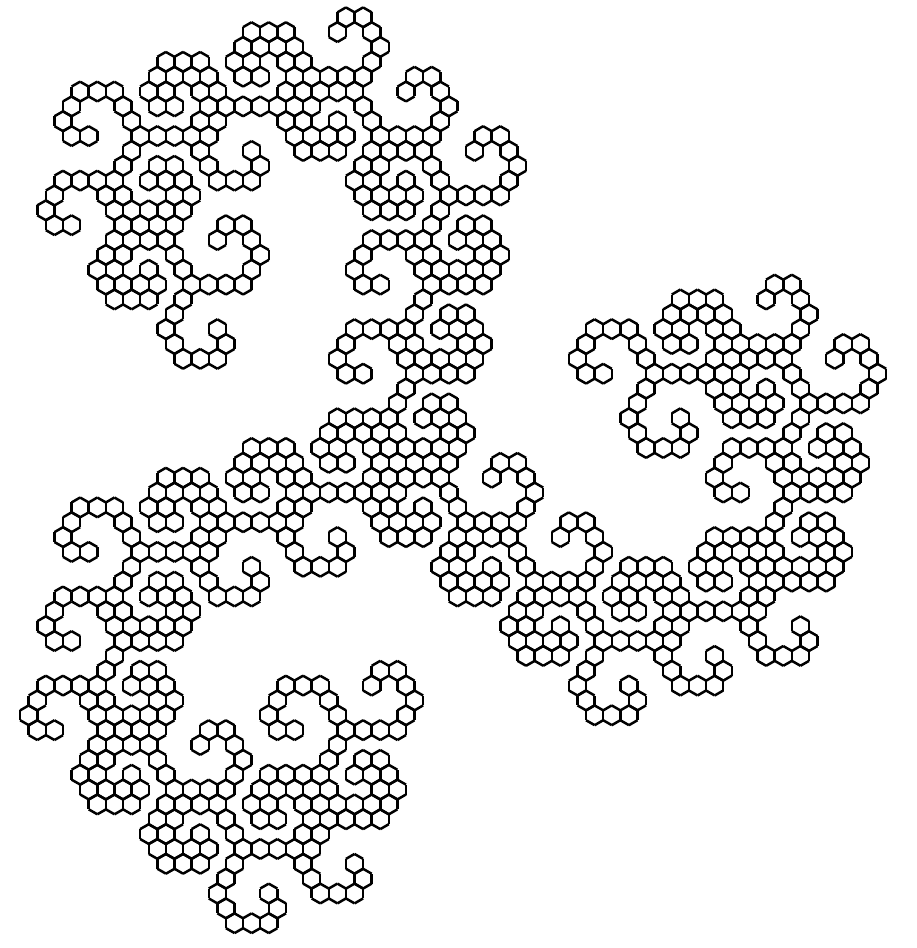}}
\end{center}
\else
\verb+{see pdf for image}+
\fi
\caption{\label{fig:r28-t-10790-tile-plus}
The tiles $\Tile{+1}$ and $\Tile{+2}$  of the curve \CID{R28-10790} on the triangular grid.}
\end{figure}
%
%%%%%%%%%%%%%%%%%%%%%%%%%%%

The tiles $\Tile{+1}$ and $\Tile{+2}$ for a curve of order 28 on the triangular grid
are shown in Figure~\ref{fig:r28-t-10790-tile-plus}.
The central hexagon is not surrounded by others on $\Tile{+1}$,
but in $\Tile{+2}$ and all tiles $\Tile{+k}$ for $k\geq{}2$.
The following \emph{appears} to hold true.
If there is a $k_0$ such that for all tiles $\Tile{k}$
with $k\geq{}k_0$ the central polygon is surrounded
by other polygons of the tile then $k_0\leq{}2$.
Here we dropped the sign of the subscript in $\Tile{k}$,
being arbitrary but fixed.

The subdivision of the tiles into smaller rotated copies of itself
corresponds to an iterated function system (IFS)
with $R$ affine maps of the form
$M_j(\vec{v}) = \vec{t}_j + \alpha\,T\,\vec{v}$ for $0\leq{}j<R$
where $\alpha$ is a scalar with
$\mathrm{abs}(\alpha)^2=1/R$, and $T$ an orthogonal matrix.
Neither $\alpha$ nor $T$ do depend on $j$,
the linear part $\alpha\,T$ of the map is a scaled rotation.

\subsection{Curves and tiles on the tri-hexagonal grid}\label{sect:tri-hex}
%% https://en.wikipedia.org/wiki/Trihexagonal_tiling

%%%%%%%%%%%%%%%%%%%%%%%%%%
% stringsubst 2 F+F+F+F+F+F  F F+F+F+F--F--F+F   + + - - | tail -1 | ./bin 6 3 0 > tmp-pic.tex && make dotex # R7-1
\begin{figure}[h!tbp]
\ifpdf
\vspace*{1mm}% layout
\begin{center}
{\includegraphics*[width=70mm, viewport={150 470 400 590}]{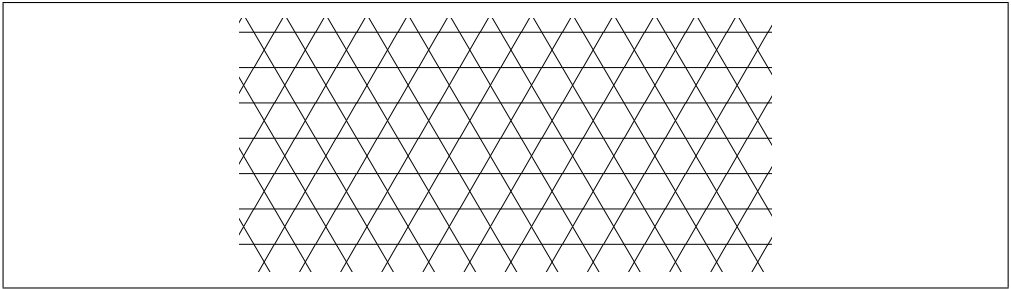}}
\end{center}
\vspace*{1mm}% layout
\else
\verb+{see pdf for image}+
\fi
\caption{\label{fig:tri-hex-grid}
%The tri-hexagonal (3.6.3.6) grid.}
The tri-hexagonal grid.}
\end{figure}
%
%%%%%%%%%%%%%%%%%%%%%%%%%%

% xxx could kill figure, see Figure~\Ref{fig:simple-grids}
Curves on the tri-hexagonal grid (see Figure~\ref{fig:tri-hex-grid}) must have
non-zero turns between adjacent edges, otherwise two dead ends remain.

\begin{figure}[h!tbp]
\ifpdf
\begin{center}
{\includegraphics*[width=90mm, viewport={60 410 490 740}]{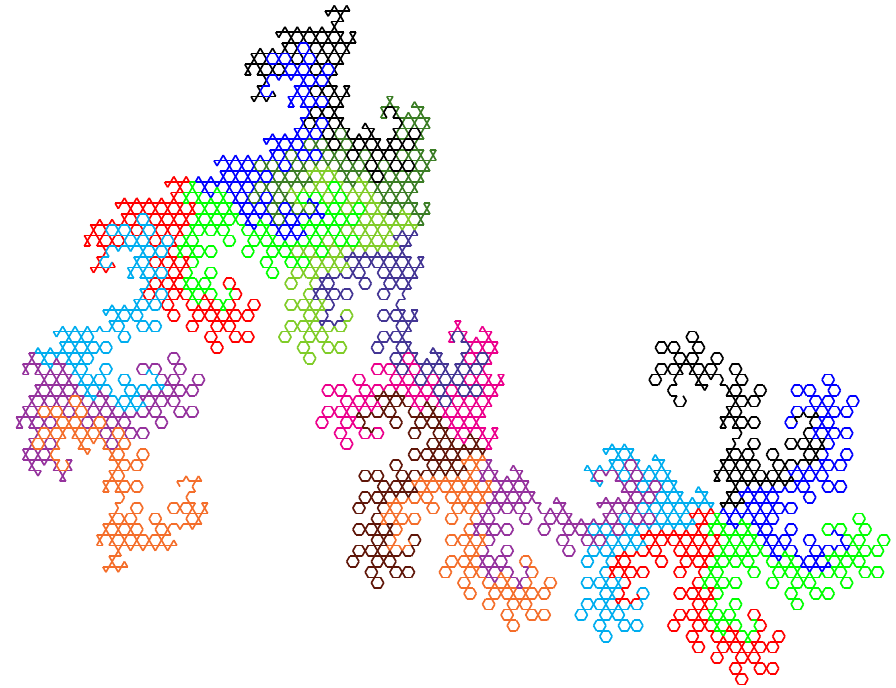}}
\end{center}
\else
\verb+{see pdf for image}+
\fi
\caption{\label{fig:r19-b-1-curve}
Third iterate of a curve of order 19 on the tri-hexagonal grid (\CID{R19-1}).
The coloration emphasizes the self-similarity.}
\end{figure}

We use turns by $60\adeg$ in the turtle graphics
and represent $120\adeg$ turns by \texttt{--} or \texttt{++}
in the L-systems.
For example, the curve of order 19 shown in Figure~\ref{fig:r19-b-1-curve}
is given by the L-system {\Lmap{F}{F+F+F+F+F--F+F+F+F--F+F+F--F--F+F+F--F+F--F}}.

%%%%%%%%%%%%%%%%%%%%%%%%%%
%% with lnth *= 2.0;  // thicker lines
% stringsubst 2 +F_+F_+F_+F_+F_--F_+F_+F_+F_--F_+F_+F_--F_--F_+F_+F_--F_+F_--F_ _ _ F F+F+F+F+F--F+F+F+F--F+F+F--F--F+F+F--F+F--F + + - - | tail -1 | ./bin 6 3 0 > tmp-pic.tex && make dotex # R19-1
%
%
%%%%%%%%%%%%%%%%%%%%%%%%%%

%%%%%%%%%%%%%%%%%%%%%%%%%%
\begin{figure}[h!tbp]
%%  ./scbin 13 6 0 c6/search-r13-b-curves.txt | ./shorten-output.pl
{\small
\begin{verbatim}
F F+F+F+F+F--F+F+F--F--F+F+F--F  R13-1  #
F F+F+F+F+F--F+F--F--F+F+F+F--F  R13-2  #
F F+F+F+F--F+F+F+F--F--F+F--F+F  R13-3  #
F F+F--F+F--F--F+F+F+F--F+F+F+F  R13-4  # ## same = 3 Z T
\end{verbatim}
}
\caption{\label{fig:output-13-b}
Descriptions of the curves of order 13 on the tri-hexagonal grid.}
\end{figure}
%
%%%%%%%%%%%%%%%%%%%%%%%%%%

Observe that all turns by $60\adeg$ are in one direction (\texttt{+})
while all turns in the other direction are by $120\adeg$ (\texttt{--}),
see Figure~\ref{fig:output-13-b} for the file giving
the L-system for the curves of order 13.

This is true for all curves on this grid:
at any step turns \texttt{++} or \texttt{-}
would give an edge ending in the middle of a hexagon.
Obviously swapping the signs would give a similar curve
whose shape is obtained by flipping the original curve,
so we only allow turns \texttt{+} and \texttt{--}.

There is no curve with an L-system beginning and ending with \texttt{F--F}
as neither tile would satisfy all conditions from section \ref{sect:conditions}.

Curves on the tri-hexagonal grid exist only
for orders $R=6k+1$ where $R$ is of the form
$x^2+x\,y+x^2$, see sequence \jjseqref{A260682} in \cite{OEIS}.
%% Odd Loeschian numbers of the form 3*k+1.
%
%% Must be of form 3k+1 (no non-turns) and odd (due to signs) ==> 6k+1
%? is(n)=(n%6==1)&&#bnfisintnorm(bnfinit(z^2+z+1), n);
%? select(n->is(n), vector(300,j,j))
% [1, 7, 13, 19, 25, 31, 37, 43, 49, 61, 67, 73, 79, 91, 97,
% 103, 109, 121, 127, 133, 139, 151, 157, 163, 169, 175, 181, 193, 199,
% 211, 217, 223, 229, 241, 247, 259, 271, 277, 283, 289]

%
The curves never have any non-trivial symmetry.

%%%%%%%%%%%%%%%%%%%%%%%%%%
% stringsubst 1 F_+F_+F_+F_+F_+F  _ _ F F+F+F+F+F--F+F+F+F--F+F+F--F--F+F+F--F+F--F   + + - - | tail -1 | ./bin 6 2 0 0 0.15 > tmp-pic.tex && make dotex # R19-1
% stringsubst 1 _F--_F--_F  _ _ F F+F+F+F+F--F+F+F+F--F+F+F--F--F+F+F--F+F--F   + + - - | tail -1 | ./bin 6 2 0 0 0.15 > tmp-pic.tex && make dotex # R19-1
%
%% with substitution rounding, render with thick lines: lnth *= 3.0; :
% stringsubst 1 F_+F_+F_+F_+F_+F _ _ F F+F+F+F+F--F+F+F+F--F+F+F--F--F+F+F--F+F--F + + - - | tail -1 | sed 's/F/FFFF/g; s/--/--t--/g; s/+/+t+/g;' | ./bin 12 3 0 0 0.15 > tmp-pic.tex && make dotex # R19-1
% stringsubst 1 _F_--F_--F_  _ _ F F+F+F+F+F--F+F+F+F--F+F+F--F--F+F+F--F+F--F + + - - | tail -1 | sed 's/F/FFFF/g; s/--/--t--/g; s/+/+t+/g;' | ./bin 12 3 0 0 0.15 > tmp-pic.tex && make dotex # R19-1
%
\begin{figure}[h!tbp]
\ifpdf
\begin{center}
{\includegraphics*[width=43mm, viewport={60 310 490 740}]{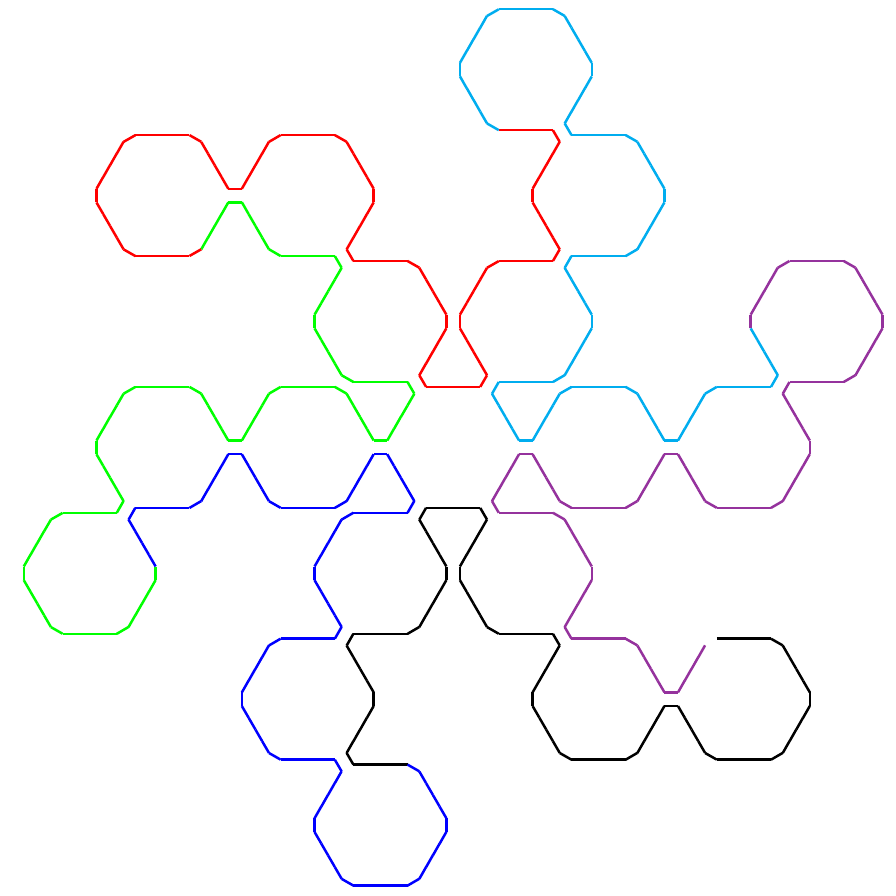}}%
\qquad% layout
{\includegraphics*[width=42mm, viewport={60 310 500 750}]{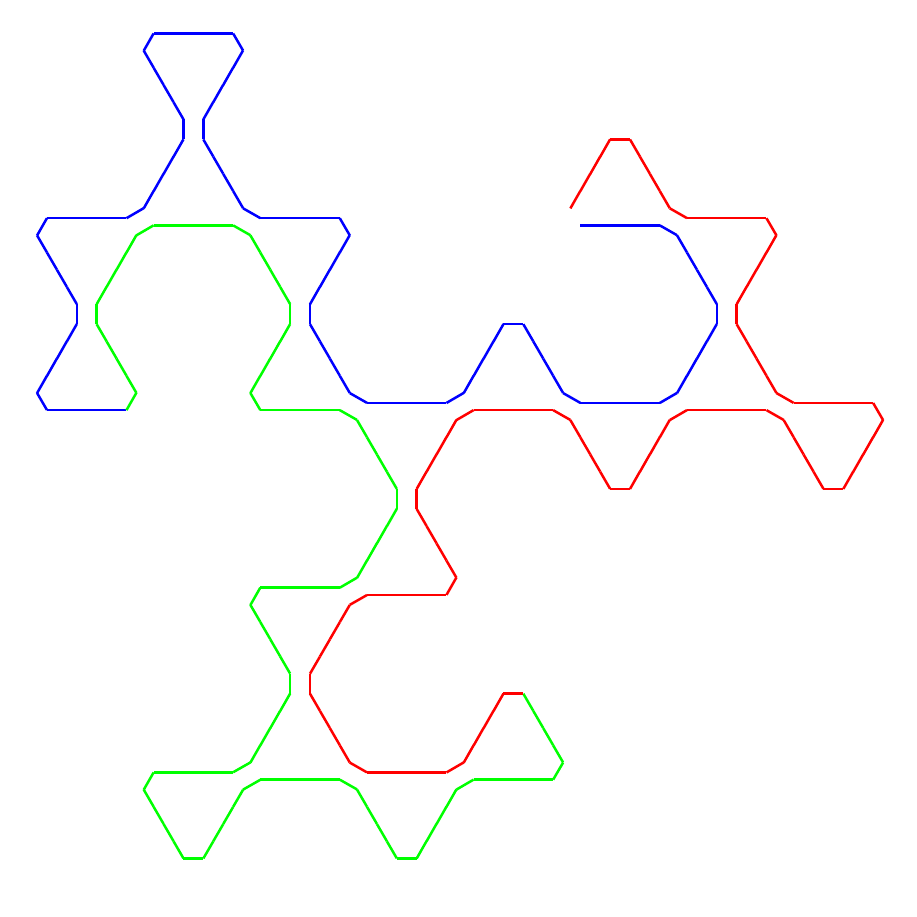}}
\end{center}
\else
\verb+{see pdf for image}+
\fi
\caption{\label{fig:r19-b-1-tiles-it1}
Tiles $\Tile{+1}$ and $\Tile{-1}$ for the curve in Figure~\ref{fig:r19-b-1-curve} (\CID{R19-1}).}
\end{figure}
%
%%%%%%%%%%%%%%%%%%%%%%%%%%

%%%%%%%%%%%%%%%%%%%%%%%%%%
% stringsubst 3 F_+F_+F_+F_+F_+F _ _ F F+F+F+F+F--F+F+F+F--F+F+F--F--F+F+F--F+F--F + + - - | tail -1 | ./bin 6 3 0 > tmp-pic.tex && make dotex # R19-1
\begin{figure}[h!tbp]
\ifpdf
\begin{center}
{\includegraphics*[width=90mm, viewport={50 340 500 740}]{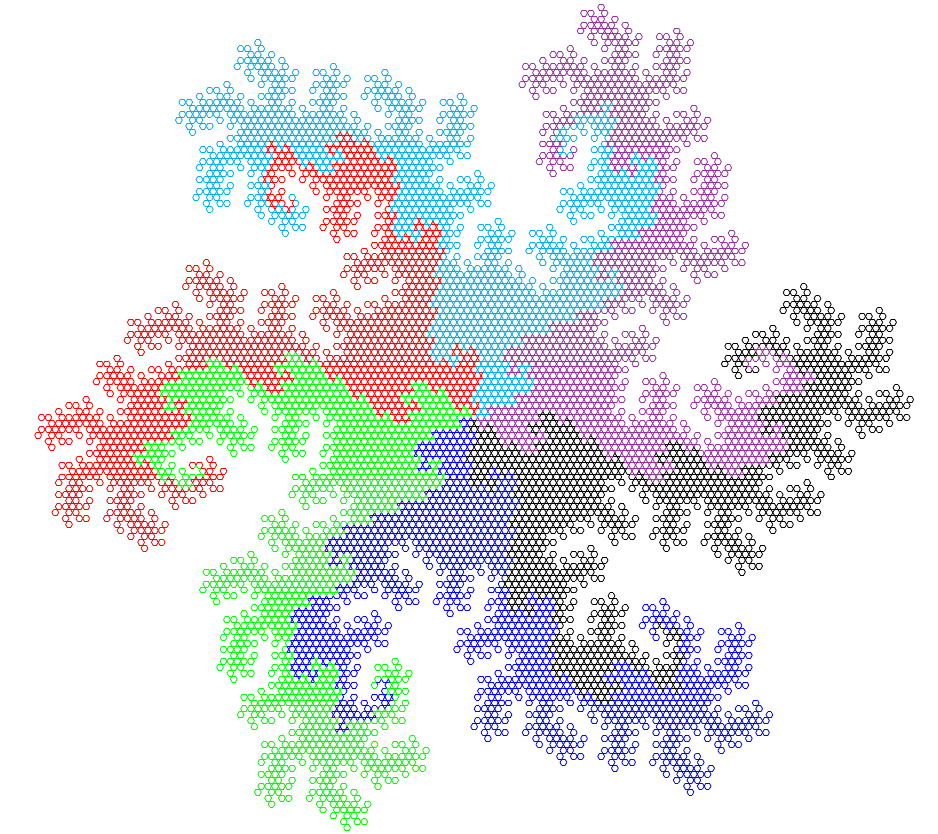}}
\end{center}
\else
\verb+{see pdf for image}+
\fi
\caption{\label{fig:r19-b-1-tile-plus}
Tile $\Tile{+3}$ for the curve in Figure~\ref{fig:r19-b-1-curve} (\CID{R19-1}).}
\end{figure}
%
%%%%%%%%%%%%%%%%%%%%%%%%%%

%%%%%%%%%%%%%%%%%%%%%%%%%%
% stringsubst 3 _F--_F--_F  _ _ F F+F+F+F+F--F+F+F+F--F+F+F--F--F+F+F--F+F--F   + + - - | tail -1 | ./bin 6 3 0 > tmp-pic.tex && make dotex # R19-1
\begin{figure}[h!tbp]
\ifpdf
\begin{center}
{\includegraphics*[width=60mm, viewport={50 280 500 750}]{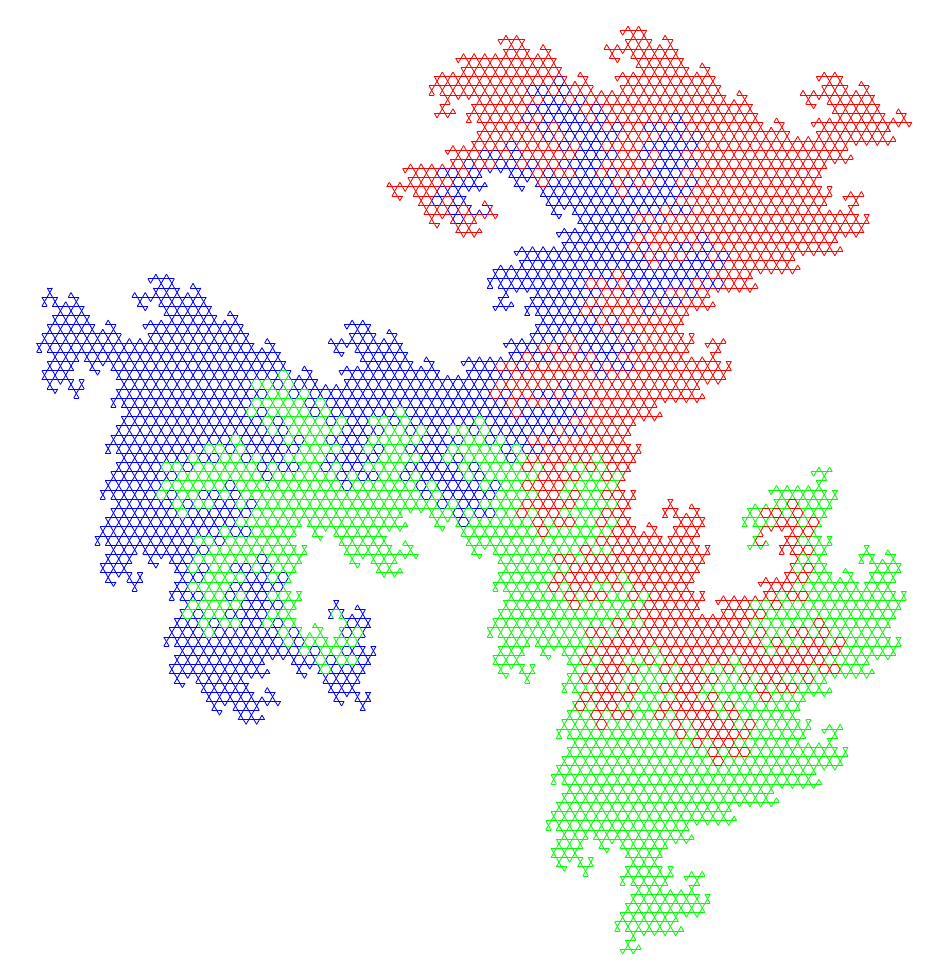}}
\end{center}
\else
\verb+{see pdf for image}+
\fi
\caption{\label{fig:r19-b-1-tile-minus}
Tile $\Tile{-3}$ for the curve in Figure~\ref{fig:r19-b-1-curve} (\CID{R19-1}).}
\end{figure}
%
%%%%%%%%%%%%%%%%%%%%%%%%%%

The two tiles $\Tile{+}$ and $\Tile{-}$ for each curve have
axioms \texttt{F+F+F+F+F+F} and \texttt{F--F--F}.
The tiles $\Tile{+1}$ and $\Tile{-1}$ for the curve mentioned are shown
in Figure~\ref{fig:r19-b-1-tiles-it1}.
All tiles $\Tile{+j}$ are centered at a hexagon and have 6-fold symmetry
while all tiles $\Tile{-j}$ are centered at a triangle and have 3-fold symmetry.
Figures~\ref{fig:r19-b-1-tile-plus} and \ref{fig:r19-b-1-tile-minus}
show the third iterates.
%$\Tile{+3}$ and $\Tile{-3}$.

%%%%%%%%%%%%%%%%%%%%%%%%%%
% stringsubst 2 +F_+F_+F_--F_+F_+F_--F_--F_+F_+F_--F_+F_--F_+F_+F_+F_--F_+F_+F_+F_+F_+F_--F_+F_--F_  _ _  F F+F+F--F+F+F--F--F+F+F--F+F--F+F+F+F--F+F+F+F+F+F--F+F--F   + + - - | tail -1 | ./bin 6 3 0 > tmp-pic.tex && make dotex # R25-11
\begin{figure}[h!tbp]
\ifpdf
\begin{center}
%% un-rotated:
%{\includegraphics*[width=128mm, viewport={70 170 490 740}]{r25-b-11-decomp.pdf}}
%% rotated:
{\includegraphics*[width=130mm, viewport={60 520 490 740}]{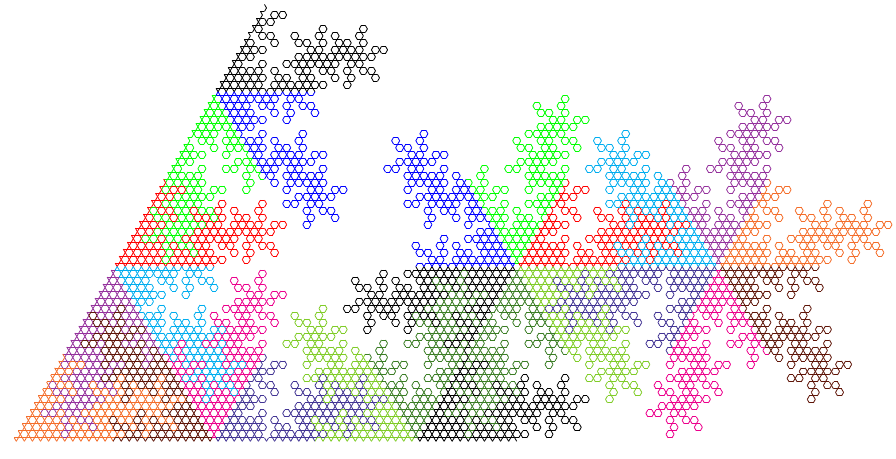}}
\end{center}
\else
\verb+{see pdf for image}+
\fi
\caption{\label{fig:r25-b-11-decomp}
Self-similarity of a curve of order 25 (\CID{R25-11}).}
\end{figure}
%
%%%%%%%%%%%%%%%%%%%%%%%%%%

The tile $\Tile{+}$ again tiles the grid with translated copies.
The tile $\Tile{-}$, however, tiles in the way the triangular grid
is tiled by triangles, half of them are rotated by $120\adeg$.
This can be seen in Figure~\ref{fig:r25-b-11-decomp}
for a curve of order 25:
%% un-rotated:
%the upper right part shows how the tiles $\Tile{+}$ match
%while in the bottom of the triangular form of the tiles $\Tile{-}$ is obvious.
%% rotated:
the right part shows how the tiles $\Tile{+}$ match
and on the left the triangular form of the tiles $\Tile{-}$ is obvious.
Curves whose tile $\Tile{-}$ has the shape of an equilateral triangle exist only
if the order is a perfect square.

The tiles $\Tile{+}$ always give integral complex numeration systems,
their shapes $\Tile{+\infty}$ are the same as the shapes of tiles with 6-fold symmetry
(in the limit) appearing for curves on the triangular grid of the same order.
The tiles $\Tile{-}$ on the tri-hexagonal grid
are never equal to any shape of a tile on the triangular grid.

%% arrangement (carousel):
% stringsubst 4 [[[_F]-_G]+++_F]++_G [ [ ] ] _ _ F F+F+F+F--F--F+F G G-G++G++G-G-G-G + + - - | tail -1 | ./bin 6 3 0 > tmp-pic.tex && make dotex

%%% Emacs:
%%% Local Variables:
%%% mode: latex
%%% MyRelDir: "."
%%% TeX-master: "arndt-curve-search.tex"
%%% dvi-file: "arndt-curve-search"
%%% makefile-dir: "./"
%%% frame-title-format: "CURVE-SEARCH (properties)"
%%% End:

\FloatBarrier

%%%%%%%%%%%%%%%%%%%%%%%%%%%%%%%%%%%%%%%%%%%%%%%%%%%%%%%%%%%%
%%%%%%%%%%%%%%%%%%%%%%%%%%%%%%%%%%%%%%%%%%%%%%%%%%%%%%%%%%%%
%\section{Curves traversing all grid points once}
\section{Plane-filling curves on all uniform grids}
%
%% xxx kludge to reset (letter at) numbering for figures:
\addtocounter{subsection}{-2}
\refstepcounter{subsection}
\addtocounter{subsection}{+1}

The curves corresponding to simple L-systems are edge-covering
(traverse each edge of the underlying grid once)
and necessarily traverse each point more than once
(hence are not point-covering).
%
%Here we observe that certain modifications of the curves
%traverse each point of some grid once.
%
Here we give methods to convert these curves
into point-covering curves on all uniform grids
and to edge-covering curves on two uniform grids.

%%%%%%%%%%%%%%%%%%%%%%%%%%
%
\begin{figure}[h!tbp]
\ifpdf
\begin{center}
{\includegraphics*[width=40mm, viewport={225 595 250 620}]{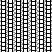}}%
\hspace*{4mm}% layout
{\includegraphics*[width=40mm, viewport={310 460 350 500}]{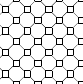}}%
\hspace*{4mm}% layout
{\includegraphics*[width=40mm, viewport={310 460 350 500}]{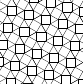}}
\end{center}
\else
\verb+{see pdf for image}+
\fi
\caption{\label{fig:some-grids-a}
The grids
$(3.3.3.4.4)=(3^3.4^2)$,
$(4.8.8)$, and
$(3.3.4.3.4)$.}
\end{figure}
%
%%%%%%%%%%%%%%%%%%%%%%%%%%

%%%%%%%%%%%%%%%%%%%%%%%%%%
%
\begin{figure}[h!tbp]
\ifpdf
\begin{center}
{\includegraphics*[width=40mm, viewport={300 450 370 520}]{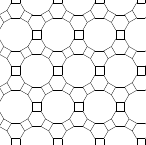}}%
\hspace*{4mm}% layout
{\includegraphics*[width=40mm, viewport={310 460 350 500}]{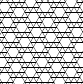}}%
\hspace*{4mm}% layout
{\includegraphics*[width=40mm, viewport={325 475 350 500}]{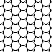}}
\end{center}
\else
\verb+{see pdf for image}+
\fi
\caption{\label{fig:some-grids-b}
The grids
$(4.6.12)$,
$(3.3.3.3.6)=(3^4.6)$, and
$(3.12.12)=(3.12^2)$.}
\end{figure}
%
%%%%%%%%%%%%%%%%%%%%%%%%%%

%% https://oeis.org/A250120
%
%% https://de.wikipedia.org/wiki/Parkettierung
%% cf. https://commons.wikimedia.org/wiki/User:Nonenmac/tessellations
%% https://en.wikipedia.org/wiki/Euclidean_tilings_by_convex_regular_polygons
%% cf. https://en.wikipedia.org/wiki/List_of_convex_uniform_tilings
%
%% "... 'vertex type' referring to the type and order of polygons surrounding each vertex in the tiling."
%% from http://probabilitysports.com/tilings.html
%
We specify the type of grid
by the symbols for the vertex type
for the corresponding uniform edge to edge tilings
of regular polygons
as in \cite[p.~63]{gruenbaum-shepard-1986},
also see
 \cite{wiki-unifom-tilings} and
 \cite{wiki-convex-tilings}:
%and
%\url{http://probabilitysports.com/tilings.html?u=0&n=1&t=1} in \cite{galebach}\xxx{Layout}
%
the symbol is the (least cyclic shift of the)
list of numbers of edges of the polygons
around any of the points.

%The grids consisting of just one type of polygons are
%% Platonic:
The symbols for the grids so far seen are
$(3.3.3.3.3.3)=(3^6)$ for the triangular grid,
$(4.4.4.4)=(4^4)$ for the square grid,
$(6.6.6)=(6^3)$ for the hexagonal grid, and
%% Archimedian:
$(3.6.3.6)$ for the tri-hexagonal grid.
We will encounter the grids
$(3.3.3.4.4)=(3^3.4^2)$,
$(4.8.8)=(4.8^2)$, and
$(3.3.4.3.4)=(3^2.4.3.4)$, shown in Figure~\ref{fig:some-grids-a},
and
$(4.6.12)$,
$(3.3.3.3.6)=(3^4.6)$ (there are two enantiomers of this grid), and
$(3.12.12)=(3.12^2)$, shown in Figure~\ref{fig:some-grids-b}.
The remaining grid $(3.4.6.4)$ is shown in section \ref{sect:3464-EC}.

We call a curve $(G)$-PC if it is point-covering on the grid $(G)$.
For example,
both the Peano and the Hilbert curve (Figure~\ref{fig:peano-hilbert}) are $(4^4)$-PC
and Gosper's flowsnake (Figure~\ref{fig:gosper-flowsnake}) is $(3^6)$-PC
(left rendering), but $(6^3)$-PC when rendered as in the right.
Similarly, we use $(G)$-EC to denote edge-covering on the grid $(G)$.

%%%%%%%%%%%%%%%%%%%%%%%%%%
%% stringsubst 4 t+Lt+t+Lt  L Lt-Lt+t+Lt-L + + - - t t | tr -d L | tail -1 | ./bin 4 2 0 0 0.15 > tmp-pic.tex && make dotex
% stringsubst 5 L  L Lt-Lt+t+Lt-L + + - - t t | tr -d L | tail -1 | ./bin 4 2 0 0 0.15 > tmp-pic.tex && make dotex
%
\begin{figure}[h!tbp]
\ifpdf
\begin{center}
{\includegraphics*[width=100mm, viewport={50 525 500 740}]{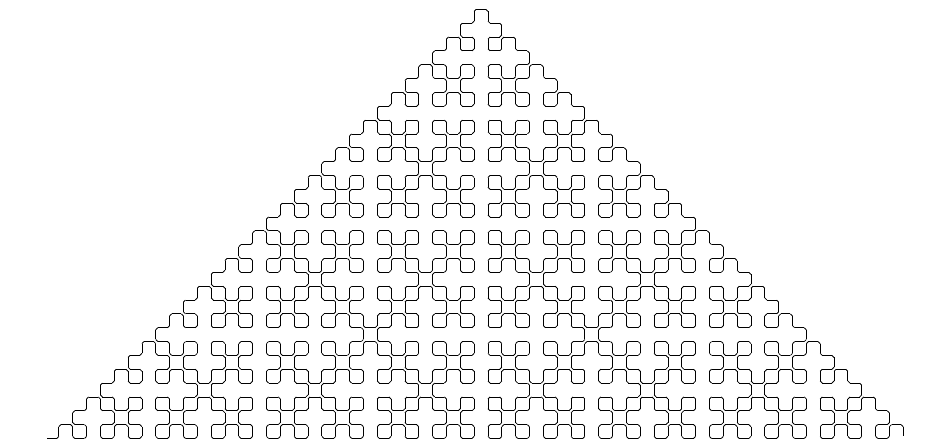}}
\end{center}
\else
\verb+{see pdf for image}+
\fi
\caption{\label{fig:not-ec-not-pc}
A plane-filling curve on the square grid that is neither EC nor PC:
all points are traversed, but some twice and
not all edges are traversed.}
\end{figure}
%
%%%%%%%%%%%%%%%%%%%%%%%%%%
% another example on the tri-hexagonal grid:
% stringsubst 3 F F F+F+F+F--F--F+F + + - - | tail -1 | sed 's/F/+F-F/g; ' | ./bin 6 3 4 0 0.15 > tmp-pic.tex && make dotex

The edge-covering and point-covering curves
are just two corner cases of the plane-filling curves
on a grid, Figure~\ref{fig:not-ec-not-pc}
shows a curve that is in neither class.

%%% Emacs:
%%% Local Variables:
%%% mode: latex
%%% MyRelDir: "."
%%% TeX-master: "arndt-curve-search.tex"
%%% dvi-file: "arndt-curve-search"
%%% makefile-dir: "./"
%%% frame-title-format: "CURVE-SEARCH (covers-intro)"
%%% End:

%\clearpage% xxx
%%%%%%%%%%%%%%%%%%%%%%%%%%
%%%%%%%%%%%%%%%%%%%%%%%%%%
\subsection{Conversions to point-covering curves}\label{sect:PC-curves}

We organize the following sections by the grids whose edge-covering curves are
used to compute point-covering curves on some uniform grids.

%%%%%%%%%%%%%%%%%%%%%%%%%%%%%%%%%%%%%%%%%%%%%%%%%%%%
%%%%%%%%%%%%%%%%%%%%%%%%%%%%%%%%%%%%%%%%%%%%%%%%%%%%
\subsubsection{Triangular grid: curves for $(6^3)$, $(3.6.3.6)$, $(3^3.4^2)$, $(3^6)$, $(3.12.12)$, $(3.4.6.4)$, $(4.6.12)$, $(3^4.6)$, and $(4^4)$}

%%%%%%%%%%%%%%%%%%%%%%%%%%
%% with lnth *= 3.0;  // thick lines
% stringsubst 2 F  F F+F-F-F+F+F-F  0 0 + + - - | tail -1 | ./bin 3 2 0 0 0.0 > tmp-pic.tex && make dotex # R7-5
% stringsubst 2 F  F F+F-F-F+F+F-F  0 0 + + - - | tail -1 | ./bin 3 2 0 0 0.1 > tmp-pic.tex && make dotex # R7-5
% stringsubst 2 F  F F+F-F-F+F+F-F  0 0 + + - - | tail -1 | ./bin 3 2 0 0 0.2 > tmp-pic.tex && make dotex # R7-5
% stringsubst 2 F  F F+F-F-F+F+F-F  0 0 + + - - | tail -1 | ./bin 3 2 0 0 0.3 > tmp-pic.tex && make dotex # R7-5
% stringsubst 2 F  F F+F-F-F+F+F-F  0 0 + + - - | tail -1 | ./bin 3 2 0 0 0.4 > tmp-pic.tex && make dotex # R7-5
% stringsubst 2 F  F F+F-F-F+F+F-F  0 0 + + - - | tail -1 | ./bin 3 2 0 0 0.5 > tmp-pic.tex && make dotex # R7-5
%
\begin{figure}[h!tbp]
\ifpdf
\begin{center}
%{\includegraphics*[width=40mm, viewport={120 420 460 730}]{r07-t-5-rnd00.pdf}}%
%{\includegraphics*[width=40mm, viewport={120 420 460 730}]{r07-t-5-rnd01.pdf}}%
%{\includegraphics*[width=40mm, viewport={120 420 460 730}]{r07-t-5-rnd02.pdf}}
%%
%{\includegraphics*[width=40mm, viewport={120 420 460 730}]{r07-t-5-rnd03.pdf}}%
%{\includegraphics*[width=40mm, viewport={120 420 460 730}]{r07-t-5-rnd04.pdf}}%
%{\includegraphics*[width=40mm, viewport={120 420 460 730}]{r07-t-5-rnd05.pdf}}
%
{\includegraphics*[width=40mm, page=1, viewport={120 420 460 730}]{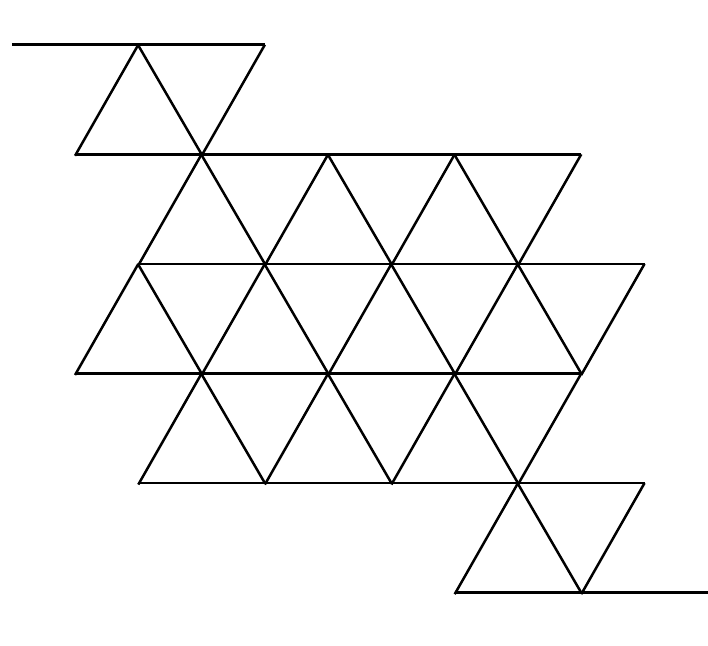}}%
{\includegraphics*[width=40mm, page=2, viewport={120 420 460 730}]{r07-t-5-rnd.pdf}}%
{\includegraphics*[width=40mm, page=3, viewport={120 420 460 730}]{r07-t-5-rnd.pdf}}
{\includegraphics*[width=40mm, page=4, viewport={120 420 460 730}]{r07-t-5-rnd.pdf}}%
{\includegraphics*[width=40mm, page=5, viewport={120 420 460 730}]{r07-t-5-rnd.pdf}}%
{\includegraphics*[width=40mm, page=6, viewport={120 420 460 730}]{r07-t-5-rnd.pdf}}
\end{center}
\else
\verb+{see pdf for image}+
\fi
\caption{\label{fig:r07-t-5-rnd}
Renderings of the curve \CID{R7-5} with map \Lmap{F}{F+F-F-F+F+F-F} on the triangular grid
for rounding parameter $e\in\{0.0,\, 0.1,\, 0.2,\, 0.3,\, 0.4,\, 0.5\}$.}
\end{figure}
%
%%%%%%%%%%%%%%%%%%%%%%%%%%%

The constructions given here only work for curves of two types.
We call curves where all turns are by a non-zero angle \jjterm{wiggly},
and curves with the same numbers of \texttt{+}, \texttt{-}, and \texttt{0}
in the production of \texttt{F} \jjterm{balanced}.
Wiggly curves exist only for odd orders,
see sequence \jjseqref{A266836} in \cite{OEIS}.
%% 1, 3, 7, 9, 13, 19, 21, 25, 27, 31, 37, 39, 43, 49, 57, 61, 63, 67, 73, 75, 79, 81, 91, 93, 97, 103, ...
%% Even orders:
%% 4, 12, 16, 28, 36, 48, 52, 64, 76, 84, 100, ...
The curve \CID{R7-5} with map \Lmap{F}{F+F-F-F+F+F-F}
whose motif is shown on the right of \Ref{fig:shape-7}
is a wiggly curve.
Unless otherwise specified our methods only work for wiggly curves.
Balanced curves only exits for orders of the form $3k+1$,
sequence \jjseqref{A202822} in \cite{OEIS}.
%% 1, 4, 7, 13, 16, 19, 25, 28, 31, 37, 43, 49, 52, 61, 64, 67, 73, 76, 79, 91, 97, 100, ...
An example is the curve \CID{R7-1} with map \Lmap{F}{F0F+F0F-F-F+F},
it is shown in Figure~\Ref{fig:iterate-3-4-5-decomp}.
%
%% None of either class for 12 * A003136:
%% 12, 36, 48, 84, 108, 144, 156, 192, 228, 252, 300, 324, 336, 372, 432, 444, 468, ...
%
In balanced curves every point is traversed once by two edges with a non-turn between them.

We have sometimes used rounded corners to make the curves visually apparent:
at each turn a fraction of $e>0$ of both edges are skipped and
a straight line is used to connect the two new end points.
Figure~\ref{fig:r07-t-5-rnd} shows renderings of a curve of order 7
on the triangular grid for $e\in\{0.0,\, 0.1,\, 0.2,\, 0.3,\, 0.4,\, 0.5\}$.

%%%%%%%%%%%%%%%%%%%%%%%%%%
% stringsubst 2 F  F F+F-F-F+F+F-F  + + - - | tail -1 | ./bin 3 3 2 > tmp-pic.tex && make dotex  # R7-5
%
% stringsubst 2 F  F F+F-F-F+F+F-F  + + - - | tail -1 | ./bin 3 3 8 > tmp-pic.tex && make dotex  # R7-5
%
\begin{figure}[h!tbp]
\ifpdf
\begin{center}
{\includegraphics*[width=63mm, viewport={150 460 470 730}]{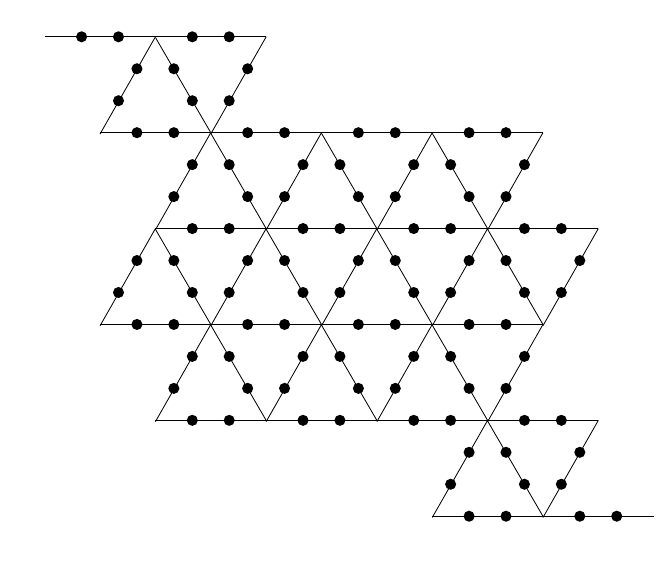}}%
{\includegraphics*[width=63mm, viewport={150 460 470 730}]{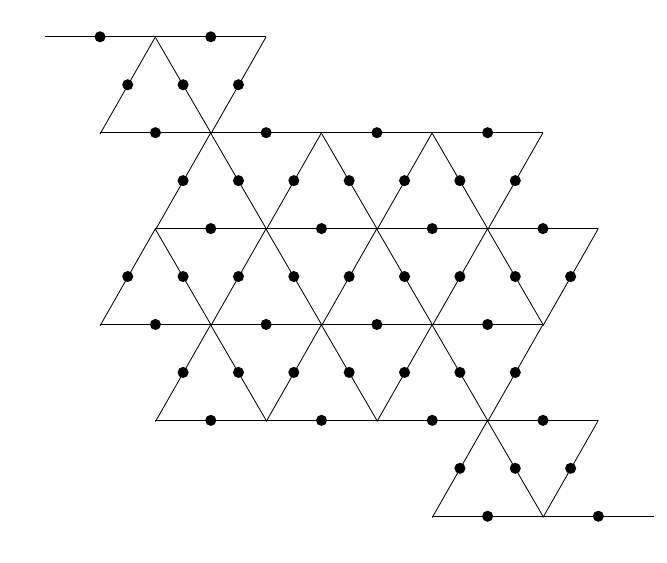}}
\end{center}
\else
\verb+{see pdf for image}+
\fi
\caption{\label{fig:r07-t-5-points}
The turning points with rounded corners for $e=1/3$ (left) and $e=1/2$ (right)
respectively give the hexagonal and tri-hexagonal grid.}
\end{figure}
%
%%%%%%%%%%%%%%%%%%%%%%%%%%%

%%%%%%%%%%%%%%%%%%%%%%%%%%
% with const char bullet[] = "{\\circle*{0.3}}";
% stringsubst 2 F  F F+F-F-F+F+F-F  0 0 + + - - | tail -1 | sed 's/+F/+F+F/g; s/-F/-F-F/g;' | ./bin 6 3 4 > tmp-pic.tex && make dotex # R7-5
%
% with const char bullet[] = "{\\circle*{0.2}}";
% stringsubst 2 F F F+F-F-F+F+F-F 0 0 + + - - | tail -1 | sed 's/F//g; s/+/+F+/g; s/-/-F-/g;' | ./bin 6 3 4 > tmp-pic.tex && make dotex # R7-5
%
\begin{figure}[h!tbp]
\ifpdf
\begin{center}
{\includegraphics*[width=55mm, viewport={110 380 510 735}]{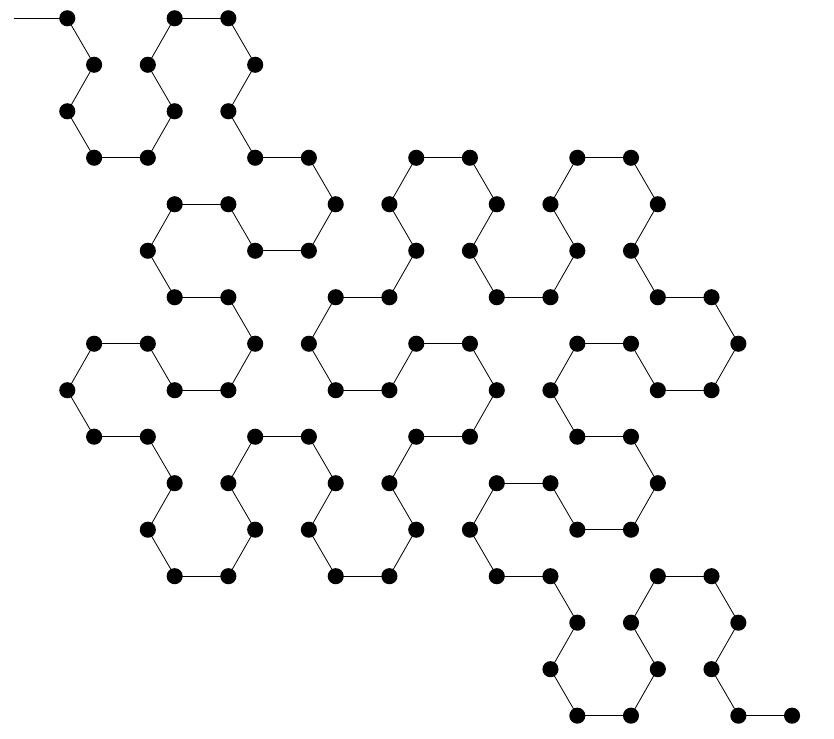}}%
{\includegraphics*[width=53mm, viewport={100 380 490 735}]{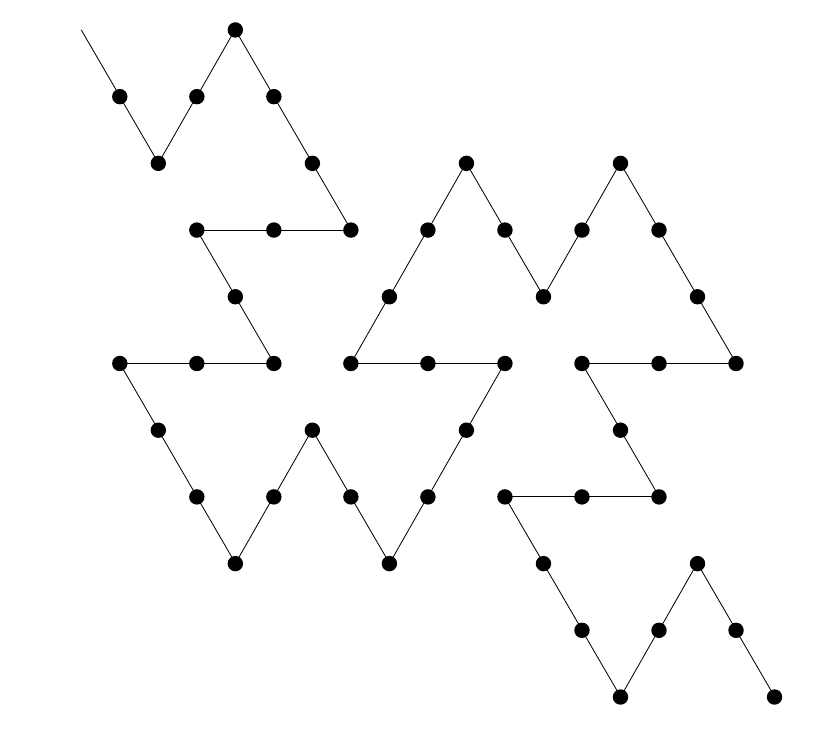}}
\end{center}
\else
\verb+{see pdf for image}+
\fi
\caption{\label{fig:r07-t-5-cover-points}
The curves obtained by post-processing
corresponding to $e=1/3$ (left), a $(6^3)$-PC curve,
and $e=1/2$, a $(3.6.3.6)$-PC curve (right).}
\end{figure}
%
%%%%%%%%%%%%%%%%%%%%%%%%%%%

The rendering in the lower left ($e=1/3$)
and lower right ($e=1/2$)
are respectively a $(6^3)$-PC curve and a $(3.6.3.6)$-PC curve.
The turning points of the rounded corners are indicated
in Figure~\ref{fig:r07-t-5-points}.
%

%% $(6^3)$-grid:
% stringsubst 3 F+F+F F F+F-F-F+F+F-F + + - - | tail -1 | sed 's/+/p/g; s/-/m/g; s/m/-F[+t]-/g; s/p/+F[-t]+/g;' | ./bin 6 3 0 > tmp-pic.tex && make dotex # R7-5

The words for these curves can be computed
using a post-processing step on the final
iterate of the L-system.
%
%%% $(6^3)$-PC:
For the $(6^3)$-PC curve,
%  $e=1/3$
replace
all \texttt{+F} by \texttt{+F+F},
all \texttt{-F} by \texttt{-F-F},
and use turns by $60\adeg$.
% Check:
% stringsubst 2 F F F+F-F-F+F+F-F 0 0 + + - - | tail -1 | sed 's/+/+F+/g; s/-/-F-/g;' | ./bin 6 3 4 > tmp-pic.tex && make dotex # R7-5
%
% Alternatively (same result):
% stringsubst 2 F  F F+F-F-F+F+F-F  0 0 + + - - | tail -1 | sed 's/+F/+F+F/g; s/-F/-F-F/g;' | ./bin 6 3 4 > tmp-pic.tex && make dotex # R7-5
%
%
%
%%% $(3.6.3.6)$-PC:
For the $(3.6.3.6)$-PC curve,
% $e=1/2$
drop all \texttt{F},
replace
all \texttt{+} by \texttt{+F+},
all \texttt{-} by \texttt{-F-},
and use turns by $60\adeg$.
%
% Check:
% stringsubst 2 F F F+F-F-F+F+F-F 0 0 + + - - | tail -1 | sed 's/F//g; s/+/+F+/g; s/-/-F-/g;' | ./bin 6 3 4 > tmp-pic.tex && make dotex # R7-5
%
% Alternatively (even simpler):
% stringsubst 2 F F F+F-F-F+F+F-F 0 0 + + - - | tail -1 | sed 's/F//g; s/+/F+/g; s/-/-F/g;' | ./bin 3 3 4 > tmp-pic.tex && make dotex # R7-5
%
% Alternatively (3 times as many points):
% stringsubst 2 F F F+F-F-F+F+F-F 0 0 + + - - | tail -1 | sed 's/F/XX/g; s/+/F+/g; s/-/-F/g;' | ./bin 3 3 4 > tmp-pic.tex && make dotex # R7-5
% Note: just keeping the F initially gives the $(3^4.6)$-PC curve
%
% Alternatively:
% stringsubst 2 F F F+F-F-F+F+F-F 0 0 + + - - | tail -1 | sed 's/F//g; s/+/+++F+/g; s/-/-F---/g;' | ./bin 12 3 4 > tmp-pic.tex && make dotex # R7-5
%
% Alternative for $(3.6.3.6)$-PC (3 times as many points):
% stringsubst 2 F F F+F-F-F+F+F-F + + - - | tail -1 | sed 's/F//g; s/+/p/g; s/-/m/g; s/p/FF++F/g; s/m/F--FF/g;' | ./bin 6 3 4 > tmp-pic.tex && make dotex
% stringsubst 2 F F F+F-F-F+F+F-F 0 0 + + - - | tail -1 | sed 's/F+/FF+F/g; s/F-/F-FF/g;' | ./bin 3 3 4 > tmp-pic.tex && make dotex
%
% stringsubst 2 F F F+F-F-F+F+F-F + + - - | tail -1 | sed 's/+/+t+/g; s/-/-t-/g; s/F/-FF+/g;' | ./bin 6 3 4 > tmp-pic.tex && make dotex
%
%
The corresponding curves with the ends of all edges marked
are shown in Figure~\ref{fig:r07-t-5-cover-points}.
%

%%%%%%%%%%%%%%%%%%%%%%%%%%
% with const char bullet[] = "{\\circle*{0.3}}";
% stringsubst 3 -F  F F+F-F-F+F+F-F  + + - - | tail -1 | tr -d F | sed 's/+/++F/g; s/-/-F-/g;' | ./bin 6 3 4 > tmp-pic.tex && make dotex  # R7-5
%
\begin{figure}[h!tbp]
\ifpdf
\begin{center}
{\includegraphics*[width=100mm, viewport={70 450 490 740}]{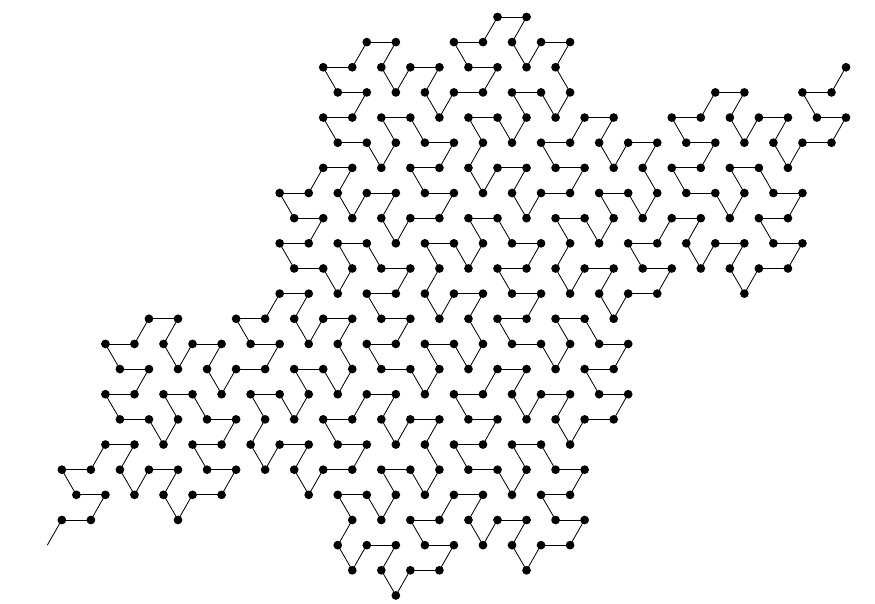}}
\end{center}
\else
\verb+{see pdf for image}+
\fi
\caption{\label{fig:333333-cover}
%A curve traversing all points of the triangular grid once.}
The $(3^6)$-PC curve from the third iterate of \CID{R7-5}.}
\end{figure}
%
%%%%%%%%%%%%%%%%%%%%%%%%%%%

%%%%%%%%%%%%%%%%%%%%%%%%%%
% with lnth *= 2.0;  // thicker lines
% with const char bullet[] = "{\\circle*{0.3}}";
% stringsubst 3 F F F0F+F0F-F-F+F 0 0 - - + + | tail -1 | sed 's/F//g; s/+/+F+/g; s/-/-F-/g; s/0/F0F/g; s/^/+/;' | ./bin 6 3 4 > tmp-pic.tex && make dotex
%
\begin{figure}[h!tbp]
\ifpdf
\begin{center}
{\includegraphics*[width=80mm, viewport={60 390 490 735}]{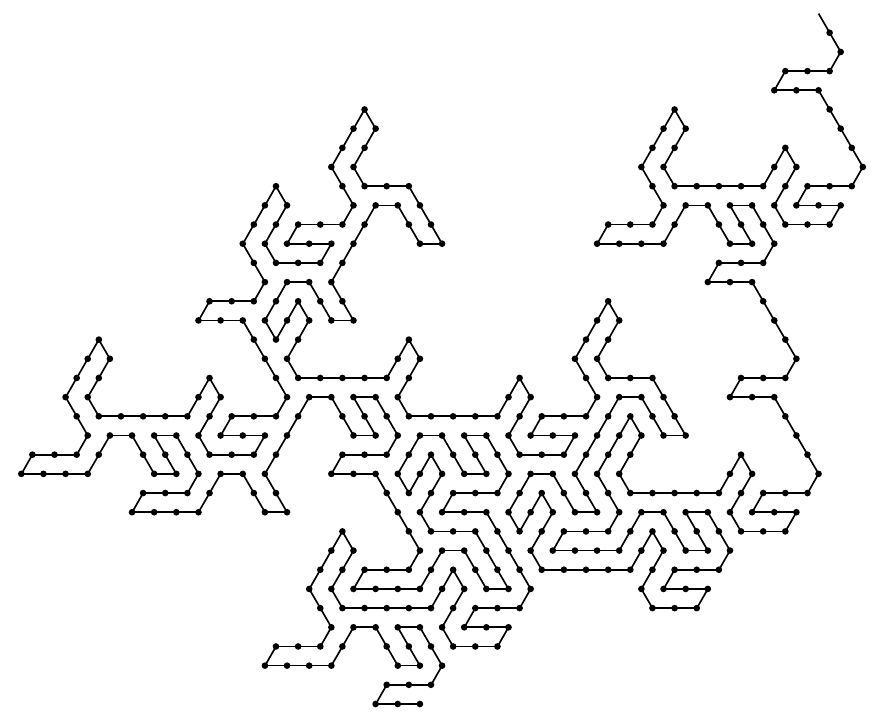}}
\end{center}
\else
\verb+{see pdf for image}+
\fi
\caption{\label{fig:333333-cover-balanced}
The $(3^6)$-PC curve from the third iterate of the balanced curve \CID{R7-1}.}
\end{figure}
%
%%%%%%%%%%%%%%%%%%%%%%%%%%%

Curves that are $(3^6)$-PC
can be obtained by post-processing as well:
drop all \texttt{F},
replace
all \texttt{+} by \texttt{++F},
all \texttt{-} by \texttt{-F-},
and use turns by $60\adeg$.
%
% Check:
% stringsubst 2 F F F+F-F-F+F+F-F 0 0 + + - - | tail -1 | tr -d F | sed 's/+/++F/g; s/-/-F-/g;' | ./bin 6 3 4 > tmp-pic.tex && make dotex # R7-5
%
% Alternatively (same but for first edge):
% stringsubst 2 F F F+F-F-F+F+F-F 0 0 + + - - | tail -1 | sed 's/+F/++F/g; s/-F/-F-/g;' | ./bin 6 3 4 > tmp-pic.tex && make dotex # R7-5
%
The third iterate of the resulting curve
is shown in Figure~\ref{fig:333333-cover}.
We can generate $(3^6)$-PC curves from balanced curves by
dropping all
\texttt{F},
replacing
all \texttt{+} by \texttt{+F+},
all \texttt{-} by \texttt{-F-}, and
all \texttt{0} by \texttt{F0F}.
Turns by $60\adeg$ are used.
The curve from \CID{R7-1} is shown in Figure~\ref{fig:333333-cover-balanced}.
% stringsubst 2 F F F0F+F+F-F-F0F 0 0 - - + + | tail -1 | sed 's/F//g; s/+/+F+/g; s/-/-F-/g; s/0/F0F/g;' | ./bin 6 3 4 > tmp-pic.tex && make dotex
% stringsubst 2 _F+_F+_F+ _ _ F F0F+F0F-F-F+F 0 0 - - + + | tail -1 | sed 's/F//g; s/+/+F+/g; s/-/-F-/g; s/0/F0F/g;' | ./bin 6 3 4 > tmp-pic.tex && make dotex
%
% Alternatively:
% stringsubst 3 F F F0F+F0F-F-F+F 0 0 + + - - | tail -1 | sed 's/F//g; s/+/+P+++/g; s/-/---M-/g; s/0/-NN+/g;' | ./bin 12 3 0 > tmp-pic.tex && make dotex #

%%%%%%%%%%%%%%%%%%%%%%%%%%
% with lnth *= 2.0;  // thicker lines
% with const char bullet[] = "{\\circle*{0.3}}";
% stringsubst 2 F F F+F-F-F+F+F-F 0 0 + + - - | tail -1 | sed 's/F//; s/+F/X/g; s/-F/Y/g; s/X/F-F++++F-F-F++++F-/g; s/Y/F+F----F+F+F----F+/g;' | ./bin 12 3 4 > tmp-pic.tex && make dotex # R7-5
%% better rotation:
% stringsubst 2 F F F+F-F-F+F+F-F 0 0 + + - - | tail -1 | sed 's/F//; s/+F/X/g; s/-F/Y/g; s/X/F-F++++F-F-F++++F-/g; s/Y/F+F----F+F+F----F+/g; s/^/-/;' | ./bin 12 3 4 > tmp-pic.tex && make dotex # R7-5
%
\begin{figure}[h!tbp]
\ifpdf
\begin{center}
{\includegraphics*[width=80mm, viewport={70 470 490 740}]{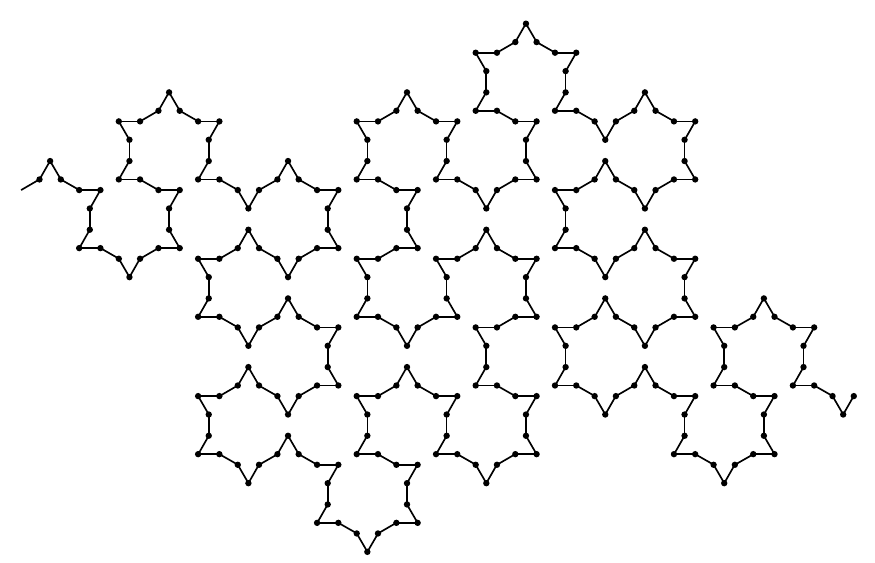}}
\end{center}
\else
\verb+{see pdf for image}+
\fi
\caption{\label{fig:31212-cover}
The $(3.12.12)$-PC curve from the second iterate of \CID{R7-5}.}
\end{figure}
%
%%%%%%%%%%%%%%%%%%%%%%%%%%%

Curves that are $(3.12.12)$-PC
like the one shown in Figure~\ref{fig:31212-cover}
can be obtained as follows.
Drop the initial \text{F},
then replace
all \texttt{+F} by \texttt{X} and
all \texttt{-F} by \texttt{Y},
then replace
all \texttt{X} by \texttt{F-F++++F-F-F++++F-} and
all \texttt{Y} by \texttt{F+F----F+F+F----F+},
use turns by $30\adeg$.
%
% Check $(3.12.12)$-PC:
% stringsubst 1 F F F+F-F-F-F+F+F+F-F + + - - | tail -1 | sed 's/F//g; s/+/p/g; s/-/m/g; s/p/F-F++++F-/g; s/m/F+F----F+/g;' | ./bin 12 3 4 > tmp-pic.tex && make dotex

% $(3.12.12)$-grid:
% stringsubst 3 F+F+F F F+F-F + + - - | tail -1 | sed 's/F//; s/+F/X/g; s/-F/Y/g; s/X/F[+t]-F[-t]++++F-F[+t]-F++++F-/g; s/Y/F[-t]+F[+t]----F+F[-t]+F----F+/g;' | ./bin 12 3 0 > tmp-pic.tex && make dotex # R7-5

%% Cannot quite adapt to curves with non-turns:
% stringsubst 4 F _ _ F F+F0F-F 0 0 + + - - | tail -1 | sed 's/F//g; s/+/X/g; s/-/Y/g; s/X/_F-F++++F-F-F++++F-/g; s/Y/_F+F----F+F+F----F+/g; s/0/_F+F----F+F-F++++F-/g;' | ./bin 12 3 4 > tmp-pic.tex && make dotex #

%%%%%%%%%%%%%%%%%%%%%%%%%%%
%% with const char bullet[] = "{\\circle*{0.3}}";
% stringsubst 2 F F F+F-F-F+F+F-F + + - - | tail -1 | sed 's/F//g; s/+/p/g; s/-/m/g; s/p/++F+++F----F+++F+++F----F+/g; s/m/--F---F++++F---F---F++++F-/g;' | ./bin 12 3 4 > tmp-pic.tex && make dotex # R7-5
%
% with const char bullet[] = "{\\circle*{0.3}}";
%  stringsubst 2 F F F+F-F-F+F+F-F + + - - | tail -1 | sed 's/F//g; s/+/p/g; s/-/m/g; s/p/+++F----F+++F+++F----F+++F/g; s/m/+++F----F-F-F----F+++F/g;' | ./bin 12 3 4 > tmp-pic.tex && make dotex # R7-5
%
\begin{figure}[h!tbp]
\ifpdf
\begin{center}
{\includegraphics*[width=63mm, viewport={80 330 490 735}]{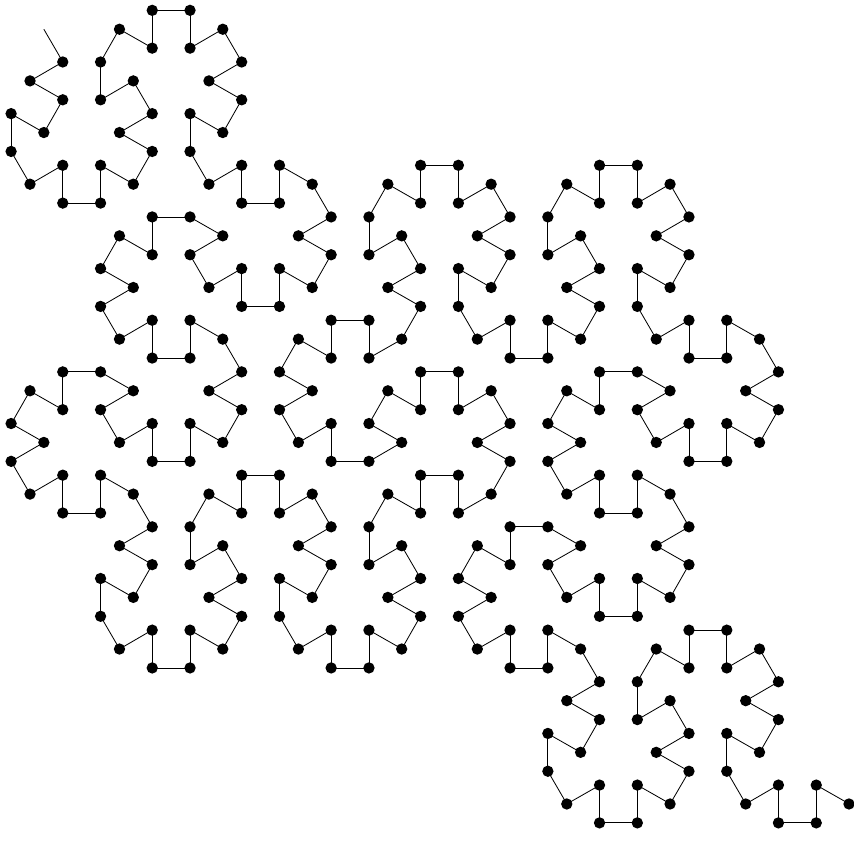}}%
{\includegraphics*[width=63mm, viewport={80 330 500 735}]{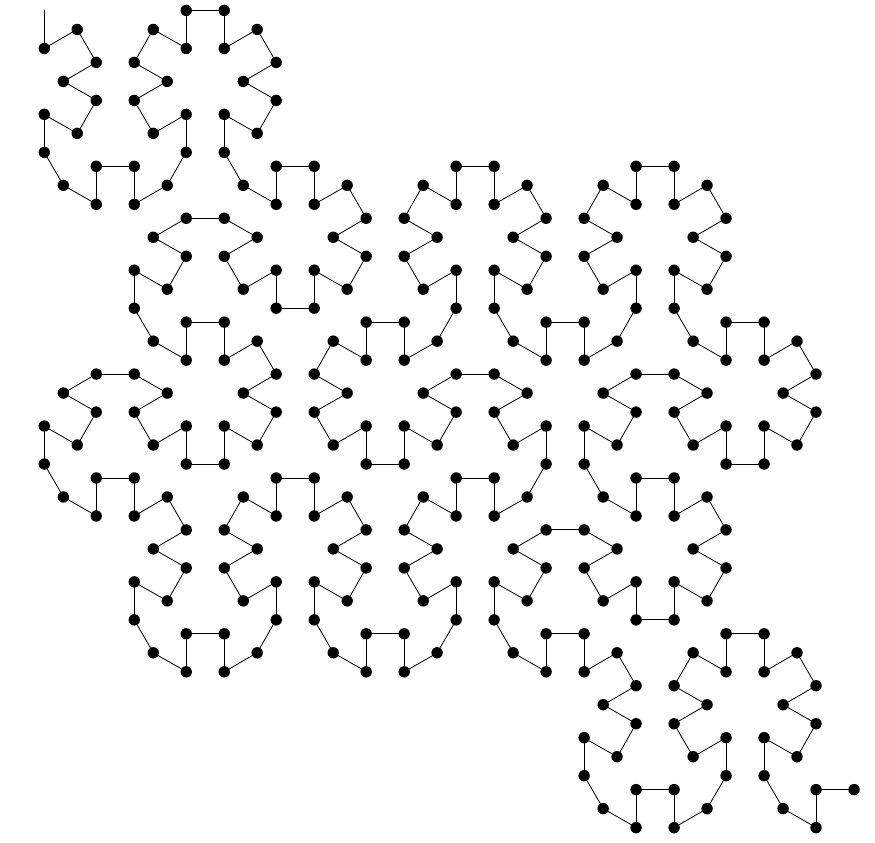}}
\end{center}
\else
\verb+{see pdf for image}+
\fi
\caption{\label{fig:3464-t-cover-both}
Two $(3.4.6.4)$-PC curves from the second iterate of \CID{R7-5}.}
\end{figure}
%
%%%%%%%%%%%%%%%%%%%%%%%%%%%

Two $(3.4.6.4)$-PC curves are shown in Figure~\ref{fig:3464-t-cover-both}.
For the curve on the left,
drop all \texttt{F},
replace
all \texttt{+} by \texttt{p} and
all \texttt{-} by \texttt{m},
then replace
all \texttt{p} by \texttt{++F+++F----F+++F+++F----F+} and
all \texttt{m} by \texttt{--F---F++++F---F---F++++F-},
use turns by $30\adeg$.
For the curve on the right,
drop all \texttt{F},
replace
all \texttt{+} by \texttt{p} and
all \texttt{-} by \texttt{m},
then replace
all \texttt{p} by \texttt{+++F----F+++F+++F----F+++F} and
all \texttt{m} by \texttt{+++F----F-F-F----F+++F},
again use turns by $30\adeg$.
%
% Check:
% stringsubst 2 F F F+F-F-F+F+F-F + + - - | tail -1 | sed 's/F//g; s/+/p/g; s/-/m/g; s/p/+++F----F+++F+++F----F+++F/g; s/m/+++F----F-F-F----F+++F/g;' | ./bin 12 3 4 > tmp-pic.tex && make dotex # R7-5

%% The $(3.4.6.4)$-grid, render with thin lines: lnth *= 1.0 / 3.0;:
% stringsubst 3 F+F+F F F+F-F-F+F+F-F + + - - | tail -1 | sed 's/F//g; s/+/p/g; s/-/m/g; s/p/+++[--y]F[+++t]----F[-t]+++F[+t]+++F[+++t]----F[-t]+++F/g; s/m/+++F[+++t]----F[+++t]-[---t]F-F[+++t]----F[-t][----t]+++F/g;' | ./bin 12 3 0 > tmp-pic.tex && make dotex # R7-5

%%%%%%%%%%%%%%%%%%%%%%%%%%%
%% with const char bullet[] = "{\\circle*{0.35}}";
%% with lnth *= 2.0;  // thicker lines
% stringsubst 3 F F F0F+F0F-F-F+F 0 0 + + - - | tail -1 | sed 's/+/p/g; s/-/m/g; s/Fp/+t++t+/g; s/Fm/--t--/g; s/F0/+t-t-t+/g; s/^/+++/g;' | ./bin 12 3 4 > tmp-pic.tex && make dotex
%
\begin{figure}[h!tbp]
\ifpdf
\begin{center}
{\includegraphics*[width=100mm, viewport={60 460 490 740}]{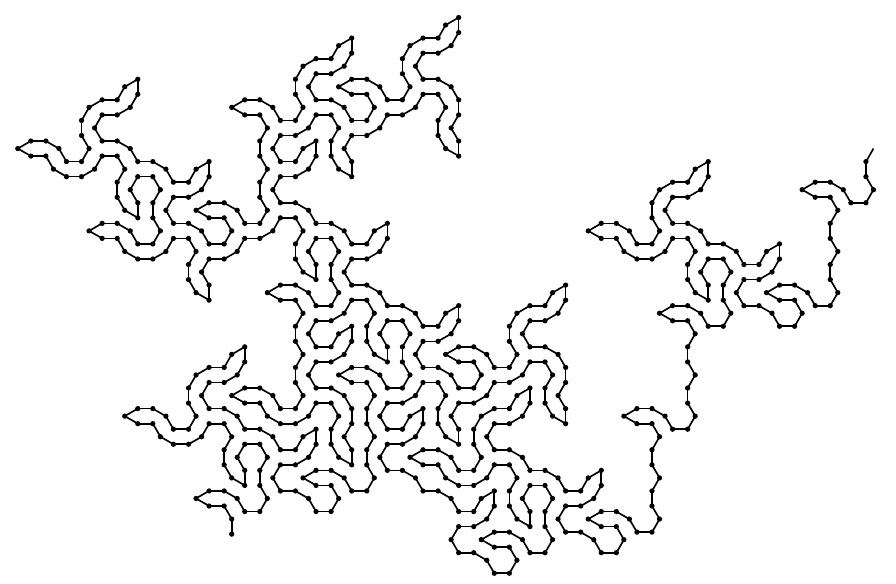}}%
\end{center}
\else
\verb+{see pdf for image}+
\fi
\caption{\label{fig:3464-t-cover-balanced}
$(3.4.6.4)$-PC curve from the third iterate of the balanced curve \CID{R7-1}.}
\end{figure}
%
%%%%%%%%%%%%%%%%%%%%%%%%%%%

For balanced curves,
$(3.4.6.4)$-PC curves are obtained by
replacing
all \texttt{+} by \texttt{p}
and all \texttt{-} by \texttt{m}, then
all \texttt{Fp} by \texttt{+F++F+},
all \texttt{Fm} by \texttt{--F--},
all \texttt{F0} by \texttt{+F-F-F+},
use turns by $30\adeg$.
The curve from \CID{R7-1} is shown in Figure~\ref{fig:3464-t-cover-balanced}.
%

%%%%%%%%%%%%%%%%%%%%%%%%%%%
% with lnth *= 2.0;  // thicker lines
% with const char bullet[] = "{\\circle*{0.4}}";
%  stringsubst 3 F F F+F-F-F+F+F-F + + - - | tail -1 | sed 's/F//g; s/+/p/g; s/-/m/g; s/p/+F+F+F+F/g; s/m/+F---F---F+F/g;' | ./bin 12 3 4 > tmp-pic.tex && make dotex # R7-5
%% better rotation:
% stringsubst 3 F F F+F-F-F+F+F-F + + - - | tail -1 | sed 's/F//g; s/+/p/g; s/-/m/g; s/p/+F+F+F+F/g; s/m/+F---F---F+F/g; s/^/--/;' | ./bin 12 3 4 > tmp-pic.tex && make dotex # R7-5
%
\begin{figure}[h!tbp]
\ifpdf
\begin{center}
{\includegraphics*[width=100mm, viewport={60 475 490 740}]{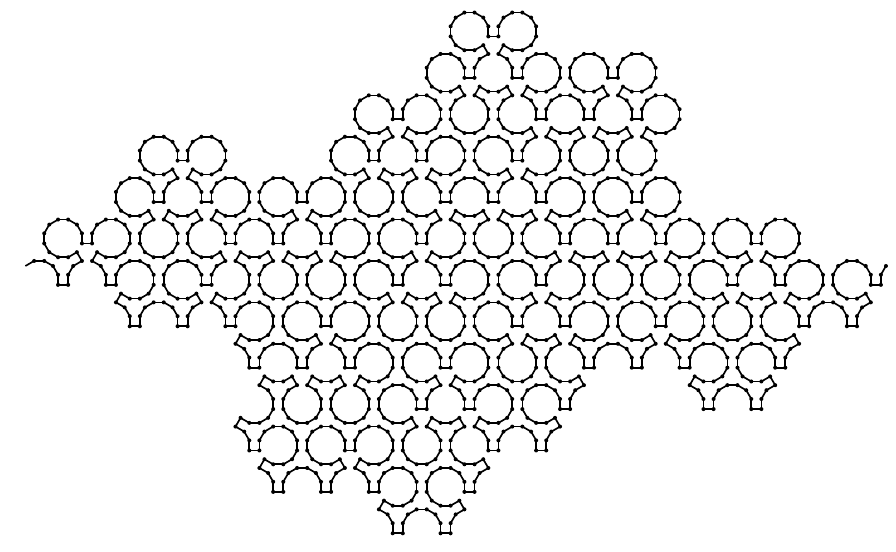}}
\end{center}
\else
\verb+{see pdf for image}+
\fi
\caption{\label{fig:4612-t-cover}
The $(4.6.12)$-PC curve from the third iterate of \CID{R7-5}.}
\end{figure}
%
%%%%%%%%%%%%%%%%%%%%%%%%%%%

For curves that are $(4.6.12)$-PC (Figure~\ref{fig:4612-t-cover})
drop all \texttt{F},
replace
all \texttt{+} by \texttt{p} and
all \texttt{-} by \texttt{m},
then replace
all \texttt{p} by \texttt{+F+F+F+F} and
all \texttt{m} by \texttt{+F---F---F+F}.
Use turns by $30\adeg$.

%% Alternative:
% stringsubst 2 F F F+F-F-F+F+F-F + + - - | tail -1 | sed 's/+/++t++/g; s/-/--t--/g; s/F/-F+F+F-/g;' | ./bin 12 3 4 0 0.3 > tmp-pic.tex && make dotex # R7-5

%% $(4.6.12)$-grid:
% stringsubst 3 F+F+F F F+F-F-F+F+F-F + + - - | tail -1 | sed 's/-/m/g; s/+/++[----t]t++/g; s/m/--t[++t]--/g; s/F/-F+F[--t]+F-/g;' | ./bin 12 3 0 > tmp-pic.tex && make dotex # R7-5

%%%%%%%%%%%%%%%%%%%%%%%%%%%
% with lnth *= 2.0;  // thicker lines
% with const char bullet[] = "{\\circle*{0.4}}";
% stringsubst 3 F F F0F+F0F-F-F+F 0 0 + + - - | tail -1 | sed 's/F+/F+t+F+t+/g; s/F-/F--t--/g; s/F0/F+t+F--t--F+t+/g;' | ./bin 12 3 4 > tmp-pic.tex && make dotex
%% better rotation:
% stringsubst 3 F F F0F+F0F-F-F+F 0 0 + + - - | tail -1 | sed 's/F+/F+t+F+t+/g; s/F-/F--t--/g; s/F0/F+t+F--t--F+t+/g; s/^/+++/g;' | ./bin 12 3 4 > tmp-pic.tex && make dotex
%
\begin{figure}[h!tbp]
\ifpdf
\begin{center}
{\includegraphics*[width=100mm, viewport={60 455 490 740}]{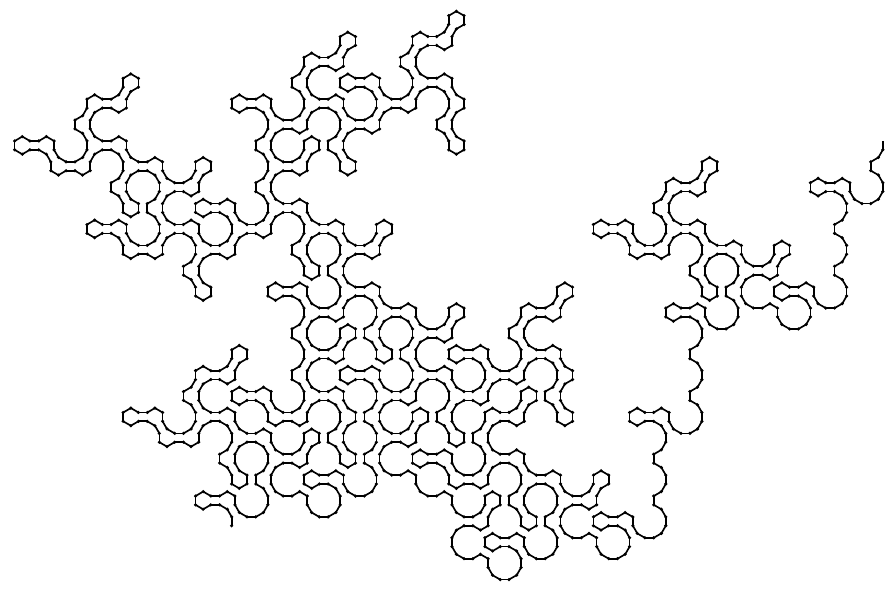}}
\end{center}
\else
\verb+{see pdf for image}+
\fi
\caption{\label{fig:4612-t-cover-balanced}
$(4.6.12)$-PC curve from the third iterate of the balanced curve \CID{R7-1}.}
\end{figure}
%
%%%%%%%%%%%%%%%%%%%%%%%%%%%

For balanced curves,
$(4.6.12)$-PC curves are obtained by replacing
all \texttt{F+} by \texttt{F+F+F+F+},
all \texttt{F-} by \texttt{F--F--}, and
all \texttt{F0} by \texttt{F+F+F--F--F+F+},
again using turns by $30\adeg$.
See figure~\ref{fig:4612-t-cover-balanced}
for the curve from the third iterate of \CID{R71-1}.

%%%%%%%%%%%%%%%%%%%%%%%%%%
% with const char bullet[] = "{\\circle*{0.3}}";
% stringsubst 2 F F F+F-F-F+F+F-F 0 0 + + - - | tail -1 | sed 's/+/F++/g; s/-/-F-/g;' | ./bin 6 3 4 > tmp-pic.tex && make dotex # R7-5
%
% stringsubst 2 F F F+F-F-F+F+F-F + + - - | tail -1 | sed 's/+/+t+/g; s/-/-t-/g; s/F/-t+/g;' | ./bin 6 3 4 > tmp-pic.tex && make dotex
%
\begin{figure}[h!tbp]
\ifpdf
\begin{center}
%{\includegraphics*[width=80mm, viewport={80 370 490 740}]{r07-t-5-33336-cover.pdf}}% third iterate
{\includegraphics*[width=64mm, viewport={80 510 520 735}]{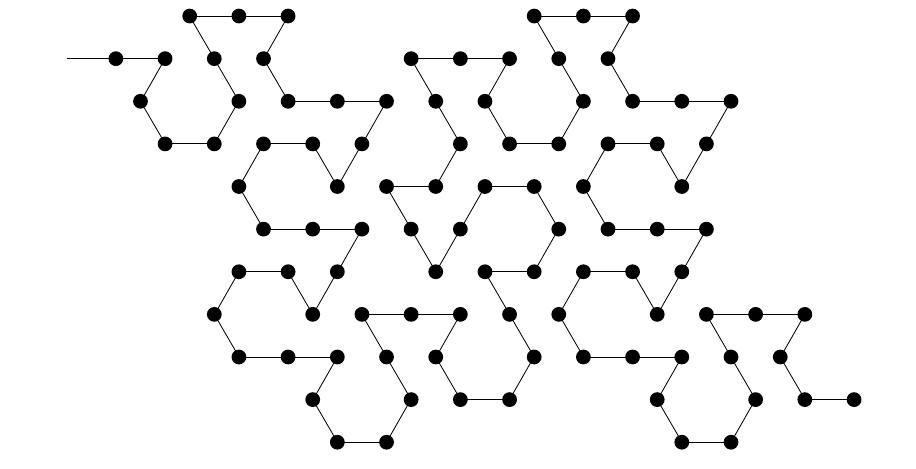}}%
{\includegraphics*[width=64mm, viewport={80 490 520 735}]{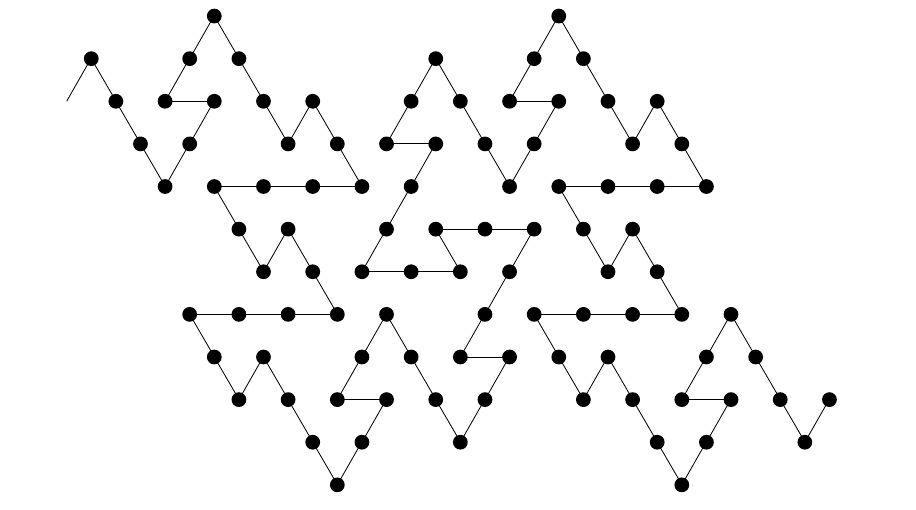}}
\end{center}
\else
\verb+{see pdf for image}+
\fi
\caption{\label{fig:r07-t-5-33336-cover}
Two $(3^4.6)$-PC curves from the second iterate of \CID{R7-5}.}
\end{figure}
%
%%%%%%%%%%%%%%%%%%%%%%%%%%%

Curves that are $(3^4.6)$-PC
are obtained by replacing
all \texttt{+} by \texttt{F++},
all \texttt{-} by \texttt{-F-},
and using turns by $60\adeg$.
Figure~\ref{fig:r07-t-5-33336-cover} (left) shows the curve from the second iterate.
For the curve for the other enantiomer of the grid,
change the first replacement to \texttt{++F}.
Replacing
all \texttt{+} by \texttt{+T+} and
all \texttt{-} by \texttt{-T-},
then
all \texttt{F} by \texttt{-T+},
drawing edges for \texttt{T}
using turns by $60\adeg$
gives the curve shown on the right of Figure~\ref{fig:r07-t-5-33336-cover}.
Change the second replacement to \texttt{+T-} for the other enantiomer.

% stringsubst 2 F F F+F-F-F+F+F-F + + - - | tail -1 | sed 's/+/+t+/g; s/-/-t-/g; s/F/-t+/g;' | ./bin 6 3 4 > tmp-pic.tex && make dotex
% stringsubst 2 F F F+F-F-F+F+F-F + + - - | tail -1 | sed 's/+/+t+/g; s/-/-t-/g; s/F/-F+/g;' | ./bin 6 3 4 > tmp-pic.tex && make dotex
%% Alternative:
% stringsubst 2 F F F+F-F-F+F+F-F + + - - | tail -1 | sed 's/+/+t+/g; s/-/-t-/g; s/F/+F-/g;' | ./bin 6 3 4 > tmp-pic.tex && make dotex # R7-5
%
% Alternatively (even simpler, compare to tri-hex cover):
% stringsubst 2 F F F+F-F-F+F+F-F 0 0 + + - - | tail -1 | sed 's/Z//g; s/+/F+/g; s/-/-F/g;' | ./bin 3 3 4 > tmp-pic.tex && make dotex # R7-5
%

%%%%%%%%%%%%%%%%%%%%%%%%%%
%% with const char bullet[] = "{\\circle*{0.2}}";
% stringsubst 2 _F_+F_+F _ _ F F+F-F-F+F+F-F 0 0 + + - - | tail -1 | sed 's/+F/+F+F/g; s/-F/-F-F/g;' |  ./bin 6 3 4 > tmp-pic.tex && make dotex # R7-5 tile-plus
%
% stringsubst 2 _F_+F_+F _ _ F F+F-F-F+F+F-F 0 0 + + - - | tail -1 | sed 's/+F/+F+F/g; s/-F/-F-F/g;' | ./hex-pc-to-33344-pc.pl | ./bin 12 3 4 > tmp-pic.tex && make dotex # R7-5 tile-plus
%
\begin{figure}[h!tbp]
\ifpdf
\begin{center}
{\includegraphics*[width=45mm, viewport={70 310 480 735}]{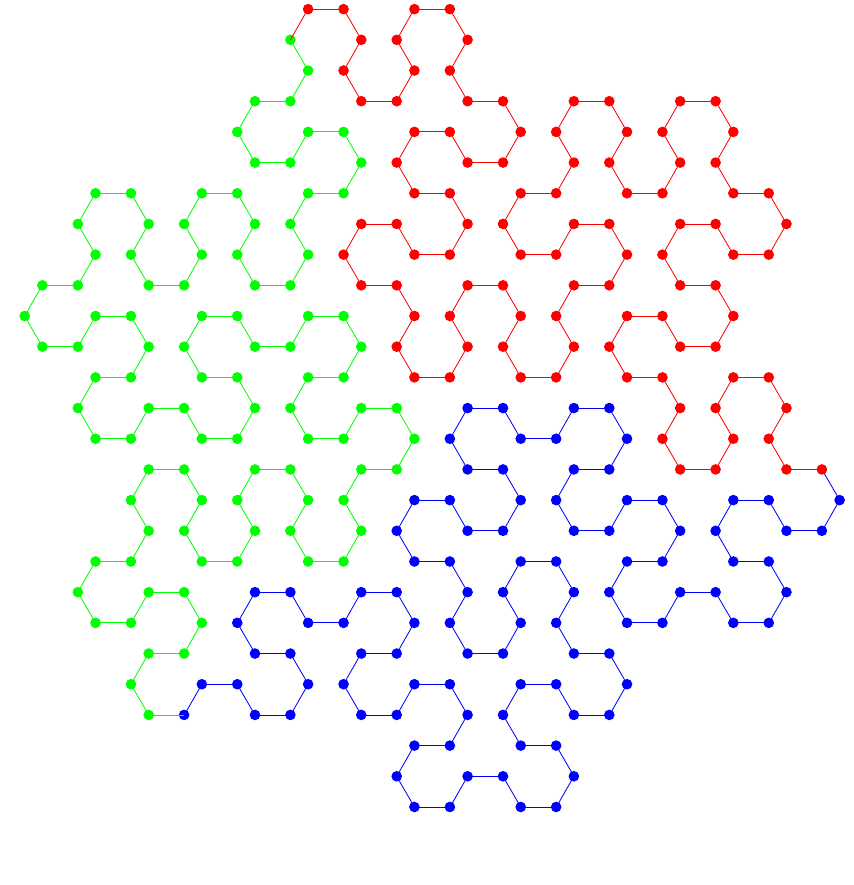}}%
{\includegraphics*[width=80mm, viewport={65 520 480 735}]{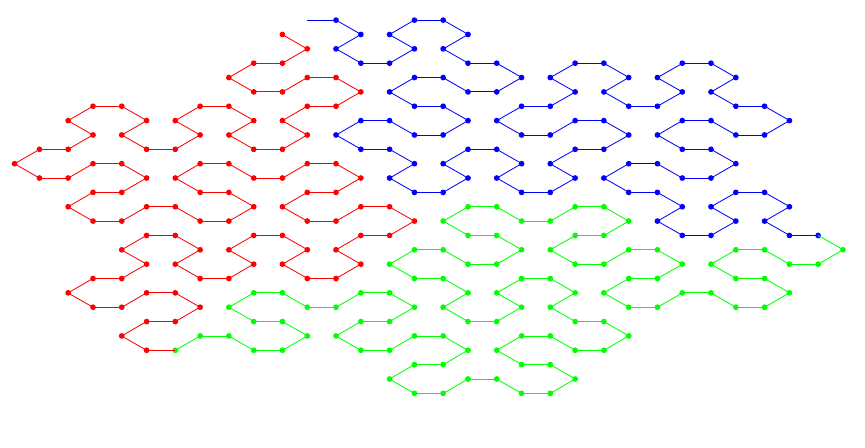}}
\end{center}
\else
\verb+{see pdf for image}+
\fi
\caption{\label{fig:33344-cover}
%The $(3^3.4^2)$-PC curve from the second iterate.}
The second iterate of the tile of \CID{R7-5} ($\Tile{+2}$, left)
and the $(3^3.4^2)$-PC curve derived from it by redirecting non-horizontal edges (right).}
\end{figure}
%
%%%%%%%%%%%%%%%%%%%%%%%%%%%

Figure~\ref{fig:33344-cover} shows the $(3^3.4^2)$-PC curve
obtained from the curve at the left
by adjusting the directions of the edges that are not horizontal.

%%%%%%%%%%%%%%%%%%%%%%%%%%
%% with    lnth *= 2.0;  // thicker lines
% with const char bullet[] = "{\\circle*{0.2}}";
%
%% shorter with new script:
% stringsubst 2 _+F_+F_+F _ _ F F+F-F-F+F+F-F + + - - | tail -1 | ./triangle-ec-to-square-pc.pl 2 | ./bin 4 3 4 0 0.15 > tmp-pic.tex && make dotex
% stringsubst 2 _+F_+F_+F _ _ F F+F-F-F+F+F-F + + - - | tail -1 | ./triangle-ec-to-square-pc.pl 2 | tr 34 45 | ./bin 6 3 4 0 0.15 > tmp-pic.tex && make dotex
%
%
\begin{figure}[h!tbp]
\ifpdf
\begin{center}
{\includegraphics*[width=60mm, viewport={40 280 470 720}]{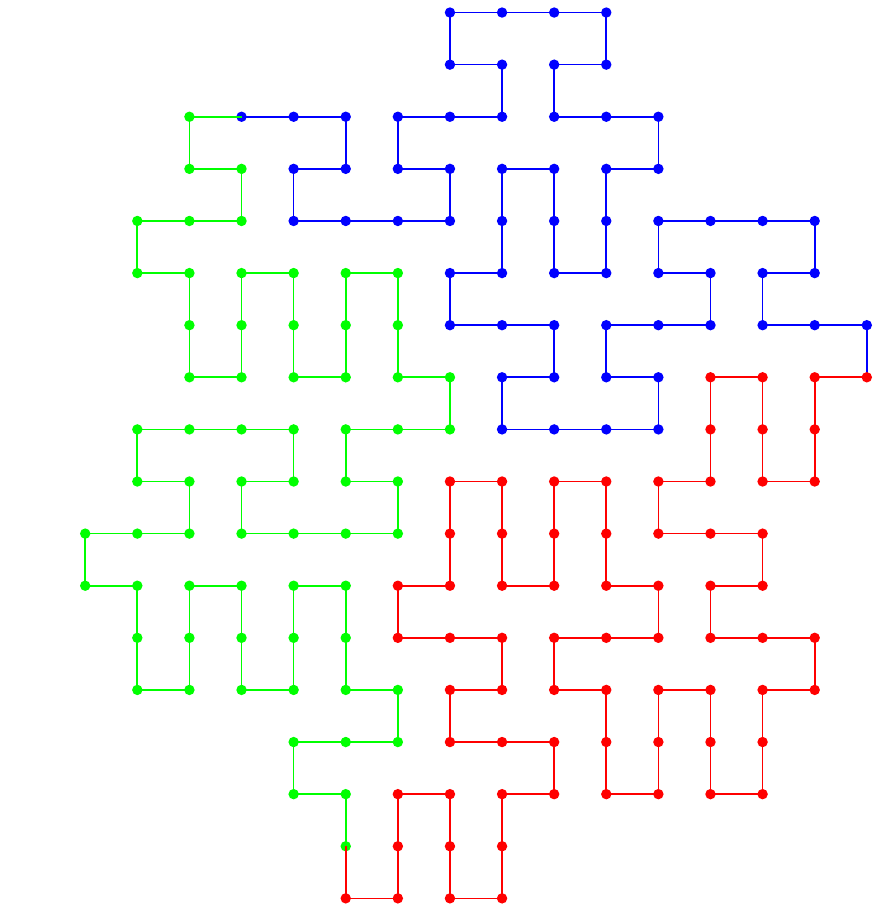}}%
{\includegraphics*[width=68mm, viewport={40 370 470 725}]{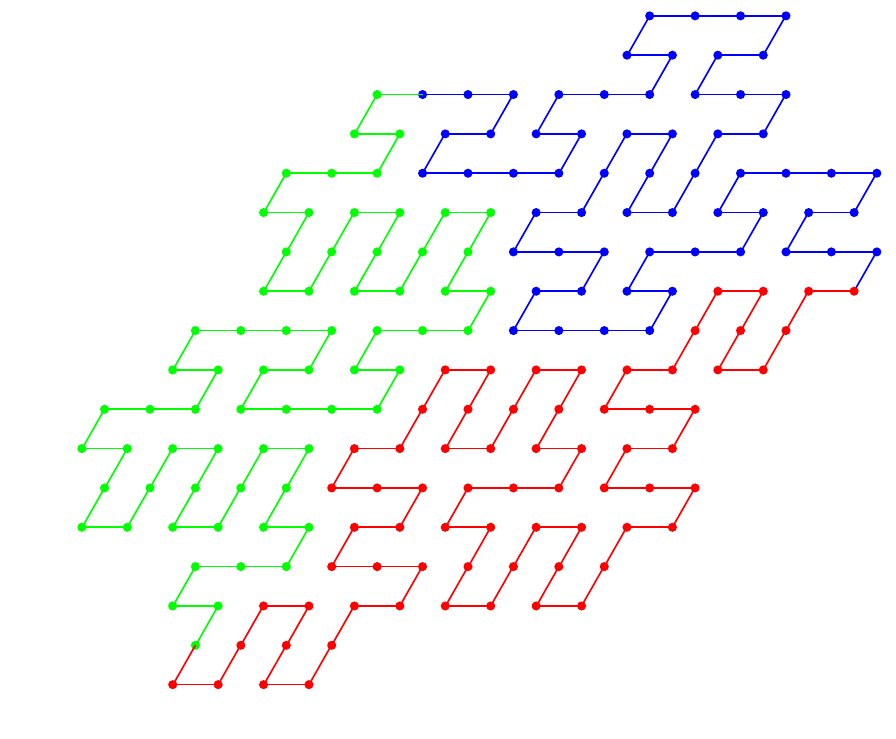}}
\end{center}
\else
\verb+{see pdf for image}+
\fi
\caption{\label{fig:r3-t-1-4444-cover}
%The $(6^3)$-PC curve from the fourth iterate of the terdragon tile ($\Tile{+4}$, left)
%and its rendering as $(4^4)$-PC curve (right).
%The $(6^3)$-PC curve from the second iterate of the tile ($\Tile{+2}$, left),
Rendering of the tile $\Tile{+2}$
as $(4^4)$-PC (left)
and $(3^6)$-PC curves (right).
}
\end{figure}
%
%%%%%%%%%%%%%%%%%%%%%%%%%%%

$(4^4)$-PC curves can be obtained by changing all non-horizontal edges
in a $(6^3)$-PC curve,
that is, mapping the directions 1, 2, \ldots, 6
shown at the right of Figure~\Ref{fig:directions}
respectively to directions 1, 2, 2, 3, 4, 4 on the square grid.
The resulting curves are stretched in the vertical direction.
Dropping the directions 3 and 5 first gives less distorted curves,
as shown at the left in Figure~\ref{fig:r3-t-1-4444-cover}.
Turning all vertical edges in the $(4^4)$-PC curve
gives the $(3^6)$-PC curve at the right.

%%%%%%%%%%%%%%%%%%%%%%%%%%
%% with    lnth *= 2.0;  // thicker lines
%
%% stringsubst 2 _F+_F+_F+ _ _ F F+F0F0F-F-F+F0F+F+F-F0F-F  0 0 + + - - | tail -1 | sed 's/F+/F+F+/g; s/F-/F--/g; s/F0/F+F--F+/g;' | ./4444-pc-filter | ./bin 12 3 0 > tmp-pic.tex && make dotex
%% stringsubst 2 _F-_F-_F- _ _ F F+F0F0F-F-F+F0F+F+F-F0F-F  0 0 + + - - | tail -1 | sed 's/F+/F+F+/g; s/F-/F--/g; s/F0/F+F--F+/g;' | ./4444-pc-filter | ./bin 12 3 0 > tmp-pic.tex && make dotex
%
\begin{figure}[h!tbp]
\ifpdf
\begin{center}
%{\includegraphics*[width=70mm, viewport={70 360 490 740}]{r13-t-15-4444-pc-tile-plus.pdf}}%
%{\includegraphics*[width=60mm, viewport={70 320 490 740}]{r13-t-15-4444-pc-tile-minus.pdf}}
{\includegraphics*[width=70mm, viewport={70 360 490 740}]{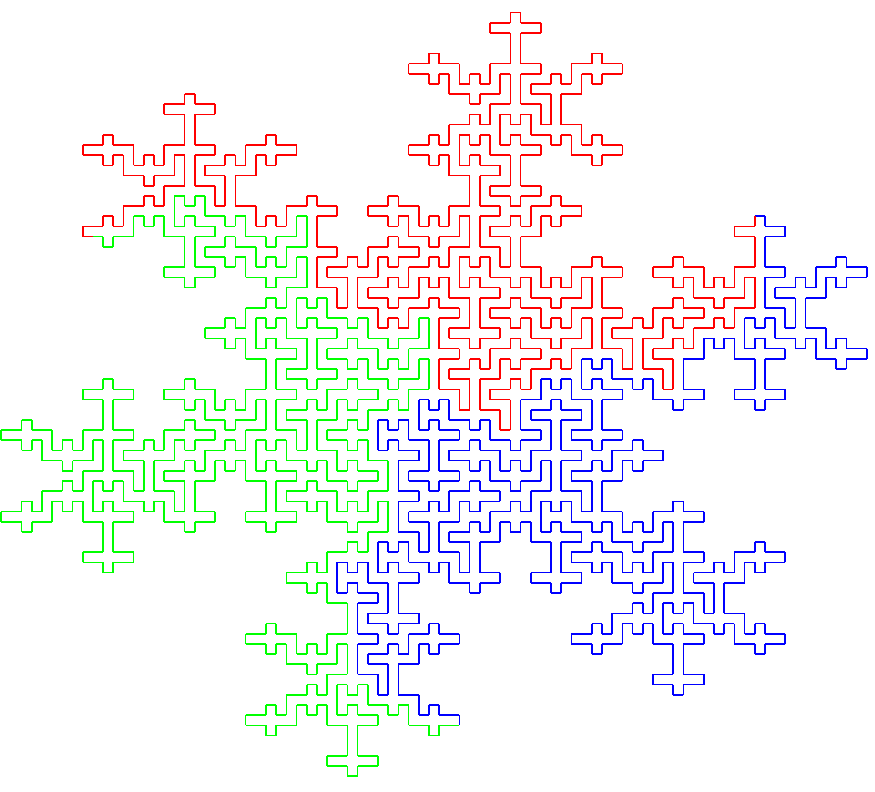}}%
{\includegraphics*[width=60mm, viewport={70 320 490 740}]{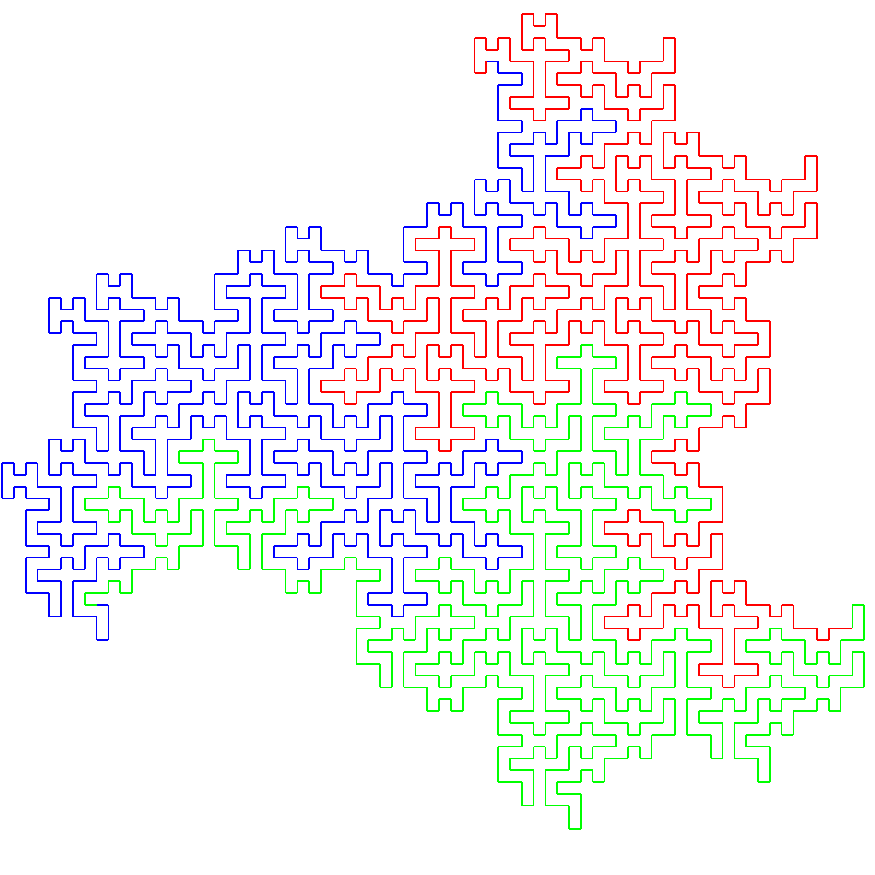}}
\end{center}
\else
\verb+{see pdf for image}+
\fi
\caption{\label{fig:r13-t-15-4444-pc-tile}
Point-covering curves on the square grid
from the tiles $\Tile{+2}$ (left) and $\Tile{-2}$ (right)
of the balanced curve \CID{R13-15} on the triangular grid.
}
\end{figure}
%
%%%%%%%%%%%%%%%%%%%%%%%%%%

%Curves of orders of the form $R=3k+1$
%that have exactly $k$ non-turns
Balanced curves
can be rendered as $(4^4)$-PC curves
such as shown in Figure~\ref{fig:r13-t-15-4444-pc-tile}
by first converting them to $(3.6.3.6)$-EC curves
as described in section \Ref{sect:3636-EC}
and then applying the conversion
used for Figure~\Ref{fig:r07-b-4444-cover}.
% that preserve the shape better

%\clearpage% xxx
%%%%%%%%%%%%%%%%%%%%%%%%%%%%%%%%%%%%%%%%%%%%%%%%%%%%
%%%%%%%%%%%%%%%%%%%%%%%%%%%%%%%%%%%%%%%%%%%%%%%%%%%%
\subsubsection{Square grid: curves for $(4.8.8)$, $(4^4)$, $(3.3.4.3.4)$, and $(3^6)$}

%%%%%%%%%%%%%%%%%%%%%%%%%%
%% with lnth *= 3.0;  // thick lines
% stringsubst 2 F  F F+F-F-F-F+F+F+F-F + + - -  | tail -1 | ./bin 4 2 0 0 0.0 > tmp-pic.tex && make dotex
% stringsubst 2 F  F F+F-F-F-F+F+F+F-F + + - -  | tail -1 | ./bin 4 2 0 0 0.1 > tmp-pic.tex && make dotex
% stringsubst 2 F  F F+F-F-F-F+F+F+F-F + + - -  | tail -1 | ./bin 4 2 0 0 0.2 > tmp-pic.tex && make dotex
% stringsubst 2 F  F F+F-F-F-F+F+F+F-F + + - -  | tail -1 | ./bin 4 2 0 0 0.3 > tmp-pic.tex && make dotex
% stringsubst 2 F  F F+F-F-F-F+F+F+F-F + + - -  | tail -1 | ./bin 4 2 0 0 0.4 > tmp-pic.tex && make dotex
% stringsubst 2 F  F F+F-F-F-F+F+F+F-F + + - -  | tail -1 | ./bin 4 2 0 0 0.5 > tmp-pic.tex && make dotex
%
\begin{figure}[h!tbp]
\ifpdf
\begin{center}
%{\includegraphics*[width=40mm, viewport={100 370 490 730}]{r09-q-rnd00.pdf}}%
%{\includegraphics*[width=40mm, viewport={100 370 490 730}]{r09-q-rnd01.pdf}}%
%{\includegraphics*[width=40mm, viewport={100 370 490 730}]{r09-q-rnd02.pdf}}
%%
%{\includegraphics*[width=40mm, viewport={100 370 490 730}]{r09-q-rnd03.pdf}}%
%{\includegraphics*[width=40mm, viewport={100 370 490 730}]{r09-q-rnd04.pdf}}%
%{\includegraphics*[width=40mm, viewport={100 370 490 730}]{r09-q-rnd05.pdf}}
%
{\includegraphics*[width=40mm, page=1, viewport={100 370 490 730}]{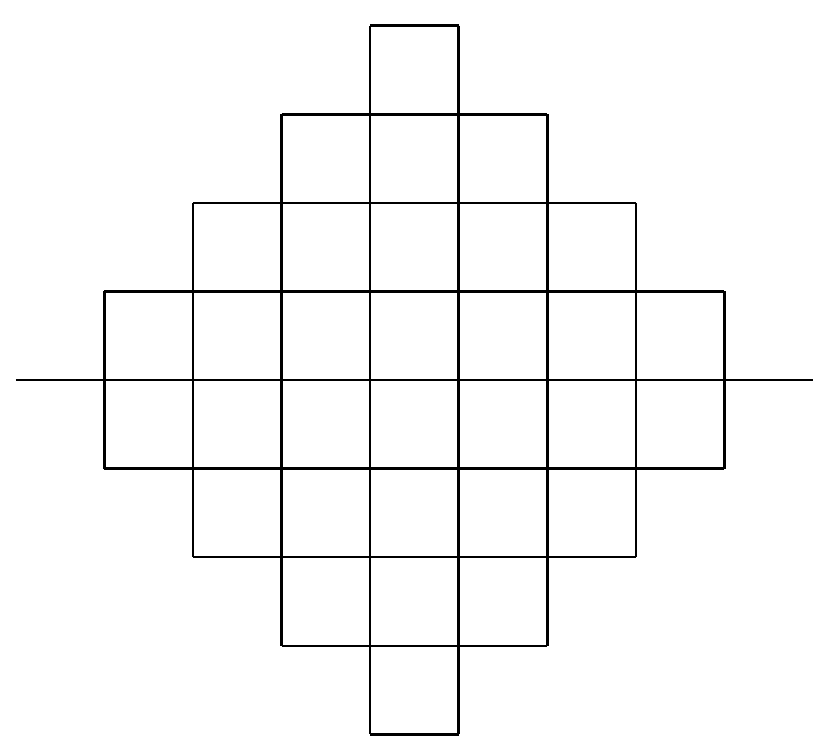}}%
{\includegraphics*[width=40mm, page=2, viewport={100 370 490 730}]{r09-q-rnd.pdf}}%
{\includegraphics*[width=40mm, page=3, viewport={100 370 490 730}]{r09-q-rnd.pdf}}
{\includegraphics*[width=40mm, page=4, viewport={100 370 490 730}]{r09-q-rnd.pdf}}%
{\includegraphics*[width=40mm, page=5, viewport={100 370 490 730}]{r09-q-rnd.pdf}}%
{\includegraphics*[width=40mm, page=6, viewport={100 370 490 730}]{r09-q-rnd.pdf}}
\end{center}
\else
\verb+{see pdf for image}+
\fi
\caption{\label{fig:r09-q-rnd}
Renderings of the curve \CID{R9-1}
with map \Lmap{F}{F+F-F-F-F+F+F+F-F}
on the square grid
for rounding parameter  $e\in\{0.0,\, 0.1,\, 0.2,\, 0.3,\, 0.4,\, 0.5\}$.}
\end{figure}
%
%%%%%%%%%%%%%%%%%%%%%%%%%%%

%%%%%%%%%%%%%%%%%%%%%%%%%%
% with const char bullet[] = "{\\circle*{0.3}}";
% stringsubst 2 F F F+F-F-F-F+F+F+F-F + + - - | tail -1 | sed 's/+F/+F+F/g; s/-F/-F-F/g;' | ./bin 8 3 4 > tmp-pic.tex && make dotex
%
% with const char bullet[] = "{\\circle*{0.2}}";
% stringsubst 2 F  F F+F-F-F-F+F+F+F-F + + - -  | tail -1 | sed 's/F//g; s/+/+F/g; s/-/F-/g;' | ./bin 4 3 4 > tmp-pic.tex && make dotex
% rotated version:
% stringsubst 2 RF R R  F F+F-F-F-F+F+F+F-F + + - -  | tail -1 | sed 's/F//g; s/+/++F/g; s/-/F--/g; s/R/-/g;' | ./bin 8 3 4 > tmp-pic.tex && make dotex
%
\begin{figure}[h!tbp]
\ifpdf
\begin{center}
{\includegraphics*[width=63mm, viewport={100 350 500 735}]{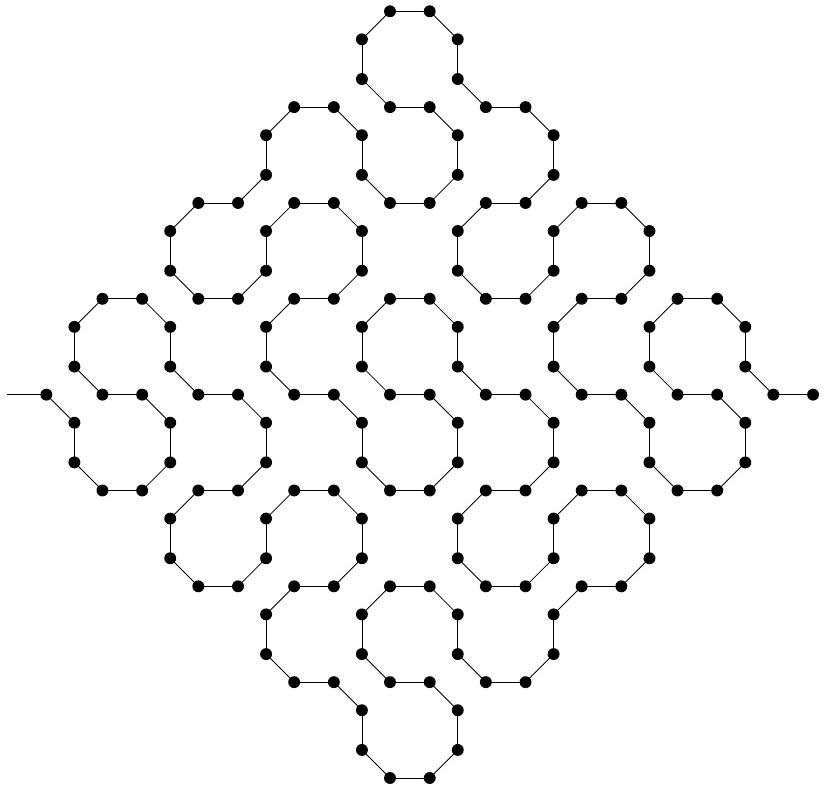}}%
{\includegraphics*[width=63mm, viewport={105 350 495 735}]{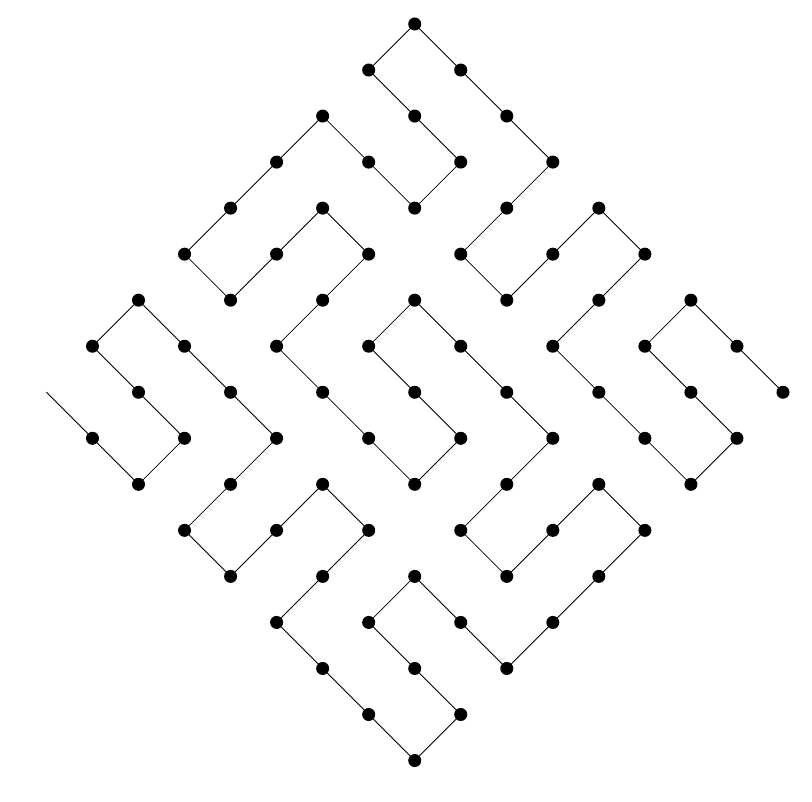}}
\end{center}
\else
\verb+{see pdf for image}+
\fi
\caption{\label{fig:r09-q-cover-points}
Points traversed by the curves corresponding to
$e=1/(2+\sqrt{2})$ $\approx{}0.292893\ldots$, a $(4.8.8)$-PC curve (left),
and $e=1/2$, a $(4^4)$-PC curve (right).}
\end{figure}
%
%%%%%%%%%%%%%%%%%%%%%%%%%%%

Figure~\ref{fig:r09-q-rnd} shows renderings
of the curve \CID{R9-1} with map \Lmap{F}{F+F-F-F-F+F+F+F-F}
for rounding parameter  $e\in\{0.0,\, 0.1,\, 0.2,\, 0.3,\, 0.4,\, 0.5\}$.
This curve is sometimes referred to
as Peano curve, see \cite[Chapter~3, pp.~31ff]{sagan}.

The curve for $e=1/3$ is \emph{almost} a $(4.8.8)$-PC curve:
%curve traversing the points of the grid (4.8.8) once:
for all edges to be of equal length,
the remaining edges must have the length of
the edges shortening the corners,
so we must have $1-2\,e=\sqrt{2}\,e$,
which gives $e=1/(2+\sqrt{2})$ $\approx{}0.292893\ldots$

For $e=1/2$ (lower right) we obtain a $(4^4)$-PC curve,
%traverses every point of the square grid once,
the points are the mid-points of the edges of the original curve.

The curves can again be obtained by post-processing steps.
For $e=1/(2+\sqrt{2})\approx{}1/3$ replace
all \texttt{+F} by \texttt{+F+F},
all \texttt{-F} by \texttt{-F-F},
and use turns by $45\adeg$.
%
% Check: round 1/0.292893.. via postprocessing:
% stringsubst 2 F F F+F-F-F-F+F+F+F-F + + - - | tail -1 | sed 's/F//g; s/+/+F+F/g; s/-/-F-F/g;' | ./bin 8 3 0 > tmp-pic.tex && make dotex
%
For $e=1/2$ drop all \texttt{F}, replace
all \texttt{+} by \texttt{+F+},
all \texttt{-} by \texttt{-F-},
and use turns by $90\adeg$.
%
% Check: round 1/2 via postprocessing:
% stringsubst 2 F  F F+F-F-F-F+F+F+F-F + + - -  | tail -1 | sed 's/F//g; s/+/+F/g; s/-/F-/g;' | ./bin 4 3 0 > tmp-pic.tex && make dotex
%
Both curves are shown in Figure~\ref{fig:r09-q-rnd}.

%%%%%%%%%%%%%%%%%%%%%%%%%%
% with const char bullet[] = "{\\circle*{0.2}}";
% stringsubst 2 F F F+F-F-F-F+F+F+F-F + + - - | tail -1 | tr -d F | sed 's/+/++F+/g; s/-/--F-/g;' | ./bin 12 3 4 > tmp-pic.tex && make dotex
% with const char bullet[] = "{\\circle*{0.2}}";
% stringsubst 2 F F F+F-F-F-F+F+F+F-F + + - - | tail -1 | sed 's/F//g; s/+/+F/g; s/-/F-/g;' | ./square-pc-to-triangle-pc.pl | ./bin 6 3 4 > tmp-pic.tex && make dotex
%
\begin{figure}[h!tbp]
\ifpdf
\begin{center}
{\includegraphics*[width=63mm, viewport={110 360 490 730}]{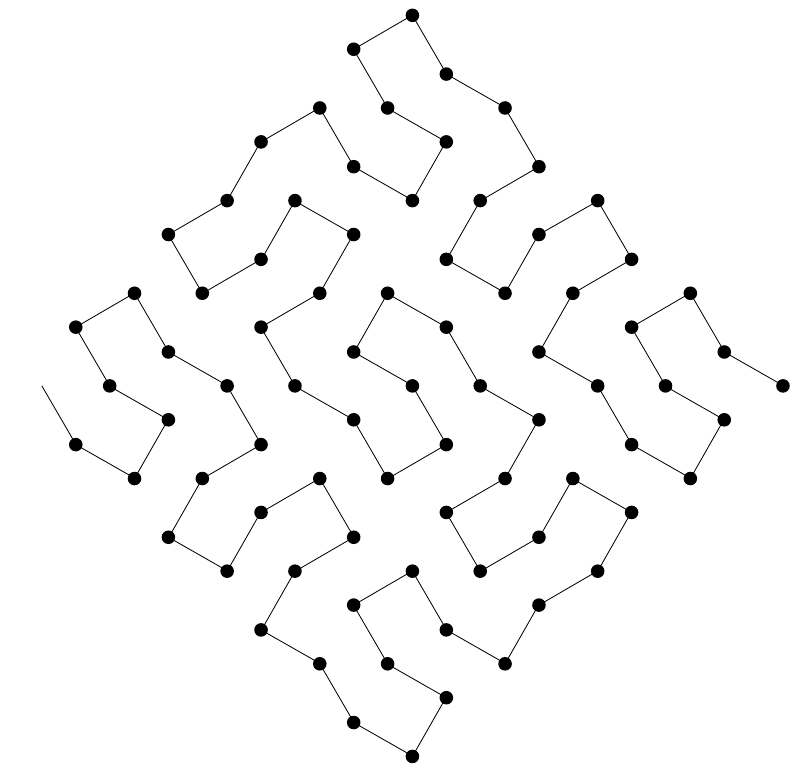}}%
{\includegraphics*[width=63mm, viewport={110 420 490 740}]{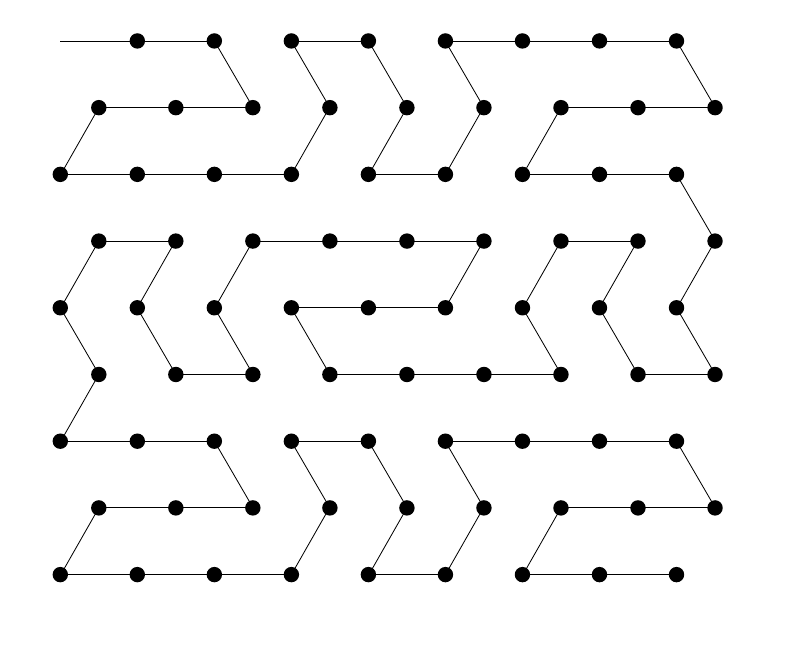}}
\end{center}
\else
\verb+{see pdf for image}+
\fi
\caption{\label{fig:r09-q-33434-cover}
The $(3.3.4.3.4)$-PC and the $(3^6)$-PC curves from the second iterate.}
\end{figure}
%
%%%%%%%%%%%%%%%%%%%%%%%%%%%

For the $(3.3.4.3.4)$-PC curve
shown on the left in Figure~\ref{fig:r09-q-33434-cover},
drop all \texttt{F}, replace
all \texttt{+} by \texttt{++F+},
all \texttt{-} by \texttt{--F-},
and use turns by $30\adeg$.
%
% Check: $(3.3.4.3.4)$
% stringsubst 2 F F F+F-F-F-F+F+F+F-F + + - - | tail -1 | tr -d F | sed 's/+/++F+/g; s/-/--F-/g;' | ./bin 12 3 4 > tmp-pic.tex && make dotex
%
% Alternatively:
% stringsubst 2 F F F+F-F-F-F+F+F+F-F + + - - | tail -1 | sed 's/F//g; s/+/++F+/g; s/-/F---/g;' | ./bin 12 3 4 0 0.15 > tmp-pic.tex && make dotex #
% stringsubst 2 F F F+F-F-F-F+F+F+F-F + + - - | tail -1 | sed 's/F//g; s/+/+F++/g; s/-/---F/g;' | ./bin 12 3 4 0 0.15 > tmp-pic.tex && make dotex #

The $(3^6)$-PC curve shown on the right in Figure~\ref{fig:r09-q-33434-cover}
results from rendering the $(4^4)$-PC curve
(right of Figure~\ref{fig:r09-q-cover-points})
such that one set of parallel edges (horizontal in the image) is kept
and the other (vertical) edges are alternatingly
turned by $\pm{}30\adeg$ against their original direction.
%
%Leaving every second row of vertical edges unchanged
%gives another $(3^3.4^2)$-PC curve.\xxx{Create image?}

%\clearpage% xxx
%%%%%%%%%%%%%%%%%%%%%%%%%%%%%%%%%%%%%%%%%%%%%%%%%%%%
%%%%%%%%%%%%%%%%%%%%%%%%%%%%%%%%%%%%%%%%%%%%%%%%%%%%
\subsubsection{Tri-hexagonal grid: curves for $(4.6.12)$, $(3.4.6.4)$, $(3^4.6)$, and $(4^4)$}\label{sect:pc-from-tri-hex}

%%%%%%%%%%%%%%%%%%%%%%%%%%
%% with lnth *= 3.0;  // thick lines
% stringsubst 2 F F F+F--F--F+F+F+F - - + + | tail -1 | ./bin 6 2 0 0 0.0 > tmp-pic.tex && make dotex # R07-b-1
% stringsubst 2 F F F+F--F--F+F+F+F - - + + | tail -1 | ./bin 6 2 0 0 0.1 > tmp-pic.tex && make dotex # R07-b-1
% stringsubst 2 F F F+F--F--F+F+F+F - - + + | tail -1 | ./bin 6 2 0 0 0.2 > tmp-pic.tex && make dotex # R07-b-1
% stringsubst 2 F F F+F--F--F+F+F+F - - + + | tail -1 | ./bin 6 2 0 0 0.3 > tmp-pic.tex && make dotex # R07-b-1
% stringsubst 2 F F F+F--F--F+F+F+F - - + + | tail -1 | ./bin 6 2 0 0 0.4 > tmp-pic.tex && make dotex # R07-b-1
% stringsubst 2 F F F+F--F--F+F+F+F - - + + | tail -1 | ./bin 6 2 0 0 0.5 > tmp-pic.tex && make dotex # R07-b-1
%
\begin{figure}[h!tbp]
\ifpdf
\begin{center}
%{\includegraphics*[width=30mm, viewport={110 370 440 730}]{r07-b-rnd00.pdf}}%
%{\includegraphics*[width=30mm, viewport={110 370 440 730}]{r07-b-rnd01.pdf}}%
%{\includegraphics*[width=30mm, viewport={110 370 440 730}]{r07-b-rnd02.pdf}}
%%
%{\includegraphics*[width=30mm, viewport={110 370 440 730}]{r07-b-rnd03.pdf}}%
%{\includegraphics*[width=30mm, viewport={110 370 440 730}]{r07-b-rnd04.pdf}}%
%{\includegraphics*[width=30mm, viewport={110 370 440 730}]{r07-b-rnd05.pdf}}
%%
{\includegraphics*[width=30mm, page=1, viewport={110 370 440 730}]{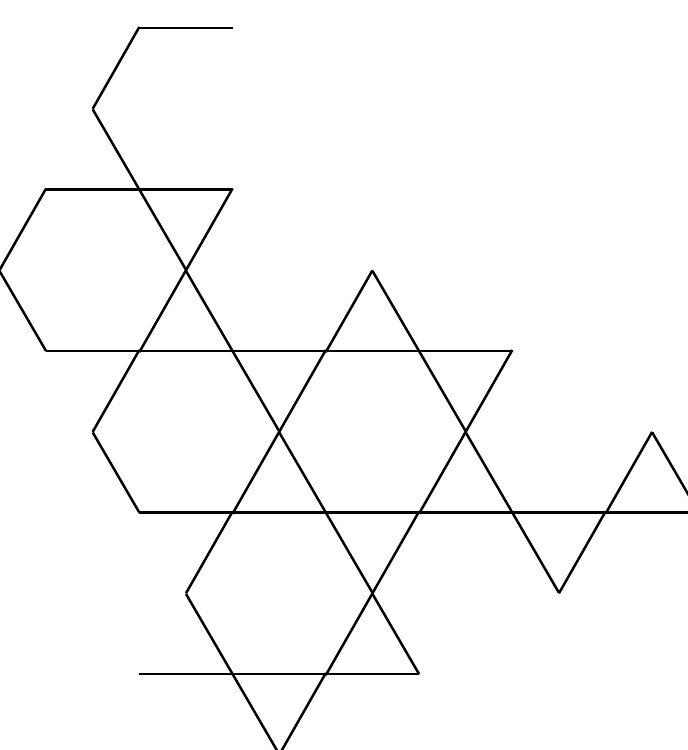}}%
\quad% layout
{\includegraphics*[width=30mm, page=2, viewport={110 370 440 730}]{r07-b-rnd.pdf}}%
\quad% layout
{\includegraphics*[width=30mm, page=3, viewport={110 370 440 730}]{r07-b-rnd.pdf}}
% line-break: xxx redo this one

{\includegraphics*[width=30mm, page=4, viewport={110 370 440 730}]{r07-b-rnd.pdf}}%
\quad% layout
{\includegraphics*[width=30mm, page=5, viewport={110 370 440 730}]{r07-b-rnd.pdf}}%
\quad% layout
{\includegraphics*[width=30mm, page=6, viewport={110 370 440 730}]{r07-b-rnd.pdf}}
\end{center}
\else
\verb+{see pdf for image}+
\fi
\caption{\label{fig:r07-b-rnd}
Renderings of the second iterate of Ventrella's curve (\CID{R7-1})
%with map \Lmap{F}{F+F--F--F+F+F+F}
%on the tri-hexagonal grid
for rounding parameter $e\in\{0.0,\, 0.1,\, 0.2,\, 0.3,\, 0.4,\, 0.5\}$.}
\end{figure}
%
%%%%%%%%%%%%%%%%%%%%%%%%%%%

%%%%%%%%%%%%%%%%%%%%%%%%%%
% with const char bullet[] = "{\\circle*{0.3}}";
%% stringsubst 2 F F F+F--F--F+F+F+F - - + + | tail -1 | sed 's/--F/+X/g; s/F/Y/g; s/Y/F+F+F+F+F+F--F--F--F--F+F+F+F+F/g; s/X/F+F+F+F+F+F--F--F--F--F+F+F--F--F/g;' | ./bin 12 3 4 > tmp-pic.tex && make dotex
%% much shorter:
% stringsubst 3 F F F+F--F--F+F+F+F - - + + | tail -1 | sed 's/F//g; s/+/p/g; s/--/m/g;  s/p/F+F+/g; s/m/F--F--/g;' | ./bin 12 3 4 > tmp-pic.tex && make dotex #
%% still better:
% stringsubst 3 F F F+F--F--F+F+F+F - - + + | tail -1 | sed 's/F//g;  s/+/F+F+/g; s/--/F--F--/g;' | ./bin 12 3 4 > tmp-pic.tex && make dotex
%
\begin{figure}[h!tbp]
\ifpdf
\begin{center}
%{\includegraphics*[width=80mm, viewport={70 290 490 740}]{r07-b-4612-cover.pdf}}
{\includegraphics*[width=80mm, viewport={70 410 490 740}]{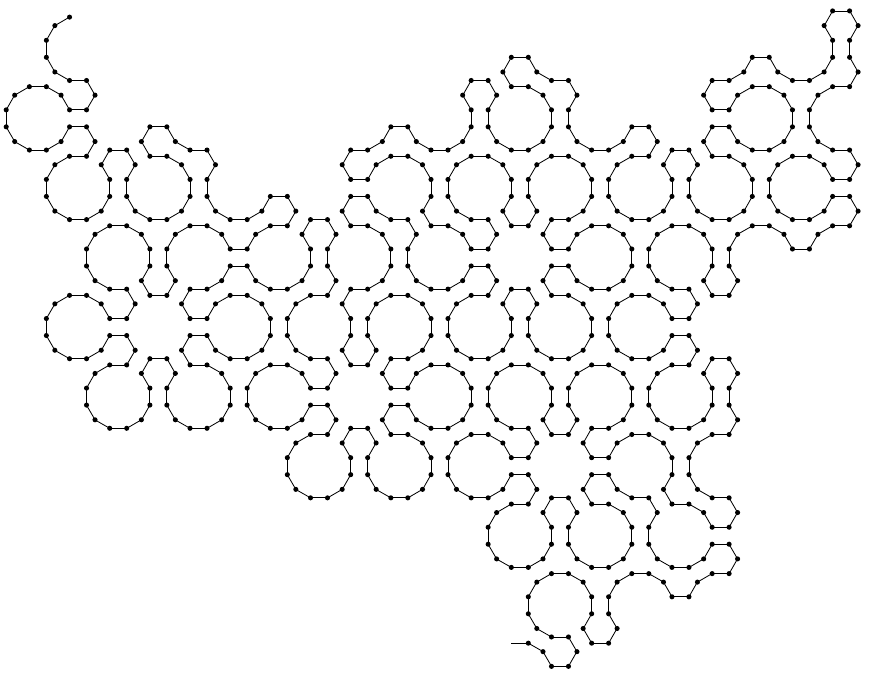}}
\end{center}
\else
\verb+{see pdf for image}+
\fi
\caption{\label{fig:r07-b-4612-cover}
A $(4.6.12)$-PC curve from the third iterate.}
\end{figure}
%
%%%%%%%%%%%%%%%%%%%%%%%%%%%

For our examples
we use Ventrella's curve (\CID{R7-1})
with the map \Lmap{F}{F+F--F--F+F+F+F}
on the tri-hexagonal grid
already seen in Figure~\Ref{fig:r7-b-curve}.
The renderings of the second iterate
for the rounding parameter $e\in\{0.0,\, 0.1,\, 0.2,\, 0.3,\, 0.4,\, 0.5\}$
are shown in Figure~\ref{fig:r07-b-rnd}.

The rendering for $e=1/3$ (lower left in Figure~\ref{fig:r07-b-rnd})
shows a distorted $(4.6.12)$-PC curve.
For the (undistorted) curve
%
%% Old construction:
%Replace all \texttt{--F} by \texttt{+X},
%then all \texttt{F} by \texttt{Y},
%then all \texttt{X} by \texttt{F+F+F+F+F+F--F--F--F--F+F+F--F--F}
%and all \texttt{Y} by \texttt{F+F+F+F+F+F--F--F--F--F+F+F+F+F}.
%Use turns by $30\adeg$.
%
drop all \texttt{F},
replace
all \texttt{+} by \texttt{F+F+}
and
all \texttt{--} by \texttt{F--F--},
use turns by $30\adeg$.
Figure~\ref{fig:r07-b-4612-cover}
shows the curve obtained from the third iterate.
%
% Check: (4.6.12) cover
% stringsubst 3 F F F+F--F--F+F+F+F - - + + | tail -1 | sed 's/F//g; s/+/p/g; s/--/m/g;  s/p/F+F+/g; s/m/F--F--/g;' | ./bin 12 3 4 > tmp-pic.tex && make dotex
%

%%%%%%%%%%%%%%%%%%%%%%%%%%
% with const char bullet[] = "{\\circle*{0.4}}";
%% stringsubst 3 F F F+F--F--F+F+F+F - - + + | tail -1 | sed 's/F//g; s/+/p/g; s/--/m/g; s/p/--F+++F+++F--/g; s/m/--FFF--/g; s/F----F/-F++F---F---F++F-/g;' | ./bin 12 3 4 > tmp-pic.tex && make dotex
%% Cleaner, also get first stroke right:
% stringsubst 3 F F F+F--F--F+F+F+F - - + + | tail -1 | sed 's/+/p/g; s/--/m/g; s/p/---F++F++F++F++F---/g; s/m/---F++F-F-F++F---/g; ' | ./bin 12 3 4 > tmp-pic.tex && make dotex
%
\begin{figure}[h!tbp]
\ifpdf
\begin{center}
%{\includegraphics*[width=80mm, viewport={70 290 490 740}]{r07-b-4612-cover-alt.pdf}}
{\includegraphics*[width=100mm, viewport={60 405 490 740}]{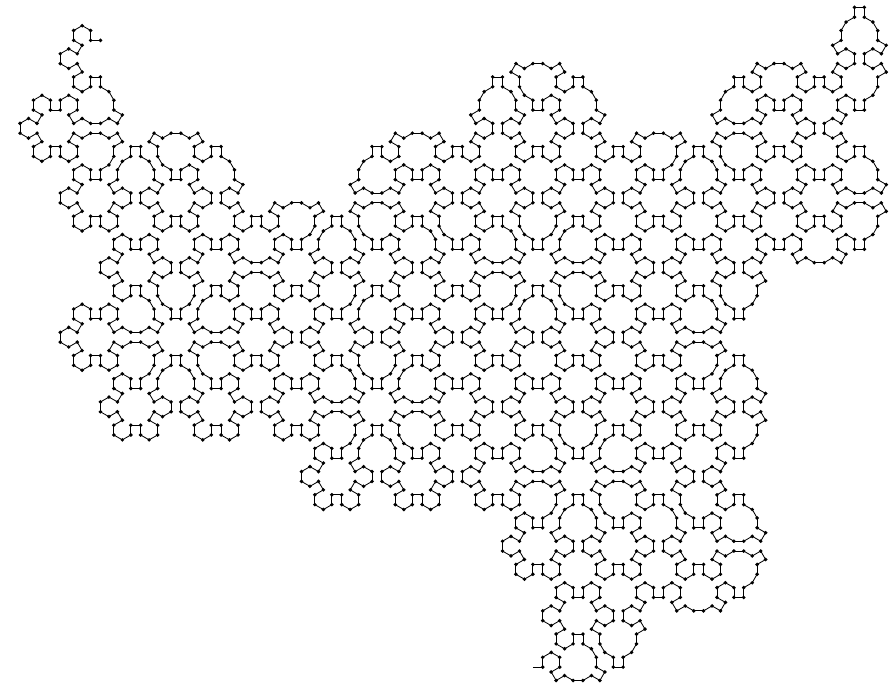}}
\end{center}
\else
\verb+{see pdf for image}+
\fi
\caption{\label{fig:r07-b-4612-cover-alt}
Another $(4.6.12)$-PC curve from the third iterate.}
\end{figure}
%
%%%%%%%%%%%%%%%%%%%%%%%%%%%

A different $(4.6.12)$-PC curve is shown in Figure~\ref{fig:r07-b-4612-cover-alt},
for the curve
%drop all \texttt{F},
replace
all \texttt{+} by \texttt{p}
and
all \texttt{--} by \texttt{m},
then
replace
all \texttt{p} by \texttt{---F++F++F++F++F---}
and
all \texttt{m} by \texttt{---F++F-F-F++F---},
use turns by $30\adeg$.

%%%%%%%%%%%%%%%%%%%%%%%%%%
% with const char bullet[] = "{\\circle*{0.3}}";
% stringsubst 3 F F F+F+F+F--F--F+F + + - - | tail -1 | sed 's/F//g; s/+/+F+/g; s/--/--F--/g;' | ./bin 12 3 4 > tmp-pic.tex && make dotex
%
\begin{figure}[h!tbp]
\ifpdf
\begin{center}
{\includegraphics*[width=80mm, viewport={70 410 490 735}]{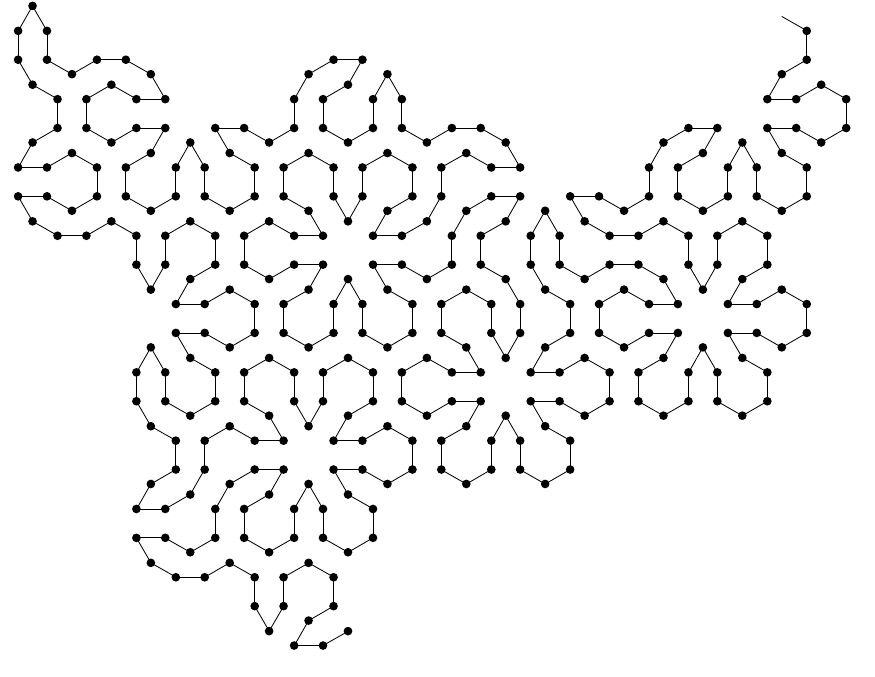}}
\end{center}
\else
\verb+{see pdf for image}+
\fi
\caption{\label{fig:r07-b-3464-cover}
%The $(3.4.6.4)$-PC curve from the second iterate.}
The $(3.4.6.4)$-PC curve from the third iterate.}
\end{figure}
%
%%%%%%%%%%%%%%%%%%%%%%%%%%%

The rendering $e=1/2$ (lower right in Figure~\ref{fig:r07-b-rnd}) similarly
suggests the following construction of a $(3.4.6.4)$-PC curve,
shown in Figure~\ref{fig:r07-b-3464-cover}.
Drop all \texttt{F},
replace
all \texttt{+} by \texttt{+F+}
and
all \texttt{--} by \texttt{--F--},
and
use turns by $30\adeg$.
% Check:
% stringsubst 3 F F F+F+F+F--F--F+F + + - - | tail -1 | sed 's/F//g; s/+/+F+/g; s/--/--F--/g;' | ./bin 12 3 4 > tmp-pic.tex && make dotex

%\FloatBarrier% xxx

%%%%%%%%%%%%%%%%%%%%%%%%%%
% with const char bullet[] = "{\\circle*{0.3}}";
% stringsubst 3 F F F+F--F--F+F+F+F - - + + | tail -1 | sed 's/F//g; s/+/p/g; s/--/m/g; s/p/--F++++F++F--/g; s/m/--FFF--/g;' | ./bin 12 3 4 > tmp-pic.tex && make dotex
%
\begin{figure}[h!tbp]
\ifpdf
\begin{center}
{\includegraphics*[width=100mm, viewport={70 390 490 740}]{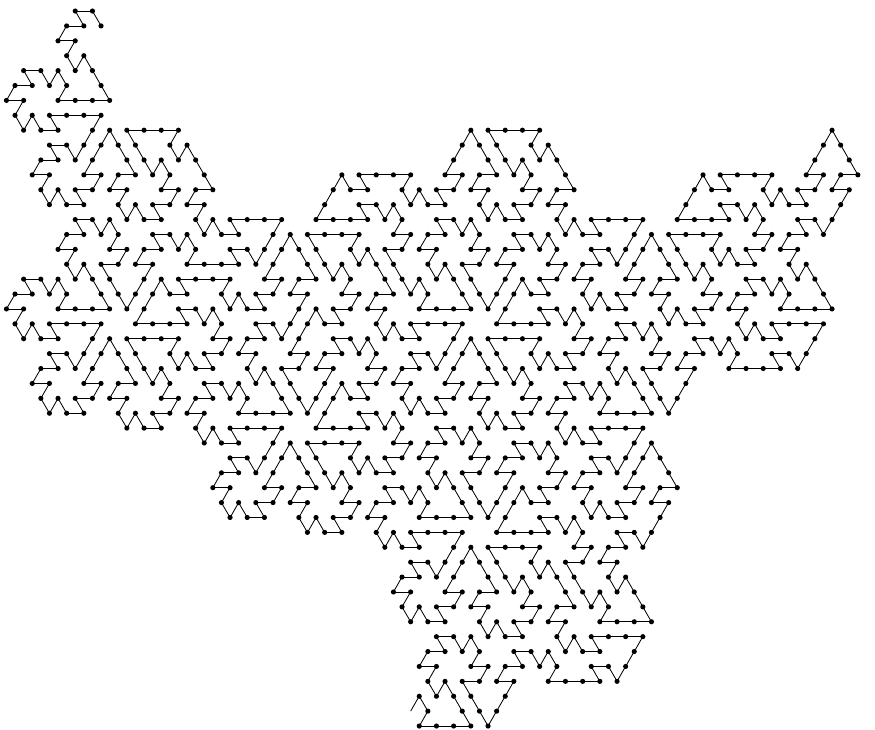}}
\end{center}
\else
\verb+{see pdf for image}+
\fi
\caption{\label{fig:r07-b-33336-cover}
A $(3^4.6)$-PC curve from the third iterate.}
\end{figure}
%
%%%%%%%%%%%%%%%%%%%%%%%%%%%

For the $(3^4.6)$-PC curve shown in Figure~\ref{fig:r07-b-33336-cover},
drop all \texttt{F},
replace
all \texttt{+} by \texttt{p}
and
all \texttt{--} by \texttt{m},
then
replace
all \texttt{p} by \texttt{--F++++F++F--}
and
all \texttt{m} by \texttt{--FFF--},
use turns by $30\adeg$.
%
% Check:
% stringsubst 3 F F F+F--F--F+F+F+F - - + + | tail -1 | sed 's/F//g; s/+/p/g; s/--/m/g; s/p/--F++++F++F--/g; s/m/--FFF--/g;' | ./bin 12 3 4 > tmp-pic.tex && make dotex

%%%%%%%%%%%%%%%%%%%%%%%%%%
%% with const char bullet[] = "{\\circle*{0.3}}";
%% with lnth *= 2.0;  // thicker lines
% stringsubst 3 F F F+F+F+F--F--F+F + + - - | tail -1 | sed 's/F//g; s/+/+F/g; s/--/-F-/g;' | ./bin 6 3 4 > tmp-pic.tex && make dotex
%% Better rotation (nope, use above):
% stringsubst 3 F F F+F+F+F--F--F+F + + - - | tail -1 | sed 's/F//g; s/+/++F/g; s/--/--F--/g; s/^/-/;' | ./bin 12 3 4 > tmp-pic.tex && make dotex
%
\begin{figure}[h!tbp]
\ifpdf
\begin{center}
{\includegraphics*[width=70mm, viewport={60 380 490 735}]{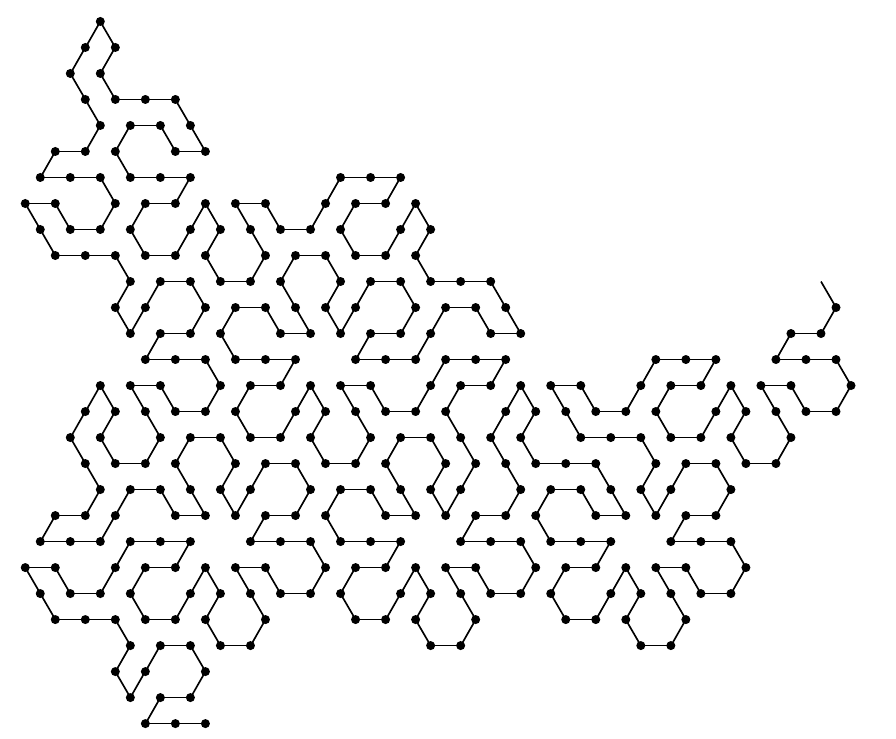}}
\end{center}
\else
\verb+{see pdf for image}+
\fi
\caption{\label{fig:r07-b-33336-cover-alt}
Another $(3^4.6)$-PC curve from the third iterate.}
\end{figure}
%
%%%%%%%%%%%%%%%%%%%%%%%%%%%

For the $(3^4.6)$-PC curve shown in Figure~\ref{fig:r07-b-33336-cover-alt},
drop all \texttt{F},
replace
all \texttt{+} by \texttt{+F}
and
all \texttt{--} by \texttt{-F-},
and
use turns by $60\adeg$.
%% Check:
% stringsubst 3 F F F+F+F+F--F--F+F + + - - | tail -1 | sed 's/F//g; s/+/+F/g; s/--/-F-/g;' | ./bin 6 3 4 > tmp-pic.tex && make dotex #
%
%% Alternative:
% stringsubst 3 F F F+F+F+F--F--F+F + + - - | tail -1 | sed 's/F--/--F/g;' | ./bin 6 3 4 > tmp-pic.tex && make dotex
%
%% Alternative:
% stringsubst 3 F F F+F+F+F--F--F+F + + - - | tail -1 | sed 's/F+/+F/g;' | ./bin 6 3 4 > tmp-pic.tex && make dotex
%
%% Alternative:
% stringsubst 3 F F F+F+F+F--F--F+F + + - - | tail -1 | sed 's/F+/++F-/g;' | ./bin 6 3 4 > tmp-pic.tex && make dotex
%
Using the replacements
 \texttt{-F} for \texttt{+}
 and
 \texttt{+F+} for \texttt{-}
gives the other enantiomer.

%\clearpage% xxx
%
%%%%%%%%%%%%%%%%%%%%%%%%%%
%% with  lnth *= 3.0;  // thick lines
% stringsubst 2 _F_+F_+F_+F_+F_+F+ _ _ F F+F+F+F--F--F+F + + - - | tail -1 | ./bin 6 2 0 0 0.15 > tmp-pic.tex && make dotex
%% rotated like the PC curve:
% stringsubst 2 _F_+F_+F_+F_+F_+F+ _ _ F F+F+F+F--F--F+F + + - - | tail -1 | sed 's/+/++/g; s/-/--/g; s/^/-/;' | ./bin 12 2 0 0 0.15 > tmp-pic.tex && make dotex
%
%% with     lnth *= 2.0;  // thicker lines
% stringsubst 2 _F_+F_+F_+F_+F_+F+ _ _ F F+F+F+F--F--F+F + + - - | tail -1 | ./4444-pc-filter | ./bin 12 3 0 > tmp-pic.tex && make dotex
%
\begin{figure}[h!tbp]
\ifpdf
\begin{center}
{\includegraphics*[width=55mm, viewport={70 300 490 740}]{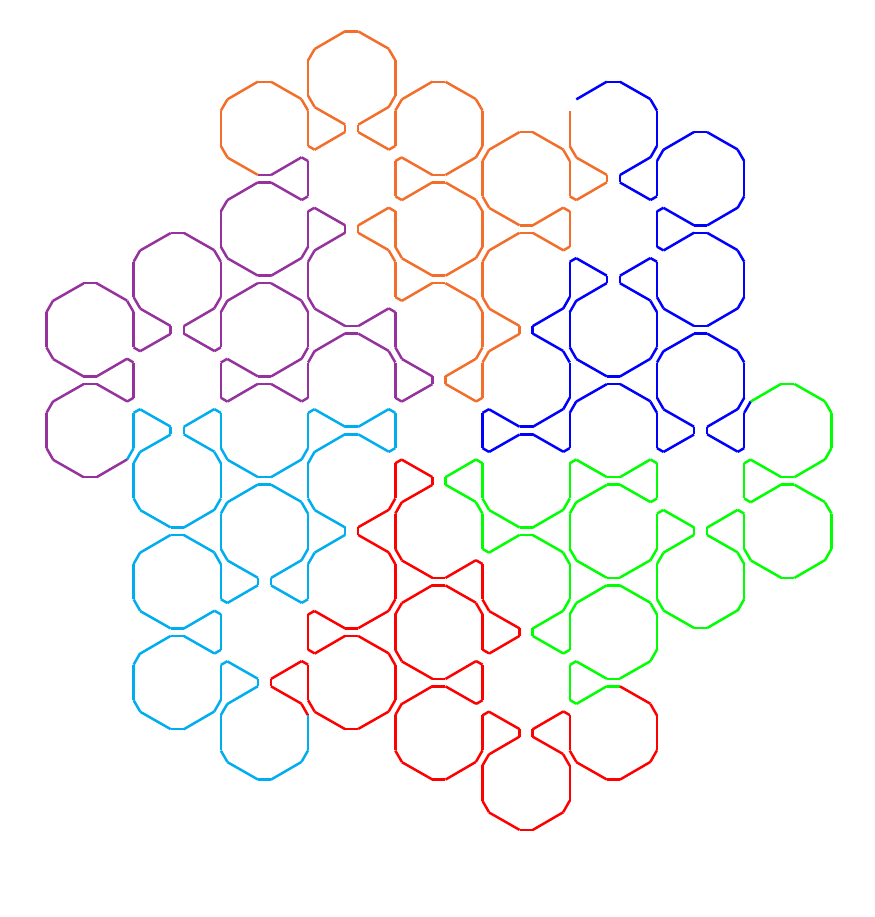}}%
{\includegraphics*[width=74mm, viewport={70 380 490 740}]{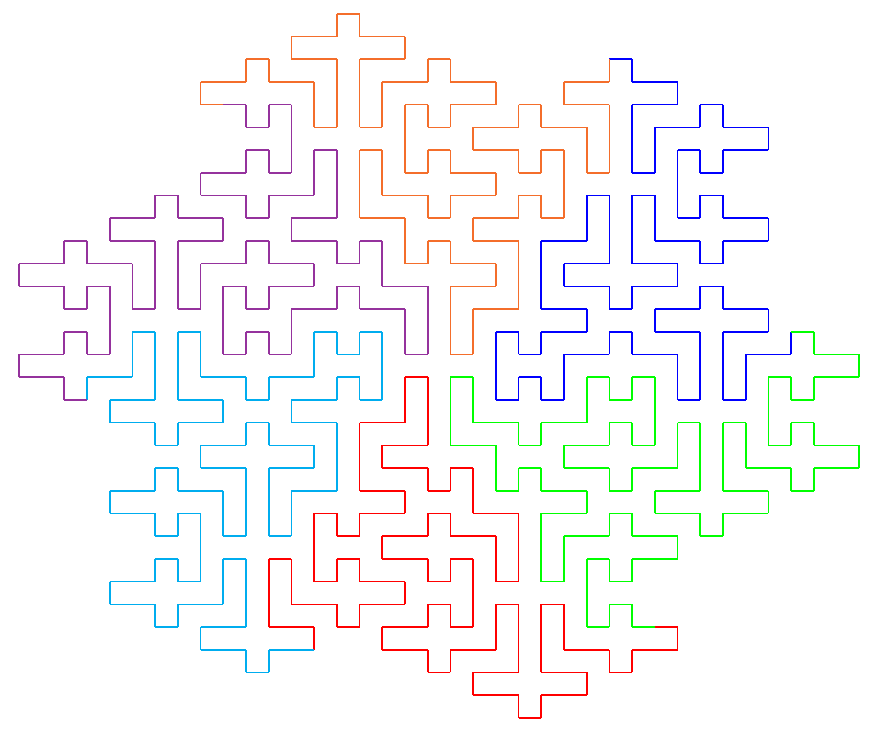}}
\end{center}
\else
\verb+{see pdf for image}+
\fi
\caption{\label{fig:r07-b-4444-cover}
A $(4^4)$-PC curve (right) from the tile $\Tile{+2}$ (left).
}
\end{figure}
%
%%%%%%%%%%%%%%%%%%%%%%%%%%%

In the $(4^4)$-PC curve shown in Figure~\ref{fig:r07-b-4444-cover}
all hexagons in the tile shown on the left
are replaced by crosses whose horizontal bars are longer than the vertical ones (right).
%% see file 4444-pc-filter.cc
%
For the conversion we keep track of the direction modulo 3 (not 6!)
which shall be held in a variable $d$.
In the following method $d$ has to be incremented or decremented
modulo 3 with every \texttt{+} or \texttt{-} in the input stream.
All letters \texttt{F} are dropped.
When a \texttt{+} is read, print
\texttt{++FF---F+++} if $d=0$,
\texttt{F+++F---F++} if $d=1$, and
\texttt{--F+++F+} if $d=2$.
When a \texttt{-} is read, print
\texttt{-F---F---F+++} if $d=0$,
\texttt{---FF---F++} if $d=1$, and
\texttt{-----FF+} if $d=2$.
Use turns by $30\adeg$, the resulting turns are all by $0$ or $\pm{}90\adeg$.

%% Check:
%void print_string_3(ulong dir, int sgn)
%{
%    dir %= 3;  // reduce!
%    if ( sgn == +1 )    {
%        switch ( dir )    {
%        case 0:        cout << "++FF---F+++";        break;
%        case 1:        cout << "F+++F---F++";        break;
%        case 2:        cout << "--F+++F+";           break;
%        }
%    }
%    else    {
%        switch ( dir )    {
%        case 0:        cout << "-F---F---F+++";        break;
%        case 1:        cout << "---FF---F++";          break;
%        case 2:        cout << "-----FF+";             break;
%        }
%    }
%}

%%%%%%%%%%%%%%%%%%%%%%%%%%
%% with  lnth *= 3.0;  // thick lines
% stringsubst 2 F+_F+_F+_F+_F+_F+ _ _  F F+F+F+F--F--F+F + + - - | tail -1 | ./trihex-ec-to-square-ec.pl 1 | ./bin 4 2 0 0 0.5 > tmp-pic.tex && make dotex
% stringsubst 2 F+_F+_F+_F+_F+_F+ _ _  F F+F+F+F--F--F+F + + - - | tail -1 | ./trihex-ec-to-square-ec.pl 2 | ./bin 4 2 0 0 0.5 > tmp-pic.tex && make dotex
% stringsubst 2 F+_F+_F+_F+_F+_F+ _ _  F F+F+F+F--F--F+F + + - - | tail -1 | ./trihex-ec-to-square-ec.pl 3 | ./bin 4 2 0 0 0.5 > tmp-pic.tex && make dotex
%
\begin{figure}[h!tbp]
\ifpdf
\begin{center}
{\includegraphics*[width=56mm, viewport={70 280 470 740}]{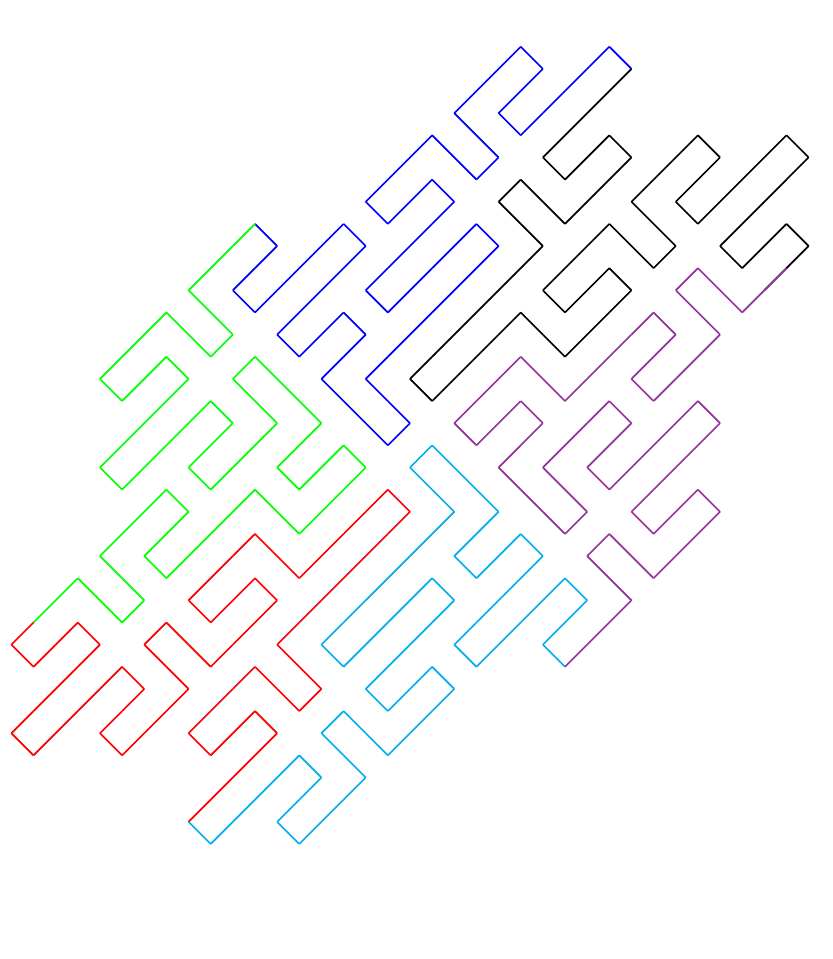}}%
%{\includegraphics*[width=72mm, viewport={70 290 470 730}]{r07-b-4444-cover-alt-b.pdf}}
%{\includegraphics*[width=72mm, viewport={40 290 490 730}]{r07-b-4444-cover-alt-b.pdf}}
{\includegraphics*[width=72mm, viewport={40 270 490 730}]{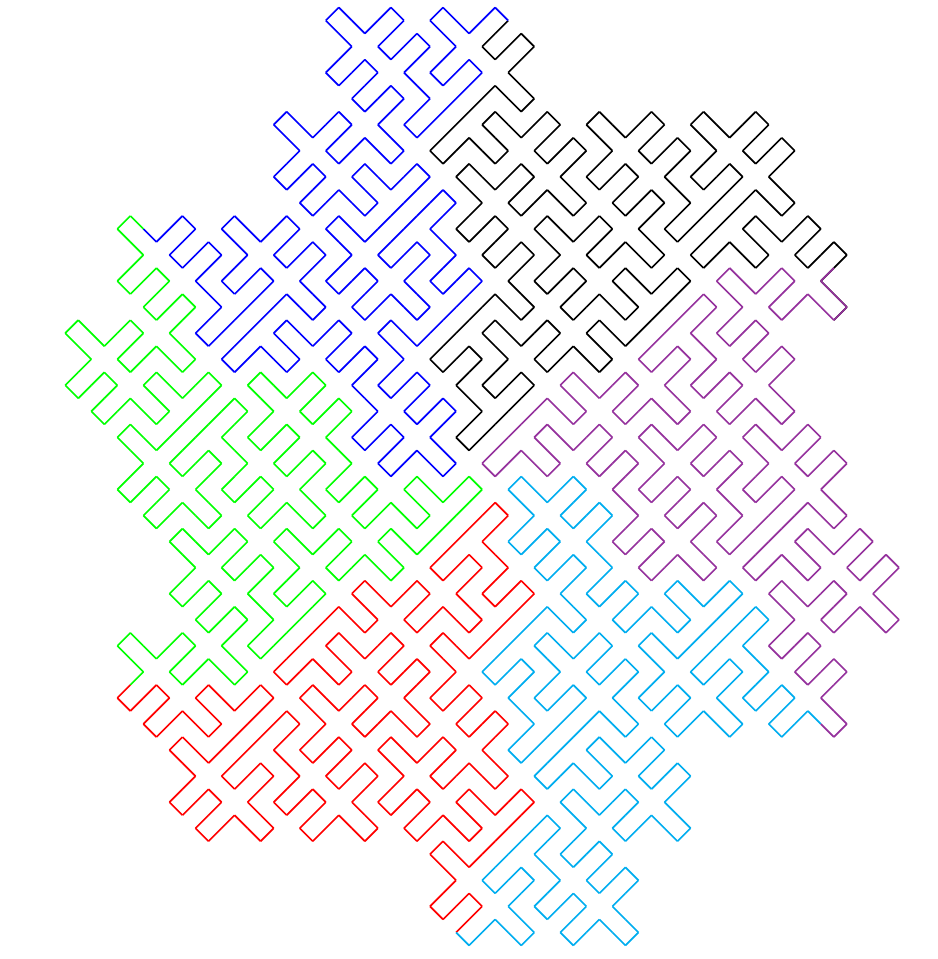}}
\end{center}
\else
\verb+{see pdf for image}+
\fi
\caption{\label{fig:r07-b-4444-alt-cover}
Two $(4^4)$-PC curves from the tile $\Tile{+2}$.
}
\end{figure}
%
%%%%%%%%%%%%%%%%%%%%%%%%%%%

The two $(4^4)$-PC curves shown in Figure~\ref{fig:r07-b-4444-alt-cover}
are renderings with rounding parameter $e=1/2$
of the $(4^4)$-EC curves computed as in section \Ref{sect:4444-EC}.

%%% Emacs:
%%% Local Variables:
%%% mode: latex
%%% MyRelDir: "."
%%% TeX-master: "arndt-curve-search.tex"
%%% dvi-file: "arndt-curve-search"
%%% makefile-dir: "./"
%%% frame-title-format: "CURVE-SEARCH (covers-points)"
%%% End:

\FloatBarrier

%%%%%%%%%%%%%%%%%%%%%%%%%%%%%%%%%%%%%%%%%%%%%%%%%%%%
%%%%%%%%%%%%%%%%%%%%%%%%%%%%%%%%%%%%%%%%%%%%%%%%%%%%
\subsection{Conversions to edge-covering curves}

The following sections are ordered by the grid of the edge-covering curves obtained.

%%%%%%%%%%%%%%%%%%%%%%%%%%%%%%%%%%%%%%%%%%%%%%%%%%%%
%%%%%%%%%%%%%%%%%%%%%%%%%%%%%%%%%%%%%%%%%%%%%%%%%%%%
\subsubsection{Edge-covering curves on $(3.4.6.4)$ from wiggly $(3^6)$-EC curves}\label{sect:3464-EC}
%% Use \texorpdfstring{$math$}{no-math} to get rid of these:
% Package hyperref Warning: Token not allowed in a PDF string (PDFDocEncoding):
% (hyperref)                removing `math shift' on input line 79.

%%%%%%%%%%%%%%%%%%%%%%%%%%
%
\begin{figure}[h!tbp]
\ifpdf
\begin{center}
{\includegraphics*[width=125mm, viewport={150 550 400 600}]{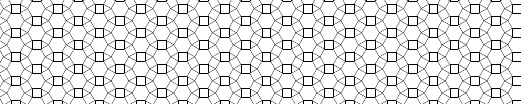}}
\end{center}
\else
\verb+{see pdf for image}+
\fi
\caption{\label{fig:3464-grid}
The $(3.4.6.4)$ grid.}
\end{figure}
%
%%%%%%%%%%%%%%%%%%%%%%%%%%

The grid $(3.4.6.4)$ (see Figure~\ref{fig:3464-grid})
has even incidences on all points, so edge-covering curves do exist.
While the grid has been left out in the search,
we give two constructions to obtain $(3.4.6.4)$-EC curves
from edge-covering curves on the triangular grid.
The grids $(3^6)$, $(3.6.3.6)$, $(4.4)$, and $(3.4.6.4)$
are the only uniform grids where edge-covering curves
can exist, so the last gap is closed.

The constructions work only for wiggly curves.
We use the terdragon (with map \Lmap{F}{F+F-F}) for our example.

%%%%%%%%%%%%%%%%%%%%%%%%%%
%
\begin{figure}[h!tbp]
%% with lnth *= 3.0;  // thick lines :
% stringsubst 3 RF  R R  F F+F-F + + - - | tail -1 | sed 's/F//g; s/-/_-tt-ttttt/g; s/+/_+tt+ttttt/g; s/R/+/;' | ./bin 6 3 0 > tmp-pic.tex && make dotex
%
% stringsubst 3 RF  R R  F F+F-F + + - - | tail -1 | sed 's/F//g; s/+/p/g; s/-/m/g; s/p/_---t++++t++++t---t---t++++t++++t---t/g; s/m/_---t++++t---t++t---t---t---t++++t++++t---t/g; s/--t/-t-tttt/g; s/++t/+t+tttt/g; s/R/++/;' | ./bin 12 3 0 > tmp-pic.tex && make dotex
%
\ifpdf
\begin{center}
{\includegraphics*[width=53mm, viewport={50 500 500 740}]{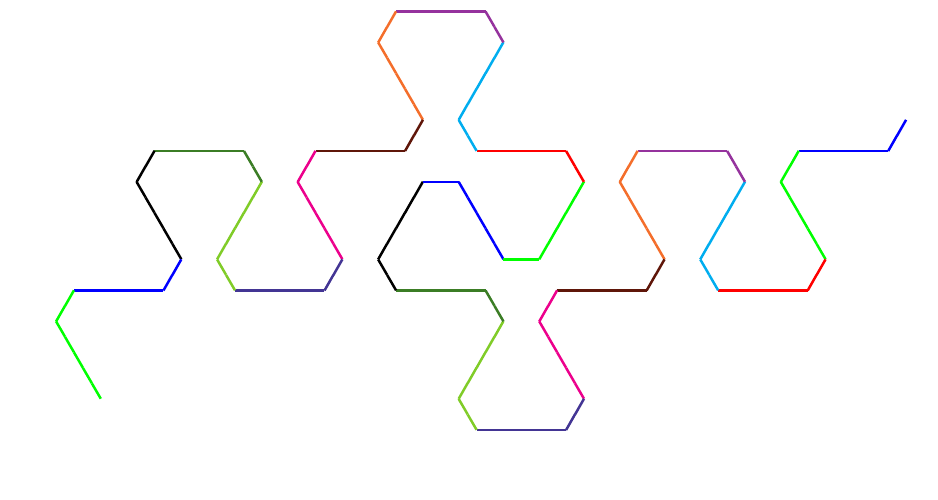}}%
{\includegraphics*[width=73mm, viewport={50 500 500 740}]{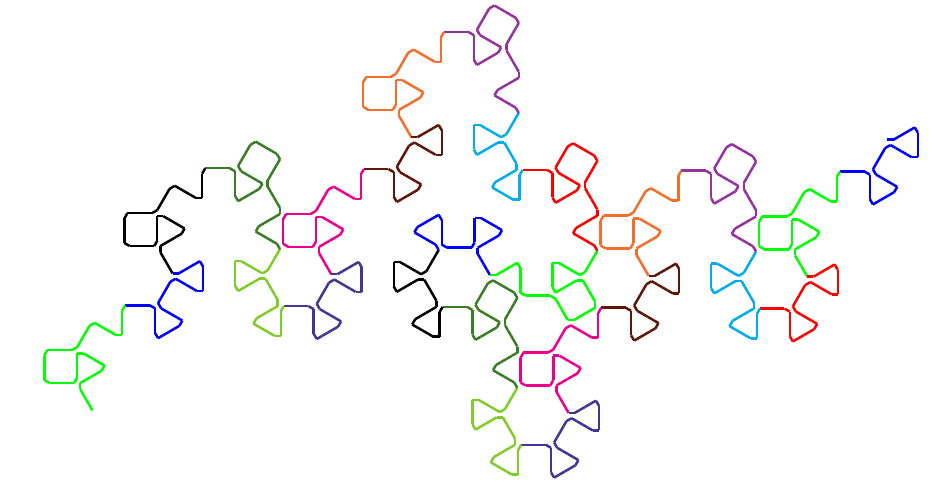}}
\end{center}
\else
\verb+{see pdf for image}+
\fi
\caption{\label{fig:r03-t-1-3464-ec-steps-0-1}
Third iterate of the terdragon with rounded corners (left)
and the corresponding curve on $(3.4.6.4)$ omitting the edges of some hexagons (right).}
\end{figure}
%
%%%%%%%%%%%%%%%%%%%%%%%%%%

\paragraph{First method.}
The initial step gives an incomplete curve.
Drop all \texttt{F},
replace all \texttt{+} by \texttt{p} and \texttt{-} by \texttt{m},
then replace
all \texttt{p} by \texttt{---F++++F++++F---F---F++++F++++F---F} and
all \texttt{m} by \texttt{---F++++F---F++F---F---F---F++++F++++F---F},
render with turns by $30\adeg$.
This is shown in Figure~\ref{fig:r03-t-1-3464-ec-steps-0-1}
where the colors help with the identification of the replacements.
Note that the edges of some hexagons are still missing.

%%%%%%%%%%%%%%%%%%%%%%%%%%
%
\begin{figure}[h!tbp]
%% with lnth *= 3.0;  // thick lines :
% stringsubst 3 F F F+F-F + + - - | tail -1 | ./3464-ec-filter | sed 's/[a-zA-Z]/t/g; s/--t/-t-tttt/g; s/++t/+t+tttt/g; s/^/++/;' | ./bin 12 3 0 > tmp-pic.tex && make dotex
%% for MAPR paper:
%% stringsubst 3 F F F+F-F + + - - | tail -1 | ./3464-ec-filter | sed 's/^/+/;' | ./bin 12 2 0 0 0.15 > tmp-pic.tex && make dotex
%% file: r03-t-1-3464-ec-step-2-monochrome.pdf
\ifpdf
\begin{center}
{\includegraphics*[width=90mm, viewport={50 500 500 740}]{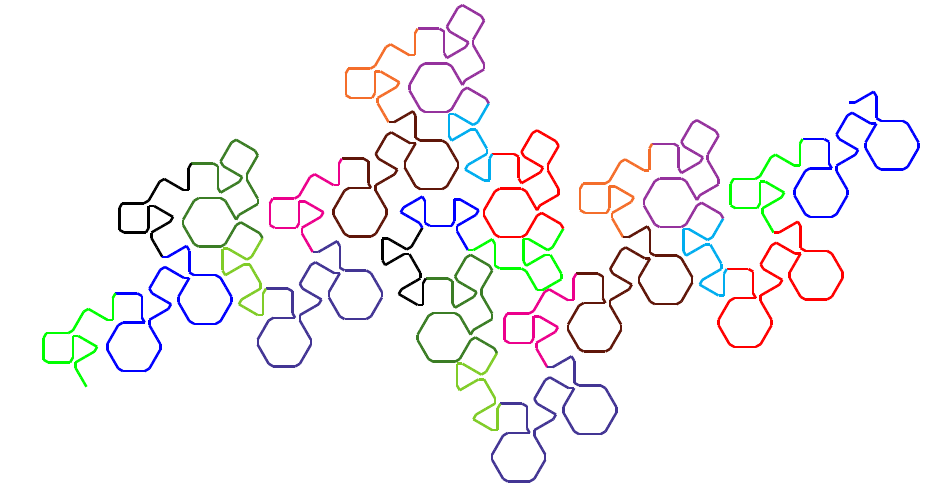}}
\end{center}
\else
\verb+{see pdf for image}+
\fi
\caption{\label{fig:r03-t-1-3464-ec-step-2}
The completed curve on $(3.4.6.4)$.}
\end{figure}
%
%%%%%%%%%%%%%%%%%%%%%%%%%%

There seems to be no way of including the hexagons by replacements
as simple as those just described.
Instead we consider the successive edges of the curves
as a walk from the left to the right in
Figure~\ref{fig:r03-t-1-3464-ec-steps-0-1},
visiting the hexagons on the right side of the current edge.
This corresponds to a tree traversal visiting all nodes of a tree
which are located where hexagons are missing.
To fill these in, we keep track of the directions modulo 3,
using representatives 0, 1, and 2.
We add two hexagons (letter \texttt{H}) by changing the replacement
for \texttt{p} to \texttt{---F++++FH++++F---F---F++++FH++++F---F}
if the direction is 0
(see the blue replacement of second edge from left in Figure~\ref{fig:r03-t-1-3464-ec-step-2})
and one hexagon by changing the replacement
for \texttt{m} to \texttt{---FH++++F---F++F---F---F---F++++F++++F---F}
if the direction is 2
(see the olive replacement of fourth edge from left in Figure~\ref{fig:r03-t-1-3464-ec-step-2}).
Finally, \texttt{H} is replaced by \texttt{---F++F++F++F++F++F-------}.

%if the original edge (left of Figure~\ref{fig:r03-t-1-3464-ec-steps-0-1})
%is horizontal and the next turn is counter clockwise\xxx{Clarify turns and dirs.}
% from left to right
% the turn is counter clockwise
%and
%add one hexagon if the original edge goes from the upper right to the lower left
% and the next turn is clockwise
%

\begin{figure}[h!tbp]
%% with lnth *= 2.0;  // thicker lines
% stringsubst 4 _-F_-F_-F  _ _ F F+F-F + + - - | tail -1 | ./3464-ec-filter | ./bin 12 2 0 0 0.20 > tmp-pic.tex && make dotex
% stringsubst 3 _-F_-F_-F  _ _ F F+F-F + + - - | tail -1 | ./3464-ec-filter | ./bin 12 2 0 0 0.20 > tmp-pic.tex && make dotex
%% rotation:
% stringsubst 3 _-F_-F_-F _ _ F F+F-F + + - - | tail -1 | ./3464-ec-filter | sed 's/^/+/;' | ./bin 12 2 0 0 0.20 > tmp-pic.tex && make dotex
%
\ifpdf
\begin{center}
%% Tile_4:
%{\includegraphics*[width=90mm, viewport={50 300 500 740}]{r03-t-1-3464-ec-tile.pdf}}
%% Tile_3:
%{\includegraphics*[width=90mm, viewport={50 330 500 730}]{r03-t-1-3464-ec-tile.pdf}}
%% rotation:
{\includegraphics*[width=80mm, viewport={50 310 500 730}]{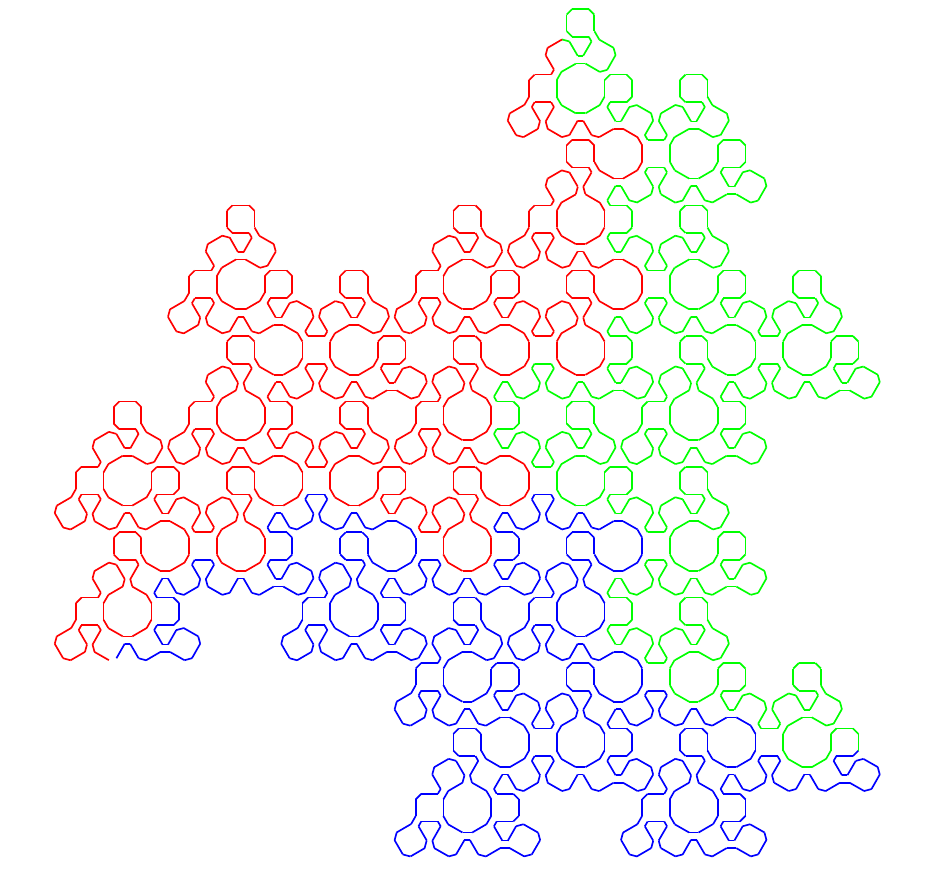}}
\end{center}
\else
\verb+{see pdf for image}+
\fi
\caption{\label{fig:r03-t-1-3464-ec-tile}
An edge-covering curve on $(3.4.6.4)$
from the tile $\Tile{-3}$ of the terdragon.}
%from the tile $\Tile{-4}$ of the terdragon.}
\end{figure}

Figure~\ref{fig:r03-t-1-3464-ec-tile}
shows the curve corresponding to the tile $\Tile{-3}$ of the terdragon.

%%%%%%%%%%%%%%%%%%%%%%%%%%
%
%
%%%%%%%%%%%%%%%%%%%%%%%%%%

%%%%%%%%%%%%%%%%%%%%%%%%%%
%
\begin{figure}[h!tbp]
%% with lnth *= 3.0;  // thick lines :
% stringsubst 3 RF  R R  F F+F-F + + - - | tail -1 | sed 's/F//g; s/-/_-tt-ttttt/g; s/+/_+tt+ttttt/g; s/R/+/;' | ./bin 6 3 0 > tmp-pic.tex && make dotex
%
% stringsubst 3 RF  R R  F F+F-F + + - - | tail -1 | sed 's/F//g; s/+/_++F++F/g; s/-/_---F---F++F/g;  s/--F/-t-tttt/g; s/++F/+t+tttt/g; s/R/+++/;' | ./bin 12 3 0 > tmp-pic.tex && make dotex
%
\ifpdf
\begin{center}
{\includegraphics*[width=40mm, viewport={50 500 500 740}]{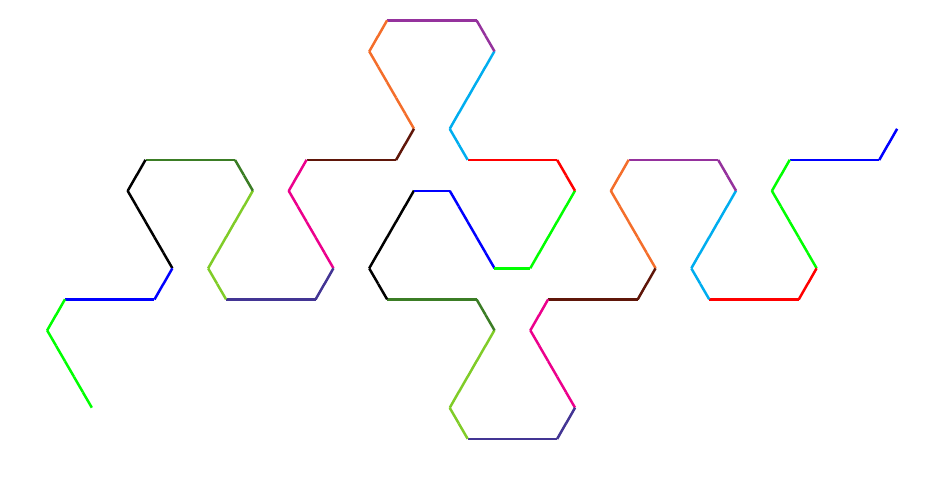}}%
{\includegraphics*[width=40mm, viewport={50 500 500 740}]{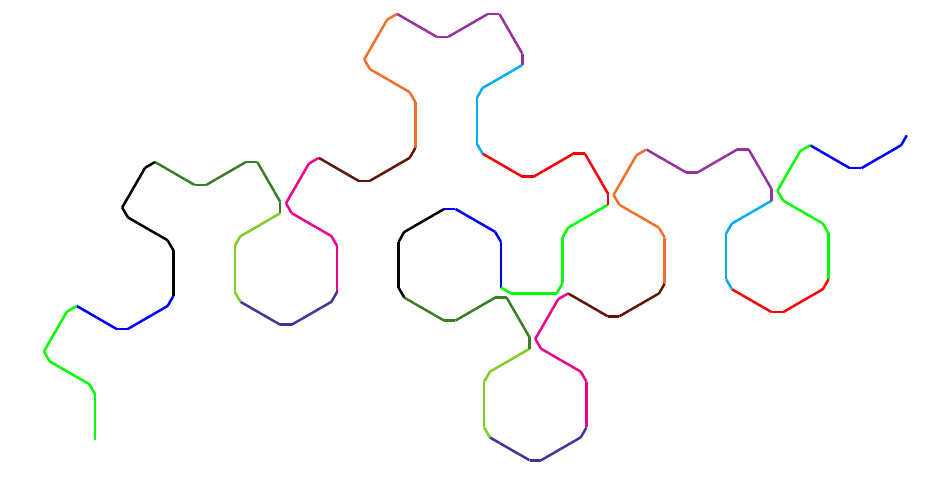}}%
{\includegraphics*[width=50mm, viewport={70 510 470 720}]{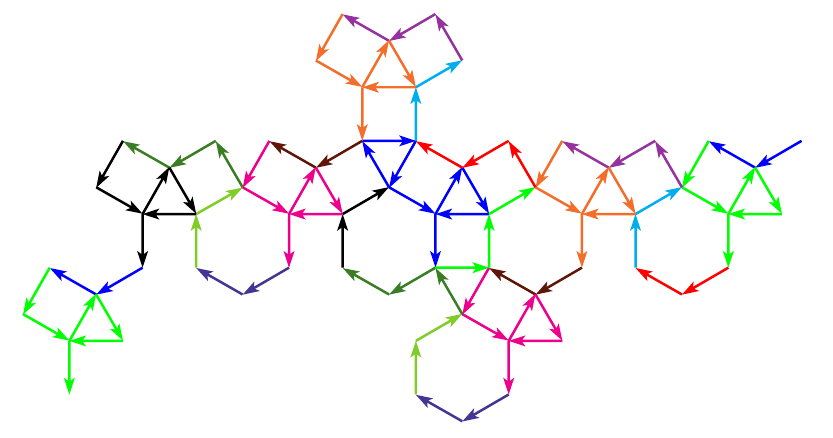}}
\end{center}
\else
\verb+{see pdf for image}+
\fi
\caption{\label{fig:r03-t-1-3464-ec-alt-steps-0-2}
Third iterate of the terdragon (left),
the corresponding curve on $(3.4.6.4)$ omitting the edges of some triangles (middle),
and the completed curve (right).}
\end{figure}
%
%%%%%%%%%%%%%%%%%%%%%%%%%%

\paragraph{Second method.}
For the incomplete curve shown in middle of Figure~\ref{fig:r03-t-1-3464-ec-alt-steps-0-2},
drop all \texttt{F},
replace
 all \texttt{+} by \texttt{++F++F}
and
 all \texttt{-} by \texttt{---F---F++F},
and use turns by $30\adeg$.
%
%For the next step we keep track of the directions modulo 3,
%using representatives 0, 1, and 2.
%
To fill in the missing triangles (letter \texttt{T}),
we change the replacement
for \texttt{-} to \texttt{---F---FTF}
if the direction is 1,
and the replacement
 for \texttt{+} to \texttt{TFTF}
if the direction is 2.
Finally, \texttt{T} is replaced by \texttt{---F++++F++++F---}.

\begin{figure}[h!tbp]
%% with lnth *= 2.0;  // thicker lines
% stringsubst 4 _-F_-F_-F  _ _ F F+F-F + + - - | tail -1 | ./3464-ec-alt-filter | ./bin 12 2 0 0 0.20 > tmp-pic.tex && make dotex
% stringsubst 3 _-F_-F_-F  _ _ F F+F-F + + - - | tail -1 | ./3464-ec-alt-filter | sed 's/^/+/g;' | ./bin 12 2 0 0 0.20 > tmp-pic.tex && make dotex
% stringsubst 4 _-F_-F_-F  _ _ F F+F-F + + - - | tail -1 | ./3464-ec-alt-filter | sed 's/^/+/g;' | ./bin 12 2 0 0 0.20 > tmp-pic.tex && make dotex
%
\ifpdf
\begin{center}
%% Tile_4:
{\includegraphics*[width=80mm, viewport={50 310 500 730}]{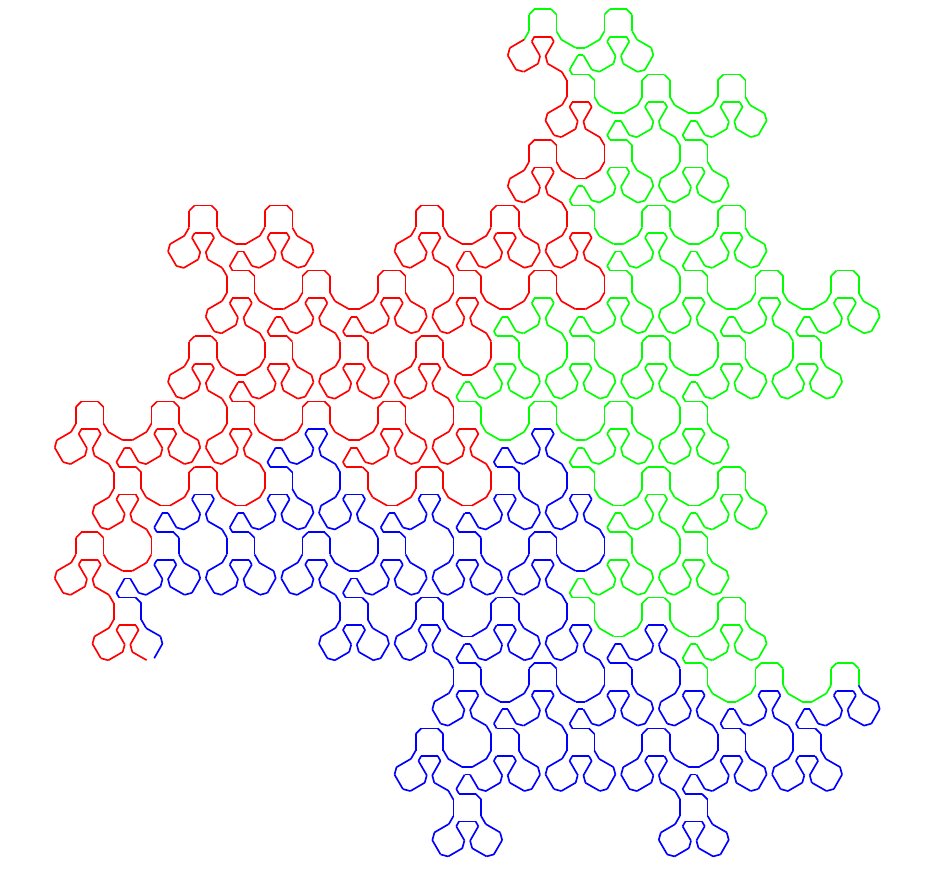}}
%% Tile_3:
%{\includegraphics*[width=70mm, viewport={60 360 480 720}]{r03-t-1-3464-ec-alt-tile.pdf}}
\end{center}
\else
\verb+{see pdf for image}+
\fi
\caption{\label{fig:r03-t-1-3464-ec-alt-tile}
An edge-covering curve on $(3.4.6.4)$
%from the tile $\Tile{-3}$ of the terdragon.}
from the tile $\Tile{-4}$ of the terdragon.}
\end{figure}

Figure~\ref{fig:r03-t-1-3464-ec-alt-tile}
shows the curve corresponding to the tile $\Tile{-4}$.

%%%%%%%%%%%%%%%%%%%%%%%%%%
%
%
%%%%%%%%%%%%%%%%%%%%%%%%%%

%%%%%%%%%%%%%%%%%%%%%%%%%%%%%%%%%%%%%%%%%%%%%%%%%%%%
%%%%%%%%%%%%%%%%%%%%%%%%%%%%%%%%%%%%%%%%%%%%%%%%%%%%
\subsubsection{Edge-covering curves on $(3.6.3.6)$}\label{sect:3636-EC}

%The resulting curves can be used as input for all conversions
%given in section \ref{sect:pc-from-tri-hex}.
%% \xxx{Give examples there?}

%%%%%%%%%%%%%%%%%%%%%%%%%%
%% with    lnth *= 3.0;  // thick lines
% stringsubst 4 _F+_F+_F+  _ _ F F+F-F + + - - | tail -1 | ./3636-ec-filter | ./bin 6 2 0 0 0.20 > tmp-pic.tex && make dotex
% stringsubst 4 _F-_F-_F-  _ _ F F+F-F + + - - | tail -1 | ./3636-ec-filter | ./bin 6 2 0 0 0.20 > tmp-pic.tex && make dotex
%
\begin{figure}[h!tbp]
\ifpdf
\begin{center}
{\includegraphics*[width=54mm, viewport={70 320 490 740}]{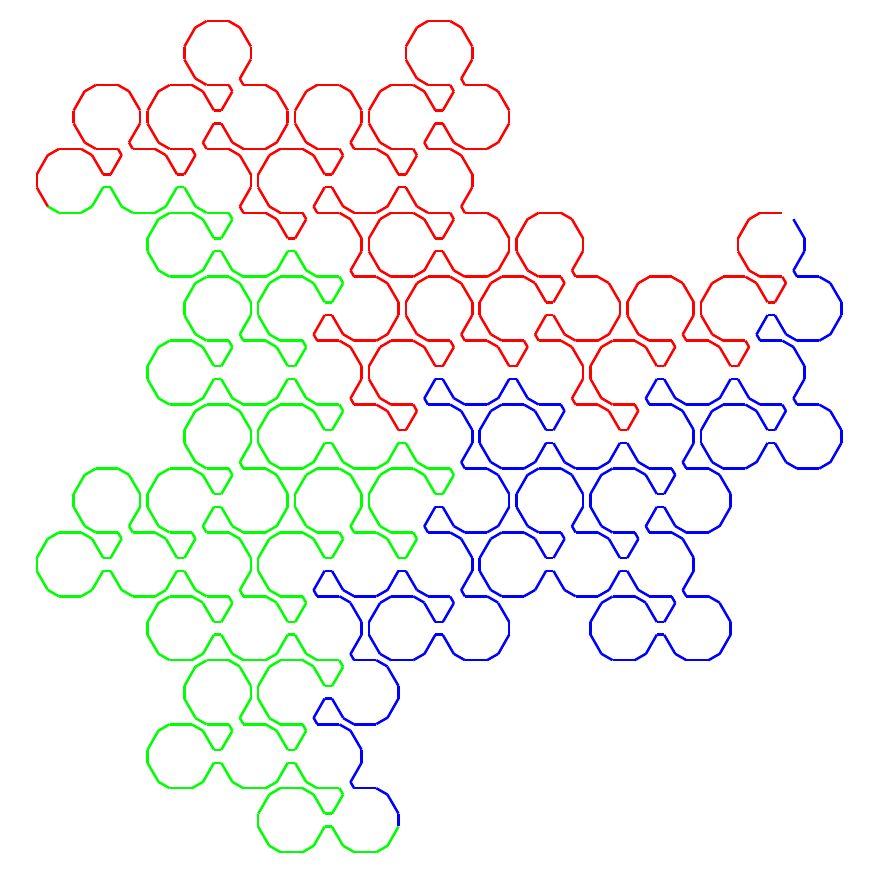}}%
{\includegraphics*[width=54mm, viewport={70 320 490 740}]{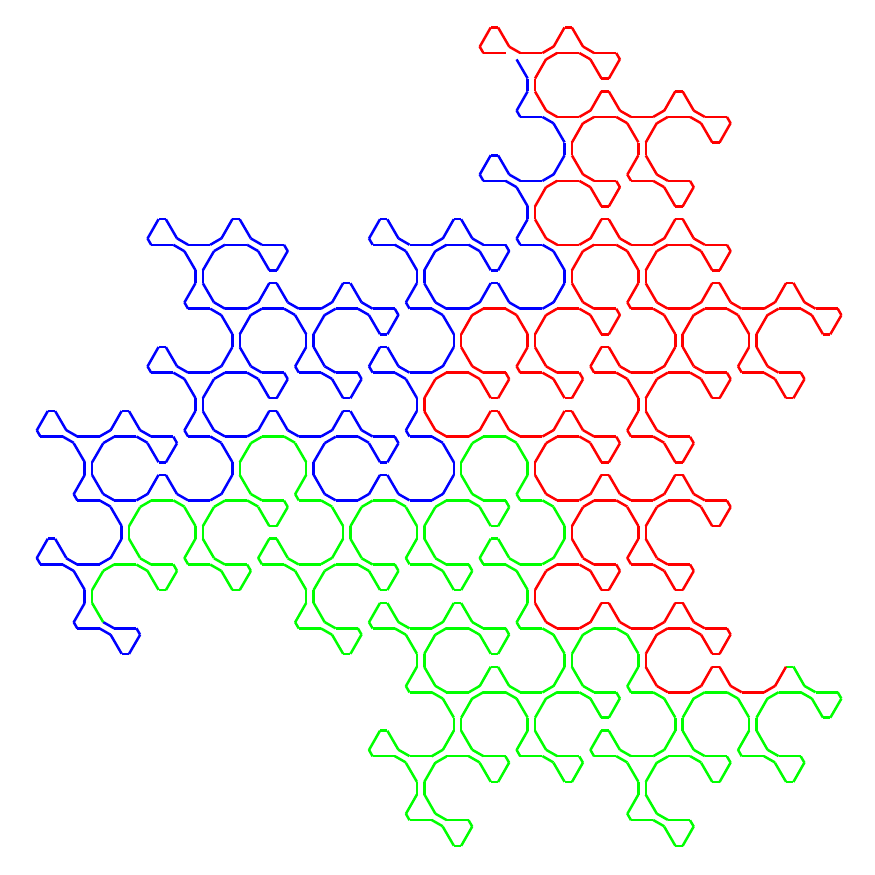}}
\end{center}
\else
\verb+{see pdf for image}+
\fi
\caption{\label{fig:r03-t-1-3636-ec-tiles}
Edge-covering curves on $(3.6.3.6)$
from the tiles $\Tile{+4}$ (left) and $\Tile{-4}$ (right) of the terdragon.}
\end{figure}
%
%%%%%%%%%%%%%%%%%%%%%%%%%%

\paragraph{Conversion from wiggly $(3^6)$-EC curves.}
Curves that are $(3.6.3.6)$-EC can be obtained from
wiggly curves on the triangular grid in a similar fashion.
Drop all \texttt{F}
and replace
 all \texttt{+} by \texttt{+F+F},
replace
 all \texttt{-} by \texttt{--F},
unless the direction is that of the initial edge,
then the replacement for \texttt{-} is \texttt{+F--F--F+F}.
Figure~\ref{fig:r03-t-1-3636-ec-tiles} shows the curves
obtained from the tiles of the terdragon.
%

%%%%%%%%%%%%%%%%%%%%%%%%%%
%% with    lnth *= 3.0;  // thick lines
% stringsubst 2 _F+_F+_F+ _ _ F F+F0F0F-F-F+F0F+F+F-F0F-F  0 0 + + - - | tail -1 | sed 's/F+/F+F+/g; s/F-/F--/g; s/F0/F+F--F+/g;' | ./bin 6 3 0 > tmp-pic.tex && make dotex
% stringsubst 2 _F-_F-_F- _ _ F F+F0F0F-F-F+F0F+F+F-F0F-F  0 0 + + - - | tail -1 | sed 's/F+/F+F+/g; s/F-/F--/g; s/F0/F+F--F+/g;' | ./bin 6 3 0 > tmp-pic.tex && make dotex
%
\begin{figure}[h!tbp]
\ifpdf
\begin{center}
{\includegraphics*[width=54mm, viewport={70 340 490 740}]{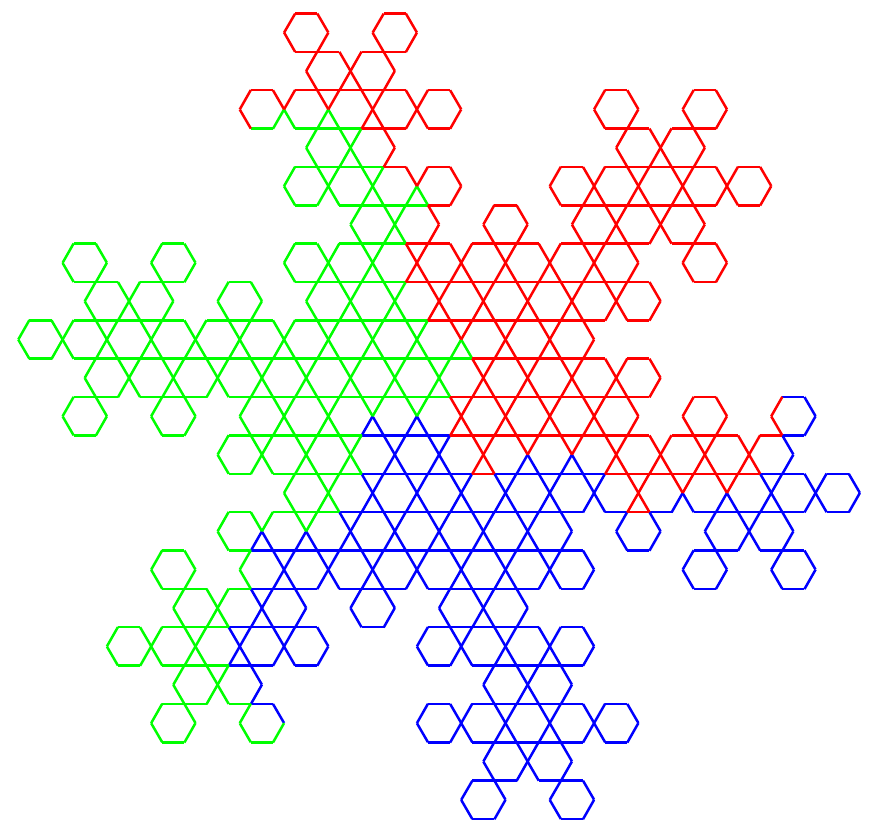}}%
{\includegraphics*[width=54mm, viewport={70 340 490 740}]{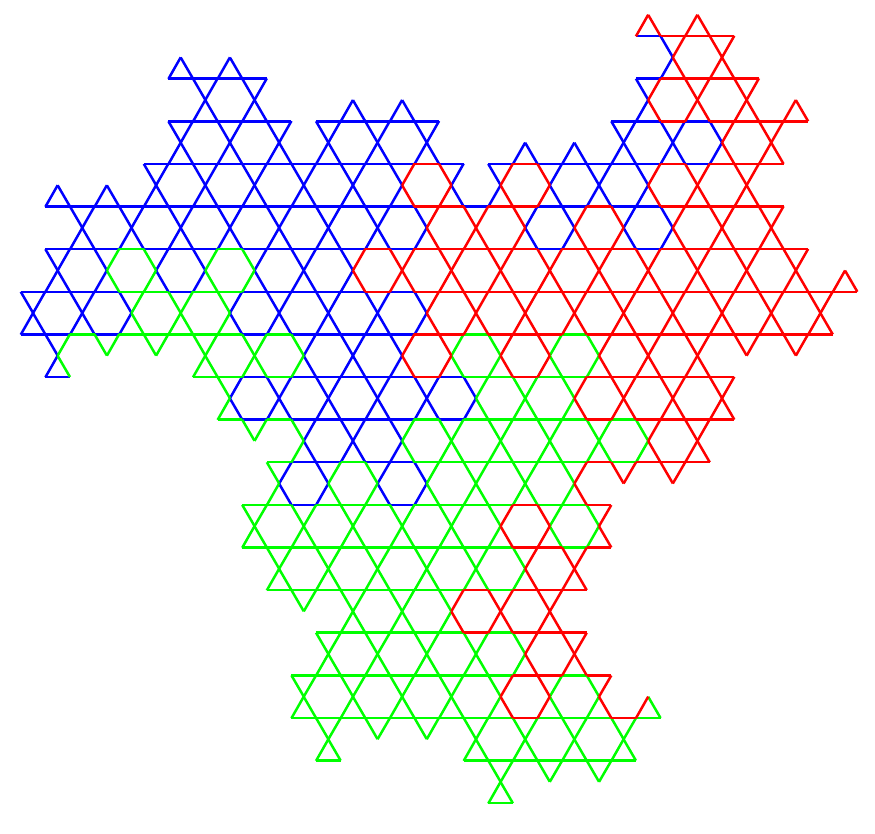}}
\end{center}
\else
\verb+{see pdf for image}+
\fi
\caption{\label{fig:r13-t-15-3636-ec-tile}
Edge-covering curves on $(3.6.3.6)$
from the tiles $\Tile{+2}$ (left) and $\Tile{-2}$ (right)
of the balanced curve \CID{R13-15} on the triangular grid.
}
\end{figure}
%
%%%%%%%%%%%%%%%%%%%%%%%%%%

\paragraph{Conversion from balanced $(3^6)$-EC curves.}
%Curves of orders of the form $R=3k+1$
%that have exactly $k$ non-turns
Balanced curves
can be rendered as $(3.6.3.6)$-EC curves
by replacing all
 \texttt{F+} by \texttt{F+F+},
 \texttt{F-} by \texttt{F--},
 and
 \texttt{F0} by \texttt{F+F--F+}.
%% sed 's/F+/F+F+/g; s/F-/F--/g; s/F0/F+F--F+/g;
%
The rendering for the curve \CID{R13-15} with rule
\Lmap{F}{F+F0F0F-F-F+F0F+F+F-F0F-F}
is shown in Figure~\ref{fig:r13-t-15-3636-ec-tile}.

%%%%%%%%%%%%%%%%%%%%%%%%%%
%% with    lnth *= 3.0;  // thick lines
%% stringsubst 3 F+_F+_F+_F+ _ _ F F+F+F-F-F + + - - | tail -1 | ./bin 4 2 0 0 0.1 > tmp-pic.tex && make dotex
%
% stringsubst 3 F+_F+_F+_F+ _ _ F F+F+F-F-F + + - - | tail -1 | ./square-ec-to-trihex-ec.pl | ./bin 6 2 0 0 0.1 > tmp-pic.tex && make dotex
%% rotation:
% stringsubst 3 F+_F+_F+_F+ _ _ F F+F+F-F-F + + - - | tail -1 | ./square-ec-to-trihex-ec.pl | tr 123456 234561 | ./bin 6 2 0 0 0.1 > tmp-pic.tex && make dotex
%
% stringsubst 3 F-_F-_F-_F- _ _ F F+F+F-F-F + + - - | tail -1 | ./square-ec-to-trihex-ec.pl | ./bin 6 2 0 0 0.1 > tmp-pic.tex && make dotex
%% rotation:
% stringsubst 3 F-_F-_F-_F- _ _ F F+F+F-F-F + + - - | tail -1 | ./square-ec-to-trihex-ec.pl | tr 123456 234561 | ./bin 6 2 0 0 0.1 > tmp-pic.tex && make dotex
%
\begin{figure}[h!tbp]
\ifpdf
\begin{center}
%{\includegraphics*[width=40mm, viewport={70 300 490 740}]{r05-q-1-tile-plus.pdf}}
{\includegraphics*[width=64mm, viewport={60 420 490 740}]{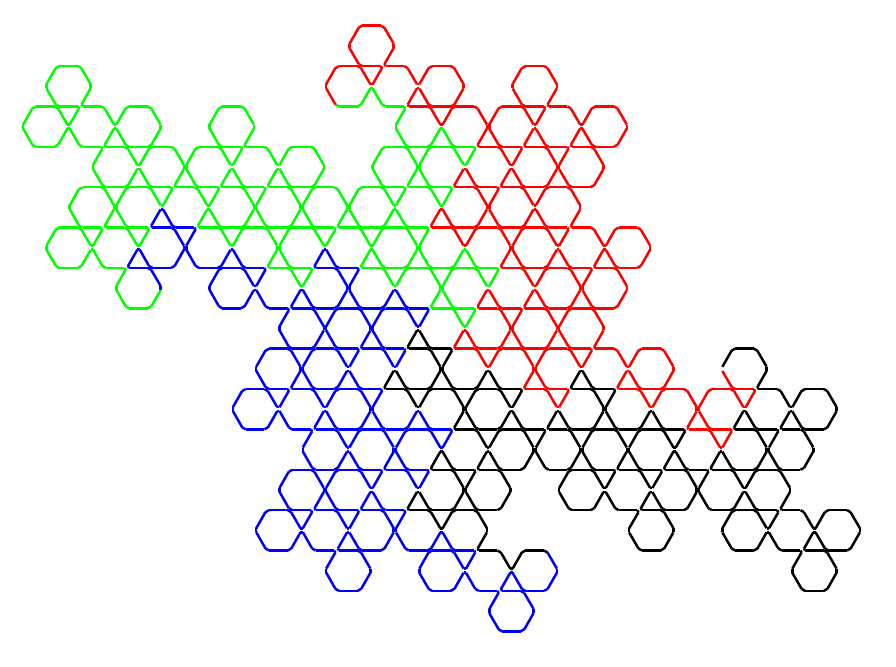}}%
{\includegraphics*[width=64mm, viewport={60 420 490 740}]{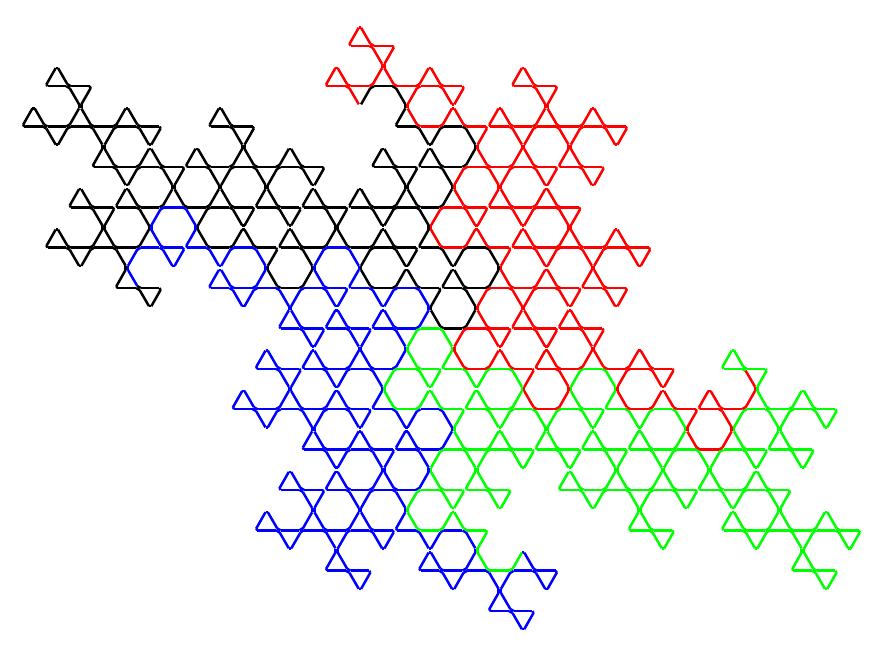}}
\end{center}
\else
\verb+{see pdf for image}+
\fi
\caption{\label{fig:r05-q-1-3636-ec-tiles}
%The tile $\Tile{+3}$ of curve \CID{R5-1} on the square grid (left)
%and edge-covering curves on $(3.6.3.6)$
%corresponding to $\Tile{+3}$ (middle) and $\Tile{-3}$ (right).
%
Edge-covering curves on $(3.6.3.6)$ corresponding to $\Tile{+3}$ (left) and $\Tile{-3}$ (right)
of the $(4^4)$-EC curve \CID{R5-1} whose tile $\Tile{+5}$ is shown in Figure~\ref{fig:r5-q-1-tile-plus}.
}
\end{figure}
%
%%%%%%%%%%%%%%%%%%%%%%%%%%

\paragraph{A conversion from $(4^4)$-EC curves.}
Figure~\ref{fig:r05-q-1-3636-ec-tiles}
shows two $(3.6.3.6)$-EC curves corresponding to the tiles
$\Tile{+3}$ (left) and $\Tile{-3}$ (right)
of the curve \CID{R5-1} on the square grid.
The conversion is surprisingly easy.
Let $k \in \{1, 2, \ldots, 6\}$ denote an edge from $0$ to
the sixth roots of unity, $k$ standing for $\exp(2\pi\,i\,(k-1)/6)$,
as shown in Figure~\ref{fig:directions}.
Horizontal edges are kept,
moves to the right and left are respectively mapped to $1$ and $4$.
Vertical edges up and down are respectively mapped to $65$ and $32$.

%\clearpage% xxx
%
%%%%%%%%%%%%%%%%%%%%%%%%%%%%%%%%%%%%%%%%%%%%%%%%%%%%
%%%%%%%%%%%%%%%%%%%%%%%%%%%%%%%%%%%%%%%%%%%%%%%%%%%%
\subsubsection{Edge-covering curves on $(4^4)$ from $(3.6.3.6)$-EC curves}\label{sect:4444-EC}

%%%%%%%%%%%%%%%%%%%%%%%%%%
%% with    lnth *= 2.0;  // thicker lines
% stringsubst 2 F+_F+_F+_F+_F+_F _ _ F F+F+F+F+F--F+F+F--F--F+F+F--F + + - - | tail -1 | ./trihex-ec-to-square-ec.pl 1 | ./bin 4 2 0 0 0.15 > tmp-pic.tex && make dotex
% stringsubst 2 F+_F+_F+_F+_F+_F _ _ F F+F+F+F+F--F+F+F--F--F+F+F--F + + - - | tail -1 | ./trihex-ec-to-square-ec.pl 2 | ./bin 4 2 0 0 0.15 > tmp-pic.tex && make dotex
% stringsubst 2 F+_F+_F+_F+_F+_F _ _ F F+F+F+F+F--F+F+F--F--F+F+F--F + + - - | tail -1 | ./trihex-ec-to-square-ec.pl 3 | ./bin 4 2 0 0 0.15 > tmp-pic.tex && make dotex
%
\begin{figure}[h!tbp]
\ifpdf
\begin{center}
{\includegraphics*[width=62mm, viewport={60 270 490 740}]{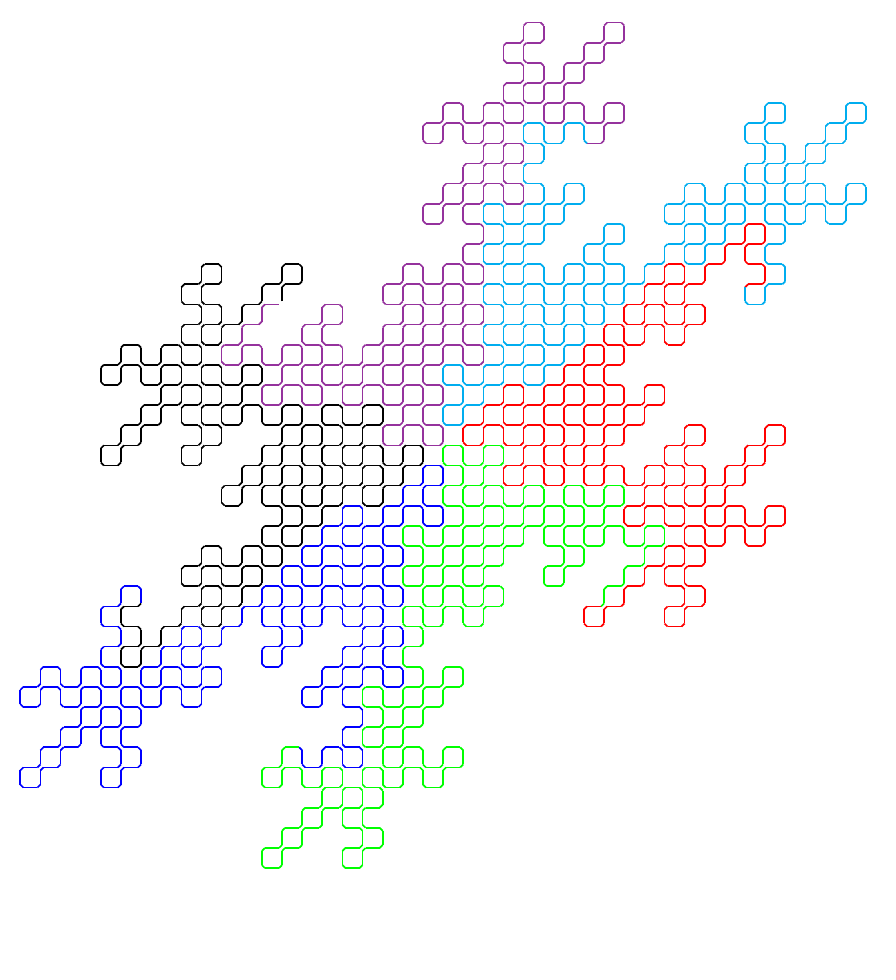}}%
%{\includegraphics*[width=66mm, viewport={60 280 490 740}]{r13-b-1-4444-ec-tile-plus-b.pdf}}
{\includegraphics*[width=66mm, viewport={60 270 490 740}]{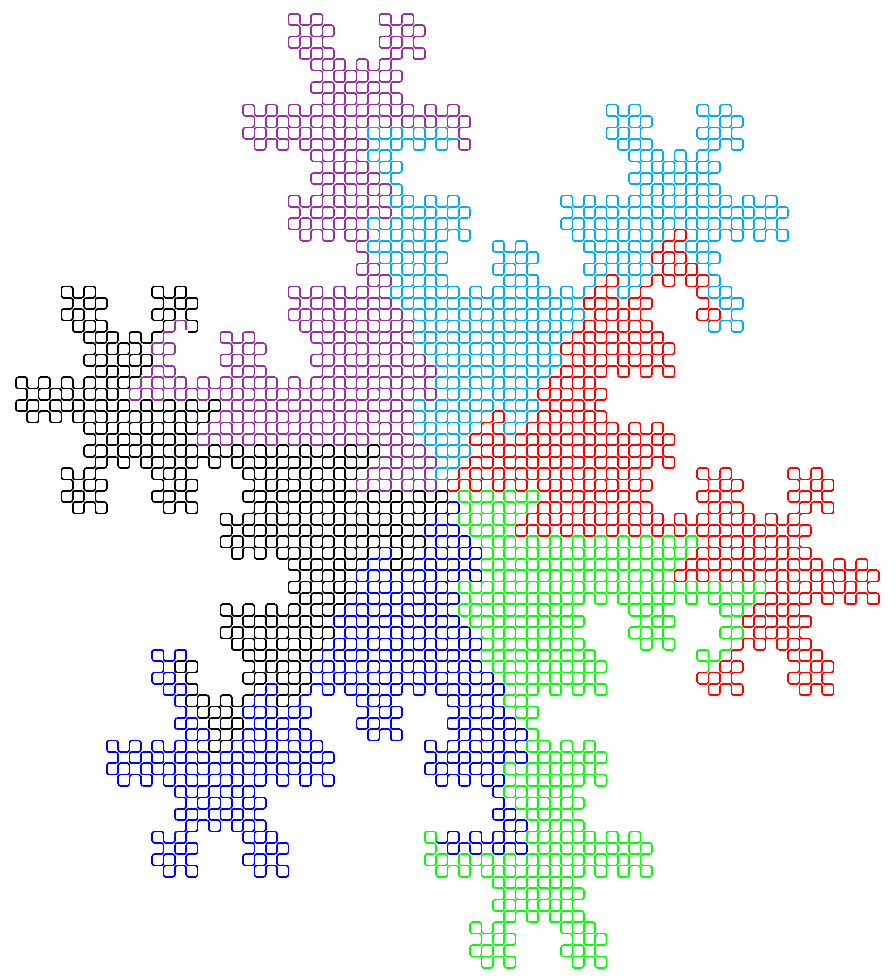}}
\end{center}
\else
\verb+{see pdf for image}+
\fi
\caption{\label{fig:r13-b-1-4444-ec-tile-plus}
Two edge-covering curves on the square grid corresponding to $\Tile{+2}$
of the curve \CID{R13-1} on the tri-hexagonal grid.
}
\end{figure}
%
%%%%%%%%%%%%%%%%%%%%%%%%%%

A simple mapping of the directions
on the tri-hexagonal grid to those of the square grid
gives the $(4^4)$-EC curve shown at the left of Figure~\ref{fig:r13-b-1-4444-ec-tile-plus}:
map the directions (as shown in Figure~\Ref{fig:directions})
1, 2, \ldots, 6 respectively to 21, 2, 3, 43, 4, and 1.
A less distorted curve
results from the mappings to 141, 2321, 2, 3, 4323, and 414.
The still less distorted curve shown
at the right of Figure~\ref{fig:r13-b-1-4444-ec-tile-plus}
uses the mappings to 12321, 2343212, 3432, 34143, 4121434, and 1214.

%%%%%%%%%%%%%%%%%%%%%%%%%%%%%%%%%%%%%%%%%%%%%%%%%%%%
%%%%%%%%%%%%%%%%%%%%%%%%%%%%%%%%%%%%%%%%%%%%%%%%%%%%
\subsubsection{Edge-covering curves on $(3^6)$ from $(3.6.3.6)$-EC curves}\label{sect:333333-EC}

%%%%%%%%%%%%%%%%%%%%%%%%%%
%% with    lnth *= 2.0;  // thicker lines
% stringsubst 2 _F--_F--_F _ _ F F+F+F+F--F--F+F + + - - | tail -1 | ./bin 6 2 1 0 0.1 > tmp-pic.tex && make dotex
% stringsubst 2 _F--_F--_F _ _ F F+F+F+F--F--F+F + + - - | tail -1 | ./trihex-ec-to-triangle-ec.pl 3 | ./bin 3 2 1 0 0.1 > tmp-pic.tex && make dotex
%
\begin{figure}[h!tbp]
\ifpdf
\begin{center}
{\includegraphics*[width=54mm, viewport={60 330 490 720}]{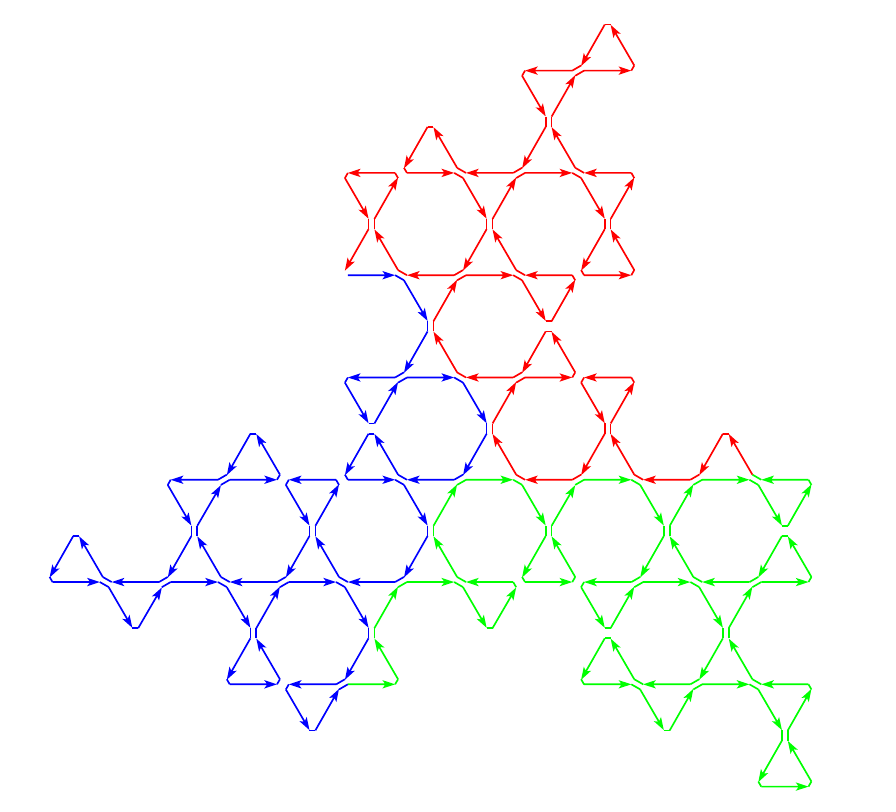}}%
{\includegraphics*[width=54mm, viewport={60 330 490 710}]{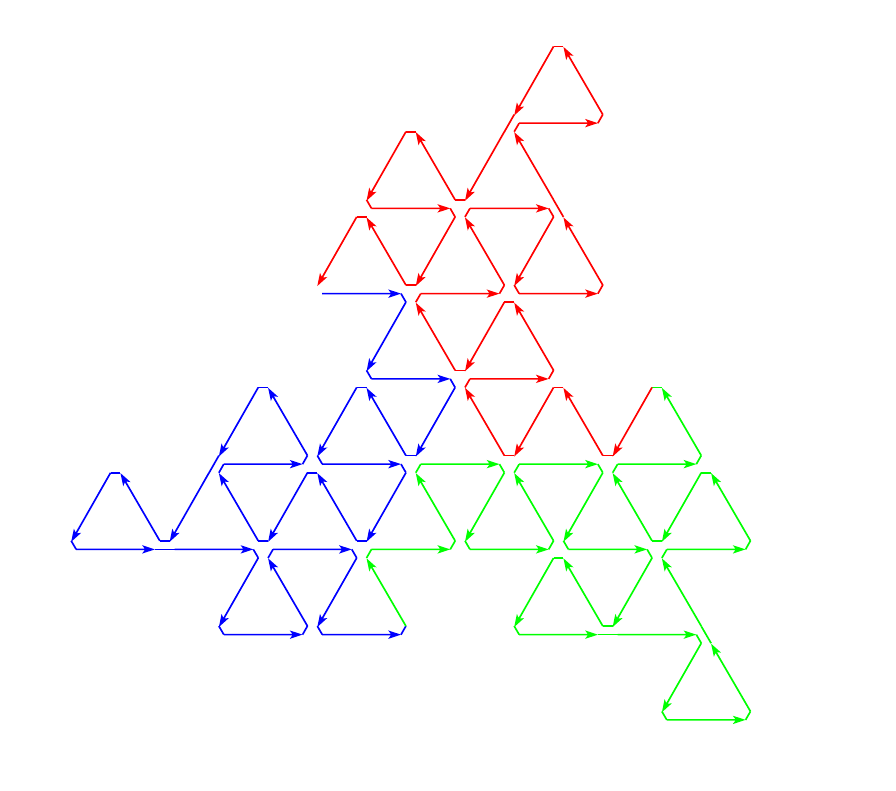}}
\end{center}
\else
\verb+{see pdf for image}+
\fi
\caption{\label{fig:r7-b-1-333333-ec-tile-minus}
The tile $\Tile{-2}$ of the curve \CID{R7-1} on the tri-hexagonal grid (left).
Simply dropping every second row of triangles gives a $(3^6)$-EC curve (right).
}
\end{figure}
%
%%%%%%%%%%%%%%%%%%%%%%%%%%

%%%%%%%%%%%%%%%%%%%%%%%%%%
%% with    lnth *= 2.0;  // thicker lines
% stringsubst 2 F+_F+_F+_F+_F+_F+ _ _   F F+F+F--F+F+F--F--F+F+F--F+F--F+F+F+F--F+F+F+F+F+F--F+F--F   + + - - | tail -1 | ./bin 6 2 0 0 0.15 > tmp-pic.tex && make dotex # R25-11
% stringsubst 2 F+_F+_F+_F+_F+_F+ _ _   F F+F+F--F+F+F--F--F+F+F--F+F--F+F+F+F--F+F+F+F+F+F--F+F--F   + + - - | tail -1 | ./trihex-ec-to-triangle-ec.pl 3 | ./bin 3 2 0 0 0.15 > tmp-pic.tex && make dotex # R25-11
%
\begin{figure}[h!tbp]
\ifpdf
\begin{center}
{\includegraphics*[width=64mm, viewport={60 310 490 740}]{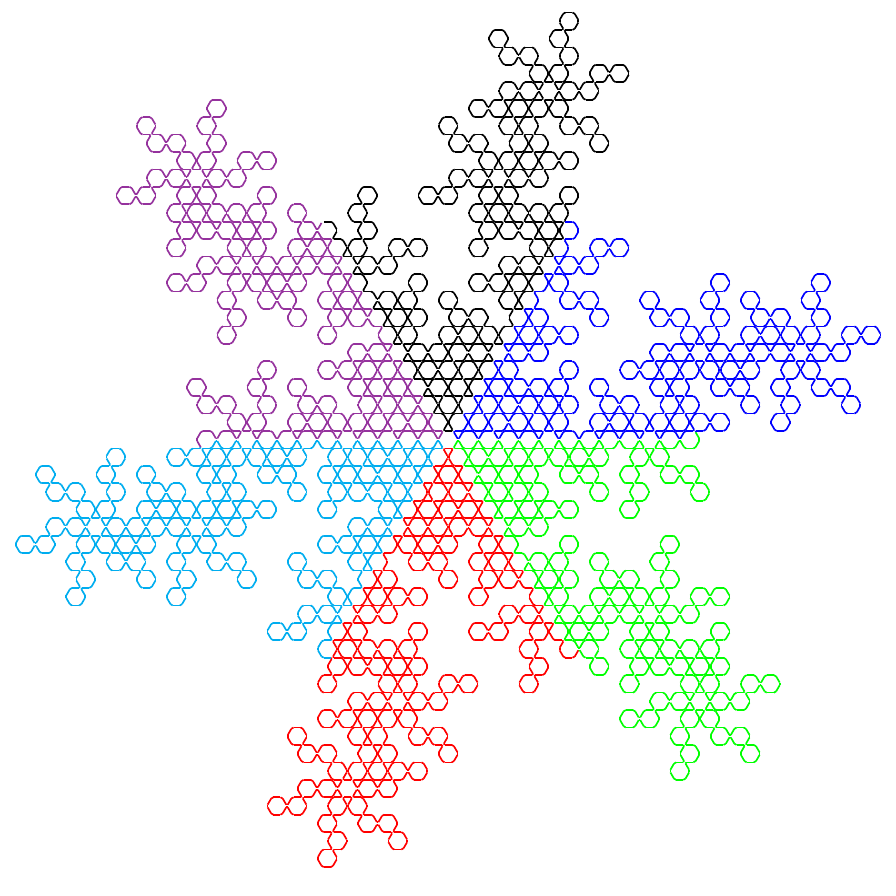}}%
{\includegraphics*[width=64mm, viewport={60 310 490 740}]{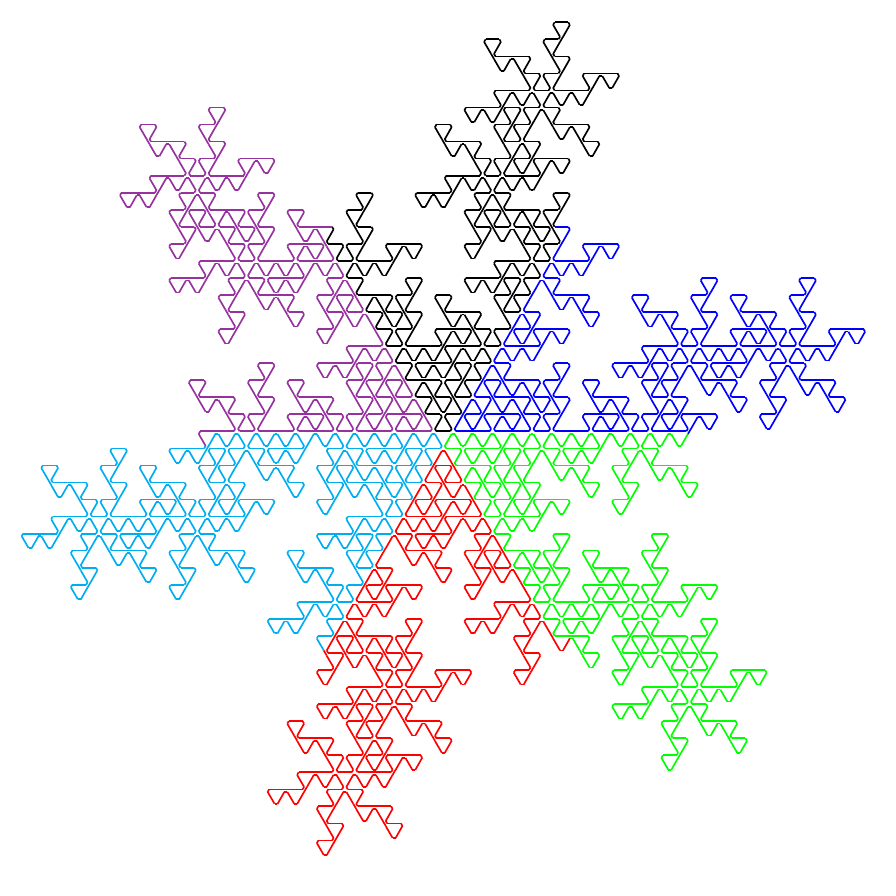}}
\end{center}
\else
\verb+{see pdf for image}+
\fi
\caption{\label{fig:r25-b-11-333333-ec-tile-plus}
The tile $\Tile{+2}$ of the curve \CID{R25-11} on the tri-hexagonal grid (left)
and the corresponding $(3^6)$-EC curve (right).
}
\end{figure}
%
%%%%%%%%%%%%%%%%%%%%%%%%%%

By simply dropping edges with every second of the six directions
in a curve on the tri-hexagonal grid a $(3^6)$-EC curve is obtained,
see Figure~\ref{fig:r7-b-1-333333-ec-tile-minus}.
There always exist curves on the triangular grid with tiles
of the same (6-fold symmetric) shape in the limit
as the shape of $\Tile{+}$ on the tri-hexagonal grid,
compare Figure~\ref{fig:r25-b-11-333333-ec-tile-plus}
and Figure~\Ref{fig:r25-t-tiles-same-shape}.
%
%% balanced?  No, rules have #+ == #- == #0 + 3 (skew-balanced)
%% Always at least one tile is 6-fold symmetric
%% See grep-rule-kludge.pl, cf. sign-balance.pl
%% Therefore the following is not clear cut:
%This construction is an inverse of that
%for Figure~\ref{fig:r13-t-15-3636-ec-tile}.

%%% Emacs:
%%% Local Variables:
%%% mode: latex
%%% MyRelDir: "."
%%% TeX-master: "arndt-curve-search.tex"
%%% dvi-file: "arndt-curve-search"
%%% makefile-dir: "./"
%%% frame-title-format: "CURVE-SEARCH (covers-edges)"
%%% End:

\FloatBarrier

%%%%%%%%%%%%%%%%%%%%%%%%%%%%%%%%%%%%%%%%%%%%%%%%%%%%%%%%%%%%
%%%%%%%%%%%%%%%%%%%%%%%%%%%%%%%%%%%%%%%%%%%%%%%%%%%%%%%%%%%%
\section{Multiplying and dividing curves}%\label{sect:kronecker}

Two curves on the same grid of orders $R_1$ and $R_2$ can be combined to obtain
curves whose order is $R_1\,R_2$
via a non-commutative product of the respective L-systems.
This construction is equivalent to
Dekking's \jjterm{folding convolution} \cite[Equation~1, pp.~21]{dekking-tiles-final}.
%\xxx{@MD: OK?}

In the other direction, any curve
can be divided into parts of equal lengths
by a division procedure of the L-system.
%
%Many\xxx{many?} curves appearing in the literature
%are special cases of these constructions.

%%%%%%%%%%%%%%%%%%%%%%%%%%%%%%%
\subsection{Products of curves}%\label{sect:}

%%%%%%%%%%%%%%%%%%%%%%%%%%
% stringsubst 6 +F_+F_0F_-F_+F_+F_0F_-F_-F_+F_0F_-F_  _ _  F G+G-G  G F+F0F-F  0 0 + + - - | tail -1 | ./bin 3 3 0 > tmp-pic.tex && make dotex # Kronecker 3 x 4 decomp
% stringsubst 6 +F_+F_-F_+F_+F_-F_0F_+F_-F_-F_+F_-F_  _ _   F G+G0G-G  G F+F-F  0 0 + + - - | tail -1 | ./bin 3 3 0 > tmp-pic.tex && make dotex # Kronecker 4 x 3 decomp
\begin{figure}[h!tbp]
\ifpdf
\begin{center}
%% for iterates 4:
%\fbox{\includegraphics*[width=52mm, viewport={70 280 490 750}]{kronecker-3-4.pdf}}%
%\fbox{\includegraphics*[width=52mm, viewport={70 280 490 750}]{kronecker-4-3.pdf}}
%% for iterates 6:
{\includegraphics*[width=63mm, viewport={65 330 490 750}]{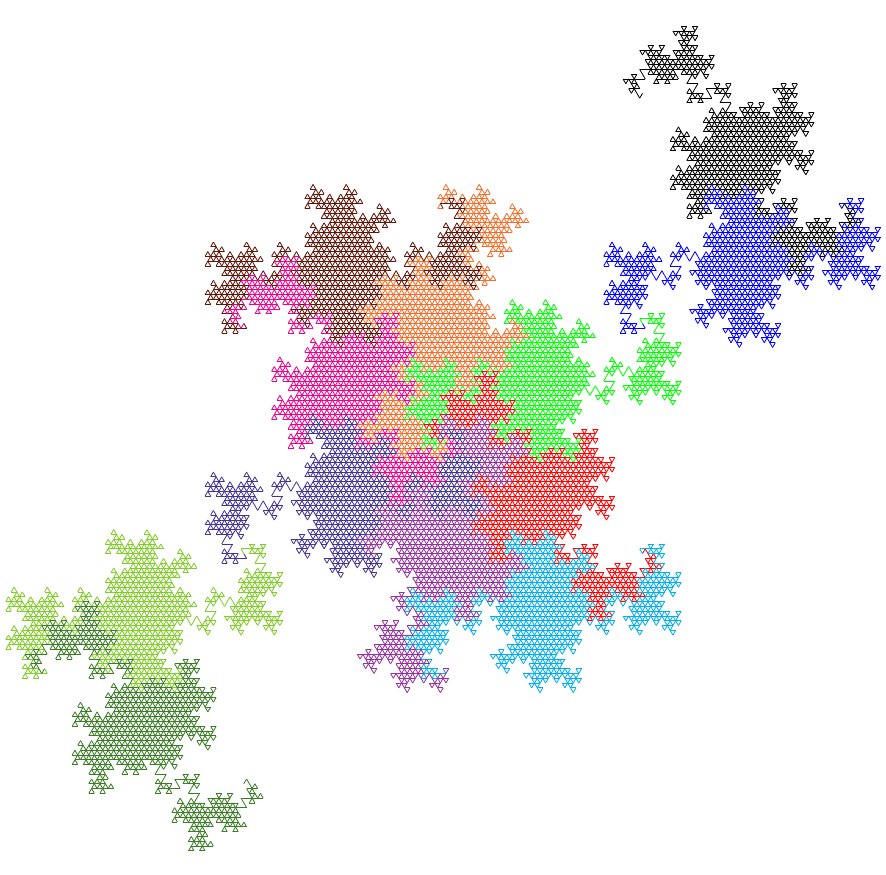}}%
{\includegraphics*[width=63mm, viewport={65 330 490 750}]{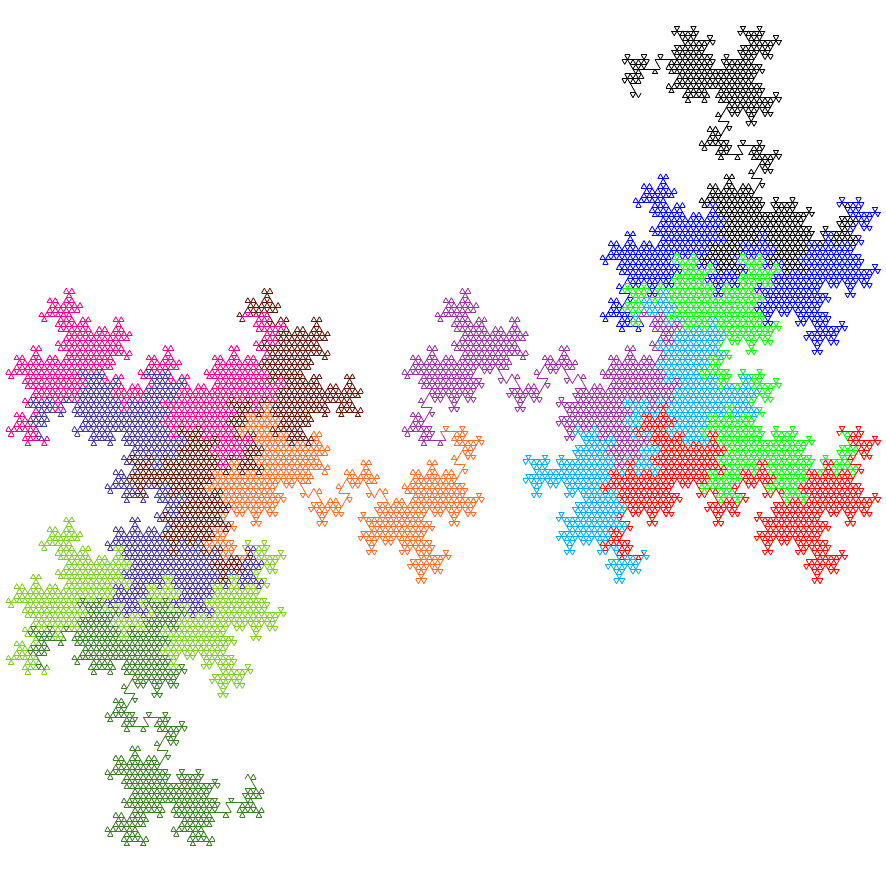}}
\end{center}
\else
\verb+{see pdf for image}+
\fi
\caption{\label{fig:kron-3-4}
Products of the terdragon (\Lmap{F}{F+F-F}) and the curve \Lmap{F}{F+F0F-F}.
The left curve is \CID{R12-10}, the right one is \CID{R12-17}.}
% R12-10 (p.6/10)
% R12-17 (p.9/10)
\end{figure}
%
%%%%%%%%%%%%%%%%%%%%%%%%%%

%%%%%%%%%%%%%%%%%%%%%%%%%%
% stringsubst 6 F_+F_0F_-F_-F_+F_0F_-F_+F_+F_0F_-F_  _ _  F G-G+G  G F+F0F-F  0 0 + + - - | tail -1 | ./bin 3 3 0 > tmp-pic.tex && make dotex # Kronecker 3m x 4
% stringsubst 6 +F_-F_+F_+F_-F_+F_0F_-F_+F_-F_-F_+F_  _ _   F G+G0G-G  G F-F+F  0 0 + + - - | tail -1 | ./bin 3 3 0 > tmp-pic.tex && make dotex # Kronecker 4 x 3m
%
\begin{figure}[h!tbp]
\ifpdf
\begin{center}
%% for iterates 4:
%{\includegraphics*[width=52mm, viewport={70 300 490 750}]{kronecker-3m-4.pdf}}%
%{\includegraphics*[width=52mm, viewport={70 180 490 750}]{kronecker-4-3m.pdf}}
%% for iterates 6:
{\includegraphics*[width=63mm, viewport={60 220 490 750}]{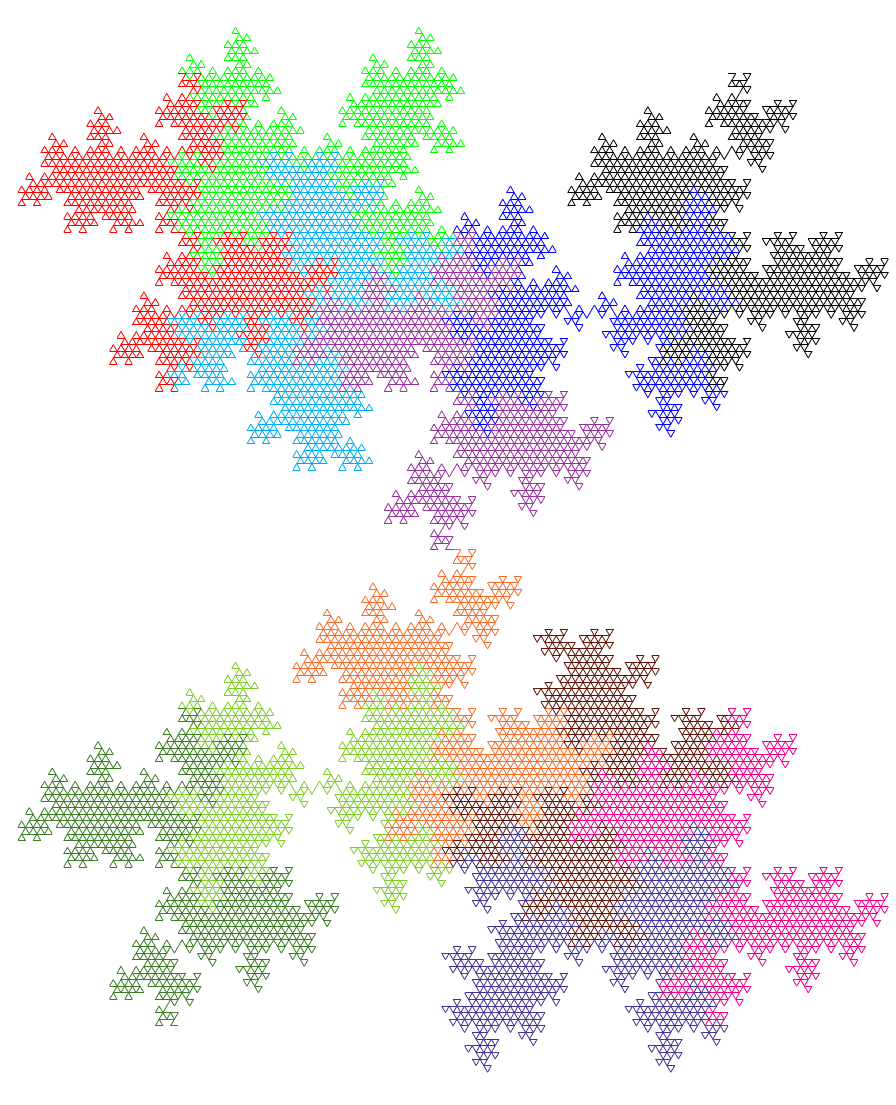}}%
{\includegraphics*[width=63mm, viewport={70 160 490 750}]{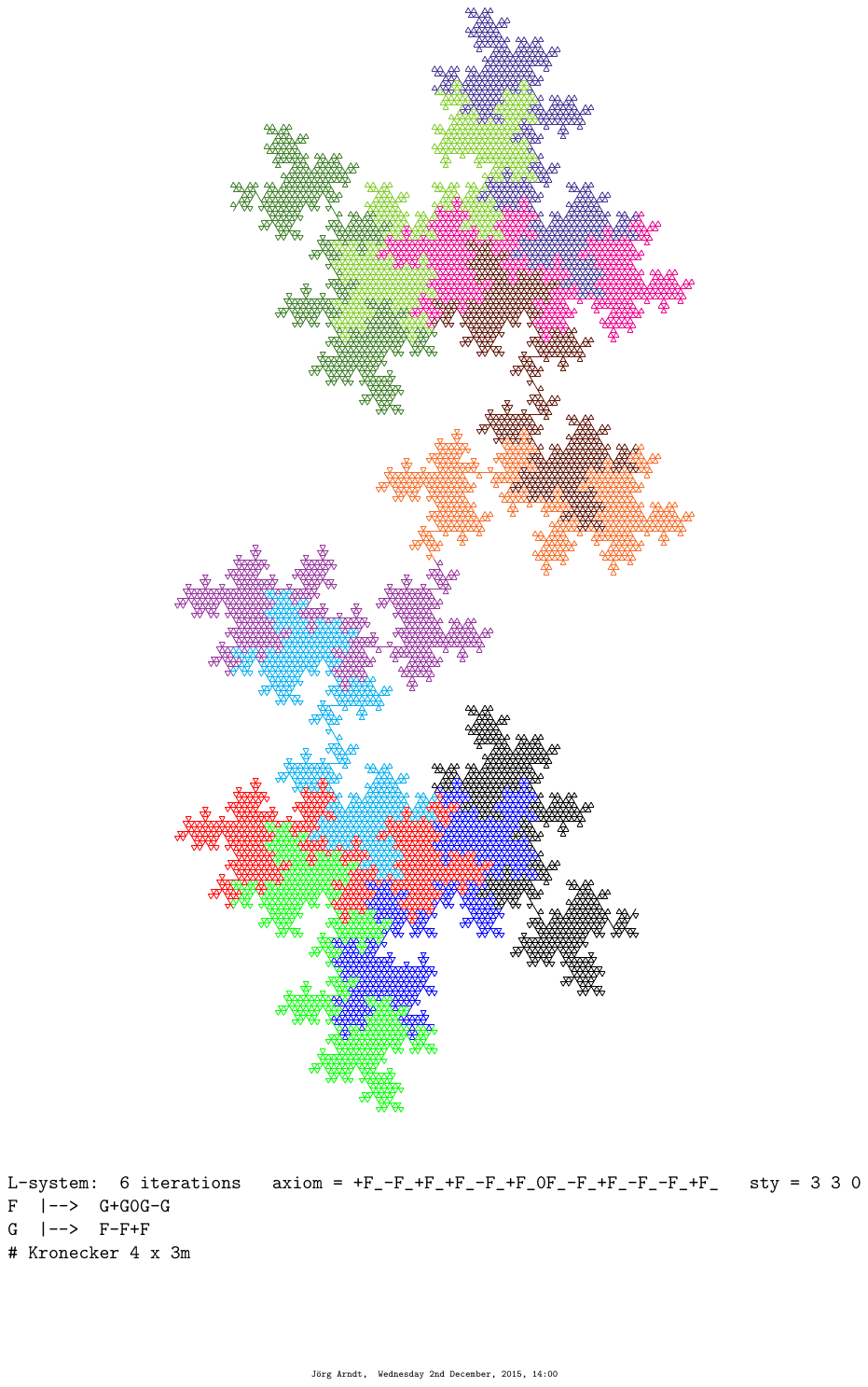}}
\end{center}
\else
\verb+{see pdf for image}+
\fi
\caption{\label{fig:kron-3m-4}
Products of the reversed terdragon (\Lmap{F}{F-F+F}) and the curve \Lmap{F}{F+F0F-F}.
The left curve is \CID{R12-13}, the right one is \CID{R12-25}.}
% R12-13 (p.8/10)
% R12-25 (p.10/10)
\end{figure}
%
%%%%%%%%%%%%%%%%%%%%%%%%%%

For two simple L-systems $L_1$ and $L_2$
with maps \Lmap{F}{f(F)} and \Lmap{F}{g(F)}
we define the \jjterm{product} $L_1 L_2$
as the L-system with map \Lmap{F}{g(f(F))}.
This product is obviously non-commutative in general.

For an example we use the curves on the triangular grid
with L-systems with maps
\Lmap{F}{F+F-F} (terdragon, $L_1$)
and
\Lmap{F}{F+F0F-F} (crab, $L_2$).
The products $L_1 L_2$ and $L_2 L_1$ respectively
give the curves shown in Figure~\ref{fig:kron-3-4}.

The map of the product $L_1 L_2$ is computed
as second iterate of \texttt{F} with maps
\Lmap{F}{G+G-G} and \Lmap{G}{F+F0F-F}.
For $L_2 L_1$ use the L-system
with maps \Lmap{F}{G+G0G-G} and \Lmap{G}{F+F-F}.
The respective maps for the simple L-systems are
\Lmap{F}{(F+F0F-F) + (F+F0F-F) - (F+F0F-F)}
and
\Lmap{F}{(F+F-F) + (F+F-F) 0 (F+F-F) - (F+F-F)}
where spaces and parentheses emphasize the structure.

Fixing the axiom \texttt{F} and using $k$ curves of pairwise different shapes we
obtain $k!$ new curves as products, one for every permutation of the maps.
The order of the product curve is the product of the orders of all curves used.

Still more curves can be obtained
by reversing the map in one
(or more) curves in the product.
Reversing the map for the terdragon
(to \Lmap{F}{F-F+F}) in the example
leads to respectively the maps
\Lmap{F}{(F+F0F-F) - (F+F0F-F) + (F+F0F-F)}
and
\Lmap{F}{(F-F+F) + (F-F+F) 0 (F-F+F) - (F-F+F)}
whose curves are shown in Figure~\ref{fig:kron-3m-4}.

In general the reversal of a map is not the same as the map with signs
swapped, so even more products can be made (this does not work for the
curves on the tri-hexagonal grid).
In general a curve is distinct from the ones obtained by reversal, swapping
signs, and doing both in the production of \texttt{F}, so even more products
exist.

%%%%%%%%%%%%%%%%%%%%%%%%%%
% stringsubst 6 F  F F+F-F  + - - + | tail -1 | ./bin 3 2 0 0 0.15 > tmp-pic.tex && make dotex # Kronecker 3 x 3m
%% better rotation:
% stringsubst 6 F F F+F-F + - - + | tail -1 | sed 's/+/++++/g; s/-/----/g; s/^/+/;' | ./bin 12 2 0 0 0.15 > tmp-pic.tex && make dotex # Kronecker 3 x 3m
%
\begin{figure}[h!tbp]
\ifpdf
\begin{center}
%% orig:
%{\includegraphics*[width=72mm, viewport={80 395 490 735}]{kronecker-3-m3.pdf}}
%% rotated:
{\includegraphics*[width=82mm, viewport={50 480 490 735}]{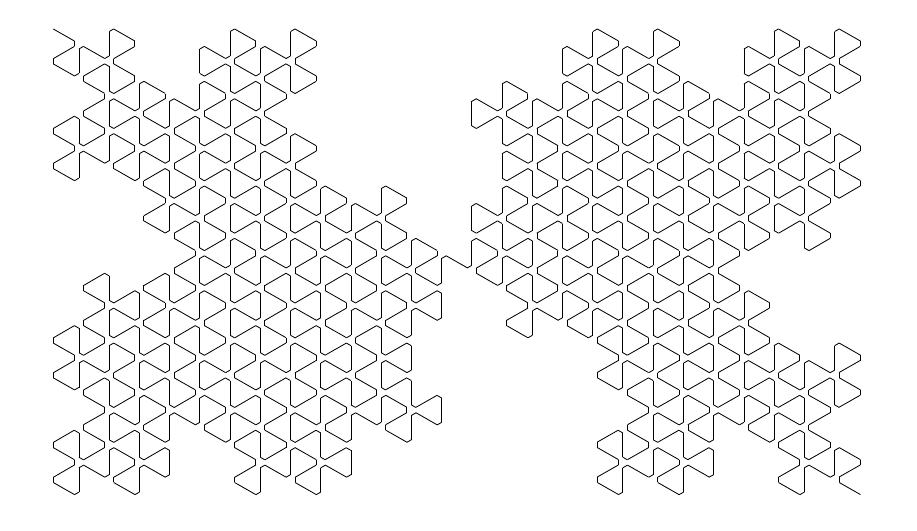}}
\end{center}
\else
\verb+{see pdf for image}+
\fi
\caption{\label{fig:kron-3-m3}
Curve (\CID{R9-8}) for the product of the maps \Lmap{F}{F+F-F} and \Lmap{F}{F-F+F}.}
\end{figure}
%
%%%%%%%%%%%%%%%%%%%%%%%%%%

%%%%%%%%%%%%%%%%%%%%%%%%%%
% stringsubst 5 _F+_F+_F-_F-_F _ _ F F+F+F-F-F + - - + | tail -1 | ./bin 4 3 0 > tmp-pic.tex && make dotex # Kronecker 5 x 5m
%% better rotation:
% stringsubst 5 _F+_F+_F-_F-_F _ _ F F+F+F-F-F + - - + | tail -1 | sed 's/+/++/g; s/-/--/g; s/^/+/;' | ./bin 8 3 0 > tmp-pic.tex && make dotex # Kronecker 5 x 5m
%
\begin{figure}[h!tbp]
\ifpdf
\begin{center}
%% orig:
%{\includegraphics*[width=92mm, viewport={60 375 490 745}]{kronecker-5-m5.pdf}}
%% rotated:
{\includegraphics*[width=100mm, viewport={60 480 490 740}]{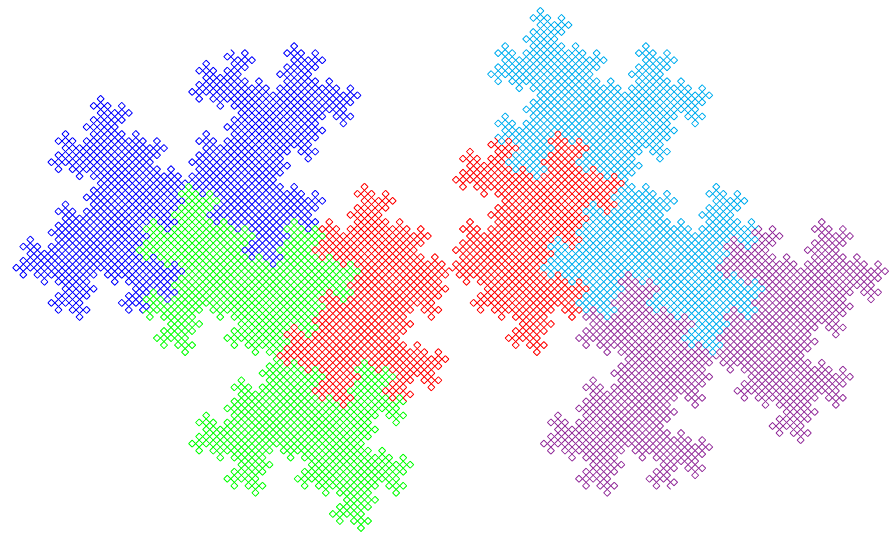}}
\end{center}
\else
\verb+{see pdf for image}+
\fi
\caption{\label{fig:kron-5-m5}
Curve (\CID{R25-45}) for the product of the maps \Lmap{F}{F+F+F-F-F} and \Lmap{F}{F-F-F+F+F},
consisting of 5 smaller copies of itself that are flipped over.}
\end{figure}
% Check:
% F |--> F+F+F-F-F-F+F+F-F-F-F+F+F-F-F+F+F+F-F-F+F+F+F-F-F
% F |--> (F+F+F-F-F) - (F+F+F-F-F) - (F+F+F-F-F) + (F+F+F-F-F) + (F+F+F-F-F)
% # R25-45 # dragon # symm-dr ## T+ = 36+ P ## T- = 44- P
%%%%%%%%%%%%%%%%%%%%%%%%%%

%%%%%%%%%%%%%%%%%%%%%%%%%%
% stringsubst 4 +F_-F_+F_+F_F_-F_F_  _ _ F F-F+F+F0F-F0F + - - +  0 0 | tail -1 | ./bin 3 3 0 > tmp-pic.tex && make dotex # Kronecker 7 x 7m (using R7-1)
%
\begin{figure}[h!tbp]
\ifpdf
\begin{center}
{\includegraphics*[width=110mm, viewport={60 450 490 745}]{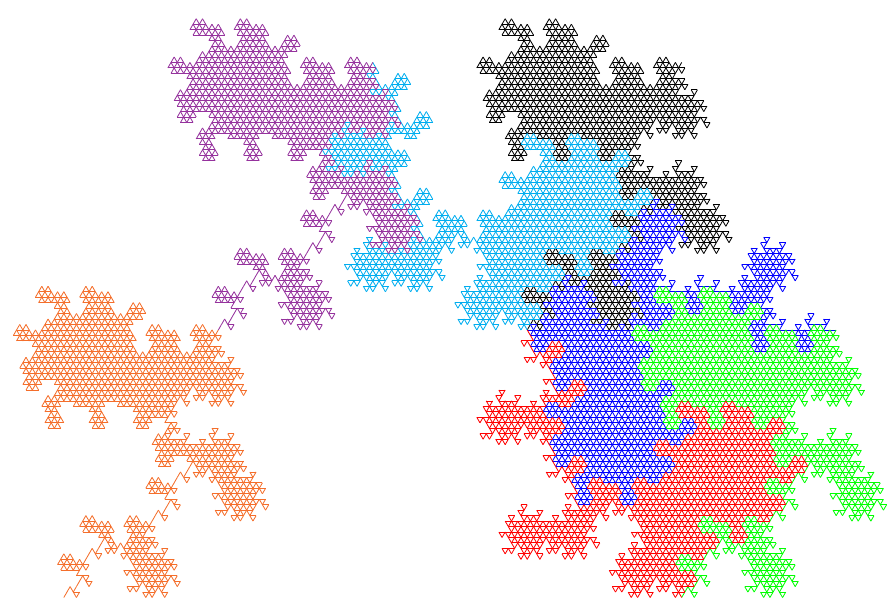}}
\end{center}
\else
\verb+{see pdf for image}+
\fi
\caption{\label{fig:kron-7-m7}
Curve for the product of the maps \Lmap{F}{F-F+F+F0F-F0F} and \Lmap{F}{F+F-F-F0F+F0F},
consisting of 7 smaller copies of itself that are flipped over.}
\end{figure}
%
%%%%%%%%%%%%%%%%%%%%%%%%%%

%%%%%%%%%%%%%%%%%%%%%%%%%%
% stringsubst 4 +G_0G_-G_0G_+G_+G_-G  _ _ F G0G-G0G+G+G-G G F+F-F-F0F+F0F 0 0 + + - - | tail -1 | ./bin 3 3 0 > tmp-pic.tex && make dotex
%
\begin{figure}[h!tbp]
\ifpdf
\begin{center}
{\includegraphics*[width=110mm, viewport={60 450 490 745}]{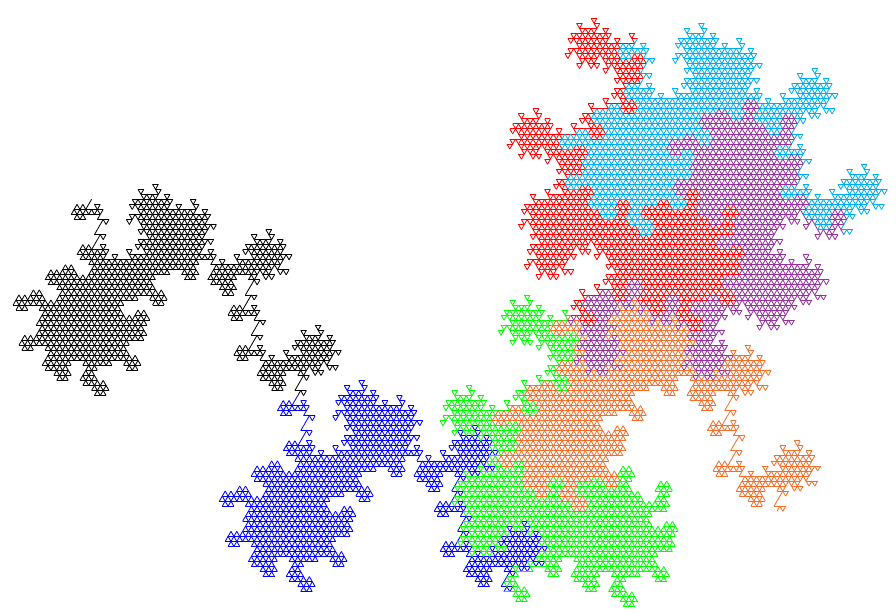}}
\end{center}
\else
\verb+{see pdf for image}+
\fi
\caption{\label{fig:kron-7-ventrella}
Curve for the product of the maps \Lmap{F}{F0F-F0F+F+F-F} and \Lmap{F}{F+F-F-F0F+F0F}.}
\end{figure}
%
%%%%%%%%%%%%%%%%%%%%%%%%%%
%% ventrella page 107: osc. signs of R7-1
%% stringsubst 4 +F_-F_+F_+F_F_-F_F_  _ _ F F-F+F+F0F-F0F + - - + | tail -1 | ./bin 3 2 0 > tmp-pic.tex && make dotex # R7-1, osc. signs

% Check: (for R7-1)
% stringsubst 2 F _ _ F F-F+F+F0F-F0F + - - + 0 0
% (F-F+F+F0F-F0F) + (F-F+F+F0F-F0F) - (F-F+F+F0F-F0F) - (F-F+F+F0F-F0F) (F-F+F+F0F-F0F) + (F-F+F+F0F-F0F) (F-F+F+F0F-F0F)
% stringsubst 4 G _ _ G F-F+F+F0F-F0F F G+G-G-G0G+G0G + + - - 0 0 | tail -1 | ./bin 3 2 0 > tmp-pic.tex && make dotex # Kronecker 7 x 7m (using R7-1)
% Different from Kronecker 7 x rev(7):
% stringsubst 4 F _ _ G F-F+F+F0F-F0F F G0G-G0G+G+G-G + + - - 0 0 | tail -1 | ./bin 3 2 0 > tmp-pic.tex && make dotex # Kronecker 7 x rev(7) (using R7-1)

A particularly simple way of multiplying an L-system with the L-system with
swapped signs is to use the maps \Lmap{+}{-} and \Lmap{-}{+}.
For example, the L-system
(\Lmap{F}{G+G-G}, \Lmap{G}{F-F+F}, \Lmap{+}{+}, \Lmap{-}{-})
is (for the even iterates) the same as
(\Lmap{F}{F-F+F}, \Lmap{+}{-}, \Lmap{-}{+}).
The resulting curve is shown in Figure~\ref{fig:kron-3-m3},
it is called the \jjterm{alternate terdragon} in \cite[Figure~21, p.~600]{davis-knuth-new},
it also appears on page 59 in \cite{ventrella}.

The curve shown at the bottom of page 84 in \cite{ventrella} can be obtained in
this way from the one at the top of the same page (the R5-dragon, also shown in
Figure~\ref{fig:r5-dragon}).  It is a curve of order $25$ (not $5$) according
to our conventions.
Figure~\ref{fig:kron-5-m5} shows the curve
decomposed into 5 smaller copies of itself that are flipped over.
All product curves of this kind have this sort of self-similarity,
Figure~\ref{fig:kron-7-m7} shows this
for the curve \CID{R7-1} on the triangular grid,
shown in Figure~\Ref{fig:iterate-3-4-5-decomp}.
% \Ref{fig:iterate-1-2-decomp} and
%
A similar (but different) curve appears on the bottom of page 107 in \cite{ventrella},
it is shown in Figure~\ref{fig:kron-7-ventrella}.
%% Ventrella's is:
% stringsubst 4 G _ _ F G+G-G-G0G+G0G G F0F-F0F+F+F-F 0 0 + + - - | tail -1 | ./bin 3 3 0 > tmp-pic.tex && make dotex
% stringsubst 4 F _ _ F G0G-G0G+G+G-G G F+F-F-F0F+F0F 0 0 + + - - | tail -1 | ./bin 3 3 0 > tmp-pic.tex && make dotex
% stringsubst 4 G_0G_-G_0G_+G_+G_-G  _ _ F G0G-G0G+G+G-G G F+F-F-F0F+F0F 0 0 + + - - | tail -1 | ./bin 3 3 0 > tmp-pic.tex && make dotex
%% Ventrella's is 7 x rev(7) with osc. signs:
%  stringsubst 5 G _ _ G F-F+F+F0F-F0F F G0G-G0G+G+G-G + - - + | tail -1 | ./bin 3 2 0 > tmp-pic.tex && make dotex # Kronecker 7 x rev(7), osc. signs (using R7-1)

% Regarding section \ref{sect:stats}: \xxx{Prod. curves not modded out}

%\xxx{Do prods. with paperfold?}
%% for every curve on square grid with paperfold, see below
%
% stringsubst 5 L L A+B R A-B A L+R+L-R-L B R+L+R-L-R + + - - | tail -1 | ./bin 4 2 0 0 0.10 > tmp-pic.tex && make dotex #
%% second gen:
% L ---> L+R+L-R-L + R+L+R-L-R
% R ---> L+R+L-R-L - R+L+R-L-R
%% so
% stringsubst 4 L L L+R+L-R-L+R+L+R-L-R R L+R+L-R-L-R+L+R-L-R + + - - | tail -1 | ./bin 4 2 0 0 0.15 > tmp-pic.tex && make dotex #
%% swap maps:
% stringsubst 4 L R L+R+L-R-L+R+L+R-L-R L L+R+L-R-L-R+L+R-L-R + + - - | tail -1 | ./bin 4 2 0 0 0.15 > tmp-pic.tex && make dotex #
%% ventrella p.160 (bottom)
%% both have the same motif and motif of tile
%
% stringsubst 5 L _ _ A L+R+L-R-L B R+L+R-L-R L A+B R A-B + + - - | tail -1 | ./bin 4 3 0 > tmp-pic.tex && make dotex #

\FloatBarrier% xxx
%%%%%%%%%%%%%%%%%%%%%%%%%%%%%%%
\subsection{Divisions of a curve}%\label{sect:}

\begin{figure}[h!tbp]
\ifpdf
\begin{center}
{\includegraphics*[width=85mm, page=1, viewport={60 490 490 740}]{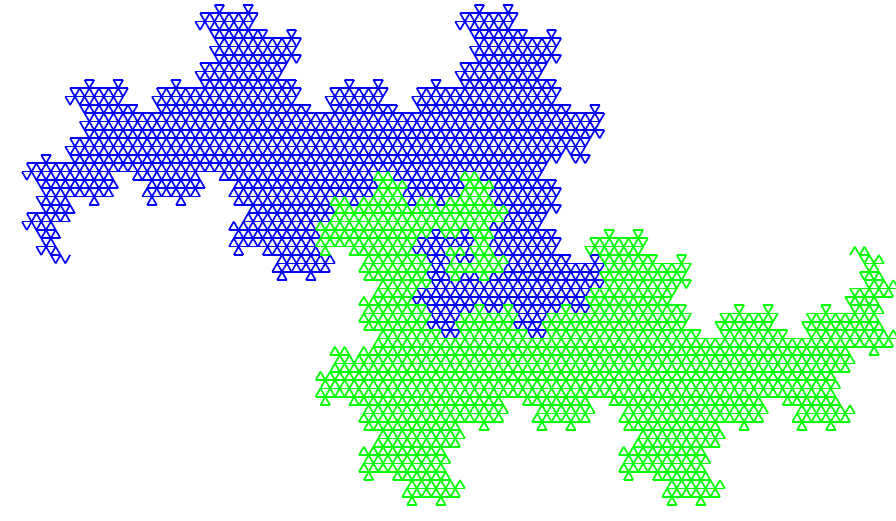}}
{\includegraphics*[width=85mm, page=2, viewport={60 490 490 740}]{terdragon-div-2-5.pdf}}
\end{center}
\else
\verb+{see pdf for image}+
\fi
\caption{\label{fig:terdragon-roots-2-3}
Dividing the terdragon (\CID{R3-1}) by 2 and 3.}
\end{figure}
%
%%%%%%%%%%%%%%%%%%%%%%%%%%

%%%%%%%%%%%%%%%%%%%%%%%%%%
%% with lnth *= 2.0;  // thicker lines
% stringsubst 8 +_X_Y  _ _  X XY+X  Y Y-XY  + + - - | tail -1 | ./bin 3 3 0 > tmp-pic.tex && make dotex # terdragon div 2
% stringsubst 8 +_X_Y_Z  _ _  X XYZ+  Y XYZ  Z -XYZ  + + - - | tail -1 | ./bin 3 3 0 > tmp-pic.tex && make dotex # terdragon div 3
% stringsubst 8 +_X_Y_Z_U  _ _  X XYZ  Y U+XY  Z ZU-X  U YZU  + + - - | tail -1 | ./bin 3 3 0 > tmp-pic.tex && make dotex # terdragon div 4
% stringsubst 8 +_X_Y_Z_U_V  _ _  X XYZ  Y UV+X  Z YZU  U V-XY  V ZUV  + + - - | tail -1 | ./bin 3 3 0 > tmp-pic.tex && make dotex # terdragon div 5
%
\begin{figure}[h!tbp]
\ifpdf
\begin{center}
{\includegraphics*[width=85mm, page=3, viewport={60 490 490 740}]{terdragon-div-2-5.pdf}}
{\includegraphics*[width=85mm, page=4, viewport={60 490 490 740}]{terdragon-div-2-5.pdf}}
\end{center}
\else
\verb+{see pdf for image}+
\fi
\caption{\label{fig:terdragon-roots-4-5}
Dividing the terdragon (\CID{R3-1}) by 4 and 5.}
\end{figure}
%
%%%%%%%%%%%%%%%%%%%%%%%%%%

The curve for a simple L-system of order $R$
can be divided into parts with equal numbers $d$ of edges as follows.
Replace each \texttt{F} in the production of \texttt{F}
by $s_1{}s_2{}\ldots{}s_d$ (where the $s_j$ are any pairwise different symbols),
then split the word obtained
into $d$ parts each containing $R$ letters $s_j$,
these are the productions of the letters $s_j$
for $1\leq{}j\leq{}d$.

The $j$th part ($1\leq{}j\leq{}d$) of the curve corresponds to the axiom $s_j$.

For example, to divide the terdragon ($R=3$) into
$d=5$ parts, first replace each \texttt{F} in
the production \texttt{F+F-F}
by \texttt{ABCDE} to obtain the word
\texttt{ABCDE+ABCDE-ABCDE},
split it into 5 parts of length 3,
used as productions, giving the maps
\Lmap{A}{ABC},
\Lmap{B}{DE+A},
\Lmap{C}{BCD},
\Lmap{D}{E-AB}, and
\Lmap{E}{CDE}.

The corresponding curves, concatenated to give the terdragon by using the axiom
%\texttt{ABCDE}, is shown in Figure~\ref{fig:terdragon-roots} (bottom).
\texttt{ABCDE}, is shown in Figure~\ref{fig:terdragon-roots-4-5} (bottom).
The division of a curve of
order $R$ by $d=R$ (here $3$, see bottom of Figure~\ref{fig:terdragon-roots-2-3})
is just another way to observe the self-similarity of the curve.
The division by two (see top of Figure~\ref{fig:terdragon-roots-2-3})
appears on page 59 in \cite{ventrella}.
%
% ventrella p.81: half crab

%%%%%%%%%%%%%%%%%%%%%%%%%%
%% with lnth *= 2.0;  // thicker lines
% stringsubst 5 +_A_B _ _  A AB0AB+AB+A  B B-AB-AB0AB  0 0 + + - - | tail -1 | ./bin 3 3 0 > tmp-pic.tex && make dotex # R7-2 div 2
%% smaller version:
% stringsubst 5 +_A_B _ _  A AB0AB+AB+A  B B-AB-AB0AB  0 0 + + - - | tail -1 | sed 's/AB/F/g; s/A_B/_F/g;' | ./bin 3 3 0 > tmp-pic.tex && make dotex # R7-2 div 2
%
%
\begin{figure}[h!tbp]
\ifpdf
\begin{center}
{\includegraphics*[width=100mm, viewport={60 370 490 740}]{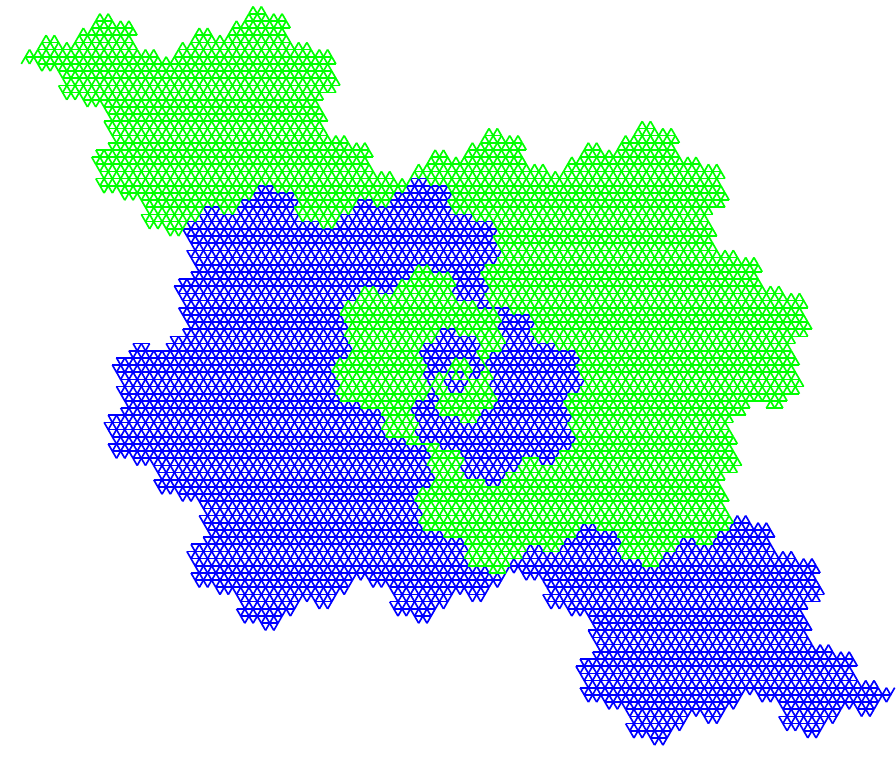}}
\end{center}
\else
\verb+{see pdf for image}+
\fi
\caption{\label{fig:r7-t-2-div-2}
%Dividing two curves of order 7 with identical shape by 2.}
Division of the curve \CID{R7-2} on the triangular grid by 2.}
\end{figure}
%
%%%%%%%%%%%%%%%%%%%%%%%%%%

%%%%%%%%%%%%%%%%%%%%%%%%%%
%% with lnth *= 2.0;  // thicker lines
% stringsubst 5 +_A_B _ _  A AB+AB-AB-A  B B+AB+AB-AB  0 0 + + - - | tail -1 | ./bin 3 3 0 > tmp-pic.tex && make dotex # R7-5 div 2
%% smaller version:
% stringsubst 5 +_A_B _ _  A AB+AB-AB-A  B B+AB+AB-AB  0 0 + + - - | tail -1 | sed 's/AB/F/g; s/A_B/_F/g;' | ./bin 3 3 0 > tmp-pic.tex && make dotex # R7-5 div 2
%
\begin{figure}[h!tbp]
\ifpdf
\begin{center}
{\includegraphics*[width=100mm, viewport={60 370 490 740}]{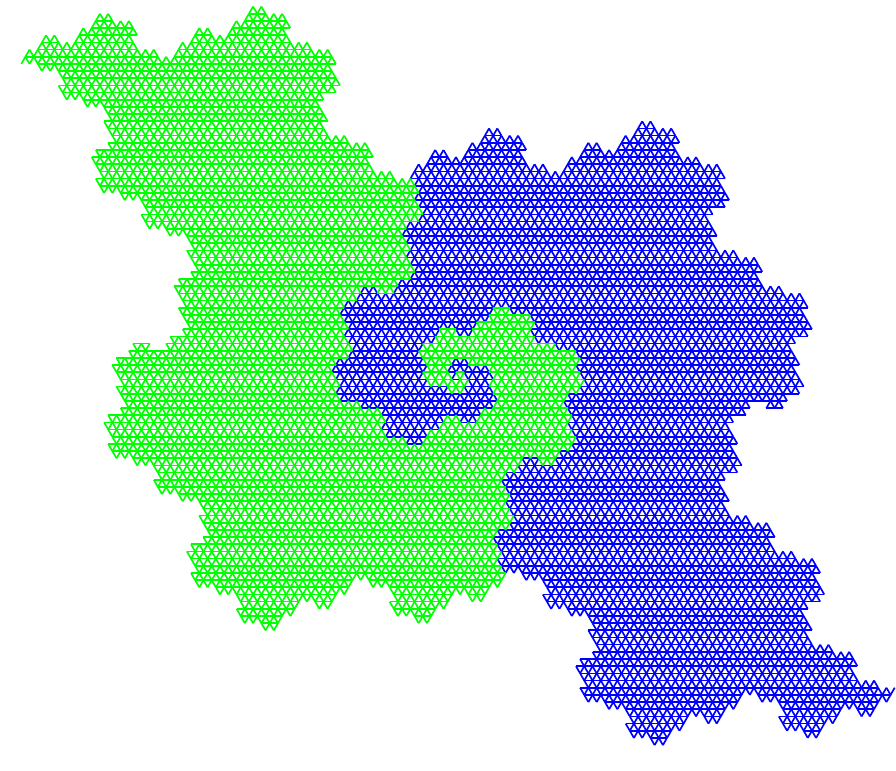}}
\end{center}
\else
\verb+{see pdf for image}+
\fi
\caption{\label{fig:r7-t-5-div-2}
%Dividing two curves of order 7 with identical shape by 2.}
Division of the curve \CID{R7-5} on the triangular grid by 2.}
\end{figure}
%
%%%%%%%%%%%%%%%%%%%%%%%%%%

The shape of the curve obtained by division is dependent on the curve,
not only on its shape.
The division by 2 of the curves \CID{R7-2} and \CID{R7-5}
whose motifs are shown in Figure~\Ref{fig:shape-7}
appear on pages 101 and 103 in \cite{ventrella},
see Figures~\ref{fig:r7-t-2-div-2} and \ref{fig:r7-t-5-div-2}.

%
%% half R7-5: ventrella p.101
% stringsubst 5 F _ _ F F+F-F-F+F+F-F  0 0 + + - - | tail -1 | ./bin 3 2 0 > tmp-pic.tex && make dotex # R7-5
% stringsubst 5 +_A_B _ _  A AB+AB-AB-A  B B+AB+AB-AB  0 0 + + - - | tail -1 | sed 's/+F/+F+F/g; s/-F/-F-F/g;' | ./bin 3 2 0 > tmp-pic.tex && make dotex # R7-5 div 2
% stringsubst 5 _A_B_C  _ _  A ABC+ABC-A  B BC-ABC+AB  C C+ABC-ABC  0 0 + + - - | tail -1 | ./bin 3 2 0 0 0.15 > tmp-pic.tex && make dotex # R7-5 div 3
% stringsubst 5 _A_B_C_D _ _  A ABCD+ABC  B D-ABCD-AB  C CD+ABCD+A  D BCD-ABCD  0 0 + + - - | tail -1 | ./bin 3 2 0 > tmp-pic.tex && make dotex # R7-5 div 4
%
%% half R7-2: ventrella p.103
% stringsubst 5 F _ _ F F0F+F+F-F-F0F  0 0 + + - - | tail -1 | ./bin 3 2 0 > tmp-pic.tex && make dotex # R7-2
% stringsubst 5 +_A_B _ _ A AB0AB+AB+A  B B-AB-AB0AB  0 0 + + - - | tail -1 | ./bin 3 2 0 > tmp-pic.tex && make dotex # R7-2 div 2

%\clearpage% xxx
%%%%%%%%%%%%%%%%%%%%%%%%%%
% stringsubst 4 _A_B _ _ A AB+A B B--AB--AB+AB+AB+AB - - + + | tail -1 | ./bin 6 3 0 > tmp-pic.tex && make dotex # R7-b
%
%% with lnth *= 2.0;  // thicker lines
%% also rotated
% stringsubst 5 R_A_B R R _ _ A AB+A B B--AB--AB+AB+AB+AB - - + + | tail -1 | sed 's/R/+/;' | ./bin 6 3 0 > tmp-pic.tex && make dotex # R7-b
%
\begin{figure}[h!tbp]
\ifpdf
\begin{center}
%{\includegraphics*[width=80mm, viewport={60 340 490 740}]{r7-b-unequal-div.pdf}}
{\includegraphics*[width=90mm, viewport={60 360 490 740}]{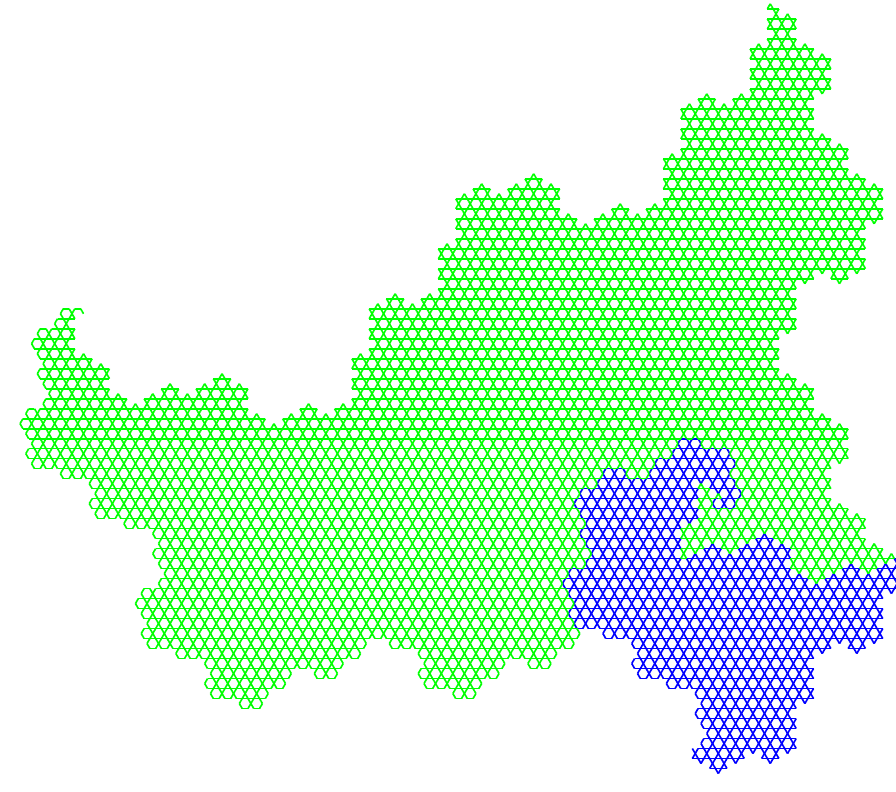}}
\end{center}
\else
\verb+{see pdf for image}+
\fi
\caption{\label{fig:r7-b-unequal-div}
Ventrella's curve split
into parts of $3/14$ and $11/14$ of the total length.}
\end{figure}
%
%%%%%%%%%%%%%%%%%%%%%%%%%%

Indeed a bit more is possible.
After substitution of the letter \texttt{F} in its own production
by a sequence of $d$ symbols, we can partition the word into maps of length $\neq{}d$.
As long as no sequence of maps
is just a cyclic permutation of the letters
such as \Lmap{B}{C}, \Lmap{C}{E}, \Lmap{E}{B}
(a fixed point such as \Lmap{B}{B} is a special case of this),
we obtain any partition of $R\,d$ into $d$ parts,
not just the partition into equal parts $R$.
For example, for Ventrella's curve with map \Lmap{F}{F+F+F+F--F--F+F}
(\CID{R7-1} on the tri-hexagonal grid see Figure~\Ref{fig:r7-b-curve}),
we can split the word \texttt{AB+AB--AB--AB+AB+AB+AB}
into the productions for the maps
\Lmap{A}{AB+A} and
\Lmap{B}{B--AB--AB+AB+AB+AB}
to partition the curve into
parts $3$ and $14-3=11$,
this is shown in Figure~\ref{fig:r7-b-unequal-div}.
%% unequal division:
% stringsubst 4 _A_B _ _ A AB+A B B--AB--AB+AB+AB+AB - - + + | tail -1 | ./bin 6 3 0 > tmp-pic.tex && make dotex # R7-b
% stringsubst 5 A_B_C _ _ A ABC+ABC--AB B C--ABC+ABC+ABC+A C BC - - + + | tail -1 | ./bin 6 3 0 0 0.1 > tmp-pic.tex && make dotex # R7-b
%
% stringsubst 7 _A_B_C _ _ A ABC+A B BC-A C BC - - + + | tail -1 | ./bin 3 3 0 0 0.1 > tmp-pic.tex && make dotex
% stringsubst 7 _A_B_C _ _ A ABC+ABC- B A C BC - - + + | tail -1 | ./bin 3 3 0 0 0.1 > tmp-pic.tex && make dotex
% no subset of production must be a cyclic permutation:
% stringsubst 7 _A_B_C _ _ A A B B C C+ABC-ABC - - + + | tail -1 | ./bin 3 3 0 0 0.1 > tmp-pic.tex && make dotex

%%%%%%%%%%%%%%%%%%%%%%%%%%
%
% stringsubst 5 _A__B _ _ A AB+CD++CD-A B B--ABAB-CD+ C -AB+CDCD++C D D+AB--AB-CD + + - - | tail -1 | ./bin 6 3 0 > tmp-pic.tex && make dotex # Gosper's "flowsnake" curve, div 2
%
% stringsubst 6 _A__B _ _ A +CDt-ABtA B B-tCD+ C -ABt+CDtC D D+tAB- + + - - t t | tr -d ABCD | tail -1 | ./bin 4 3 0 > tmp-pic.tex && make dotex # Hilbert curve, div
%
\begin{figure}[h!tbp]
\ifpdf
\begin{center}
{\includegraphics*[width=63mm, viewport={60 310 495 740}]{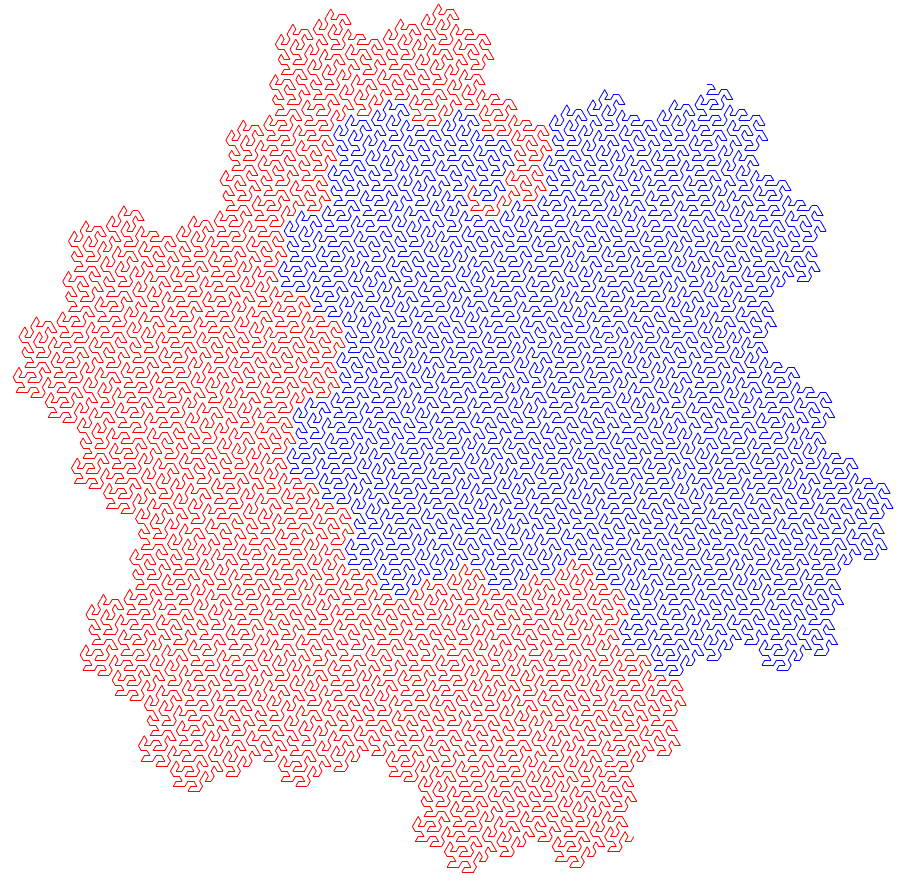}}%
{\includegraphics*[width=63mm, viewport={60 310 495 750}]{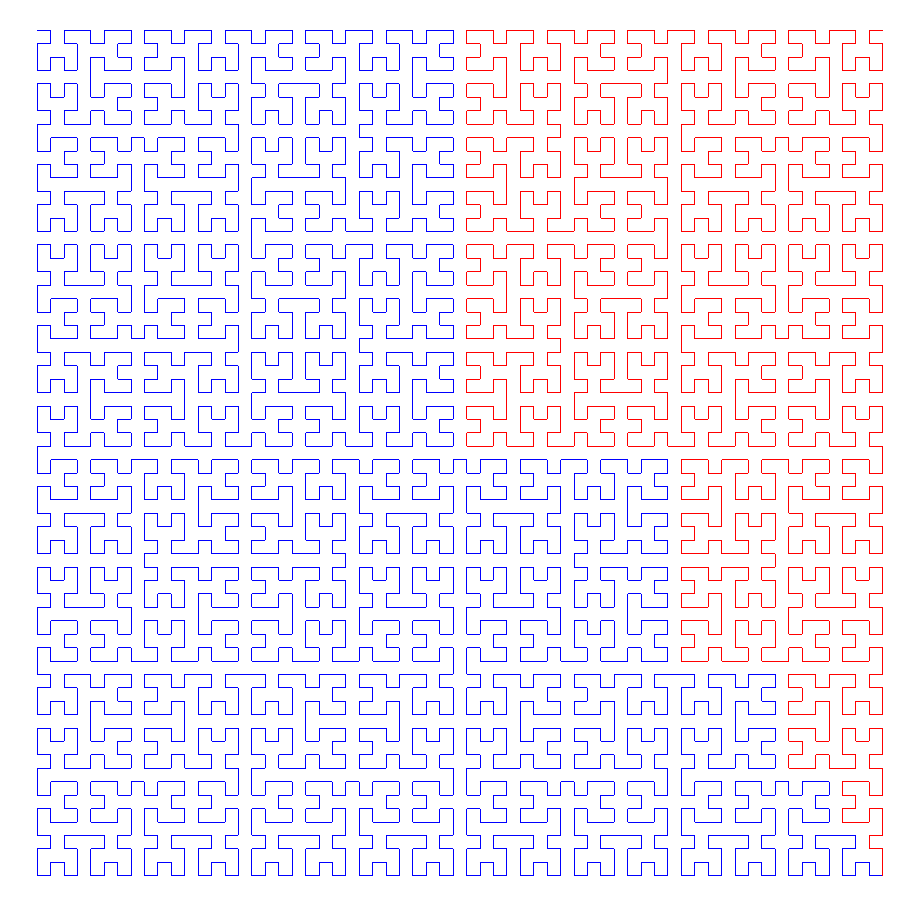}}
\end{center}
\else
\verb+{see pdf for image}+
\fi
\caption{\label{fig:gosper-hilbert-div}
Gosper's flowsnake split into halves (left) and
the Hilbert curve split into parts of $2/3$ and $1/3$ of the total length (right).}
\end{figure}
%
%%%%%%%%%%%%%%%%%%%%%%%%%%

Splitting curves with non-simple L-systems can similarly be achieved.
For Gosper's curve with maps
\Lmap{L}{L+R++R-L--LL-R+} and
\Lmap{R}{-L+RR++R+L--L-R},
replace all
\texttt{L} by \texttt{AB}
and all
\texttt{R} by \texttt{CD}
to obtain the two words
\texttt{AB+CD++CD-AB--ABAB-CD+}
and
\texttt{-AB+CDCD++CD+AB--AB-CD}.
Divide both evenly, giving the maps
\Lmap{A}{AB+CD++CD-A},
\Lmap{B}{B--ABAB-CD+},
\Lmap{C}{-AB+CDCD++C}, and
\Lmap{D}{D+AB--AB-CD}.
The split into halves using the axiom \texttt{AB}
is shown at the left in Figure~\ref{fig:gosper-hilbert-div}.

For the Hilbert curve
with the maps
\Lmap{L}{+Rt-LtL-tR+} and
\Lmap{R}{-Lt+RtR+tL-}
(where finally only the letters \texttt{t} are used for drawing edges),
we obtain the two words
\texttt{+CDt-ABtAB-tCD+}
and
\texttt{-ABt+CDtCD+tAB-}.
Note that the letters \texttt{t} are kept unmodified.
Here we divide unevenly to make the picture more interesting, using the maps
\Lmap{A}{+CDt-ABtA},
\Lmap{B}{B-tCD+},
\Lmap{C}{-ABt+CDtC}, and
\Lmap{D}{D+tAB-}.
The split resulting from the axiom \texttt{AB}
% \xxx{below is over-simplified}
%gives parts $2/3$ and $1/3$
%as there are two letters \texttt{t}
%in the rules for \texttt{A} and \texttt{C}
%and one in the rules for \texttt{B} and \texttt{D}.
is shown at the right in Figure~\ref{fig:gosper-hilbert-div}.

% stringsubst 3 L   L L+R++R-L--LL-R+   R -L+RR++R+L--L-R   + +   - - | tail -1 | ./bin 6 3 0 > tmp-pic.tex && make dotex # Gosper's "flowsnake" curve
% stringsubst 4 _A__B _ _ A AB+CD++CD-A B B--ABAB-CD+ C -AB+CDCD++C D D+AB--AB-CD + + - - | tail -1 | ./bin 6 3 0 > tmp-pic.tex && make dotex
%
% stringsubst 4 L  L +Rt-LtL-tR+  R -Lt+RtR+tL-  + + - - t t | tr -d LR | tail -1 | ./bin 4 3 0 > tmp-pic.tex && make dotex # Hilbert curve
% stringsubst 6 _A__B _ _ A +CDt-ABtA B B-tCD+ C -ABt+CDtC D D+tAB- + + - - t t | tr -d ABCD | tail -1 | ./bin 4 3 0 > tmp-pic.tex && make dotex # Hilbert curve, div

%%% Emacs:
%%% Local Variables:
%%% mode: latex
%%% MyRelDir: "."
%%% TeX-master: "arndt-curve-search.tex"
%%% dvi-file: "arndt-curve-search"
%%% makefile-dir: "./"
%%% frame-title-format: "CURVE-SEARCH (products)"
%%% End:

%\FloatBarrier

\clearpage% xxx

%%%%%%%%%%%%%%%%%%%%%%%%%%%%%%%%%%%%%%%%%%%%%%%%%%%%%%%%%%%%%%
%%%%%%%%%%%%%%%%%%%%%%%%%%%%%%%%%%%%%%%%%%%%%%%%%%%%%%%%%%%%%%
\section{Eye of the beholder}

%% attic/all-r12-t-decomp.txt
%% attic/all-r17-q-decomp.txt
%% attic/all-r19-b-decomp.txt

We show the fourth iterates for all
shapes of curves for the triangular, square, and tri-hexagonal grid
of respectively orders $R=12$, $R=17$, and $R=19$,
using coloration to make the self-similarity apparent.

%%%%%%%%%%%%%%%%%%%%%%%%%%%%%%%%%%%%%%%%%%%%%%%%%%%%%%%%%%%%%%
%\subsection{Triangular grid}

% triangular grid: R12, 10 total, minus 4 shown as products ==> 6

%% ./mk-decomp-rule.pl < c3/search-r12-curves.txt | grep -v same | sed 's/ 3 F/ 2 F/;' > r12-t-decomp.txt
%% edit r12-t-decomp.txt (iter=2, and rotations)
%% zsh omfg.zsh r12-t-decomp.txt
%% cp -auvi all-omfg.pdf all-r12-t-decomp.pdf

%%%%%%%%%%%%%%%%%%%%%%%%%%
% ./scbin 12 3 0 c3/search-r12-curves.txt | grep -v same | ./shorten-output.pl
\begin{figure}[h!tbp]
{\small
\begin{verbatim}
  01:  F F0F+F0F-F+F-F-F0F+F+F-F  R12-1  #
  02:  F F0F+F+F-F0F-F-F+F0F+F-F  R12-3  #
  03:  F F0F+F+F-F-F0F+F+F-F-F0F  R12-4  # # symm-dr
  04:  F F0F+F-F+F0F+F-F-F0F-F+F  R12-6  #
  05:  F F+F0F+F-F+F0F-F+F-F0F-F  R12-9  # # symm-dr
  06:  F F+F0F+F-F-F0F+F+F-F0F-F  R12-10  # # symm-dr
  07:  F F+F0F-F+F+F-F+F-F-F+F-F  R12-12  #
  08:  F F+F0F-F-F+F0F-F+F+F0F-F  R12-13  # # symm-dr
  09:  F F+F-F+F+F-F0F+F-F-F+F-F  R12-17  # # symm-dr
  10:  F F+F-F-F+F-F0F+F-F+F+F-F  R12-25  # # symm-dr
\end{verbatim}
}
\caption{\label{fig:output-12-t-no-same}
The curves of order 12 on triangular grid without those duplicating a prior shape.}
\end{figure}
%
%%%%%%%%%%%%%%%%%%%%%%%%%

For the triangular grid we skip four curves already seen,
\CID{R12-10} and \CID{R12-17}
appear in Figure~\ref{fig:kron-3-4},
\CID{R12-13} and \CID{R12-25}
in Figure~\ref{fig:kron-3m-4}.

%%%%%%%%%%%%%%%%%%%%%%%%%%%%%%%%%%%%%%%%%%%%%%%%%%%%%%%%%%%%%%
%\subsection{Square grid}

% square grid: R17 ==> 6

%%%%%%%%%%%%%%%%%%%%%%%%%%
% ./scbin 17 4 0 c4/search-r17-q-curves.txt | grep -v same | ./shorten-output.pl
\begin{figure}[h!tbp]
{\small
\begin{verbatim}
  01:  F F+F+F-F+F-F-F-F+F-F+F+F+F-F+F-F-F  R17-1  # # symm-dr
  02:  F F+F+F-F-F+F+F+F-F+F-F-F-F+F+F-F-F  R17-2  # # symm-dr
  03:  F F+F+F-F-F+F+F+F-F-F+F+F-F-F-F+F-F  R17-3  #
  04:  F F+F-F+F+F+F-F-F+F-F+F+F-F-F-F+F-F  R17-6  # # symm-dr
  05:  F F+F-F+F+F-F+F+F-F-F-F+F+F-F-F-F+F  R17-10  #
  06:  F F+F-F+F+F-F+F-F+F+F-F-F-F+F-F-F+F  R17-11  #
\end{verbatim}
}
\caption{\label{fig:output-17-q-no-same}
The curves of order 17 on the square grid without those duplicating a prior shape.}
\end{figure}
%
%%%%%%%%%%%%%%%%%%%%%%%%%

%For the square grid the curves of order $R=17$ are shown.
%

%% ./mk-decomp-rule.pl < c4/search-r17-q-curves.txt | grep -v same | sed 's/ 3 F/ 2 F/;' > r17-q-decomp.txt
%% edit r17-q-decomp.txt (iter=2, and rotations)
%% zsh omfg.zsh r17-q-decomp.txt
%% cp -auvi all-omfg.pdf all-r17-q-decomp.pdf

%%%%%%%%%%%%%%%%%%%%%%%%%%%%%%%%%%%%%%%%%%%%%%%%%%%%%%%%%%%%%%
%\subsection{Tri-hexagonal grid}

%
% tri-hexagonal grid: R19, 7 total, minus 1 shown, ==> 6
%%%%%%%%%%%%%%%%%%%%%%%%%%
%  ./scbin 19 6 0 c6/search-r19-b-curves.txt | grep -v same | ./shorten-output.pl
\begin{figure}[h!tbp]
{\small
\begin{verbatim}
  01:  F F+F+F+F+F--F+F+F+F--F+F+F--F--F+F+F--F+F--F  R19-1  #
  02:  F F+F+F+F+F--F--F+F+F--F--F+F+F--F+F+F+F+F--F  R19-2  #
  03:  F F+F+F+F--F+F--F--F+F+F+F--F--F+F+F+F+F+F--F  R19-3  #
  04:  F F+F+F+F--F--F+F+F--F--F+F+F--F+F+F+F+F--F+F  R19-5  #
  05:  F F+F+F--F+F+F+F--F--F+F--F+F--F+F+F+F+F+F--F  R19-6  #
  06:  F F+F+F--F+F--F--F+F+F+F--F+F--F+F+F+F+F+F--F  R19-7  #
  07:  F F+F--F+F+F+F--F--F+F--F+F+F--F+F+F+F+F+F--F  R19-14  #
\end{verbatim}
}
\caption{\label{fig:output-19-b-no-same}
The curves of order 19 on the tri-hexagonal grid without those duplicating a prior shape.}
\end{figure}
%
%%%%%%%%%%%%%%%%%%%%%%%%%%

For the tri-hexagonal
% grid the curves of order $R=19$ are shown.
the third iterate of curve \CID{R19-1} already appeared in Figure~\ref{fig:r19-b-1-curve}.

%% ./mk-decomp-rule.pl < c6/search-r19-b-curves.txt | grep -v same | sed 's/ 3 F/ 2 F/;' > r19-b-decomp.txt
%% edit r19-b-decomp.txt (iter=2, and rotations)
%% zsh omfg.zsh r19-b-decomp.txt
%% cp -auvi all-omfg.pdf all-r19-b-decomp.pdf

%\FloatBarrier
\clearpage% xxx
%%%%%%%%%%%%%%%%%%%%%%%%%%%%%%%%%%%%%%%%%%%%%%%%%%%%
%% TRIANGULAR grid
\subsection{Triangular grid}
%\setcounter{subsection}{1}
%\stepcounter{subsection}
\quad\\% xxx layout

%%%%%%%%%%%%%%%%%%%%%%%%%%
%
\begin{figure}[h!tbp]
\ifpdf
\begin{center}
%\fbox{\includegraphics*[width=90mm, page=1, viewport={60 290 490 730}]{all-r12-t-decomp.pdf}}
{\includegraphics*[width=130mm, page=1, viewport={65 340 490 740}]{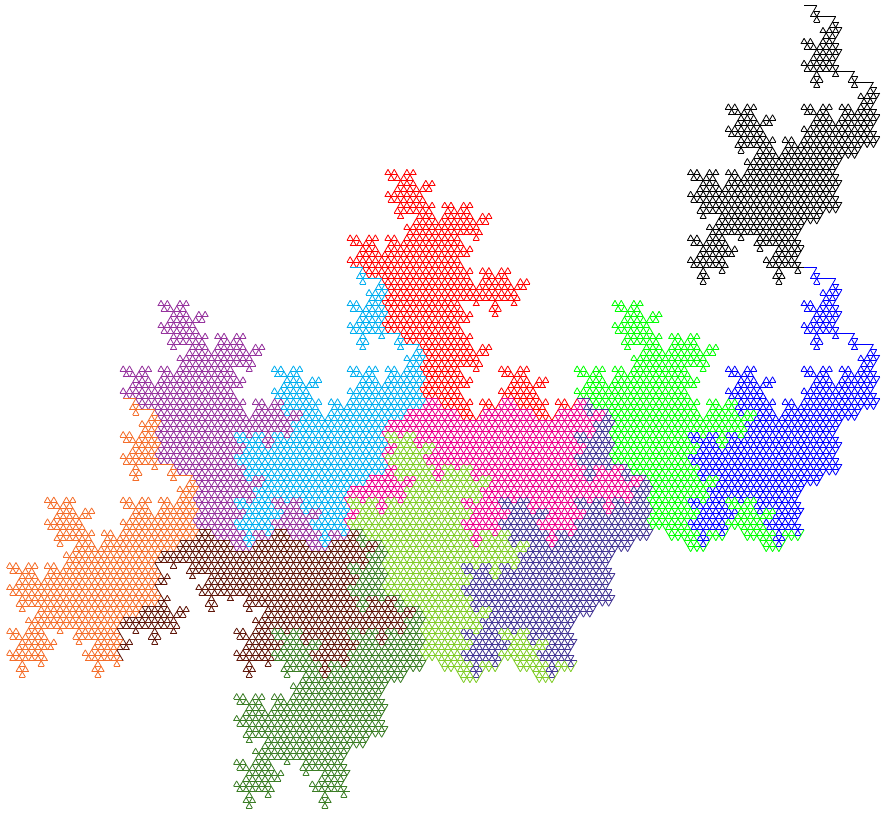}}
\end{center}
\else
\verb+{see pdf for image}+
\fi
\caption{\label{fig:R12-1}
\CID{R12-1} on the triangular grid.}
\end{figure}
%
%%%%%%%%%%%%%%%%%%%%%%%%%%

%%%%%%%%%%%%%%%%%%%%%%%%%%
%
\begin{figure}[h!tbp]
\ifpdf
\begin{center}
%\fbox{\includegraphics*[width=90mm, page=2, viewport={60 400 490 740}]{all-r12-t-decomp.pdf}}
{\includegraphics*[width=130mm, page=2, viewport={60 210 490 740}]{all-r12-t-decomp.pdf}}
\end{center}
\else
\verb+{see pdf for image}+
\fi
\caption{\label{fig:R12-3}
\CID{R12-3} on the triangular grid.}
\end{figure}
%
%%%%%%%%%%%%%%%%%%%%%%%%%%

%%%%%%%%%%%%%%%%%%%%%%%%%%
%
\begin{figure}[h!tbp]
\ifpdf
\begin{center}
%\fbox{\includegraphics*[width=90mm, page=3, viewport={60 420 490 740}]{all-r12-t-decomp.pdf}}
{\includegraphics*[width=130mm, page=3, viewport={60 490 490 740}]{all-r12-t-decomp.pdf}}
\end{center}
\else
\verb+{see pdf for image}+
\fi
\caption{\label{fig:R12-4}
\CID{R12-4} on the triangular grid.}
\end{figure}
%
%%%%%%%%%%%%%%%%%%%%%%%%%%

%%%%%%%%%%%%%%%%%%%%%%%%%%
%
\begin{figure}[h!tbp]
\ifpdf
\begin{center}
%\fbox{\includegraphics*[width=90mm, page=4, viewport={60 260 490 730}]{all-r12-t-decomp.pdf}}
{\includegraphics*[width=130mm, page=4, viewport={60 360 490 740}]{all-r12-t-decomp.pdf}}
\end{center}
\else
\verb+{see pdf for image}+
\fi
\caption{\label{fig:R12-6}
\CID{R12-6} on the triangular grid.}
\end{figure}
%
%%%%%%%%%%%%%%%%%%%%%%%%%%

%%%%%%%%%%%%%%%%%%%%%%%%%%
%
\begin{figure}[h!tbp]
\ifpdf
\begin{center}
%\fbox{\includegraphics*[width=90mm, page=5, viewport={60 340 490 740}]{all-r12-t-decomp.pdf}}
{\includegraphics*[width=130mm, page=5, viewport={60 460 490 740}]{all-r12-t-decomp.pdf}}
\end{center}
\else
\verb+{see pdf for image}+
\fi
\caption{\label{fig:R12-9}
\CID{R12-9} on the triangular grid.}
\end{figure}
%
%%%%%%%%%%%%%%%%%%%%%%%%%%
%% skip 6

%%%%%%%%%%%%%%%%%%%%%%%%%%
%
\begin{figure}[h!tbp]
\ifpdf
\begin{center}
%\fbox{\includegraphics*[width=90mm, page=7, viewport={60 300 490 740}]{all-r12-t-decomp.pdf}}
{\includegraphics*[width=130mm, page=7, viewport={60 430 490 740}]{all-r12-t-decomp.pdf}}
\end{center}
\else
\verb+{see pdf for image}+
\fi
\caption{\label{fig:R12-12}
\CID{R12-12} on the triangular grid.}
\end{figure}
%
%%%%%%%%%%%%%%%%%%%%%%%%%%
%% skip 8
%% skip 9
%% skip 10

%\FloatBarrier
\clearpage% xxx
%%%%%%%%%%%%%%%%%%%%%%%%%%%%%%%%%%%%%%%%%%%%%%%%%%%%
\subsection{Square grid}
%% SQUARE grid
%\stepcounter{subsection}
\quad\\% xxx layout

%%%%%%%%%%%%%%%%%%%%%%%%%%
%
\begin{figure}[h!tbp]
\ifpdf
\begin{center}
%\fbox{\includegraphics*[width=120mm, viewport={60 440 490 740}]{r17-q-all-1.pdf}}
%\fbox{\includegraphics*[width=120mm, page=1, viewport={60 440 490 740}]{all-r17-q-decomp.pdf}}
{\includegraphics*[width=130mm, page=1, viewport={60 290 490 740}]{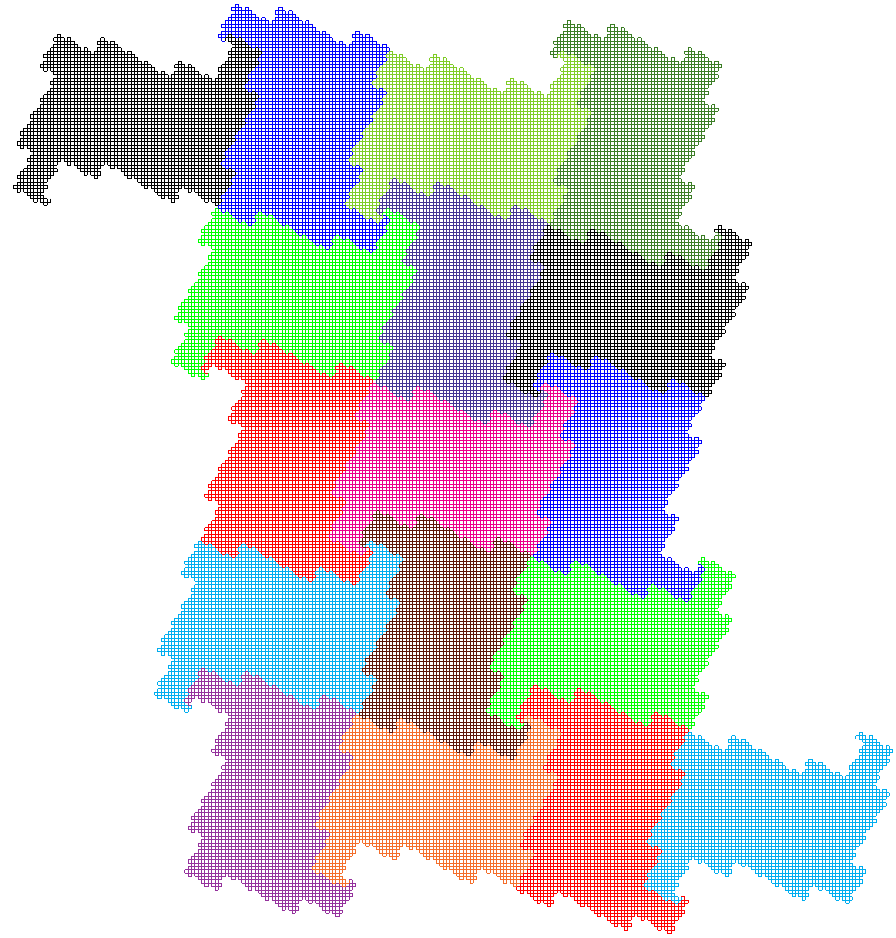}}
\end{center}
\else
\verb+{see pdf for image}+
\fi
\caption{\label{fig:R17-1}
\CID{R17-1} on the square grid.}
\end{figure}
%
%%%%%%%%%%%%%%%%%%%%%%%%%%

%%%%%%%%%%%%%%%%%%%%%%%%%%
%
\begin{figure}[h!tbp]
\ifpdf
\begin{center}
%\fbox{\includegraphics*[width=120mm, viewport={60 420 490 740}]{r17-q-all-2.pdf}}
%\fbox{\includegraphics*[width=120mm, page=2, viewport={60 420 490 740}]{all-r17-q-decomp.pdf}}
{\includegraphics*[width=130mm, page=2, viewport={60 160 490 740}]{all-r17-q-decomp.pdf}}
\end{center}
\else
\verb+{see pdf for image}+
\fi
\caption{\label{fig:R17-2}
\CID{R17-2} on the square grid.}
\end{figure}
%
%%%%%%%%%%%%%%%%%%%%%%%%%%

%%%%%%%%%%%%%%%%%%%%%%%%%%
%
\begin{figure}[h!tbp]
\ifpdf
\begin{center}
%\fbox{\includegraphics*[width=120mm, viewport={60 440 490 740}]{r17-q-all-3.pdf}}
%\fbox{\includegraphics*[width=120mm, page=3, viewport={60 440 490 740}]{all-r17-q-decomp.pdf}}
{\includegraphics*[width=130mm, page=3, viewport={60 160 490 740}]{all-r17-q-decomp.pdf}}
\end{center}
\else
\verb+{see pdf for image}+
\fi
\caption{\label{fig:R17-3}
\CID{R17-3} on the square grid.}
\end{figure}
%
%%%%%%%%%%%%%%%%%%%%%%%%%%

%%%%%%%%%%%%%%%%%%%%%%%%%%
%
\begin{figure}[h!tbp]
\ifpdf
\begin{center}
%\fbox{\includegraphics*[width=120mm, viewport={60 450 490 740}]{r17-q-all-4.pdf}}
%\fbox{\includegraphics*[width=120mm, page=4, viewport={60 450 490 740}]{all-r17-q-decomp.pdf}}
{\includegraphics*[width=130mm, page=4, viewport={60 160 490 740}]{all-r17-q-decomp.pdf}}
\end{center}
\else
\verb+{see pdf for image}+
\fi
\caption{\label{fig:R17-6}
\CID{R17-6} on the square grid.}
\end{figure}
%
%%%%%%%%%%%%%%%%%%%%%%%%%%

%%%%%%%%%%%%%%%%%%%%%%%%%%
%
\begin{figure}[h!tbp]
\ifpdf
\begin{center}
%\fbox{\includegraphics*[width=120mm, viewport={60 480 490 740}]{r17-q-all-5.pdf}}
%\fbox{\includegraphics*[width=120mm, page=5, viewport={60 480 490 740}]{all-r17-q-decomp.pdf}}
{\includegraphics*[width=130mm, page=5, viewport={60 160 490 740}]{all-r17-q-decomp.pdf}}
\end{center}
\else
\verb+{see pdf for image}+
\fi
\caption{\label{fig:R17-10}
\CID{R17-10} on the square grid.}
\end{figure}
%
%%%%%%%%%%%%%%%%%%%%%%%%%%

%%%%%%%%%%%%%%%%%%%%%%%%%%
%
\begin{figure}[h!tbp]
\ifpdf
\begin{center}
%\fbox{\includegraphics*[width=120mm, viewport={60 440 490 740}]{r17-q-all-6.pdf}}
%\fbox{\includegraphics*[width=120mm, page=6, viewport={60 450 490 740}]{all-r17-q-decomp.pdf}}
{\includegraphics*[width=130mm, page=6, viewport={60 160 490 740}]{all-r17-q-decomp.pdf}}
\end{center}
\else
\verb+{see pdf for image}+
\fi
\caption{\label{fig:R17-11}
\CID{R17-11} on the square grid.}
\end{figure}
%
%%%%%%%%%%%%%%%%%%%%%%%%%%

%\FloatBarrier
\clearpage% xxx
%%%%%%%%%%%%%%%%%%%%%%%%%%%%%%%%%%%%%%%%%%%%%%%%%%%%
\subsection{Tri-Hexagonal grid}
%% TRI-HEXAGONAL grid
%\stepcounter{subsection}
\quad\\% xxx layout

%%%%%%%%%%%%%%%%%%%%%%%%%%
%
\begin{figure}[h!tbp]
\ifpdf
\begin{center}
%\fbox{\includegraphics*[width=90mm, page=1, viewport={60 340 490 740}]{all-r19-b-decomp.pdf}}
{\includegraphics*[width=130mm, page=1, viewport={60 235 490 740}]{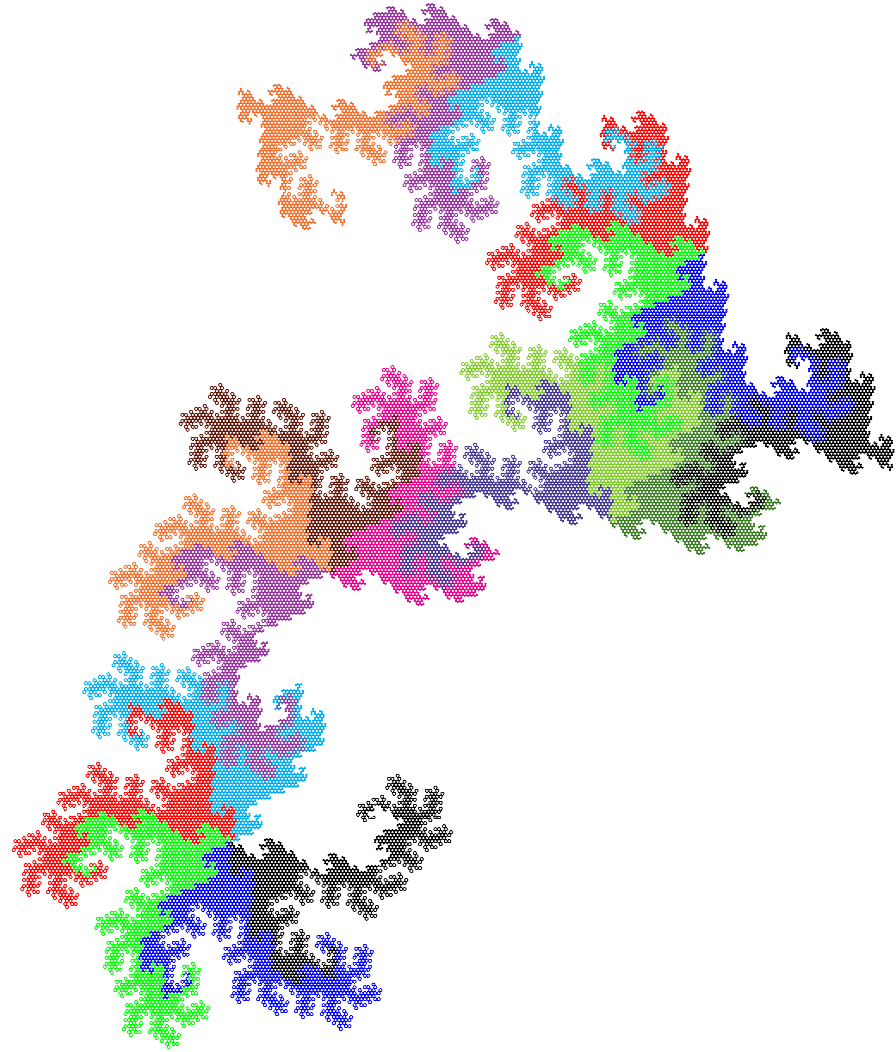}}
\end{center}
\else
\verb+{see pdf for image}+
\fi
\caption{\label{fig:R19-1}
\CID{R19-1} on the tri-hexagonal grid.}
\end{figure}
%
%%%%%%%%%%%%%%%%%%%%%%%%%%

%%%%%%%%%%%%%%%%%%%%%%%%%%
%
\begin{figure}[h!tbp]
\ifpdf
\begin{center}
%\fbox{\includegraphics*[width=90mm, page=2, viewport={60 340 490 740}]{all-r19-b-decomp.pdf}}
{\includegraphics*[width=130mm, page=2, viewport={60 280 490 740}]{all-r19-b-decomp.pdf}}
\end{center}
\else
\verb+{see pdf for image}+
\fi
\caption{\label{fig:R19-2}
\CID{R19-2} on the tri-hexagonal grid.}
\end{figure}
%
%%%%%%%%%%%%%%%%%%%%%%%%%%

%%%%%%%%%%%%%%%%%%%%%%%%%%
%
\begin{figure}[h!tbp]
\ifpdf
\begin{center}
%\fbox{\includegraphics*[width=90mm, page=3, viewport={60 370 490 740}]{all-r19-b-decomp.pdf}}
{\includegraphics*[width=130mm, page=3, viewport={60 370 490 740}]{all-r19-b-decomp.pdf}}
\end{center}
\else
\verb+{see pdf for image}+
\fi
\caption{\label{fig:R19-3}
\CID{R19-3} on the tri-hexagonal grid.}
\end{figure}
%
%%%%%%%%%%%%%%%%%%%%%%%%%%

%%%%%%%%%%%%%%%%%%%%%%%%%%
%
\begin{figure}[h!tbp]
\ifpdf
\begin{center}
%\fbox{\includegraphics*[width=90mm, page=4, viewport={60 420 490 740}]{all-r19-b-decomp.pdf}}
{\includegraphics*[width=130mm, page=4, viewport={60 395 490 740}]{all-r19-b-decomp.pdf}}
\end{center}
\else
\verb+{see pdf for image}+
\fi
\caption{\label{fig:R19-5}
\CID{R19-5} on the tri-hexagonal grid.}
\end{figure}
%
%%%%%%%%%%%%%%%%%%%%%%%%%%

%%%%%%%%%%%%%%%%%%%%%%%%%%
%
\begin{figure}[h!tbp]
\ifpdf
\begin{center}
%\fbox{\includegraphics*[width=90mm, page=5, viewport={60 320 490 740}]{all-r19-b-decomp.pdf}}
{\includegraphics*[width=130mm, page=5, viewport={60 370 490 740}]{all-r19-b-decomp.pdf}}
\end{center}
\else
\verb+{see pdf for image}+
\fi
\caption{\label{fig:R19-6}
\CID{R19-6} on the tri-hexagonal grid.}
\end{figure}
%
%%%%%%%%%%%%%%%%%%%%%%%%%%

%%%%%%%%%%%%%%%%%%%%%%%%%%
%
\begin{figure}[h!tbp]
\ifpdf
\begin{center}
%\fbox{\includegraphics*[width=90mm, page=6, viewport={60 330 490 740}]{all-r19-b-decomp.pdf}}
{\includegraphics*[width=130mm, page=6, viewport={60 420 490 740}]{all-r19-b-decomp.pdf}}
\end{center}
\else
\verb+{see pdf for image}+
\fi
\caption{\label{fig:R19-7}
\CID{R19-7} on the tri-hexagonal grid.}
\end{figure}
%
%%%%%%%%%%%%%%%%%%%%%%%%%%

%%%%%%%%%%%%%%%%%%%%%%%%%%
%
\begin{figure}[h!tbp]
\ifpdf
\begin{center}
%\fbox{\includegraphics*[width=90mm, page=7, viewport={60 370 490 740}]{all-r19-b-decomp.pdf}}
{\includegraphics*[width=130mm, page=7, viewport={60 420 490 740}]{all-r19-b-decomp.pdf}}
\end{center}
\else
\verb+{see pdf for image}+
\fi
\caption{\label{fig:R19-14}
\CID{R19-14} on the tri-hexagonal grid.}
\end{figure}
%
%%%%%%%%%%%%%%%%%%%%%%%%%%

%%% Emacs:
%%% Local Variables:
%%% mode: latex
%%% MyRelDir: "."
%%% TeX-master: "arndt-curve-search.tex"
%%% dvi-file: "arndt-curve-search"
%%% makefile-dir: "./"
%%% frame-title-format: "CURVE-SEARCH (beholder)"
%%% End:

\FloatBarrier

%
% xxx thin curves and mantas
% xxx borders
% xxx 3D and beyond
%

%%%%%%%%%%%%%%%%%%%%%%%%%%%%%%%%%%%%%%%%%%%%%%%%%%%%%%%%%%%%%%
%%%%%%%%%%%%%%%%%%%%%%%%%%%%%%%%%%%%%%%%%%%%%%%%%%%%%%%%%%%%%%
\section*{Acknowledgments}
%
% see file thanks.txt
%
It is my pleasure to thank the following people
for their support, criticism, corrections, and improvements.
\begin{quote}
Michael Beeler,
Michel Dekking,
Brian Galebach,
Bill Gosper,
Julia Handl,
Michael Lontke,
Tarkan Nein,
Edith Parzefall,
Andreas Pazureck,
James Propp,
Kevin Ryde,
%% Subject: uniform grids E curve
%% Message-ID: <877fadioau.fsf@blah.blah>
%
Jeffrey Shallit,
Marcus Veitenthal,
and
Allan Wechsler.
\end{quote}

%\clearpage% xxx
%%% bibliography

%%%%%%%%%%%%%%%%%%%%%%%%%%%%%%%%%%%%%%%%%%%%%%%%%%%%%%%%%%%%
\newcommand{\jjbibtitle}[1]{{\small\bfseries #1}}
\newcommand{\bdate}[1]{(#1)}

%%% Emacs:
%%% Local Variables:
%%% mode: latex
%%% MyRelDir: "."
%%% dvi-file: "arndt-curve-search"
%%% makefile-dir: "./"
%%% frame-title-format: "CURVE-SEARCH (bib)"
%%% End:

\end{document}